\documentclass{article} 
\usepackage{iclr2027_conference,times}

\usepackage{amsmath,amsfonts,bm}

\def\1{\bm{1}}

\def\ra{{\textnormal{a}}}

\DeclareMathAlphabet{\mathsfit}{\encodingdefault}{\sfdefault}{m}{sl}
\SetMathAlphabet{\mathsfit}{bold}{\encodingdefault}{\sfdefault}{bx}{n}

\newcommand{\R}{\mathbb{R}}

\DeclareMathOperator*{\argmax}{arg\,max}
\DeclareMathOperator*{\argmin}{arg\,min}

\usepackage{lipsum}
\usepackage{amsfonts}
\usepackage{graphicx}
\usepackage{epstopdf}
\usepackage{algorithmic}
\usepackage{microtype}
\usepackage{subcaption}
\usepackage{booktabs}
\usepackage{amsfonts,mathrsfs, amsthm}
\usepackage{algorithm}
\usepackage{verbatim}
\usepackage{acronym,color}
\usepackage{pifont}
\usepackage[normalem]{ulem}
\usepackage{yhmath}
\usepackage{natbib}
\usepackage{amsmath}
\usepackage{amssymb}
\usepackage{mathtools}
\usepackage{enumerate}
\usepackage{caption}
\usepackage{xcolor}

\usepackage{hyperref}
\usepackage{url}

\def\la{\langle}
\def\ra{\rangle}

\def\cF{{\cal F}}

\def\1{{\bf 1}}

\def\d{\delta}

\def\R{\mathbb{R}}

\newcommand{\EXP}[1]{\mathbb{E}\!\left[#1\right] }

\def\be{\begin{equation}}
\def\ee{\end{equation}}

\newcommand{\abs}[1]{\left|#1\right|}

\newtheorem{theorem}{Theorem}[section]

\newtheorem{lemma}[theorem]{Lemma}
\newtheorem{corollary}[theorem]{Corollary}

\newtheorem{assumption}[theorem]{Assumption}
\newtheorem{remark}[theorem]{Remark}

\title{Universal Extra-Gradient Method for Stochastic Variational Inequalities}

\author{
Abhishek Chakraborty$^{1}$, Anton Rodomanov$^{2}$, Sebastian U.~Stich$^{2}$, and Angelia Nedi\'c$^{1}$ \\[1ex]
$^{1}$School of Electrical, Computer and Energy Engineering, Arizona State University\\
$^{2}$CISPA Helmholtz Center for Information Security\\[1ex]
\texttt{achakr61@asu.edu},
\texttt{anton.rodomanov@cispa.de},\\
\texttt{sebastian.stich@cispa.de},
\texttt{angelia.nedic@asu.edu}
}

\iclrpreprintcopy 
\begin{document}

\maketitle

\begin{abstract}
    We study the stochastic Extra-Gradient Method (EGM) for solving a stochastic variational inequality (VI) over a compact convex set. The main assumption is $(\delta, L)$-continuity of the VI mapping, which includes H\"older mappings and their linear combinations. We develop new adaptive, line-search-free step-size rules for stochastic EGM. These rules can be interpreted as enhanced AdaGrad-type updates and require no prior knowledge of either the continuity modulus of the VI mapping or the variance of the stochastic oracle. The resulting methods are universal and attain the best complexity bound among the admissible continuity regimes. We prove convergence guarantees in terms of an expected residual function that upper bounds the relevant accuracy measures for several important problem classes, including the dual gap for monotone VIs, and the primal-dual gap for saddle-point and optimization problems. For H\"older-continuous mappings, our bounds recover the optimal dimension-independent complexity of first-order methods.
\end{abstract}

\section{Introduction}\label{sec_intro}

Variational inequalities (VIs) provide a unified framework for a wide range of problems in optimization, games, and machine learning. Convex minimization, convex-concave saddle-point problems, monotone inclusions, and equilibrium problems in games can all be formulated as VIs; see, e.g., \citep{facchinei2003finite}. This generality is especially relevant to modern machine learning, where adversarial and robust training~\citep{goodfellow2020generative,madry2018towards}, multi-agent reinforcement learning, and constrained Markov decision processes~\citep{jin2020efficiently} naturally lead to saddle-point or game-theoretic formulations rather than the minimization of a single objective. 

A VI is characterized by an underlying mapping, and methods that use the mapping evaluations  only are attractive because they avoid Jacobian evaluations, which can be computationally expensive. However, simple projection methods may fail to converge even for monotone mappings. 
The extragradient method (EGM)~\citep{korpelevich1976extragradient} 
overcomes this difficulty by evaluating the mapping at an intermediate point before performing the ``main" update. This extra-gradient correction ensures convergence under suitable step-size choices and has led to several extensions, including the Mirror-Prox method in deterministic~\citep{nemirovski2004prox} and stochastic settings~\citep{juditsky2011solving}.

A practical limitation of classical stochastic EGM is its dependence on problem parameters. Standard convergence guarantees typically require step sizes tuned to quantities such as the continuity modulus of the mapping, the Lipschitz or H\"older constants, the H\"older exponent, the diameter of the constraint set, or the variance of the stochastic oracle. In applications, these quantities are often unknown, expensive to estimate, or available only through conservative worst-case bounds. Moreover, a single problem may satisfy several continuity descriptions, and the most favorable one may not be apparent in advance. A step size tuned to the wrong regime can therefore lose much of the practical advantage promised by the theory.


 Universal methods remove this dependence by adapting automatically to the regularity of the problem. In optimization, this idea goes back to Nesterov's universal gradient methods~\citep{nesterov2015universal} and has since led to a broad literature on adaptive and line-search-based schemes for optimization~\citep{yurtsever2015universal, guminov2019universal, kamzolov2021universal, ito2023parameter}. For VIs, line-search variants have also been developed~\citep{stonyakin2021inexact,stonyakin2022generalized}, but they typically require additional mapping evaluations or repeated local tests. While such tests are suitable in deterministic problems, they become delicate in stochastic settings because noisy mapping evaluations can make acceptance conditions unreliable. Moreover, a line-search-free adaptive method is preferable when stochastic oracle calls are costly.


Apart from the line-search techniques, other adaptive algorithms for VIs have also been studied. One line of work develops Golden Ratio Algorithms for deterministic monotone VIs, including an adaptive variant for locally Lipschitz-continuous mappings~\citep{malitsky2020golden}. This variant
estimates the local curvature from consecutive iterates and therefore requires neither prior knowledge of a Lipschitz constant nor a line search. Although its initial step size can be arbitrary in theory, this choice influences the subsequent step sizes in practice, and the original formulation also employs a prescribed upper bound on the step sizes. A method with a similar flavor was studied under Bregman continuity in~\citep{antonakopoulos2019adaptive}. It was shown to be universal by establishing rates for the bounded-variation and Lipschitz-continuous regimes. However, the method uses a ``shrink ratio" hyperparameter, whose value can be sensitive in practice. 

Another line of fully adaptive and universal methods uses AdaGrad- or AdaGrad-Norm-type updates. In~\citep{antonakopoulos2021sifting}, universal methods for stochastic VIs were developed using variants of AdaGrad-Norm stepsizes. However, the selected stepsizes are always smaller than $1$ because of suboptimal choices of the associated constants. Moreover, the analysis assumes that the stochastic mapping is uniformly bounded almost surely, which is a strong assumption. A similar bounded-mapping assumption was used in~\citep{bach2019universal}, where better constants allow the AdaGrad-Norm stepsizes to exceed $1$. However, the method requires an arbitrary initial estimate of the mapping-norm bound, treated as a hyperparameter, as well as knowledge of the diameter of the underlying bounded constraint set. Its convergence rate also incurs an additional $\sqrt{\log T}$ factor in the stochastic term, making it suboptimal. 
Later work~\citep{ene2021adaptive} developed an AdaGrad method for EGM that requires no hyperparameter other than an upper bound on the $\ell_\infty$-diameter of the underlying constraint set, which may be unknown in practice. Although this method is universal and addresses an important gap, its setting differs from ours since it uses coordinate-wise AdaGrad updates and an $\ell_\infty$-based analysis. Moreover, its stochastic rate depends on the problem dimension and incurs an additional $\sqrt{\log T}$ factor in the stochastic term. Inspired by work~\citep{rodomanov2024universal} in optimization, an attempt to close this gap by developing an adaptive, line-search-free, and universal stepsize for stochastic VIs was made in~\citep{klimza2024universal}. The work establishes rates for H\"older-continuous mappings but requires knowledge of the diameter of the underlying constraint set; some steps of its analysis appear to be incomplete, as discussed in~Appendix~\ref{app_klimza}.

Universality has no single formal definition, and its meaning varies across the literature. In the stochastic VI literature, a method is often called universal if it attains optimal convergence rates for Lipschitz-continuous mappings and mappings with bounded variance, or, more generally, for H\"older-continuous mappings. 
A broader notion of universality is introduced in~\citep{nesterov2025universal} for optimization through a Global Curvature Bound (GCB) that acts as a modulus of smoothness. H\"older-continuous functions and their linear combinations arise as special cases of this framework. Subsequently, \citep{doikov2025universal} extended GCB approach to universal higher-order methods for deterministic VIs using line search. The modulus of smoothness is similar in spirit to $(\delta,L)$-smoothness in optimization~\citep{rodomanov2024universality}. We extend this notion to mappings requiring that the constant $L$ depends on  $\delta$, $\delta>0$, (see Section~\ref{sec_problem}). To the best of our knowledge, no universal method for stochastic VIs exists that avoids all problem parameters and hyperparameters, requires no line search, and adapts to a broader class of mappings encompassing several continuity regimes.

\textbf{Contributions:} Our main contributions are summarized below:

(1) We introduce the class of $(\delta,L)$-continuous VI mappings, an approximate-continuity condition tailored to variational inequalities. This class includes H\"older-continuous mappings, linear combinations of mappings with different continuity regimes, and mappings arising from inexact computations. This framework provides a unified way to analyze stochastic VIs beyond a single prescribed class.

(2) We develop adaptive, line-search-free step-size rules for the stochastic Extra-Gradient Method. The proposed rules are derived from balance equations that control the key error terms in the EGM analysis. They can be viewed as enhanced AdaGrad-type updates and do not require prior knowledge of the continuity modulus, H\"older exponent, Lipschitz constant, stochastic-oracle variance, or the exact diameter of the constraint set.


(3) We prove a unified convergence guarantee for the proposed methods in terms of an expected residual function. This residual function upper bounds natural accuracy measures across several important settings, including the dual gap for monotone VIs, and the primal-dual gap for saddle-point and optimization problems. For H\"older-continuous mappings, our bounds recover the optimal dimension-independent complexity of first-order methods.

(4) We further develop distance-adaptive variants that reduce sensitivity to the initial diameter estimate. These variants preserve the same optimal convergence rates while improving the dependence of the constants on the quality of the initial estimate. Numerical experiments on optimization and min-max problems illustrate the practical effectiveness of the proposed adaptive rules.

\paragraph{Notation.} 
For any vectors $x, y \in \mathbb{R}^n$, we denote $\langle x, y \rangle$ as their standard inner product. We denote the weighted Euclidean norm of a vector $x \in \mathbb{R}^n$ induced by a fixed positive definite matrix $B \in \R^{n \times n}$ as $\|x\| = \sqrt{\langle B x, x \rangle}$. The dual norm is defined as $\|g\|_* = \max_{\|x\| \leq 1} \la g,x \ra = \sqrt{\la g, B^{-1} g \ra}$. For the standard Euclidean norm, we have $B = I$, where  $I$ is the identity matrix. For a scalar $a \in \mathbb{R}$, we use the notation $[a]_+ = \max\{0, a\}$. We use $\mathbb{E}[\cdot]$ for the expectation mapping.

\section{Problem Formulation}\label{sec_problem}

We consider the following Variational Inequality (VI) problem:
\begin{align}
    \text{Find $x^* \in X$}\colon \quad \langle F(x^*),\, x - x^* \rangle \geq 0 \qquad \forall x \in X , \label{VI_prob}
\end{align}
where $X \subseteq \R^n$ is a nonempty closed convex set, and $F \colon X \to \R^n$ is a monotone mapping, i.e., 
\begin{align}
    \langle F(x) - F(y), x - y \rangle \geq 0\qquad \forall x, y \in X. \label{monotone_map}
\end{align}
We assume that the constraint set $X$ is easy to project onto, so that projection onto $X$ can be performed efficiently. We use the following assumption on the constraint set $X$.
\begin{assumption}\label{asum_set}
    The constraint set $X$ in problem~\eqref{VI_prob} is bounded, i.e., there exists a constant $D > 0$ such that $\max_{x,y \in X} \|x-y\| \leq D$.
\end{assumption}
A point $x^*$ satisfying inequality~\eqref{VI_prob} is usually called a \emph{strong solution} of the VI corresponding to $(X, F)$. In this regard, we would like to define also the concept of a weak (or Minty) solution:
\begin{align}
    \text{Find $x_* \in X$}\colon \quad \langle F(x),\, x - x_* \rangle \geq 0 \qquad \forall x \in X . \label{weak_soln}
\end{align}
Assumption~\ref{asum_set} ensures that the Minty solution exists by~\citep[Theorem~4]{lassonde2015convex}; this result does not require continuity of $F$.
If $F$ is monotone, each strong solution~\eqref{VI_prob} of the VI is also a weak solution. Conversely, if $F$ is continuous, each weak solution is also a strong one~\citep{minty1962monotone}, \citep[Lemma 1.5, page 84]{kinderlehrer2000introduction}. Therefore, when $F$ is both monotone and continuous, the two solution concepts coincide.

\subsection{$(\delta, L)$-continuity of Mapping}\label{delta_L_cont_subsec}

We first state the notion of $(\delta,L)$-smoothness in optimization, which has been studied in~\citep{devolder2014first,nesterov2015universal,rodomanov2024universality,savchuk2025adaptive}, among others. Consider the optimization problem $\min_{x\in X_1} f(x)$, where $f$ is convex and $X_1$ is a convex set. The function $f$ is said to be $(\delta,L)$-smooth over the set $X_1$ if
\begin{align}
0
\leq f(x)-\left(f(y)+\langle\nabla f(y),x-y\rangle\right)
\leq \frac{L}{2}\|x-y\|^2+\delta
\qquad \forall x,y\in X_1.
\label{opt_deltaL}
\end{align}
Work~\citep{nesterov2025universal} generalized this definition using a modulus of smoothness. This notion is similar in spirit to~\eqref{opt_deltaL}, which can be imposed for every $\delta_1>0$, with some $L\equiv L(\delta_1)$ depending on $\delta_1$. However for VIs, only mapping evaluations are available, and there is no corresponding notion of objective-function values. Therefore, the approximate-smoothness condition from optimization must be adapted to define approximate continuity for VI mappings. Previous works on VIs~\citep{juditsky2011solving,chakraborty2025popov} consider mappings satisfying 
\[\|F(x) - F(y)\|_* \leq L\|x-y\| + \delta \qquad\forall x, y \in X.\] 

In this work, we consider a reparameterized version of the preceding relation and say that $F$ is $(\delta,L)$-continuous on $X$ if, for every $\delta>0$, there exists $L(\delta)>0$ such that
\begin{align}
    \|F(x) - F(y)\|_*^2 \leq L^2(\delta) \|x - y\|^2 + 4 L(\delta) \delta \qquad \forall x, y \in X. \label{eqasum_deltaL}
\end{align}
Definition~\eqref{eqasum_deltaL} of $(\delta,L)$-continuous mapping $F$ is analogous to the one used in the optimization literature\footnote{If $F = \nabla f$ and $X \subset X_1$, then relation~\eqref{opt_deltaL} implies relation~\eqref{eqasum_deltaL}~\citep[Corollary~1]{devolder2014first}.}. In what follows, we assume that $F$ is $(\delta,L)$-continuous mapping.

\begin{assumption}\label{asum_map_growth}
    The mapping $F$ is $(\delta,L)$-continuous over the set $X$ (cf.\ relation~\eqref{eqasum_deltaL}).
\end{assumption}

Assumption~\ref{asum_map_growth} may arise in settings where the mapping $F$ is computed using an inexact oracle (see~\citep{devolder2014first} for examples in optimization). An important example of $(\delta,L)$-continuous mapping is given by a H\"older-continuous mapping.

\begin{lemma}\label{lem_holder_to_delta_L}
    Let $F$ be a H\"older continuous mapping, i.e., there exists $\nu \in [0,1]$ and $H_\nu > 0$ such that
    \begin{align}
        \|F(x) - F(y)\|_* \leq H_\nu \|x - y\|^\nu  \qquad \forall x, y \in X. \label{holder_growth_reln}
    \end{align}
    Then, $F$ is $(\delta,L)$-continuous with 
    \begin{align*}
        L(\delta) = \left(\frac{1-\nu}{4\delta}\right)^{\frac{1-\nu}{1+\nu}}\nu^{\frac{\nu}{1+\nu}} H_\nu^{\frac{2}{1+\nu}},
    \end{align*}
    where we use the convention $0^0 = 1$.
\end{lemma}

The $(\delta,L)$-continuity properties are preserved under the following calculus rules: (i)~scalar multiplication, (ii)~summation, and (iii)~composition with an affine mapping (see Lemma~\ref{lem_calculus}). In particular, the sum of H\"older mappings with different H\"older exponents is also $(\delta,L)$-continuous, even though the resulting mapping may not be H\"older~\citep[Example~1]{nesterov2025universal}, where we do not take into account Assumption~\ref{asum_set}. Thus, the class of $(\delta,L)$-continuous mappings is wider than that of H\"older-continuous ones.

\begin{remark}\label{rem_vi_opt_comp}
    In the context of smooth optimization, we have $F(\cdot) = \nabla f(\cdot)$, a different $(\delta, L)$ condition than that of relation~\eqref{opt_deltaL} follows from H\"older-continuity of the gradient with $\tilde{L}(\delta) = \left(\frac{1 - \nu}{2 (1 + \nu) \delta}\right)^{\frac{1 - \nu}{1 + \nu}} H_{\nu}^{\frac{2}{1 + \nu}}$ \citep[Remark F.1]{rodomanov2024universality}.
    For a H\"older-continuous mapping $F$, Lemma~\ref{lem_holder_to_delta_L} derives $L(\delta)$ which satisfies $L(\delta)\leq \tilde{L}(\delta)$, as shown in Lemma~\ref{lem_L_delta}.
\end{remark}



%
%

\subsection{Residual Function}\label{subsec_residual}

Let $(x_t)_{t=1}^{k}$ be a sequence of iterates generated by an algorithm up to iteration $k \geq 1$, where $x_t \in X$ for all $1 \leq t \leq k$, and let $(\lambda_t)_{t=1}^{k}$ be a sequence of nonnegative weights. We introduce the accuracy certificate $(x_t, \lambda_t)_{t=1}^k$ as the main object of interest. Define $\Lambda_k = \sum_{t=1}^k \lambda_t$ and assume that $\Lambda_k > 0$ for all $k \geq 1$. We consider the residual function $\epsilon_k$ associated with the accuracy certificate $(x_t, \lambda_t)_{t=1}^k$ over the set $X$, defined as~\citep{nesterov2009primal, nemirovski2010accuracy, rodomanov2023subgradient},
\begin{align}
    \epsilon_k \equiv \epsilon_k((x_t, \lambda_t)_{t=1}^k) := \max_{x \in X} \frac{1}{\Lambda_k} \sum_{t = 1}^{k} \lambda_t \la F(x_t), x_t - x \ra. \label{residual_val_eq}
\end{align}
In view of relation~\eqref{weak_soln}, we have $\epsilon_k \geq 0$. Our goal is to develop a universal and parameter-free algorithm that produces an accuracy certificate whose residual is bounded by a prescribed tolerance $\varepsilon$. Note that this residual is a weighted average of the well-known gap function~\citep[equation (1.5.15)]{facchinei2003finite} studied in the variational inequality literature. Hence, it serves as a meaningful measure for analyzing the convergence of an algorithm. Below, we provide several examples illustrating the applicability of this residual, similar to those studied in~\citep{rodomanov2023subgradient}.

i) \textbf{Convex minimization problem:} Let us consider the problem $f^* = \min_{x \in X} f(x)$, where $f : X \rightarrow \mathbb{R}$ is a differentiable convex function, $X$ is a convex constraint set, and let $F = \nabla f$. Then, for $\widehat{x}_k = \frac{1}{ \Lambda_k} \sum_{t = 1}^k \lambda_t x_t$, we have
    \begin{align*}
        f(\widehat{x}_k) - f^* \leq \epsilon_k .
    \end{align*}
    Therefore, establishing a convergence rate for the quantity $\epsilon_k$ directly implies the same convergence rate for the function values evaluated at the averaged iterates. Moreover, the problem admits the dual formulation $f^* = \max_{s \in \mathrm{dom}\, f_*} \left[ - f_*(s) - \xi_X(-s) \right]$, where $f_*$ is the conjugate of $f$ and $\xi_X(-s) := \max_{x \in X} \langle -s, x \rangle$. In this setting, $\epsilon_k$ is an upper bound for the corresponding primal-dual gap~\citep{rodomanov2023subgradient}, i.e., for all $k \geq 1$ and $\widehat s_k := \frac{1}{\Lambda_k} \sum_{t=1}^k \lambda_t \nabla f(x_t)$,
    \begin{align*}
        f(\widehat x_k) + f_*(\widehat s_k) + \xi_X(- \widehat s_k) \leq \epsilon_k .
    \end{align*}

ii) \textbf{Convex-concave min-max problem:} For convex and closed sets $U \subseteq \mathbb{R}^{n_1}$ and $V \subseteq \mathbb{R}^{n_2}$, consider the saddle-point problem $\min_{u \in U} \max_{v \in V} f(u,v)$, where the function $f: U \times V \rightarrow \mathbb{R}$ is convex in $u$ for every fixed $v \in V$ and concave in $v$ for every fixed $u \in U$. With $F(u,v) = (\nabla_u f(u,v),\; -\nabla_v f(u,v) )$, and $\widehat u_k = \frac{1}{\Lambda_k} \sum_{t=1}^k \lambda_t u_t$ and $\widehat v_k = \frac{1}{\Lambda_k} \sum_{t=1}^k \lambda_t v_t$, we have
    \begin{align*}
        \max_{v \in V} f(\widehat u_k,v) - \min_{u \in U} f(u, \widehat v_k) \leq \epsilon_k .
    \end{align*}

iii) \textbf{Monotone Variational Inequality:}
    Consider VI problem~\eqref{VI_prob}, with a monotone $F$ (cf.\ relation~\eqref{monotone_map}), and its  dual gap function~\citep{auslender1974resolution, marcotte1998weak, nemirovski2004prox}
    \begin{align*}
        \bar G(y) =\max_{x \in X} \la F(x), y - x \ra.
    \end{align*}
    Then, for $\widehat{x}_k = \frac{1}{\Lambda_k} \sum_{t=1}^k \lambda_t x_t$, we have for all $k \geq 1$,
    \begin{align*}
        \bar G(\widehat x_k) \leq \epsilon_k.
    \end{align*}
    The dual gap function characterizes a weak (or Minty) solution (cf.\ relation~\eqref{weak_soln}). It is always non-negative and takes zero value only when $y = x_*$ is a weak solution. 

    %
    %
    %
    %



\section{Algorithm for Stochastic Variational Inequalities}\label{sec_algo}

We consider a stochastic VI with the mapping $F(x) = \EXP{\widehat F(x,\omega)}$, $x \in X$, where $\omega$ is a random variable. 
 We use the following assumption on the stochastic realizations of the mapping $\widehat F(x,\omega)$.
\begin{assumption}\label{asum_stochastic}
    The mapping $\widehat F(x,\omega)$ is unbiased and, for some $\sigma > 0$, has a bounded variance
    \begin{align*}
        \EXP{\widehat F(x,\omega)} = F(x) \quad \text{and} \quad \EXP{\|\widehat F(x,\omega) - F(x)\|_*^2} \leq \sigma^2 \qquad\forall x\in X.
    \end{align*}
\end{assumption}

For solving stochastic VI problem~\eqref{VI_prob}, we consider the stochastic Extra-Gradient Method (EGM)~\citep{korpelevich1976extragradient}, with iterates given by
for $k \geq 0$:
\begin{align}
    &x_{k+1} = \argmin_{x \in X} \left\{ \langle \widehat F(v_k, \omega_k),\, x - v_k \rangle + \frac{H_{k}}{2} \|x - v_k\|^2 \right\}, \label{Kor_algo} \\
    &v_{k+1} = \argmin_{x \in X} \left\{ \langle \widehat F(x_{k+1}, \omega_{k+1}),\, x - x_{k+1} \rangle + \frac{H_{k}}{2} \|x - v_k\|^2 \right\}, \nonumber
\end{align}
where $H_0 = 0$, $H_k \geq 0$ for $k \geq 1$, and $(\omega_k)_{k \geq 0}$ is a sequence of stochastic samples. 
The method starts with some random initial point $v_0 \in X$.
%
%
Under Assumptions~\ref{asum_set} and \ref{asum_stochastic}, for the method~\eqref{Kor_algo}, we obtain the following upper bound on the expected residual function for any $k \geq 0$ (cf.\ Lemma~\ref{lem_generic_reln}):
\begin{align}
    \mathbb{E}[\epsilon_T] &\leq \frac{1}{T}\sum_{k=0}^{T-1} \mathbb{E}[\Gamma_k] +
    \frac{D \sigma}{\sqrt{T}} \quad \text{where} \quad \Gamma_k \equiv \Delta_k(H_{k+1}) + \frac{3}{2} (H_{k+1} - H_{k}) D^2. \label{eq1_lem_generic_reln_dup}
\end{align}
The quantity $\Delta_k(H) \equiv \xi_k - \frac{H}{2} r_k^2$ with $r_k^2 = \|v_k - x_{k+1}\|^2 + \|v_{k+1} - x_{k+1}\|^2$ and
\begin{align*}
    \xi_k &= \langle \widehat F(x_{k+1}, \omega_{k+1}), x_{k+1} - v_{k+1} \rangle - \frac{H_k}{2} \|v_k - v_{k+1}\|^2 + \frac{H_k}{2} r_k^2 \nonumber\\
    &\leq \langle \widehat F(x_{k + 1}, \omega_{k+1}) - \widehat F(v_k, \omega_k), x_{k + 1} - v_{k + 1} \rangle \equiv \xi_k'.
\end{align*}
The relation $\Delta_k(H_{k+1}) = \xi_k - \frac{H_{k+1}}{2} r_k^2$ is established using the optimality condition for the $v_k$-update and is valid for any sequence $(x_k)_{k\geq 0}$. 
The optimality conditions for both the $v_k$- and $x_k$-updates imply
$\xi_k\le \xi'_k$. Relation~\eqref{eq1_lem_generic_reln_dup} is the key inequality guiding the choice of the sequence $(H_k)_{k\geq 0}$. In order to obtain the convergence rate with respect to the expected residual function, we require an upper bound on the quantity $\Gamma_k$. Its second quantity can be controlled via summation, and we require the upper bound on $\Delta_k(H_{k+1})$, which motivates the use of balance equations for updating $H_{k+1}$. In particular,  by controlling either
$\Delta_k(H_{k+1}) = \xi_k - \frac{H_{k+1}}{2} r_k^2$ or its upper bound $\xi_k' - \frac{H_{k+1}}{2} r_k^2$, we will obtain update rules for $H_{k+1}$, which we will describe next. 

Let $\widehat D > 0$ be an arbitrary estimate of the true diameter $D$ of the set $X$. We update $(H_{k+1})_{k \geq 0}$ through a balance equation that controls $\Gamma_k$ and yields, for all $k\geq0$, the upper bound (cf.\ Lemma~\ref{lem_algo1_2_gen})
\begin{align}
    \EXP{\Gamma_k} \leq \frac{\EXP{\left[\min\{H_{k+1}, \bar H\} - H_k \right]_+}}{2} (c_1 D^2 + c_2 \widehat D^2) + \frac{4 L \delta + 2 \sigma^2}{\bar H \widehat D^2} (3 D^2 + \widehat D^2) , \label{gam_k_reln}
\end{align}
where $\bar H \geq \sqrt{2} \bar L(\delta)$, for $\delta > 0$, is a parameter that is a byproduct of the analysis and not required to implement the methods. We choose $\bar H$ optimally to obtain sharper convergence rates. Next, we present several examples that satisfy the preceding relation with different constants $c_1$, $c_2$, and $\bar L(\delta)$.

To control either $\Delta_k(H_{k+1})$ in relation~\eqref{eq1_lem_generic_reln_dup} or its upper bound
$\xi_k'-\frac{H_{k+1}}{2}r_k^2$, we choose $H_{k+1}$ to satisfy one of the following balance equations: for all $k\geq0$,
\begin{align}
    \frac{H_{k+1} - H_k}{2}  \widehat D^2 &= [\Delta_k(H_{k + 1})]_+ \quad &\text{(Rule~1),} \label{baleqn1} \\
    \frac{H_{k+1} - H_k}{2}  \widehat D^2 &= \left[\xi_k' - \frac{H_{k+1}}{2} r_k^2 \right]_+ \quad &\text{(Rule~2).} \label{baleqn2}
\end{align}
The balance equations~\eqref{baleqn1} and \eqref{baleqn2} yield relation~\eqref{gam_k_reln} with constants $c_1 = 3$, $c_2 = 1$, $\bar L(\delta) = L(\delta)$. The explicit update and algorithmic details for Rules 1 and 2 are provided in Appendix~\ref{app_subsec_algos}. In the context of optimization, a similar balance equation~\eqref{baleqn2} appears with $\widehat D^2 = 2D^2$ in~\citep{rodomanov2024universal, rodomanov2024universality}, where the knowledge of the set diameter $D$ (cf.\ Assumption~\ref{asum_set}) is used. However, this requirement can be relaxed by modifying the analysis in~\citep{rodomanov2024universal, rodomanov2024universality}, which allows the use of an arbitrary $\widehat D > 0$, instead. 
Along this line, an update rule simpler than the classic AdaGrad-norm can be obtained via the following balance equation: for all $k \geq 0$, 
\begin{align}
    (H_{k+1}^2 - H_k^2) \widehat D^2 = \|\widehat F(x_{k + 1}, \omega_{k+1}) - \widehat F(v_k, \omega_k) \|_*^2 \qquad \text{(Rule~3),}  \label{balance_adagd_norm}
\end{align}
which also satisfies relation~\eqref{gam_k_reln} with the constants $c_1 = 3$, $c_2 = 2$, and $\bar L(\delta) = L(\delta) \sqrt{1+ 3D^2/ \widehat D^2}$.


\section{Convergence Rate Analysis}\label{sec_convergence}

In this section, 
a unified convergence rate analysis is provided for the iterates of the stochastic EGM in~\eqref{Kor_algo} with update rules given by
balance equations~\eqref{baleqn1}--\eqref{balance_adagd_norm}.

\begin{theorem}\label{thm_universal}
    Let Assumptions~\ref{asum_set}, \ref{asum_map_growth}, and \ref{asum_stochastic} hold. Let $\widehat D > 0$ and $\d>0$ be arbitrary, and consider update Rules~1--3 (balance equations~\eqref{baleqn1}--\eqref{balance_adagd_norm}). Then, the following relation holds for the expected residual function value 
    \begin{align*}
        \mathbb{E}[\epsilon_T]
        \leq
        \frac{\bar{L}(\delta)\bar{D}_1^2}{\sqrt{2}\,T} +
        \frac{2\sqrt{2L(\delta)\delta}}{\sqrt{T}}
        \left[ \frac{\bar{D}_1\bar{D}_2}{\widehat{D}} \right] + \frac{\sigma}{\sqrt{T}} \left[D + \frac{2\bar{D}_1\bar{D}_2}{\widehat{D}} \right] \qquad \forall T \geq 1,
    \end{align*}
    where $\epsilon_T$ is the residual function (cf.\ relation~\eqref{residual_val_eq}), $D$ is the diameter of the set $X$, $L(\delta)$ is the constant characterizing $(\delta,L)$-continuity of $F$ (cf.~\eqref{eqasum_deltaL}),
    $\bar D_1^2 = c_1 D^2 + c_2 \widehat D^2$ and $\bar D_2^2 = 3D^2 + \widehat D^2$. The constants $c_1,c_2,$ and $\bar L(\delta)$ depend on Rules~1--3, as follows:
    $c_1 = 3$, $c_2 = 1$, $\bar L(\delta) = L(\delta)$ for Rules~1 and 2, and $c_1 = 3$, $c_2 = 2$, $\bar L(\delta) = L(\delta) \sqrt{1 + 3 D^2/\widehat D^2}$ for Rule~3.
\end{theorem}
Theorem~\ref{thm_universal} shows that the 
expected residual function is bounded by $(\delta,L)$-continuity related term ${\bar L(\delta)\bar D_1^2}/{T} + {\sqrt{L(\delta)\delta}\,\bar D_1\bar D_2}/{(\widehat D\sqrt{T})}$ and the variance-based term $\frac{\sigma}{\sqrt{T}}
\left[{2\bar D_1\bar D_2}/{\widehat D}+D\right]$, both decaying at the rate of $O(1/\sqrt{T})$. Our method is universal as it achieves the desired rate for stochastic VIs under Assumption~\ref{asum_map_growth} without knowledge of any problem-dependent parameters. In particular, the error due to $(\delta,L)$-continuity decays in the order of 
$O(1/\sqrt{T})$, rather than producing the non-vanishing
$O(\delta)$ error floor observed in the optimization literature~\citep{devolder2014first,rodomanov2024universality}, where a different approximate-smoothness condition is imposed on the objective function (cf.~\eqref{opt_deltaL}). 
The bound in Theorem~\ref{thm_universal} can be optimized over $\delta>0$, since it holds  for any $\delta>0$.

For a 
H\"older continuous mapping, using 
Lemma~\ref{lem_holder_to_delta_L}, we have the following
corollary of Theorem~\ref{thm_universal}.

\begin{corollary}\label{cor_holder}
    Let Assumptions~\ref{asum_set} and \ref{asum_stochastic} hold, and let the mapping $F$ be H\"older continuous (cf.\ relation~\eqref{holder_growth_reln}) for some $\nu \in [0,1]$ and  $H_\nu>0$. Then, for any $\widehat D > 0$, 
    the expected residual function (cf.\ relation~\eqref{residual_val_eq}) satisfies the following inequality
    \begin{align*}
        \mathbb{E}[\epsilon_T]
        \leq \frac{2^{\frac{1}{2}-\nu}H_\nu}
        {(1-\nu)^{\frac{1-\nu}{2}}\nu^{\frac{\nu}{2}}}
        \left[\frac{\bar D_2}{\widehat D}\right]^{1+\nu(a-1)}
        \left[\frac{\bar D_1^2}{T}\right]^{\frac{1+\nu}{2}}
        +
        \frac{\sigma}{\sqrt{T}}
        \left[\frac{2\bar D_1\bar D_2}{\widehat D}+D
        \right] \qquad \forall T \geq 1,
    \end{align*}
    where $\bar D_1^2 = c_1 D^2 + c_2 \widehat D^2$ and $\bar D_2^2 = 3D^2 + \widehat D^2$, with the constants $c_1 = 3$, $c_2 = 1$, $a = 0$ for Rules~1 and 2, and $c_1 = 3$, $c_2 = 2$, $a = 1$ for Rule~3, and with the convention $0^0=1$.
\end{corollary}
Hence, the balance equations~\eqref{baleqn1}--\eqref{balance_adagd_norm} prescribe rules that achieve the optimal convergence rate for H\"older-continuous mappings. The constant $\widehat D$ in both Theorem~\ref{thm_universal} and Corollary~\ref{cor_holder} may be chosen arbitrarily,
although its relation to $D$ affects the constants in the bound through the ratio $\bar D_2/ \widehat D = \sqrt{1 + 3D^2/ \widehat D^2}$. 


\section{Algorithm with Distance Adaptation}\label{sec_dist_ada}
In this section, to obtain better constants in the convergence rate and eliminate the ratio $\bar D_2 / \widehat D$ appearing in Theorem~\ref{thm_universal} and Corollary~\ref{cor_holder}, we consider the EGM algorithm with distance adaptation. Specifically,
an alternative to the use of a constant $\widehat D>0$,  the distance adaptation~\citep{carmon2022making, ivgi2023dog, khaled2023dowg} provides iterative estimates $D_k$ of the set $X$ diameter $D$. Starting with some initial $D_0 = \widehat D>0$,  the adaptive distance update is given by
\begin{align}
    D_{k+1} = \max \{ D_k, \|v_{k+1} - v_0\|, r_k \} \qquad \forall k \geq 0, \label{adapt_dia_up}
\end{align}
where $r_k=\sqrt{\|v_k - x_{k+1}\|^2 + \|v_{k+1} - x_{k+1}\|^2}$.  

In this case, the iterates of the Distance-Adaptive EGM (DA EGM) are generated as: for all $k\geq 0$,
\begin{align*}
    &x_{k+1} = \argmin_{x \in X} \left\{ D_k \langle \widehat F(v_k, \omega_k),\, x - v_k \rangle + \frac{H_{k}}{2} \|x - v_k\|^2 \right\}, \nonumber \\
    &v_{k+1} = \argmin_{x \in X} \left\{ D_k \langle \widehat F(x_{k+1}, \omega_{k+1}),\, x - x_{k+1} \rangle + \frac{H_{k}}{2} \|x - v_k\|^2 \right\} . \nonumber
\end{align*}
starting with a random initial point $v_0\in X$. Hence, both EGM (cf.\ relation~\eqref{Kor_algo}) and DA EGM are explicitly the same algorithm but with the effective step sizes as $1/H_k$ and $D_k/H_k$ respectively for all $k \geq 1$. At each iteration $k$, after performing the DA EGM updates, we update the diameter $D_{k+1}$ according to~\eqref{adapt_dia_up} and then use an appropriate error balance equation similar to those for Rules~1 and~2 (cf.\ relations~\eqref{baleqn1} and~\eqref{baleqn2}), but with an additional scaling factor involving distances. We refer to such balance equations as ``Methods~1 and~2" which are presented below, i.e., for all $k \geq 0$, 
\begin{align}
    \frac{H_{k+1} - H_k}{2} D_{k+1}^2 &= \left[\widehat \Delta_k(H_{k+1}) \right]_+ = \left[\widehat \xi_k - \frac{H_{k+1}}{2} r_k^2 \right]_+ \qquad &\text{(Method 1)} , \label{bal_ada1}\\
    \frac{H_{k+1} - H_k}{2} D_{k+1}^2 &= \left[\widehat \xi_k' - \frac{H_{k+1}}{2} r_k^2 \right]_+ \qquad &\text{(Method 2)} \label{bal_ada2},
\end{align}
\begin{align*}
    \hspace{-1cm} \text{where} \quad \widehat \xi_k &\equiv D_k \langle \widehat F(x_{k+1}, \omega_{k+1}), x_{k+1} - v_{k+1} \rangle - \frac{H_k}{2} \|v_k - v_{k+1}\|^2 + \frac{H_k}{2} r_k^2 \\
    & \leq D_k \la \widehat F(x_{k+1},\omega_{k+1}) - \widehat F(v_k,\omega_k), x_{k+1} - v_{k+1} \ra \equiv \widehat \xi_k' .
\end{align*}
If $D_k = \widehat D$ for all $k \geq 0$, then Methods~1 and 2 reduce to Rules~1 and 2 respectively. A detailed discussion about the algorithm, balance equations, and its related analysis is presented in Appendix~\ref{app_dist_ada}. Under Assumption~\ref{asum_map_growth}, both Methods~1 and 2 are universal because for any $T \geq 1$,
\begin{align*}
    \EXP{\epsilon_{\tau(T)}}
        \leq \rho^{1/T}\log\left( e\rho \right)
        \left[\frac{3\sqrt{2}L(\delta) D^2}{T}
        +
        \frac{
        12D\sqrt{L(\delta) \delta}
        +(6\sqrt{2}+5)\sigma D
        }{\sqrt{T}}
        \right]  \quad \text{(Theorem~\ref{thm_distada})},
\end{align*}
where $\tau(T) := \argmin\limits_{1 \leq t \leq T} \frac{D_t}{\sum_{k=0}^{t-1} D_k}$ and $\rho = \max \left\{1, \tfrac{D}{\widehat D}\right\}$. For H\"older continuous mapping, we have
\begin{align*}
    \EXP{\epsilon_{\tau(T)}}
    &\leq \rho^{1/T} \log\left( e \rho \right)
    \left[
    \frac{
    3 \times 2^{1-\nu/2}H_\nu D^{1+\nu}
    }{
    (1-\nu)^{(1-\nu)/2}T^{(1+\nu)/2}
    }
    +
    (6\sqrt{2}+5)\frac{\sigma D}{\sqrt{T}}
    \right]  \qquad \text{(Corollary~\ref{cor_adadist_holder})} ,
\end{align*}
with the convention $0^0 = 1$. In comparison to Theorem~\ref{thm_universal} and Corollary~\ref{cor_holder}, the factor $D/\widehat D$ does not appear in the preceding results for distance-adaptive methods. 
If $\widehat D<< D$, then the non-distance-adaptive bounds may contain a large factor proportional to $D/\widehat D$. In contrast, in the distance-adaptive methods the factor $\rho^{1/T}\log(e\rho)$ is milder, approaching $\log(e\rho)$ as $T$ increases.


\section{Experiments}\label{sec_experiments}

We evaluate our algorithms on both optimization and VI problems. For optimization, we consider (i) H\"older-smooth polyhedral feasibility,
(ii) logistic regression, and (iii) softmax problems. We compare our methods EGM with Rules~1 and 2 (aka. Algorithm~\ref{algo_params_free_vi}), EGM with Rule~3 (aka. Algorithm~\ref{algo_adagd_norm}), and Distance-Adaptive EGM (DA EGM) Methods~1 and 2 (aka. Algorithm~\ref{algo_dist_ada_vi}), with two variants of the Universal
Gradient Method (UGM), corresponding to Equations~(3.2) and~(3.3) of
\citet{rodomanov2024universality}, which are denoted by UGM with AdaGrad and UGM,
respectively. We also consider projected AdaGrad-Norm
\citep{ward2020adagrad} and DoG~\citep{ivgi2023dog} methods. Note that EGM with Rules~1--3 are provided in relations~\eqref{baleqn1}--\eqref{balance_adagd_norm}, whereas DA EGM Methods~1 and 2 are provided in relations~\eqref{bal_ada1} and \eqref{bal_ada2}. For the VI experiments, we consider (i) saddle-point formulation of a constrained $\ell_2$-regression problem and compare our methods with EGM where step sizes are selected using grid-search, aka. EGM (grid-tuned), and also with Generalized Extra-Gradient (GEG)~\citep{antonakopoulos2021sifting} and (ii) bilinear saddle point problem, where we compare our methods with GEG. For all the experiments, we also study sensitivity of the selection of $\widehat D$ for EGM Rules~1--3 and $D_0$ for DA EGM Methods~1 and 2 in Appendix~\ref{app_experiments}.

In this section, we present experiments on H\"older-smooth polyhedral feasibility and both the VI saddle-point problems. Additional experiments, implementation details, and results for all problems are provided in Appendix~\ref{app_experiments}. All experiments were conducted on a 2021 MacBook Pro equipped with an Apple M1 Pro chip and 16~GB of RAM.

\begin{figure*}[htbp]
\centering
\begin{subfigure}{0.32\textwidth}
    \centering
    \includegraphics[width=\linewidth]{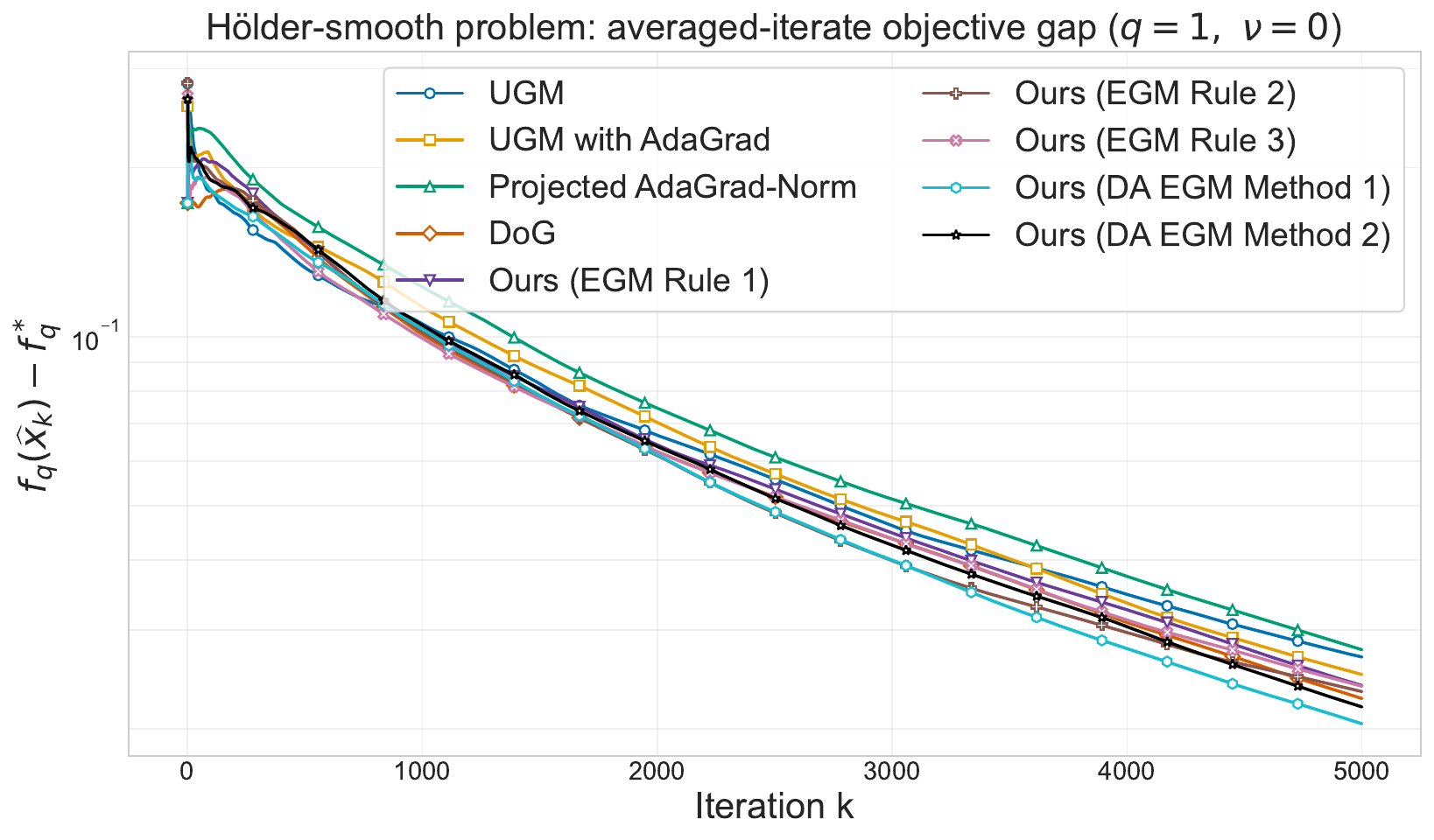}
\end{subfigure}
\hfill
\begin{subfigure}{0.32\textwidth}
    \centering
    \includegraphics[width=\linewidth]{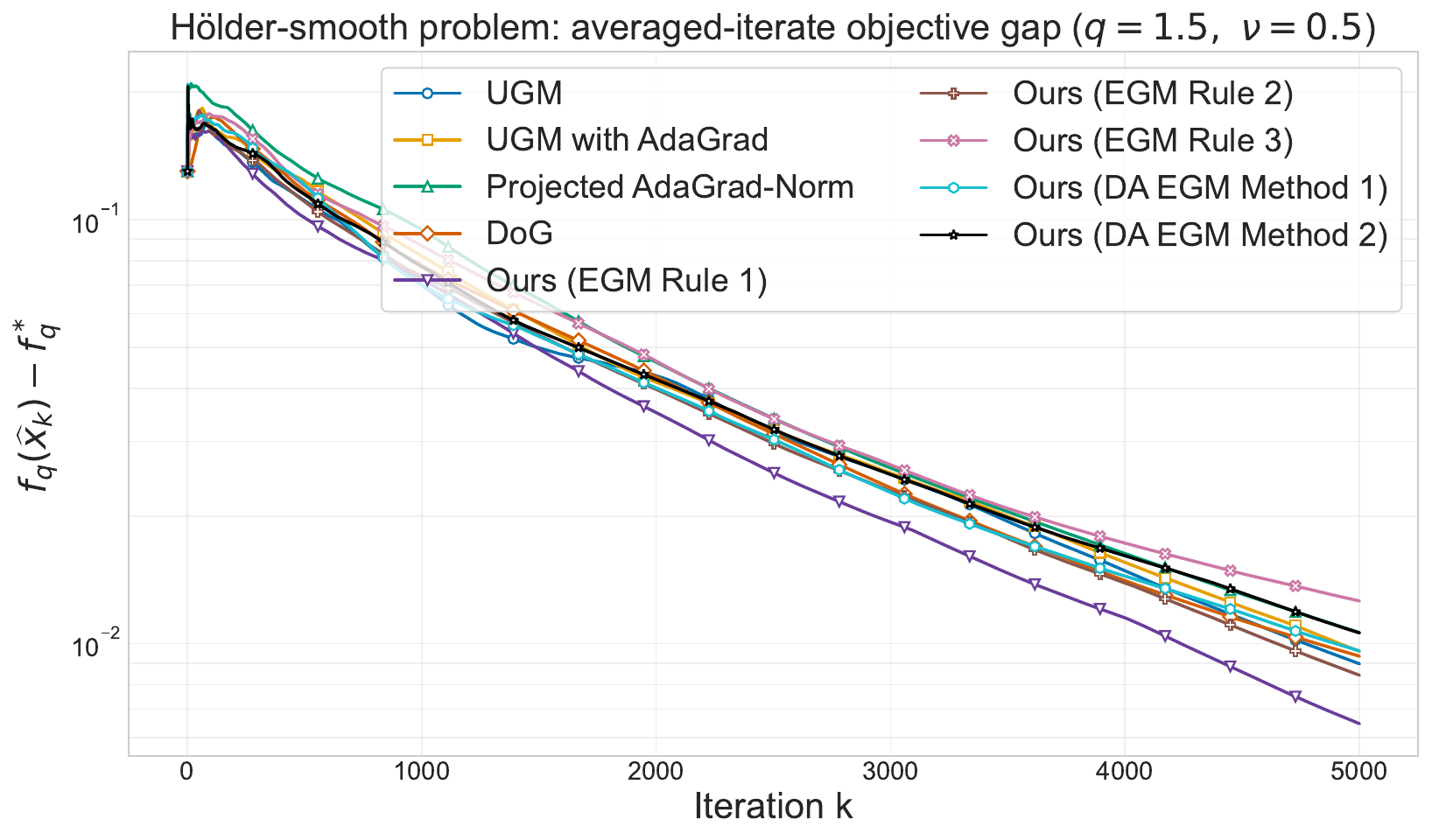}
\end{subfigure}
\hfill
\begin{subfigure}{0.32\textwidth}
    \centering
    \includegraphics[width=\linewidth]{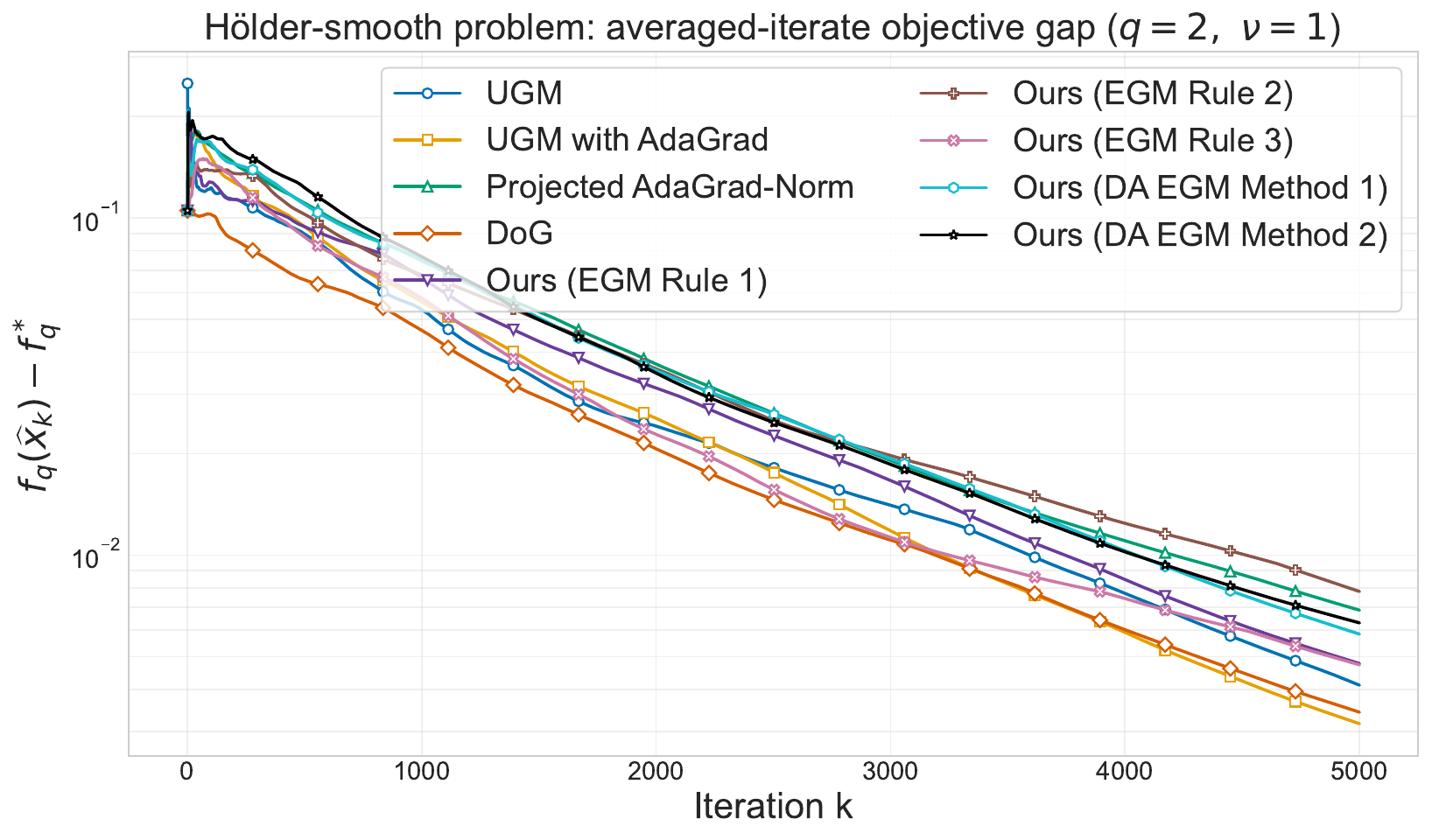}
\end{subfigure}
\caption{Results for the stochastic H\"older-smooth polyhedral feasibility problem. Averaged-iterate objective function gaps are presented for $q=1$, $q=1.5$, and $q=2$, respectively.}
\label{fig_holder_avgitr}
\end{figure*}

\paragraph{H\"older-smooth polyhedral feasibility problem.}
We consider the constrained optimization problem
\begin{align}
    f_q^\star
    :=
    \min_{\|x\|_2 \leq R} f_q(x) \quad \text{where} \quad 
    f_q(x)
    :=
    \frac{1}{n}
    \sum_{i=1}^{n}
    \left[\langle a_i,x\rangle-b_i\right]_+^q, \label{eq:holder_experiment}
\end{align}
where $a_i\in\mathbb{R}^d$, $b_i\in\mathbb{R}$, and
$q\in[1,2]$. This objective belongs to the H\"older-smooth
function class with H\"older exponent $\nu=q-1\in[0,1]$. Thus, $q=1$ corresponds to the non-smooth endpoint, $q=2$ corresponds to the Lipschitz-gradient case, and values $q\in(1,2)$ describe intermediate levels of smoothness. Problem~\eqref{eq:holder_experiment} can also be viewed as a penalty formulation of the polyhedral feasibility problem $\langle a_i,x\rangle\leq b_i$ for all $i=1,\ldots,n$. In particular, $f_q(x)=0$ if and only if $x$ satisfies all inequalities. The problem data are generated following the procedure in \citep[Appendix~E]{moshtaghifar2026dada}. The construction ensures that $x^*$
lies in the Euclidean ball and satisfies every inequality constraint.
Consequently, $f_q(x^*)=f_q^*=0$ for each tested value of $q$. The algorithms are run stochastically with batch size $1$: every oracle call samples an index out of $\{1, \ldots, n\}$ uniformly at random with replacement and evaluates its (sub)gradient. Figure~\ref{fig_holder_avgitr} shows that the methods exhibit comparable empirical performance. For $q=1$ and $q=1.5$, our DA EGM Method~1 and EGM Rule~1 perform slightly better, respectively. Whereas for $q = 2$, UGM with AdaGrad and DoG perform slightly better.

\begin{figure*}[t]
\centering

\begin{subfigure}{0.48\textwidth}
    \centering
    \includegraphics[width=\linewidth]{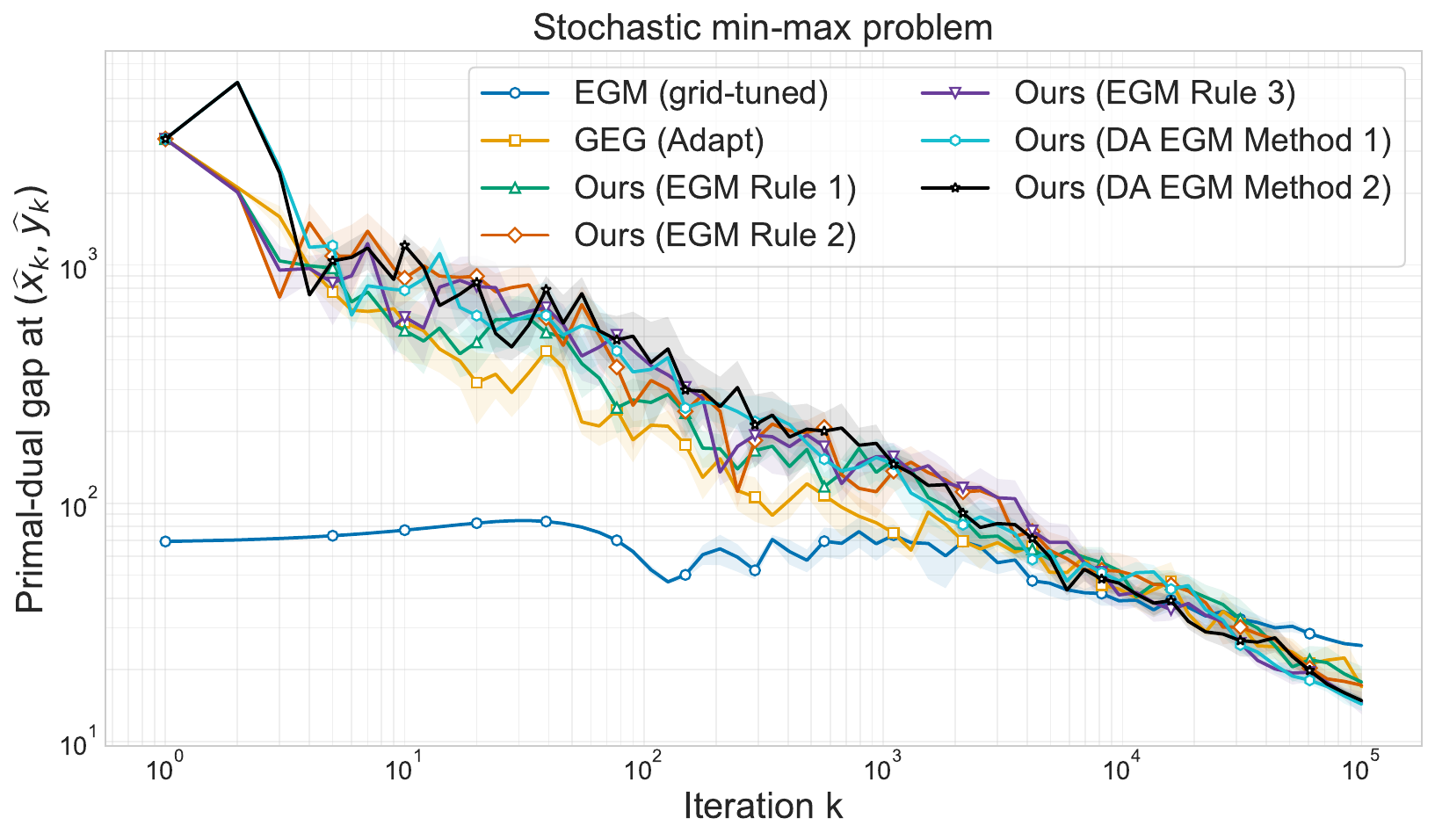}
\end{subfigure}
\hfill
\begin{subfigure}{0.48\textwidth}
    \centering
    \includegraphics[width=\linewidth]{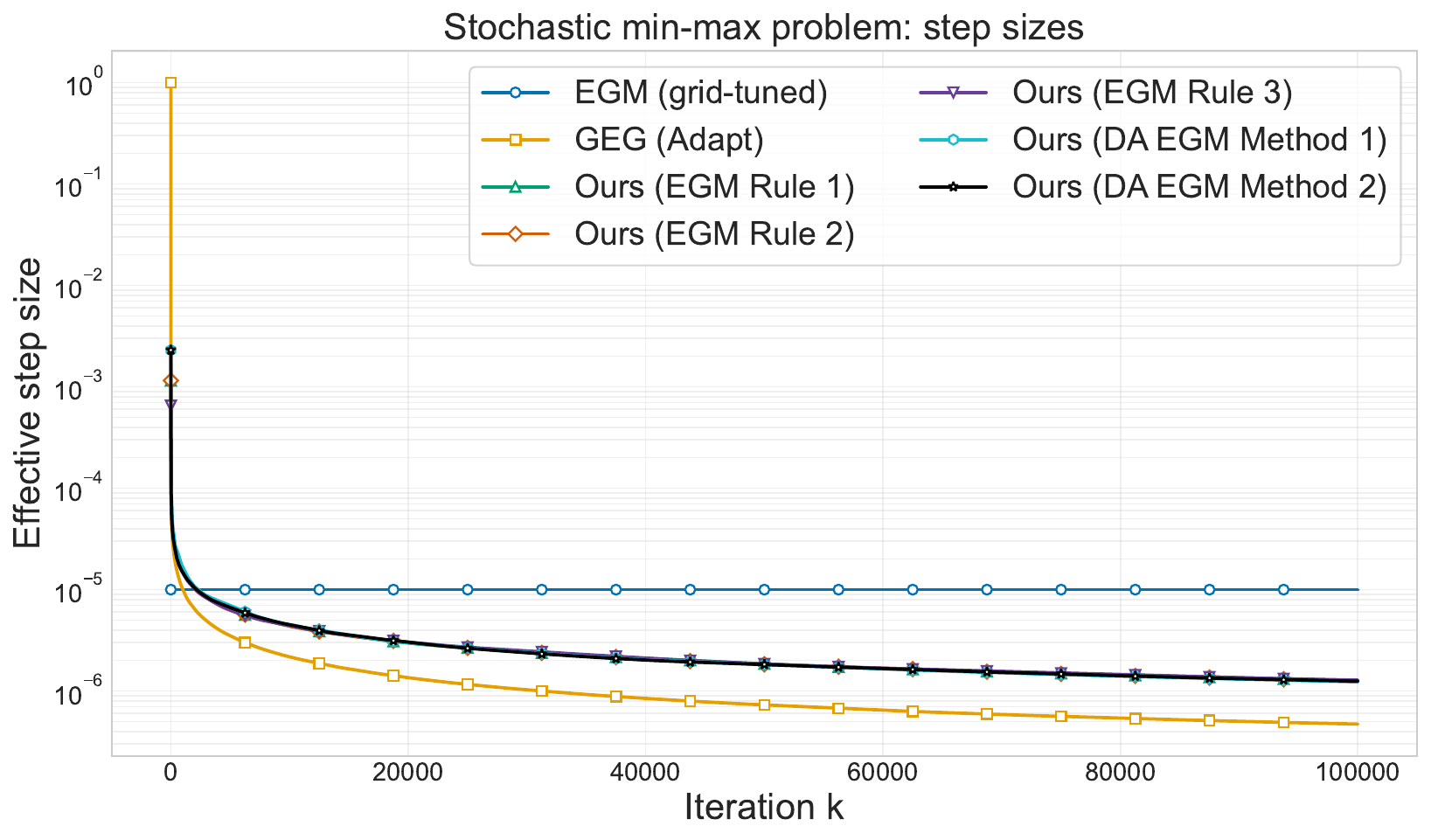}
\end{subfigure}

\caption{Results for the min-max dual formulation of constrained $\ell_2$-regression problem.}
\label{fig_minmax}
\end{figure*}



\paragraph{Min-max dual formulation of a constrained $\ell_2$-regression problem.}

We consider the problem
\begin{align*}
    \min_{\|x\|_2\leq 1}\|Ax-b\|_2
    =
    \min_{\|x\|_2\leq 1}
    \max_{\|y\|_2\leq 1}
    \langle Ax-b,y\rangle .
\end{align*}
The dual reformulation is a min-max problem with the associated VI mapping given by
\begin{align*}
    F([x,y]) = \begin{bmatrix}
        A^\top y \\
        b-Ax
    \end{bmatrix}.
\end{align*}
If the dimension of matrix $A$ is large, then computing the matrix-vector operation while computing $F([x,y])$ can be costly. We employ the randomized matrix-vector oracle described in~\citep[Section~3.6.6]{nesterov2013firstorder} to construct the stochastic mapping oracle $\widehat F([x,y], \omega)$, where the randomness arises from sampling a column of the matrix. Further details are provided in Appendix~\ref{app_experiments}. 

Figure~\ref{fig_minmax} reports the primal-dual gap at the averaged iterate and the step sizes generated by the various methods. EGM (grid-tuned) uses a step size selected via grid search and therefore performs well during the initial iterations. However, later on, the performance of all methods becomes comparable. DA-EGM Methods~1 and~2, as well as EGM with Rule~3, perform slightly better in the final iterations.


\paragraph{Bilinear saddle-point problem.}
We consider the following nonsmooth saddle-point problem:
\begin{equation}
    \min_{x\in[-1,1]}
    \max_{y\in[-1,1]}
    \Phi(x,y)
    \quad \text{where} \quad 
    \Phi(x,y)
    =
    Mxy+\beta |x|-\beta |y|. \nonumber
\end{equation}
%
The corresponding VI mapping of the problem is
\begin{align*}
    F([x,y])
    =
    \begin{bmatrix}
        \partial_x\Phi(x,y)\\
        -\partial_y\Phi(x,y)
    \end{bmatrix}
    =
    \begin{bmatrix}
        My+\beta\operatorname{sign}(x)\\
        -Mx+\beta\operatorname{sign}(y)
    \end{bmatrix}.
\end{align*}
The mapping consists of a combination of a smooth and a non-smooth component and it satisfies $(\delta,L)$-continuity (cf. Assumption~\ref{asum_map_growth}) with parameter
\begin{align*}
    L_{F}(\delta) = M+\frac{2\beta^2}{\delta} ,
\end{align*}
which is derived in detail towards the end of Appendix~\ref{app_experiments}. For our experiments, we consider stochastic mapping denoted as
\begin{align*}
    \widehat F([x,y],\omega) = F([x,y]) + \omega \quad \text{where} \quad \omega = \begin{bmatrix}
        \omega_1\\
        \omega_2
    \end{bmatrix}
    \sim
    \mathcal N(0,\sigma^2I_2) .
\end{align*}
Hence, Assumption~\ref{asum_stochastic} is also satisfied. We report the residual function values for our experiment. If we substitute the value of $L_{F}(\delta)$ in Theorems~\ref{thm_universal} or \ref{thm_distada} and then optimize the right hand side with respect to $\delta$ (similar to Corollaries~\ref{cor_holder} or \ref{cor_adadist_holder}), we will obtain the universal rate
\begin{align*}
    \mathbb{E} \left[\epsilon_T\right] \leq O \left(\frac{M}{T} + \frac{\beta+\sigma}{\sqrt{T}} \right) 
\end{align*}
up to multiplicative factors depending on absolute constants, the diameter of the feasible set, and the initial diameter estimate used by our algorithms. For Theorem~\ref{thm_distada}, the left-hand side of the preceding relation is $\mathbb{E} \left[\epsilon_{\tau(T)}\right]$. To the best of our knowledge, no existing method in the stochastic VI literature achieves such universal rates without requiring knowledge of the problem parameters.


\begin{figure*}[t]
\centering

\begin{subfigure}{0.48\textwidth}
    \centering
    \includegraphics[width=\linewidth]{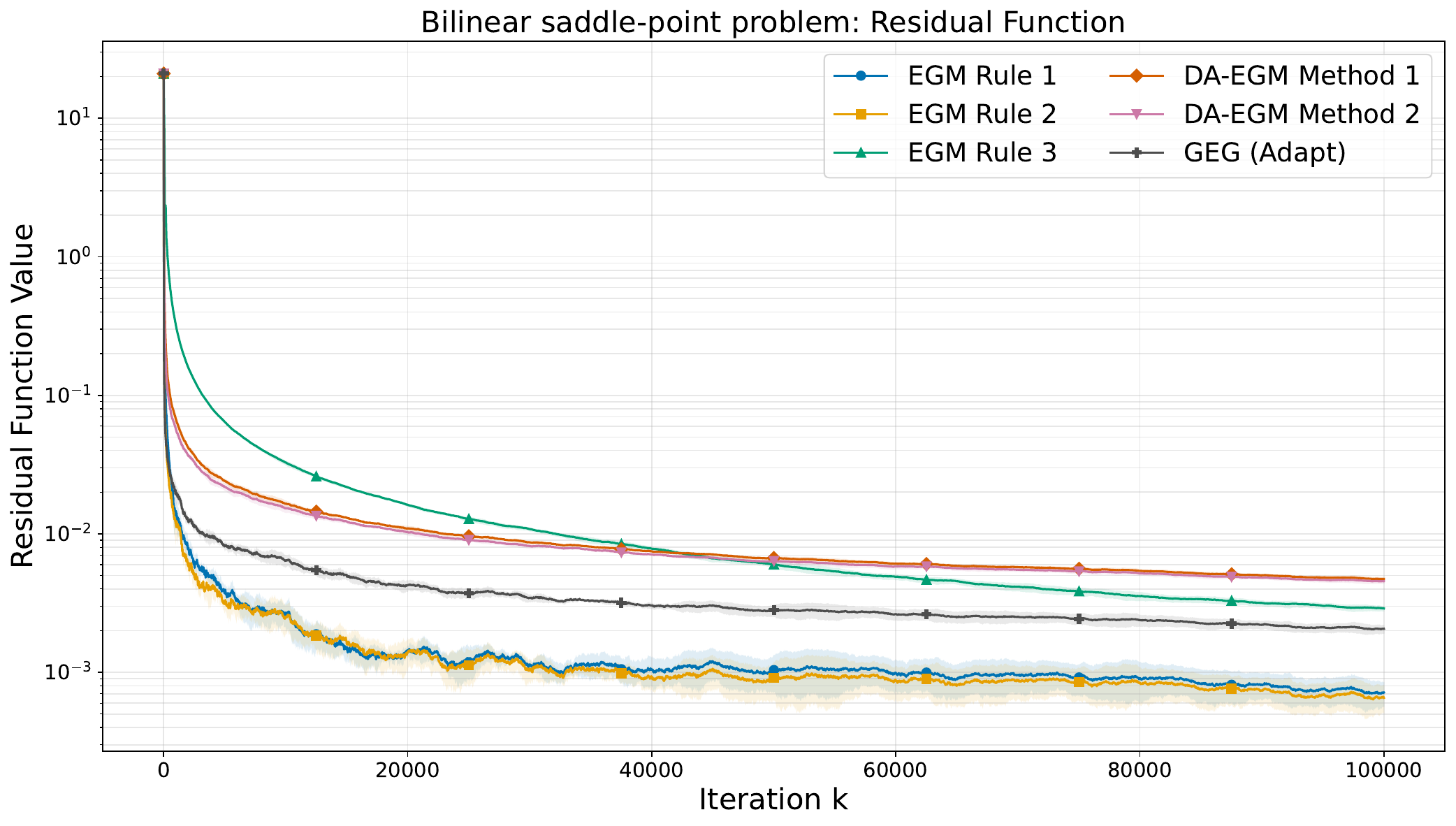}
\end{subfigure}
\hfill
\begin{subfigure}{0.48\textwidth}
    \centering
    \includegraphics[width=\linewidth]{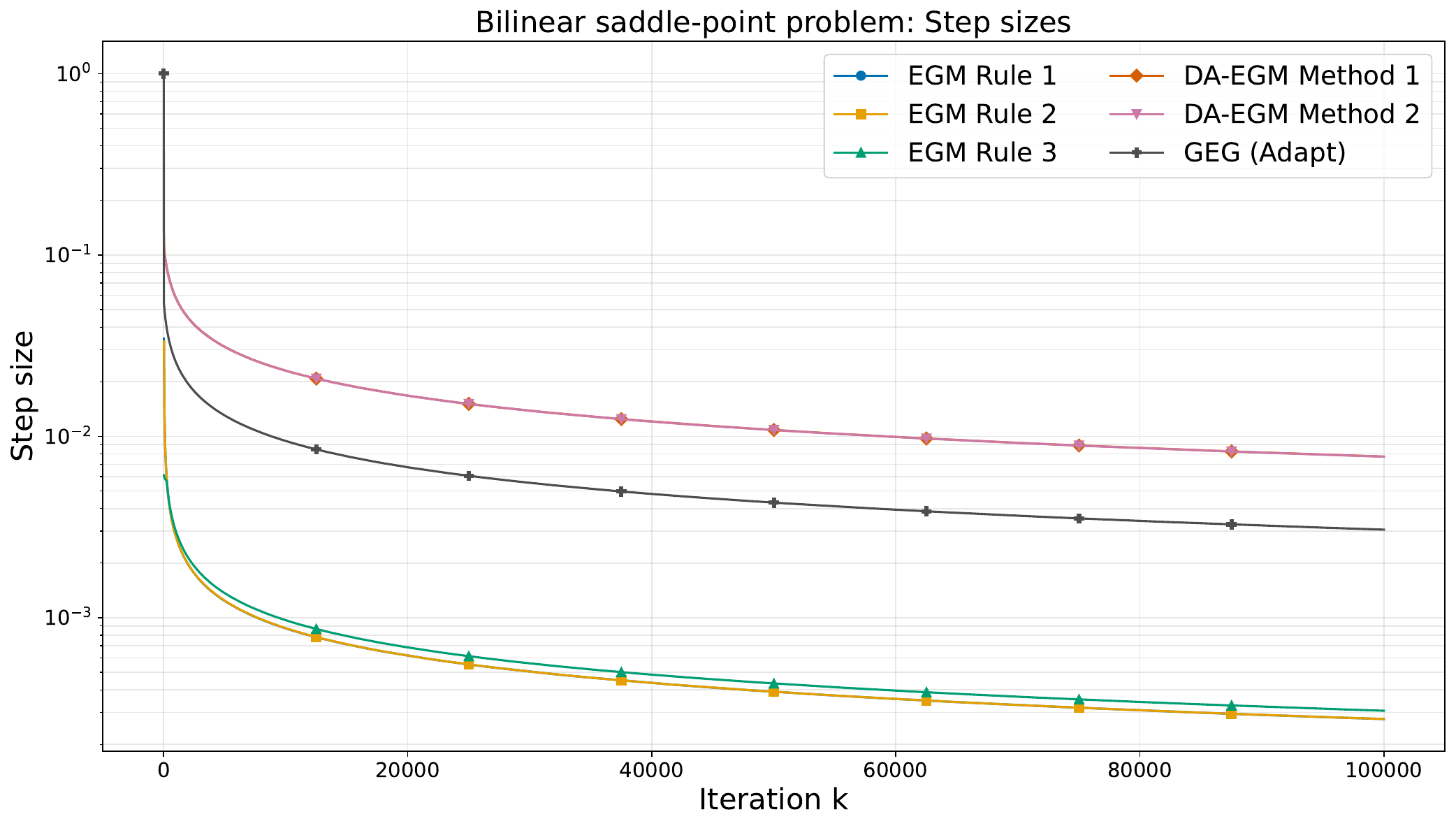}
\end{subfigure}

\caption{Results for the bilinear saddle-point problem.}
\label{fig_saddle_bilinear}
\end{figure*}


Figure~\ref{fig_saddle_bilinear} shows the performance of our EGM and DA EGM step-size update strategies, along with that of GEG. We compare the methods in terms of the residual function value and also report their step sizes over the iterations. The performance of all the methods is comparable, with EGM Rules~1 and~2 performing slightly better. Details of the initialization and residual function calculation are provided in Appendix~\ref{app_experiments}.



\section{Conclusion}

We developed universal extra-gradient methods for stochastic VIs with approximately continuous mappings, a property we refer to as $(\delta,L)$-continuity. We established calculus rules for $(\delta,L)$-continuous mappings and showed that H\"older-continuous mappings and their linear combinations belong to this class. We derived several line-search-free step-size rules in the form of balance equations that adapt to both the continuity of the mapping and the stochastic-oracle variance without requiring prior knowledge of either quantity. Our analysis establishes optimal convergence guarantees through a unified residual function encompassing optimization, saddle-point, and monotone VI problems. For H\"older-continuous mappings, the resulting rates recover the optimal dimension-independent complexity. We also presented distance-adaptive variants that reduce sensitivity to the initial diameter estimate. These variants preserve the same optimal convergence rates while improving some of the constants relative to their non-distance-adaptive counterparts, whose bounds depend more strongly on the initial estimate $\widehat D$. Numerical experiments on optimization and min-max problems demonstrate the practical effectiveness of the proposed methods across different continuity and stochastic noise regimes. Future work may extend these methods to non-Euclidean geometries and to unbounded constraint sets.




\subsection*{AI use statement}


In this work, we used generative AI tools to assist with implementing the experimental methods, including writing, debugging, and editing portions of the experimental code. We also used generative AI tools to correct grammar and improve the readability of author-written text.

We did not use generative AI tools to develop the theoretical or conceptual framework, formulate mathematical claims, assist with proofs, propose or refine hypotheses, design the research methodology or experiments, generate synthetic datasets, analyze or interpret experimental results, or draw the scientific conclusions of this work. The remaining tasks requiring disclosure were not applicable to this work.

All AI-assisted code was manually reviewed and tested by the authors before being used to produce the reported results. All AI-assisted language edits were reviewed to ensure that they preserved the authors' intended technical meaning and claims. We take responsibility for the final content of this work, including all text, claims, results, and artifacts produced with the aid of generative AI.

\bibliography{iclr2027_conference}
\bibliographystyle{iclr2027_conference}


\appendix

\section{Comparison with \citep{klimza2024universal}}\label{app_klimza}

\citet[Algorithm~1]{klimza2024universal} consider H\"older-continuous $F$ and the balance equation
\begin{align*}
    \frac{H_{k+1} - H_k}{2}  D^2 = \left[\langle \widehat F(x_{k+1}, \omega_{k+1}), x_{k+1} - v_{k+1} \rangle - \frac{H_{k+1}}{2} \|v_k - v_{k+1}\|^2 \right]_+ ,
\end{align*}
which is equivalent to the explicit update
\begin{align*}
    H_{k+1} = H_k + \frac{2 \left[ \langle \widehat F(x_{k+1}, \omega_{k+1}), x_{k+1} - v_{k+1} \rangle - \frac{H_{k}}{2} \|v_k - v_{k+1}\|^2 \right]_+}{D^2 + \|v_k - v_{k+1}\|^2} .
\end{align*}
We point out three steps of the analysis that we were unable to verify.

\emph{(i) Lower bound in \citet[Equation~33]{klimza2024universal}.} The lower bound $D^2 + \tfrac{1}{2} \|v_k - v_{k+1}\|^2 - \tfrac{1}{2} \|v_k + v_{k+1} - 2x_{k+1}\|^2 \geq \tfrac{D^2}{2}$ is claimed using compactness. However, this does not follow from boundedness alone. For instance, if $v_{k+1} = v_k$, the left-hand side equals $D^2 - 2\|x_{k+1} - v_k\|^2$, which is smaller than $D^2/2$ whenever $\|x_{k+1} - v_k\| > D/2$.

\emph{(ii) Order of maximum and expectation in Equations (35)--(37).} The analysis exchanges $\max_{x \in X}$ and $\mathbb{E}$. Since $\mathbb{E}[\max_{x} \,\cdot\,] \geq \max_{x} \mathbb{E}[\,\cdot\,]$, the resulting bound controls only the weaker quantity and not the expected gap. We avoid this by using an auxiliary sequence $(h_k)_{k\geq 0}$ (Lemma~\ref{lem_stoc_error}),
which allows us to take $\max_{x \in X}$ first and then the total expectation.

\emph{(iii) Corner cases.} As noted in the updated arXiv version of \citet{klimza2024universal}, some steps of the proof do not apply when the H\"older exponent is $\nu = 0$ or $\nu = 1$. The guarantees therefore do not cover the bounded-variation and Lipschitz regimes, both of which our analysis includes.

\section{Auxiliary Results}

\begin{lemma}\citep[Lemma A.6]{rodomanov2024universality} \label{lem_sum_H_k}
    Let $(H_k)_{k=0}^\infty$ be a non-negative non-decreasing sequence and let $\bar{H} \geq 0$ be a constant. Then, for any $T \geq 1$, we have
    \begin{align*}
        \sum_{k=0}^{T-1} \left[\min\{\bar H, H_{k+1}\} - H_k \right]_+ = \left[\min\{\bar H, H_{T}\} - H_0 \right]_+ .
    \end{align*}
\end{lemma}


\begin{lemma}\label{lem_aux_seq}
Let $(\gamma_k)_{k\geq 0}$ be a nonnegative sequence such that
\begin{align}
    \gamma_0\leq V \quad \text{and} \quad \gamma_{k+1}-\gamma_k
    \leq
    a\left[V-\gamma_{k+1}\right]_+ + e_k
    \qquad \forall k\geq 0 ,
    \label{eq:gamma-recursion}
\end{align}
where $V, a > 0$ and $e_k \geq 0$ for all $k \geq 0$. Then, for any $k\geq 1$,
\begin{equation}
    \gamma_k
    \leq
    V+\sum_{i=0}^{k-1}e_i. \label{eq:gamma-bound}
\end{equation}
\end{lemma}

\begin{proof}
We prove relation~\eqref{eq:gamma-bound} by induction. For $k = 0$, it is indeed satisfied by our assumption $\gamma_0 \leq V$. Now suppose relation~\eqref{eq:gamma-bound} has been proved for some index $k \geq 0$, and let us prove it for the next index $k' = k + 1$. We can assume that $\gamma_{k + 1} > V$ since otherwise the induction step is trivial. Then, combining relations~\eqref{eq:gamma-recursion} and \eqref{eq:gamma-bound}, we obtain
\begin{align*}
    \gamma_{k+1} \leq \gamma_k+e_k \leq
    V+\sum_{i=0}^{k-1}e_i+e_k =V+\sum_{i=0}^{k}e_i.
\end{align*}
Thus, the induction is now complete.
\end{proof}


\begin{lemma}\citep[Lemma 30]{liu2023stochastic}, \citep[Lemma 2]{moshtaghifar2026dada}, \citep[Lemma B.3]{chakraborty2026randomized} \label{lem_adadist_bound}
    Let the sequence $(D_t)_{t \geq 0}$ be positive and non-decreasing. Then for any $T \geq 1$, we obtain
    \begin{align*}
        \min_{1 \leq t \leq T} \frac{D_t}{\sum_{k=0}^{t-1} D_k} \leq \frac{\left( \frac{D_T}{D_0} \right)^{\frac{1}{T}} \log \frac{eD_T}{D_0}}{T} .
    \end{align*}
\end{lemma}


\begin{lemma}~\citep[Theorem 4.4.4]{durrett2019probability}\label{lem_doob}
    If the sequence $(B_t)_{t = 1}^T$ is a martingale, then by Doob's $L^2$ inequality we obtain the following for any $T \geq 1$:
    \begin{align*}
        \EXP{\max_{1 \leq t \leq T} |B_t|^2} \leq 4 \EXP{|B_T|^2} .
    \end{align*}
\end{lemma}


\section{Proofs for Section~\ref{sec_problem}}

\subsection{Proof of Lemma~\ref{lem_holder_to_delta_L}}

\begin{proof}
    Let us fix arbitrary $x, y \in X$ and denote $r = \| x - y \|$. Suppose that $\nu < 1$. Then, for any $L > 0$, we can estimate
    \begin{align*}
        \|F(x) - F(y)\|_*^2 &\leq H_\nu^2 r^{2\nu} = L^2 r^2 + (H_\nu^2 r^{2\nu} - L^2 r^2) \leq L^2 r^2 + \Delta(L) ,
    \end{align*}
    where
    \begin{align*}
        \Delta(L) &= \max_{r \geq 0} \{H_\nu^2 r^{2\nu} - L^2 r^2\} = \max_{\tau \geq 0} \left\{ H_\nu^2 \tau^{\nu} - L^2 \tau \right\} = (1 - \nu) \left(\frac{\nu}{L^2} \right)^{\frac{\nu}{1 - \nu}} H_{\nu}^{\frac{2}{1 - \nu}}, 
    \end{align*}
    where the final identity follows from obtaining the maximum value of the concave function. Let us fix an arbitrary $\delta > 0$ and choose $L$ such that $\Delta(L) = 4L \delta$. This is $L = L(\delta)$ as defined in the statement of Lemma~\ref{lem_holder_to_delta_L}. The claim now follows.

    If $\nu = 1$, then $\| F(x) - F(y) \|_*^2 \leq H_1^2 r^2$ and $L(\delta) = H_1$, and the claim is again trivial.
\end{proof}


\subsection{Proof of the Claim in Remark~\ref{rem_vi_opt_comp}}

\begin{lemma}\label{lem_L_delta}
    For any $\nu \in [0,1]$, $H_\nu > 0$, and $\delta > 0$, let $L(\delta)$ and $\tilde{L}(\delta)$ be defined as
    \begin{align*}
        L(\delta) = \left(\frac{1-\nu}{4\delta}\right)^{\frac{1-\nu}{1+\nu}}\nu^{\frac{\nu}{1+\nu}} H_\nu^{\frac{2}{1+\nu}} , \quad \tilde{L}(\delta) = \left(\frac{1 - \nu}{2 (1 + \nu) \delta}\right)^{\frac{1 - \nu}{1 + \nu}} H_{\nu}^{\frac{2}{1 + \nu}} 
    \end{align*}
    with the convention $0^0 = 1$. Then, 
    \begin{align*}
        L(\delta) \leq \tilde{L}(\delta) .
    \end{align*}
\end{lemma}
\begin{proof}
    Using the definitions of $L(\delta)$ and $\tilde{L}(\delta)$, we have
    \begin{align*}
        \frac{L(\delta)}{\tilde{L}(\delta)} = \left( \frac{1+\nu}{2} \right)^{\frac{1-\nu}{1+\nu}} \nu^{\frac{\nu}{1+\nu}} .
    \end{align*}
    The maximum value on the right-hand side is attained at $\nu = 1$. In order to show that, we apply log on both sides (since log is a monotone function) to obtain
    \begin{align}
        (1+ \nu) \ln\left(\frac{L(\delta)}{\tilde{L}(\delta)} \right) = (1 - \nu)\ln\left(\frac{1+\nu}{2}\right) + \nu \ln \nu := \phi(\nu) . \label{eqn_phi}
    \end{align}
    Applying double derivative on $\phi(\nu)$ with respect to $\nu$, we obtain
    \begin{align*}
        \nabla_\nu^2 \phi(\nu) = \frac{1 - \nu}{\nu (1 + \nu)^2} \geq 0 \quad \text{for $\nu \in (0,1]$}.
    \end{align*}
    Hence, $\phi(\nu)$ is a convex function with respect to $\nu$. Since the function $\phi(\nu)$ is convex, the maximum value will be attained at one of the end points. We see $\lim_{\nu \rightarrow 0+} \phi(\nu) = \ln(0.5) + \lim_{\nu \rightarrow 0+} \nu \ln \nu = \ln(0.5) + \lim_{\nu \rightarrow 0+} \frac{1/\nu}{-1/\nu^2} = \ln(0.5) < 0$ and $\phi(1) = 0$. Hence, we obtain $\phi(\nu) \leq 0$ for any $\nu \in [0,1]$, which when used in relation~\eqref{eqn_phi} yields $L(\delta) \leq \tilde{L}(\delta)$.
\end{proof}


\subsection{Calculus properties of $(\delta, L)$-continuity for VI}

\begin{lemma}\label{lem_calculus}
Suppose that, for every $\delta>0$, the mapping $F$ admits a
constant $L_F(\delta)>0$ such that
\[
    \|F(x)-F(y)\|_*^2
    \leq
    L_F^2(\delta) \|x-y\|^2
    +4L_F(\delta)\delta.
\]
Then the following properties hold:\\
    (i) (Scalar multiplication)
    If $c>0$, then an admissible continuity modulus for $cF$ is
    \[
        L_{cF}(\delta)
        :=
        cL_F\left(\frac{\delta}{c}\right).
    \]

    (ii) (Composition with an affine mapping)
    Let $G(x):=F(Ax+b)$, where $A \in \R^{m \times n}$ and $b \in \R^m$. Then an admissible continuity
    modulus for $G$ is
    \[
        L_G(\delta)
        :=
        \|A\|L_F\bigl(\|A\|\delta\bigr).
    \]

    (iii) (Summation) Suppose that $F_1$ and $F_2$ admit continuity moduli
    $L_1(\delta)$ and $L_2(\delta)$, respectively. Then, for every
    $\delta_1,\delta_2>0$ satisfying $\delta_1+\delta_2=\delta$, the mapping $F = F_1+F_2$ has the continuity modulus
    \begin{align*}
        L_{F}(\delta) = L_{F_1} (\delta_1) + L_{F_2} (\delta_2) .
    \end{align*}
\end{lemma}

\begin{proof}
\emph{Part (i):}
Let us fix $\delta>0$ and apply the continuity condition for $F$ with
the parameter $\delta/c$. Then
\begin{align*}
    \|cF(x)-cF(y)\|_*^2
    &=
    c^2\|F(x)-F(y)\|_*^2
    \\
    &\leq
    c^2L_F^2\left(\frac{\delta}{c}\right)
    \|x-y\|^2
    +
    4c^2L_F\left(\frac{\delta}{c}\right)
    \frac{\delta}{c}
    \\
    &=
    L_{cF}^2(\delta)\|x-y\|^2
    +4L_{cF}(\delta)\delta,
\end{align*}
where
\[
    L_{cF}(\delta)
    =
    cL_F\left(\frac{\delta}{c}\right).
\]

\emph{Part (ii):}
Let us fix $\delta>0$ and apply the continuity condition for $F$ with
the parameter $\|A\|\delta$. Then
\begin{align*}
    \|G(x)-G(y)\|_*^2
    &=
    \|F(Ax+b)-F(Ay+b)\|_*^2
    \\
    &\leq
    L_F^2(\|A\|\delta) \|A(x-y)\|^2
    +
    4L_F(\|A\|\delta)\|A\|\delta
    \\
    &\leq
    \bigl(\|A\|L_F(\|A\|\delta)\bigr)^2
    \|x-y\|^2 + 4\bigl(\|A\|L_F(\|A\|\delta)\bigr)\delta.
\end{align*}
Hence,
\[
    L_G(\delta)
    =
    \|A\|L_F(\|A\|\delta).
\]

\emph{Part (iii):} For simplicity, let us define
\begin{align*}
    \Delta_1 = F_1(x) - F_1(y) \quad \text{and} \quad \Delta_2 = F_2(x) - F_2(y) .
\end{align*}
Then we see that
\begin{align}
    &\|F(x) - F(y)\|_*^2 = \|\Delta_1 + \Delta_2\|_*^2 \leq (1+\tau)\|\Delta_1\|_*^2 + \left(1+\frac{1}{\tau}\right)
    \|\Delta_2\|_*^2 \nonumber\\
    & \leq \left( (1+\tau) L_{F_1}^2(\delta_1) + \left(1+ \frac{1}{\tau} \right) L_{F_2}^2 (\delta_2) \right) \|x-y\|^2 + 4 \left( (1+\tau) L_{F_1}(\delta_1) \delta_1 + \left(1+ \frac{1}{\tau} \right) L_{F_2} (\delta_2) \delta_2 \right) , \label{eq-est-wtau}
\end{align}
where the last two relations are obtained from Young’s inequality for any constant $\tau > 0$ and then the $(\delta, L)$-continuity of the mappings. The optimal value $\tau^*>0$ that minimizes the coefficient of the first term in~\eqref{eq-est-wtau} is given by
\begin{align*}
    \tau^* = \frac{L_{F_2}(\delta_2)}{L_{F_1}(\delta_1)} .
\end{align*}
Substituting the value of $\tau^*$ back in equation~\eqref{eq-est-wtau} yields
\begin{align*}
    \|F(x) - F(y)\|_*^2 \leq (L_{F_1}(\delta_1) + L_{F_2}(\delta_2))^2 \|x-y\|^2 + 4 (L_{F_1}(\delta_1) + L_{F_2}(\delta_2)) (\delta_1 + \delta_2) ,
\end{align*}
which implies that $F$ is $(\delta_1 + \delta_2, L_{F_1}(\delta_1) + L_{F_2}(\delta_2))$ continuous.
\end{proof}



\begin{lemma}\label{lem_cont_multiple}
    Let $F_i$ is H\"older continuous with parameters $(\nu_i, H_{\nu_i})$ for all $i \in \abs{\mathcal{I}_1}$ and let $\bar F_j$ is $(\delta_j, L_j)$-continuous for all $j \in \abs{\mathcal{I}_2}$, where $\abs{\mathcal{I}_1} \geq 0$ and $\abs{\mathcal{I}_2} \geq 0$, then $F = \sum_{i \in \abs{\mathcal{I}_1}} F_i + \sum_{j \in \abs{\mathcal{I}_2}} \bar F_j$ is $(\sum_{i \in \abs{\mathcal{I}_1}} \delta_i + \sum_{j \in \abs{\mathcal{I}_2}} \delta_j, \sum_{i \in \abs{\mathcal{I}_1}} L_i + \sum_{j \in \abs{\mathcal{I}_2}} L_j)$-continuous, where $\delta_i > 0$ is an arbitrary constant and
    \begin{align*}
        L_i \equiv L_i(\delta_i) = \left(\frac{1-\nu_i}{4\delta_i}\right)^{\frac{1-\nu_i}{1+\nu_i}}\nu_i^{\frac{\nu_i}{1+\nu_i}} H_{\nu_i}^{\frac{2}{1+\nu_i}} \quad \forall i \in \abs{\mathcal{I}_1} ,
    \end{align*}
    with the convention $0^0=1$.
\end{lemma}

\begin{proof}
    If $F_i$ is H\"older continuous with parameters $(\nu_i, H_{\nu_i})$ for any $i \in \abs{\mathcal{I}_1}$, then by Lemma~\ref{lem_holder_to_delta_L}, we can fix any arbitrary $\delta_i > 0$ and obtain the expression for $L_i(\delta_i)$ as provided in the lemma. Now, by applying Lemma~\ref{lem_calculus}(iii) recursively, we obtain the desired result.
\end{proof}


\section{Proofs for Section~\ref{sec_algo}}

For the algorithmic analysis, we encode the stochastic history generated by EGM (cf. relation~\eqref{Kor_algo}) up to iteration $k$ in the sigma-algebra $\cF_k$, defined as follows:
\begin{align}
    \cF_k = \{v_0\} \cup \{\omega_t \mid 0 \leq t \leq k \} \qquad \forall k \geq 0. \label{sigma_field}
\end{align}
For notational simplicity, we define the quantities
\begin{align}
    & Y_k = \widehat F(v_k, \omega_k) - F(v_k) \quad \text{and} \quad Z_k = \widehat F(x_k, \omega_k) - F(x_k)\qquad \forall k\ge0 . \label{stoc_errors}
\end{align}

\subsection{Generic Lemma for EGM}\label{app_gen_egm}

\begin{lemma}\label{lem_generic_reln}

Under Assumptions~\ref{asum_set} and \ref{asum_stochastic}, and for any $k \geq 0$, EGM algorithm satisfies
\begin{align*}
    \mathbb{E}[\epsilon_T] &\leq \frac{1}{T}\sum_{k=0}^{T-1} \mathbb{E}[\Gamma_k] +
    \frac{D \sigma}{\sqrt{T}} \quad \text{where} \quad \Gamma_k \equiv \Delta_k(H_{k+1}) + \frac{3}{2} (H_{k+1} - H_{k}) D^2,
\end{align*}
$\Delta_k(H) \equiv \xi_k - \frac{H}{2} r_k^2$ with $r_k^2 = \|v_k - x_{k+1}\|^2 + \|v_{k+1} - x_{k+1}\|^2$ and
\begin{align}
    \xi_k &= \langle \widehat F(x_{k+1}, \omega_{k+1}), x_{k+1} - v_{k+1} \rangle - \frac{H_k}{2} \|v_k - v_{k+1}\|^2 + \frac{H_k}{2} r_k^2 \nonumber\\
    &\leq \langle \widehat F(x_{k + 1}, \omega_{k+1}) - \widehat F(v_k, \omega_k), x_{k + 1} - v_{k + 1} \rangle \equiv \xi_k'. \label{xi_xi_prime_upbd}
\end{align}
Moreover,
\begin{align}
    \Delta_k(H) &\leq \xi_k' - \frac{H}{2} r_k^2 \leq \frac{1}{2 H} \| \widehat F(x_{k + 1}, \omega_{k+1}) - \widehat F(v_k, \omega_k) \|_*^2 - \frac{H}{2} \| x_{k + 1} - v_k \|^2 \equiv \bar{\Delta}_k(H). \label{eq3_lem_generic_reln}
\end{align}
\end{lemma}

\begin{proof}
    Let $x \in X$ and $k \geq 0$ be arbitrary. From the definition of $v_{k + 1}$ and using $H_k$-strong convexity of the update, it follows that
    \begin{align*}
        &\langle \widehat F(x_{k+1}, \omega_{k+1}),\, v_{k+1} - x_{k+1} \rangle
        + \frac{H_{k}}{2}\|v_k - v_{k+1}\|^2 + \frac{H_{k}}{2}\|x - v_{k+1}\|^2 \\
        &\leq \langle \widehat F(x_{k+1}, \omega_{k+1}),\, x - x_{k+1} \rangle 
        + \frac{H_{k}}{2}\|v_k - x\|^2 .
    \end{align*}
    Rearranging terms along with algebraic manipulations, we obtain
    \begin{align}
        \langle \widehat F(x_{k+1}, \omega_{k+1}), x_{k+1} - x \rangle \leq \xi_k &- \frac{H_{k+1}}{2} r_k^2 + \frac{H_{k+1} - H_{k}}{2} [r_k^2 + \|x-v_{k+1}\|^2] \nonumber\\
        &+ \frac{H_{k}}{2}\|x - v_k\|^2 - \frac{H_{k+1}}{2}\|x - v_{k+1}\|^2 , \label{main_eqn_lem_gen}
    \end{align}
    where the quantity $r_k^2 = \|v_k - x_{k+1}\|^2 + \|v_{k+1} - x_{k+1}\|^2$ and
    \begin{align}
        \xi_k = \langle \widehat F(x_{k+1}, \omega_{k+1}), x_{k+1} - v_{k+1} \rangle - \frac{H_k}{2} \|v_k - v_{k+1}\|^2 + \frac{H_k}{2} r_k^2 . \label{xi_eqn}
    \end{align}
    For notational simplicity, we define
    \begin{align*}
        \Delta_k(H) = \xi_k - \frac{H}{2} r_k^2.
    \end{align*}
    Hence, applying the preceding relation and Assumption~\ref{asum_set} to relation~\eqref{main_eqn_lem_gen}, we obtain
    \begin{align}
        \langle \widehat F(x_{k+1}, \omega_{k+1}), x_{k+1} - x \rangle &\leq \Gamma_k + \frac{H_{k}}{2}\|x - v_k\|^2 - \frac{H_{k+1}}{2}\|x - v_{k+1}\|^2 , \label{eq1_reln}
    \end{align}
    where the quantity
    \begin{align}
        \Gamma_k \equiv \Delta_k(H_{k+1}) + \frac{3}{2} (H_{k+1} - H_{k}) D^2 . \label{gammak}
    \end{align}
    The left-hand side of the preceding relation can be written as
    \begin{align}
        \langle \widehat F(x_{k+1}, \omega_{k+1}), x_{k+1} - x \rangle = \langle F(x_{k+1}), x_{k+1} - x \rangle + \langle Z_{k+1}, x_{k+1} - h_k \rangle + \langle Z_{k+1}, h_k - x \rangle , \label{split_eqn}
    \end{align}
    where $Z_{k+1} = \widehat F(x_{k+1}, \omega_{k+1}) - F(x_{k+1})$ (cf.\ relation~\eqref{stoc_errors}) and with $h_0 \in X$, $h_{k+1}$ for any $k \geq 0$ is defined as
    \begin{align*}
        h_{k+1} = \argmin_{x \in X} \left\{ \langle - Z_{k+1}, x - h_k \rangle + \frac{\eta}{2} \|x - h_k\|^2 \right\} ,
    \end{align*}
    where $\eta \geq 0$ is an arbitrary constant. Substituting relation~\eqref{split_eqn} to relation~\eqref{eq1_reln} we obtain
    \begin{align*}
        \langle F(x_{k+1}), x_{k+1} - x \rangle &\leq \Gamma_k + \langle Z_{k+1}, h_k - x_{k+1} \rangle + \langle Z_{k+1}, x - h_k \rangle \\
        & \qquad + \frac{H_{k}}{2}\|x - v_k\|^2 - \frac{H_{k+1}}{2}\|x - v_{k+1}\|^2 .
    \end{align*}
    Applying Lemma~\ref{lem_stoc_error} on the third term of the right-hand side of the preceding relation yields
    \begin{align}
        \langle F(x_{k+1}), x_{k+1} - x \rangle &\leq \Gamma_k + \langle Z_{k+1}, h_k - x_{k+1} \rangle  + \frac{\|Z_{k+1}\|_*^2}{2 \eta} + \frac{\eta}{2} [\|x - h_k\|^2 - \|x - h_{k+1}\|^2] \nonumber\\
        & \hspace{3cm}+ \frac{H_{k}}{2}\|x - v_k\|^2 - \frac{H_{k+1}}{2}\|x - v_{k+1}\|^2 , \nonumber
    \end{align}
    We sum both sides of the preceding relation from $k = 0$ to $T-1$ and divide both sides by $T$ to obtain
    \begin{align*}
        \frac{1}{T}\sum_{k=0}^{T-1} \langle F(x_{k+1}),x_{k+1}-x\rangle
        &\leq \frac{1}{T}\sum_{k=0}^{T-1} \Gamma_k + \frac{1}{T}\sum_{k=0}^{T-1}
        \langle Z_{k+1},h_k-x_{k+1}\rangle +
        \frac{1}{2\eta T}
        \sum_{k=0}^{T-1}\|Z_{k+1}\|_*^2 \\
        & + \frac{\eta}{2 T} [\|x - h_0\|^2 - \|x - h_T\|^2] +
        \frac{H_0}{2T}\|x-v_0\|^2
        -\frac{H_T}{2T}\|x-v_T\|^2 .
    \end{align*}
    Recall that we choose $H_0 = 0$. Next, we take maximum with respect to $x \in X$ on both sides of the preceding relation and apply Assumption~\ref{asum_set}, yielding
    \begin{align*}
        \epsilon_T &\leq \frac{1}{T}\sum_{k=0}^{T-1} \Gamma_k + \frac{1}{T}\sum_{k=0}^{T-1}
        \langle Z_{k+1},h_k-x_{k+1}\rangle +
        \frac{1}{2\eta T}
        \sum_{k=0}^{T-1}\|Z_{k+1}\|_*^2 + \frac{\eta D^2}{2T} ,
    \end{align*}
    where $\epsilon_T \equiv \epsilon_T((x_t, 1)_{t=1}^T)$ is the residual function given in relation~\eqref{residual_val_eq}. Next, we take total expectation on both sides of the preceding relation to obtain
    \begin{align*}
        \mathbb{E}[\epsilon_T] &\leq \frac{1}{T}\sum_{k=0}^{T-1} \mathbb{E}[\Gamma_k] + \frac{1}{T}\sum_{k=0}^{T-1}
        \mathbb{E} \left[\langle Z_{k+1},h_k-x_{k+1}\rangle \right] +
        \frac{1}{2\eta T}
        \sum_{k=0}^{T-1} \mathbb{E} \left[\|Z_{k+1}\|_*^2 \right] + \frac{\eta D^2}{2T} .
    \end{align*}
    Next, we apply Assumption~\ref{asum_stochastic}. The second quantity on the right-hand side can be simplified as
    \begin{align*}
        \sum_{k=0}^{T-1} \EXP{\langle Z_{k+1}, h_k - x_{k+1} \rangle} = \sum_{k=0}^{T-1} \EXP{\langle \EXP{Z_{k+1} \mid \cF_k}, h_k - x_{k+1} \rangle} = 0 ,
    \end{align*}
    where the sigma algebra $\cF_k$ is defined in relation~\eqref{sigma_field}. Hence, we obtain
    \begin{align*}
        \mathbb{E}[\epsilon_T] &\leq \frac{1}{T}\sum_{k=0}^{T-1} \mathbb{E}[\Gamma_k] +
        \frac{\sigma^2}{2\eta}
        + \frac{\eta D^2}{2T} .
    \end{align*}
    Minimizing the right hand side of the preceding relation with respect to $\eta > 0$ yields the optimal
    \begin{align*}
        \eta^* = \frac{\sigma \sqrt{T}}{D} ,
    \end{align*}
    which when substituted back yields the first relation of Lemma~\ref{lem_generic_reln}.

    Next, we prove that $\xi_k \leq \xi_k'$. We have
    \begin{align}
        \xi_k = &\langle \widehat F(x_{k+1}, \omega_{k+1}) - \widehat F(v_k, \omega_k), x_{k+1} - v_{k+1} \rangle - \frac{H_{k}}{2}\|v_k - v_{k+1}\|^2 \nonumber\\
        &+ \langle \widehat F(v_k, \omega_k), x_{k+1} - v_{k+1} \rangle + \frac{H_k}{2} r_k^2 . \label{xi_split}
    \end{align}
    Since the iterate $x_{k+1}$ is the solution of $H_k$-strongly convex function minimization over $X$, it follows that for any $y \in X$,
    \begin{align*}
        &\langle \widehat F(v_k, \omega_k), x_{k+1} - v_k \rangle + \frac{H_{k}}{2} \|v_k - x_{k+1}\|^2 + \frac{H_{k}}{2}\|y - x_{k+1}\|^2 \\
        &\leq \langle \widehat F(v_k, \omega_k), y - v_k \rangle + \frac{H_{k}}{2}\|y - v_k\|^2 ,
    \end{align*}
    Substituting $y = v_{k+1}$ in the preceding relation and rearranging terms yields the following upper bound
    \begin{align*}
        \la \widehat F(v_k, \omega_k) , x_{k+1} - v_{k+1} \ra \leq \frac{H_{k}}{2}\|v_k - v_{k+1}\|^2 - \frac{H_{k}}{2} r_k^2 ,
    \end{align*}
    which when substituted back in relation~\eqref{xi_split} yields relation~\eqref{xi_xi_prime_upbd}.

    Finally, the upper bound in relation~\eqref{eq3_lem_generic_reln} is obtained using Young's inequality. This concludes the proof.
\end{proof}


\subsection{Lemma to control the stochastic error}

We present a lemma which is a modification of \citep[Lemma 3 and Corollary 2]{juditsky2011solving} and \citep[Lemma 5]{chakraborty2025popov} that controls the stochastic error quantity required for our analysis.

\begin{lemma}\label{lem_stoc_error} 
    Let $h_0 \in X$ and $h_{k+1}$ for all $k \geq 0$ be updated as 
    \begin{align*}
        h_{k+1} = \argmin_{x \in X} \left\{ \langle - Z_{k+1}, x - h_k \rangle + \frac{\eta}{2} \|x - h_k\|^2 \right\} ,
    \end{align*}
    where $\eta \geq 0$ is a constant and $Z_{k+1} := \widehat F(x_{k+1}, \omega_{k+1}) - F(x_{k+1})$ for any $\omega_{k+1} \in \Xi$. Then, for any $k \geq 0$ and any $x \in X$, we have
    \begin{align*}
        \langle Z_{k+1}, x - h_k \rangle &\leq \frac{\|Z_{k+1}\|_*^2}{2 \eta} + \frac{\eta}{2} [\|x - h_k\|^2 - \|x - h_{k+1}\|^2] .
    \end{align*}
\end{lemma}
\begin{proof}
    Using the $H_k$-strong convexity of the function in the update rule for $h_{k+1}$, we obtain for any $k \geq 0$ and any $x \in X$,
    \begin{align*}
        &\la -Z_{k+1} , h_{k+1} - h_k \ra + \frac{\eta}{2} \|h_{k+1} - h_k\|^2 + \frac{\eta}{2} \|x - h_{k+1}\|^2 \\
        &\leq \la -Z_{k+1} , x - h_k \ra + \frac{\eta}{2} \|x - h_k\|^2 .
    \end{align*}
    Simplifying the preceding relation yields
    \begin{align*}
        \langle Z_{k+1}, x - h_k \rangle &\leq \la Z_{k+1}, h_{k+1} - h_k \ra - \frac{\eta}{2} \|h_{k+1} - h_k\|^2 + \frac{\eta}{2} [\|x - h_k\|^2 - \|x - h_{k+1}\|^2] .
    \end{align*}
    Applying Young's inequality to the first quantity on the right-hand side yields the desired result.
\end{proof}


\subsection{Algorithmic updates for various examples}\label{app_subsec_algos}

Here, we derive the explicit updates for $H_{k+1}$ under different balance equation examples. We first present a lemma that facilitates these derivations.
\begin{lemma}\label{lem_balance_eqns} \citep[Lemma E.1]{rodomanov2024universal}
    Let $\Omega > 0$ be an arbitrary constant and let $(\zeta_k)_{k \geq 0}$ be a sequence of non-negative scalars. Let $(H_k)_{k \geq 0}$ be a sequence of non-decreasing scalars and $(\beta_k)_{k \geq 0}$ be a scalar sequence. Then, for all $k \geq 0$, the balance equation
    \begin{align*}
        \frac{H_{k+1} - H_k}{2} \Omega = \left[ \beta_k - \frac{H_{k+1}}{2} \zeta_k \right]_+
    \end{align*}
    is equivalent to the balance equation
    \begin{align*}
        \frac{H_{k+1} - H_k}{2} (\Omega + \zeta_k)= \left[ \beta_k - \frac{H_{k}}{2} \zeta_k \right]_+ .
    \end{align*}
\end{lemma}

We now present the explicit update rules for the sequence ${H_k}$ corresponding to each balance-equation example separately, followed by the complete algorithm.

\paragraph{Rule~1:} By Lemma~\ref{lem_balance_eqns}(i), the balance equation~\eqref{baleqn1} has a unique solution given by
\begin{align}
    H_{k+1} &= H_k + \frac{2 [\Delta_k(H_k)]_+}{R_{k}^2} = H_k + \frac{2 \left[\xi_k - \frac{H_k}{2} r_k^2 \right]_+}{R_{k}^2} \label{bal_up1}\\
    &=  H_k + \frac{2 \left[\langle \widehat F(x_{k+1}, \omega_{k+1}), x_{k+1} - v_{k+1} \rangle - \frac{H_k}{2} \|v_k - v_{k+1}\|^2 \right]_+}{R_{k}^2}, \nonumber
\end{align}
where the quantity $R_k$ is defined as
\begin{align}
    R_{k}^2 = \widehat D^2 + r_k^2 = \widehat D^2 + \| v_k - x_{k + 1} \|^2 + \| v_{k + 1} - x_{k + 1} \|^2. \label{ada_rad_quantity}
\end{align}

\paragraph{Rule~2:} Here, for the balance equation~\eqref{baleqn2}, we obtain a unique solution for $H_{k+1}$ for all $k \geq 0$, using Lemma~\ref{lem_balance_eqns}(i) as follows:
\begin{align}
    H_{k+1} = H_k + \frac{2 \left[\xi_k' - \frac{H_{k}}{2} r_k^2 \right]_+}{R_k^2} , \label{bal_up2}
\end{align}
where the quantity $R_k$ is defined in relation~\eqref{ada_rad_quantity}.

We present both Rules~1 and 2 in a combined way in Algorithm~\ref{algo_params_free_vi}.
%
%
%
%
\begin{algorithm}
\caption{Line search-free EGM using $\xi_k$ or $\xi_k'$}
\label{algo_params_free_vi}
\begin{algorithmic}[1]
\STATE \textbf{Input:} $v_0 \in X$, $H_0 = 0 $, $\widehat D > 0$ 
\FOR{$k = 0, 1, 2, \dots, T-1$}
    \STATE Compute $x_{k+1}$ and $v_{k+1}$ according to EGM (relation~\eqref{Kor_algo})
    \STATE Compute $R_{k}$ from relation~\eqref{ada_rad_quantity}
    \STATE $H_{k+1} = H_{k} + \frac{2 \left[\bar \xi_k - \frac{H_k}{2} r_k^2 \right]_+}{R_{k}^2}$ where $\bar \xi_k = \xi_k$ or $\xi_k'$ in Lemma~\ref{lem_generic_reln}
\ENDFOR
\end{algorithmic}
\end{algorithm}
%
%
%
%
If one compares the $H_{k+1}$ update rules in Rules 1 (cf.\ relation~\eqref{bal_up1}) and 2 (cf.\ relation~\eqref{bal_up2}), we observe that since $\xi_k \leq \xi_k'$, it follows that
\[
H_k + \frac{2 \left[\xi_k - \frac{H_k}{2} r_k^2 \right]_+}{R_{k}^2} \leq H_k + \frac{2 \left[\xi_k' - \frac{H_k}{2} r_k^2 \right]_+}{R_{k}^2}.
\]
This indicates that the balance equation in Rule~1 may yield larger step sizes $(\propto 1/H_{k+1})$ compared to Rule~2, which can be beneficial in practice. 


\paragraph{Rule~3:} The errors in Lemma~\ref{lem_generic_reln} can also be controlled using a classical AdaGrad-norm-style balance equation. In particular, large variations in $\widehat F$ should naturally lead to more conservative steps, while stable regions permit more aggressive updates. The classical AdaGrad-norm mechanism captures this intuition, as reflected in the balance equation~\eqref{balance_adagd_norm}. The complete method is presented in Algorithm~\ref{algo_adagd_norm}.
%
%
%
%
\begin{algorithm}
\caption{Line search-free EGM using classic AdaGrad-norm}
\label{algo_adagd_norm}
\begin{algorithmic}[1]
\STATE \textbf{Input:} $v_0 \in X$, $H_0 = 0$, $\widehat D > 0$
\FOR{$k = 0, 1, 2, \dots, T-1$}
    \STATE Compute $x_{k+1}$ and $v_{k+1}$ according to EGM (relation~\eqref{Kor_algo})
    \STATE $H_{k+1} = \sqrt{H_k^2 + \frac{\| \widehat F(x_{k+1}, \omega_{k+1}) - \widehat F(v_k, \omega_k) \|_*^2}{\widehat D^2}}$
\ENDFOR
\end{algorithmic}
\end{algorithm}
%
%
%
%

Since $H_k \leq H_{k+1}$, the following relation holds for Rule~3, for all $k \geq 0$:
\begin{align}
    \frac{\|\widehat F(x_{k + 1}, \omega_{k+1}) - \widehat F(v_k, \omega_k)\|_*^2}{2H_{k+1}} \leq (H_{k+1} - H_k) \widehat D^2 \leq \frac{\|\widehat F(x_{k + 1}, \omega_{k+1}) - \widehat F(v_k, \omega_k)\|_*^2}{H_{k+1}} . \label{adagd_ineq}
\end{align}
Relation~\eqref{adagd_ineq} is essential for controlling the error quantities arising in the analysis of the method.


\subsection{upper bound on $\Gamma_k$}

For the convenience of notation, we state $L(\delta) = L$ and $\bar L(\delta) = \bar L$ everywhere in the analysis.

\begin{lemma}\label{lem_algo1_2_gen}
    Let Assumptions~\ref{asum_set}, \ref{asum_map_growth}, and \ref{asum_stochastic} hold. Let $\widehat D > 0$ be a chosen constant. Define $\bar L = L$ for Algorithm~\ref{algo_params_free_vi} and $\bar L = L \sqrt{1 + 3 D^2/\widehat D^2}$ for Algorithm~\ref{algo_adagd_norm}. Then the iterates of Rules~1 and 2 (a.k.a. Algorithm~\ref{algo_params_free_vi}) and Rule~3 (a.k.a. Algorithm~\ref{algo_adagd_norm}) satisfy the following for any constant $\bar{H} \geq \sqrt{2} \bar L$, any $x \in X$, and any $k \geq 0$:
    \begin{align*}
        \EXP{\Gamma_k} \leq \frac{\EXP{\left[\min\{H_{k+1}, \bar H\} - H_k \right]_+}}{2} (c_1 D^2 + c_2 \widehat D^2) + \frac{4 L \delta + 2 \sigma^2}{\bar H \widehat D^2} (3 D^2 + \widehat D^2) ,
    \end{align*}
    where $c_1 = 3$, $c_2 = 1$ for Algorithm~\ref{algo_params_free_vi} and $c_1 = 3$, $c_2 = 2$ for Algorithm~\ref{algo_adagd_norm}.
\end{lemma}
\begin{proof}
    We consider the quantity $\Gamma_k$ (cf.\ relation~\eqref{eq1_lem_generic_reln_dup}) defined as
    \begin{align}
        \Gamma_k \equiv \Delta_k(H_{k+1}) + \frac{3}{2} (H_{k+1} - H_{k}) D^2 , \nonumber 
    \end{align}
    where $\Delta_k(H)$ is defined in Lemma~\ref{lem_generic_reln}.
    We choose a constant $\bar H > 0$ which can be chosen arbitrarily. We consider two cases. For each case, we analyze each example of Algorithm~\ref{algo_params_free_vi} and Algorithm~\ref{algo_adagd_norm} separately.

    \textit{Case $H_{k+1} \leq \bar H$:} For Rule~1 (cf.\ relation~\eqref{baleqn1}), we obtain
    \begin{align}
        \Gamma_k &\leq [\Delta_k(H_{k+1})]_+ + \frac{3}{2} (H_{k+1} - H_{k}) D^2 = \frac{H_{k + 1} - H_k}{2} (3D^2 + \widehat D^2) \nonumber\\
        & = \frac{\left[\min\{H_{k+1}, \bar H\} - H_k \right]_+}{2} (3D^2 + \widehat D^2) . \label{reln1_barH_large}
    \end{align}
    For Rule~2 (cf.\ relation~\eqref{baleqn2}), we apply the fact that $\Delta_k(H_{k+1}) \leq \xi_k' - \frac{H_{k+1}}{2} r_k^2$, where $r_k^2$ is defined in Lemma~\ref{lem_generic_reln}, and obtain
    \begin{align*}
        \Gamma_k &\leq \left[\xi_k' - \frac{H_{k+1}}{2} r_k^2 \right]_+ + \frac{3}{2} (H_{k+1} - H_{k}) D^2 = \frac{H_{k + 1} - H_k}{2} (3D^2 + \widehat D^2) .
    \end{align*}
    Following similar lines of argument as above, we obtain relation~\eqref{reln1_barH_large}.

    Next, we analyze Rule~3 (a.k.a. Algorithm~\ref{algo_adagd_norm}). Application of the upper bounds of relation~\eqref{eq3_lem_generic_reln} in Lemma~\ref{lem_generic_reln} to relation~\eqref{gammak} yields
    \begin{align}
        \Gamma_k \leq \bar \Delta_k(H_{k+1}) + \frac{3}{2} (H_{k+1} - H_k) D^2 , \label{gm_k_upbd}
    \end{align}
    which using the definition of $\bar \Delta_k(H_{k+1})$ can be algebraically simplified as
    \begin{align}
        \Gamma_k &\leq \frac{\| \widehat F(x_{k + 1}, \omega_{k+1}) - \widehat F(v_k, \omega_k) \|_*^2}{2H_{k+1}} - \frac{ H_{k+1}}{2} \|v_k - x_{k+1}\|^2 + \frac{3}{2} (H_{k + 1} - H_k) D^2 \nonumber\\
        &\leq \frac{H_{k+1}-H_k}{2} (3D^2 + 2 \widehat D^2) - \frac{ H_{k+1}}{2} \|v_k - x_{k+1}\|^2  \label{adagd_gammak2} \\
        & \leq \frac{\left[\min\{H_{k+1}, \bar H\} - H_k \right]_+}{2} (3D^2 + 2 \widehat D^2) . \label{analysis1algo_adagd}
    \end{align}
    where relation~\eqref{adagd_gammak2} is obtained by applying relation~\eqref{adagd_ineq}.

    \textit{Case $H_{k+1} \geq \bar H$:} For this case, relation~\eqref{gm_k_upbd} is satisfied. Since $H_{k+1} \geq \bar H$, hence $\bar \Delta_k(H_{k+1}) \leq \bar \Delta_k(\bar H)$ and we obtain
    \begin{align}
        \EXP{\Gamma_k} \leq \bar \Delta_k(\bar H) + \frac{3}{2} (H_{k + 1} - H_k) D^2 = \bar \Delta_k(\bar H) + \left[\frac{H_{k + 1} - H_k}{2} \widehat D^2 \right] \frac{3D^2}{\widehat D^2} . \label{gamma_upbd_deltabar}
    \end{align}
    Applying the upper bound relation~\eqref{eq3_lem_generic_reln} of Lemma~\ref{lem_generic_reln} and the fact that $H_{k+1} \geq \bar H$, for both Rules~1 and~2 of Algorithm~\ref{algo_params_free_vi}, we have  
    \begin{align*}
        \frac{H_{k+1} - H_k}{2} \widehat D^2 \leq [\bar \Delta_k(H_{k+1})]_+ \leq [\bar \Delta_k(\bar H)]_+ .
    \end{align*}
    Hence, substituting the preceding relation back to relation~\eqref{gamma_upbd_deltabar} yields
    \begin{align}
        \Gamma_k \leq [\bar \Delta_k(\bar H)]_+ \left[ 1 + \frac{3 D^2}{\widehat D^2} \right] . \label{gamma_upbd_deltabar2}
    \end{align}
    Next, we show the upper bound of $\bar \Delta_k(\bar H)$. We have
    \begin{align*}
        \bar \Delta_k(\bar H) & = \frac{1}{2 \bar H} \| \widehat F(x_{k + 1}, \omega_{k+1}) - \widehat F(v_k, \omega_k) \|_*^2 - \frac{\bar H}{2} \| x_{k + 1} - v_k \|^2 \\
        & = \frac{1}{2 \bar H} \| [F(x_{k + 1}) - F(v_k)] + [Z_{k+1} - Y_k] \|_*^2 - \frac{\bar H}{2} \| x_{k + 1} - v_k \|^2 \\
        & \leq \frac{1}{\bar H} \|F(x_{k + 1}) - F(v_k)\|_*^2 - \frac{\bar H}{2} \| x_{k + 1} - v_k \|^2 + \frac{1}{\bar H} \|Z_{k+1} - Y_k\|_*^2 ,
    \end{align*}
    where the stochastic errors $Y_k$ and $Z_k$ are defined in relation~\eqref{gam_k_reln}. Applying Assumption~\ref{asum_map_growth} to the preceding relation, we obtain
    \begin{align*}
        \bar \Delta_k(\bar H) \leq \left[\frac{L^2}{\bar H} - \frac{\bar H}{2} \right] \| x_{k + 1} - v_k \|^2 + \frac{4 L \delta}{\bar H} + \frac{1}{\bar H} \|Z_{k+1} - Y_k\|_*^2 .
    \end{align*}
    With the choice of $\bar H \geq \sqrt{2} L$, we obtain
    \begin{align*}
        \bar \Delta_k(\bar H) \leq \frac{4 L \delta}{\bar H} + \frac{1}{\bar H} \|Z_{k+1} - Y_k\|_*^2 .
    \end{align*}
    Substituting the preceding relation to relation~\eqref{gamma_upbd_deltabar2}, we obtain
    \begin{align}
        \Gamma_k \leq \frac{4 L \delta + \|Z_{k+1} - Y_k\|_*^2}{\bar H \widehat D^2} (3 D^2 + \widehat D^2) . \label{analysis3algo1}
    \end{align}

    For Rule~3 (Algorithm~\ref{algo_adagd_norm}), we apply relation~\eqref{adagd_ineq} to relation~\eqref{gamma_upbd_deltabar} to obtain the upper bound
    \begin{align*}
        \Gamma_k &\leq \bar \Delta_k(\bar H) + \frac{\|\widehat F(x_{k + 1}, \omega_{k+1}) - \widehat F(v_k, \omega_k)\|_*^2}{2 H_{k+1}} \frac{3D^2}{\widehat D^2} \\
        & \leq \frac{1}{2 \bar H} \left[1 + \frac{3 D^2}{\widehat D^2} \right] \|\widehat F(x_{k + 1}, \omega_{k+1}) - \widehat F(v_k, \omega_k)\|_*^2 - \frac{\bar H}{2} \|v_k - x_{k+1}\|^2 \\
        & \leq \left[ \frac{L^2}{\bar H} \left[1 + \frac{3 D^2}{\widehat D^2} \right] - \frac{\bar H}{2} \right] \|x_{k+1} - v_k\|^2 + \frac{4 L \delta + \|Z_{k+1} - Y_k\|_*^2}{\bar H \widehat D^2} (3 D^2 + \widehat D^2) .
    \end{align*}
    We choose $\bar H \geq \sqrt{2} L \left[\frac{\sqrt{3D^2 + \widehat D^2}}{\widehat D}\right]$, and see that for Rule~3 also relation~\eqref{analysis3algo1} is satisfied. For notational simplicity, we define the constant $\bar D_2^2 = 3 D^2 + \widehat D^2$. The constants $D_1$ and $\bar L$ are defined as
    \begin{align*}
        D_1^2 = \begin{cases}
            \bar D_2^2 \quad &\text{for Algorithm~\ref{algo_params_free_vi},} \\
            \bar D_2^2 + \widehat D^2 \quad &\text{for Algorithm~\ref{algo_adagd_norm}.}
        \end{cases} 
        \qquad \text{and }\qquad
        \bar L = \begin{cases}
            L \quad &\text{for Algorithm~\ref{algo_params_free_vi},} \\
            L \bar D_2/\widehat D \quad &\text{for Algorithm~\ref{algo_adagd_norm}.}
        \end{cases} 
    \end{align*}
    Now, we combine  both cases $H_{k+1} \leq \bar H$ and $H_{k+1} \geq \bar H$ for both Rules 1 and 2 of Algorithm~\ref{algo_params_free_vi} (cf.\ relations~\eqref{reln1_barH_large} and \eqref{analysis3algo1}) and Rule~3 of Algorithm~\ref{algo_adagd_norm} (cf.\ relations~\eqref{analysis1algo_adagd} and \eqref{analysis3algo1}) to obtain the following for any $\bar H \geq \sqrt{2} \bar L$, i.e.,
    \begin{align*}
        \Gamma_k &\leq \max \left\{ \frac{\left[\min\{H_{k+1}, \bar H\} - H_k \right]_+}{2} \bar D_1^2, \frac{4 L \delta + \|Z_{k+1} - Y_k\|_*^2}{\bar H \widehat D^2} \bar D_2^2 \right\} \\
        &\leq \frac{\left[\min\{H_{k+1}, \bar H\} - H_k \right]_+}{2} \bar D_1^2 + \frac{4 L \delta + \|Z_{k+1} - Y_k\|_*^2}{\bar H \widehat D^2} \bar D_2^2 .
    \end{align*}
    Next, we take expectation on both sides of the preceding relation and obtain
    \begin{align*}
        \EXP{\Gamma_k} \leq \frac{\EXP{\left[\min\{H_{k+1}, \bar H\} - H_k \right]_+}}{2} \bar D_1^2 + \frac{4 L \delta + \EXP{\|Z_{k+1} - Y_k\|_*^2}}{\bar H \widehat D^2} \bar D_2^2 .
    \end{align*}
    The second quantity on the right hand side can be simplified using Assumption~\ref{asum_stochastic} yielding
    \begin{align}
        \EXP{\|Z_{k+1} - Y_k\|_*^2} &= \EXP{\|Z_{k+1}\|_*^2 + \|Y_k\|_*^2 - 2 \langle \EXP{Z_{k+1} \mid \cF_k}, Y_k \rangle } \nonumber\\
        &= \EXP{\|Z_{k+1}\|_*^2 + \|Y_k\|_*^2} \leq 2 \sigma^2. \label{reln_stoc_split}
    \end{align}
    Hence, we obtain the desired relation of the lemma.
    \begin{align*}
        \EXP{\Gamma_k} \leq \frac{\EXP{\left[\min\{H_{k+1}, \bar H\} - H_k \right]_+}}{2} \bar D_1^2 + \frac{4 L \delta + 2 \sigma^2}{\bar H \widehat D^2} \bar D_2^2 .
    \end{align*}
    Finally, substituting the values of $\bar D_1$ and $\bar D_2$ concludes the proof.
\end{proof}
%
%


\section{Proofs for Section~\ref{sec_convergence}}

Throughout this section, for any fixed $\delta>0$, we write $L(\delta)=L$ and $\bar L(\delta)=\bar L$ for notational simplicity.

\subsection{Proof of Theorem~\ref{thm_universal}}

\begin{proof}
    Let $x \in X$ and $T \geq 1$ be arbitrary. Let us denote
    \begin{align}
        \bar D_1^2 = c_1 D^2 + c_2 \widehat D^2 \qquad \text{and} \qquad \bar D_2^2 = 3D^2 + \widehat D^2 . \label{notation_dia}
    \end{align}
    Applying relation~\eqref{gam_k_reln} to relation~\eqref{eq1_lem_generic_reln_dup} and use of Lemma~\ref{lem_sum_H_k} yields
    \begin{align*}
        \mathbb{E}[\epsilon_T] &\leq \frac{\bar D_1^2}{2 T} \EXP{\left[\min\{H_T, \bar H\} - H_0 \right]_+} + \frac{4 L \delta + 2\sigma^2}{\bar H \widehat D^2} \bar D_2^2 + \frac{D \sigma}{\sqrt{T}} .
    \end{align*}
    Recall that we choose $H_0 = 0$ for all the algorithms. Hence,
    \begin{align*}
        \left[\min\{H_{T}, \bar H\} - H_0 \right]_+ = \min\{H_{T}, \bar H\} \leq \bar H .
    \end{align*}
    So we obtain
    \begin{align}
        \mathbb{E}[\epsilon_T] \leq \frac{\bar H \bar D_1^2}{2 T} + \frac{4 L \delta + 2 \sigma^2}{\bar H \widehat D^2} \bar D_2^2 + \frac{D \sigma}{\sqrt{T}} . \label{eqn_thm_main1}
    \end{align}
    For simplicity, let us denote
    \begin{align}
        \theta = \frac{4 L \delta + 2 \sigma^2}{\widehat D^2} \bar D_2^2 . \label{theta_reln}
    \end{align}
    We minimize the right-hand side of the relation~\eqref{eqn_thm_main1} with respect to $\bar H \geq \sqrt{2} \bar L$ and obtain the optimal
    \begin{align*}
        \bar H^* = \max\left\{\sqrt{2} \bar{L}, \frac{\sqrt{2T\theta}}{\bar{D}_1} \right\}.
    \end{align*}
    Substituting $\bar H = \bar H^*$ back to relation~\eqref{eqn_thm_main1}, we obtain
    \begin{align*}
    \mathbb{E}[\epsilon_T] &\leq
    \max\left\{ \frac{\bar{L}\bar{D}_1^2}{\sqrt{2}T}, \sqrt{\frac{\bar{D}_1^2\theta}{2T}} \right\} + \min\left\{ \frac{\theta}{\sqrt{2}\bar{L}}, \sqrt{\frac{\bar{D}_1^2\theta}{2T}} \right\} + \frac{D\sigma}{\sqrt{T}}
    \nonumber\\
    &\leq
    \frac{\bar{L}\bar{D}_1^2}{\sqrt{2}T} + 2\sqrt{\frac{\bar{D}_1^2\theta}{2T}} + \frac{D\sigma}{\sqrt{T}}
    \nonumber\\
    &= \frac{\bar{L}\bar{D}_1^2}{\sqrt{2}T} + \sqrt{\frac{2\bar{D}_1^2\theta}{T}} + \frac{D\sigma}{\sqrt{T}}.
\end{align*}
Substituting the value of $\theta$ from relation~\eqref{theta_reln} and applying the sub-additivity of square-root function yields the final relation of the theorem.
\end{proof}
%
%

\subsection{Proof of Corollary~\ref{cor_holder}}

\begin{proof}
    We start from the relation of Theorem~\ref{thm_universal} and use $\bar L = L \left(\frac{\bar D_2}{\widehat D}\right)^a = L(\delta) \left(\frac{\bar D_2}{\widehat D}\right)^a$, where $a = 0$ for Algorithm~\ref{algo_params_free_vi} and $a = 1$ for Algorithm~\ref{algo_adagd_norm}, and the quantity $L(\delta)$ is derived in Lemma~\ref{lem_holder_to_delta_L}, and obtain
    \begin{align}
        \mathbb{E}[\epsilon_T]
        &\leq
        \frac{\bar D_1^2}{\sqrt{2}T}
        \left(\frac{\bar D_2}{\widehat D}\right)^a
        \left(\frac{1-\nu}{4}\right)^{\frac{1-\nu}{1+\nu}}
        \nu^{\frac{\nu}{1+\nu}}
        H_\nu^{\frac{2}{1+\nu}}
        \delta^{-\frac{1-\nu}{1+\nu}}
        \nonumber \\
        &\quad+
        \sqrt{\frac{\bar D_1^2}{T}}
        \left[
        \frac{2\sqrt{2}\,\bar D_2}{\widehat D}
        \left(\frac{1-\nu}{4}\right)^{\frac{1-\nu}{2(1+\nu)}}
        \nu^{\frac{\nu}{2(1+\nu)}}
        H_\nu^{\frac{1}{1+\nu}}
        \delta^{\frac{\nu}{1+\nu}}
        +
        \frac{2\bar D_2}{\widehat D}\sigma
        \right]
        +
        \frac{D\sigma}{\sqrt{T}} \label{eps_holder1}
    \end{align}
    for all $\nu \in [0,1]$ with the convention that $0^0 = 1$, where $D_1^2 = c_1 D^2 + c_2 \widehat D^2$ and $D_2^2 = 3D^2 + \widehat D^2$. Since $L(\delta)$ in Lemma~\ref{lem_holder_to_delta_L} holds for any $\delta > 0$, we minimize relation~\eqref{eps_holder1} with respect to $\delta$ and obtain the optimal for any $\nu \in (0,1)$ as
    \begin{align*}
        \delta^*
        =
        \left[
        \frac{1-\nu}{4\nu}
        \left(\frac{\widehat D}{\bar D_2}\right)^{1-a}
        \sqrt{\frac{\bar D_1^2}{T}}
        \right]^{1+\nu}
        \left(\frac{1-\nu}{4}\right)^{\frac{1-\nu}{2}}
        \nu^{\frac{\nu}{2}}
        H_\nu.
    \end{align*}
    Substituting $\delta = \delta^*$ back to relation~\eqref{eps_holder1} yields the following for any $\nu \in (0,1)$:
    \begin{align}
        \mathbb{E}[\epsilon_T]
        &\leq
        \frac{2^{\frac{1}{2}-\nu}H_\nu}
        {(1-\nu)^{\frac{1-\nu}{2}}\nu^{\frac{\nu}{2}}}
        \left(\frac{\bar D_2}{\widehat D}\right)^{1+\nu(a-1)}
        \left(\frac{\bar D_1^2}{T}\right)^{\frac{1+\nu}{2}}
        +
        \sqrt{\frac{\bar D_1^2}{T}}\,
        \frac{2\bar D_2}{\widehat D}\sigma
        +
        \frac{D\sigma}{\sqrt{T}} \nonumber
        \\
        &=
        \frac{2^{\frac{1}{2}-\nu}H_\nu}
        {(1-\nu)^{\frac{1-\nu}{2}}\nu^{\frac{\nu}{2}}}
        \left(\frac{\bar D_2}{\widehat D}\right)^{1+\nu(a-1)}
        \left(\frac{\bar D_1^2}{T}\right)^{\frac{1+\nu}{2}}
        +
        \frac{\sigma}{\sqrt{T}}
        \left(
        \frac{2\bar D_1\bar D_2}{\widehat D}+D
        \right). \label{eps_holder_main}
    \end{align}
    Next, we consider the corner cases. When $\nu = 0$, relation~\eqref{eps_holder1} yields
    \begin{align*}
        \mathbb{E}[\epsilon_T]
        &\leq
        \inf_{\delta>0}
        \left\{
        \frac{\bar D_1^2}{4\sqrt{2}T}
        \left(\frac{\bar D_2}{\widehat D}\right)^a
        H_0^2\delta^{-1}
        +
        \sqrt{\frac{\bar D_1^2}{T}}
        \left[
        \frac{\sqrt{2}\,\bar D_2}{\widehat D}H_0
        +
        \frac{2\bar D_2}{\widehat D}\sigma
        \right]
        \right\}
        +
        \frac{D\sigma}{\sqrt{T}} \nonumber
        \\
        &=
        \lim_{\delta\to\infty}
        \left\{
        \frac{\bar D_1^2}{4\sqrt{2}T}
        \left(\frac{\bar D_2}{\widehat D}\right)^a
        H_0^2\delta^{-1}
        \right\}
        +
        \sqrt{\frac{\bar D_1^2}{T}}
        \left[
        \frac{\sqrt{2}\,\bar D_2}{\widehat D}H_0
        +
        \frac{2\bar D_2}{\widehat D}\sigma
        \right]
        +
        \frac{D\sigma}{\sqrt{T}} \nonumber
        \\
        &=
        \sqrt{\frac{\bar D_1^2}{T}}
        \left[
        \frac{\sqrt{2}\,\bar D_2}{\widehat D}H_0
        +
        \frac{2\bar D_2}{\widehat D}\sigma
        \right]
        +
        \frac{D\sigma}{\sqrt{T}} \nonumber
        \\
        &=
        \frac{\sqrt{2}\,H_0\bar D_1\bar D_2}
        {\widehat D\sqrt{T}}
        +
        \frac{\sigma}{\sqrt{T}}
        \left(
        \frac{2\bar D_1\bar D_2}{\widehat D}+D
        \right).
    \end{align*}
    On the other hand, when $\nu = 1$, relation~\eqref{eps_holder1} simplifies to
    \begin{align*}
        \mathbb{E}[\epsilon_T]
        &\leq
        \inf_{\delta>0}
        \left\{
        \frac{\bar D_1^2}{\sqrt{2}T}
        \left(\frac{\bar D_2}{\widehat D}\right)^a H_1
        +
        \sqrt{\frac{\bar D_1^2}{T}}
        \left[
        \frac{2\sqrt{2}\,\bar D_2}{\widehat D}
        \sqrt{\delta H_1}
        +
        \frac{2\bar D_2}{\widehat D}\sigma
        \right]
        \right\}
        +
        \frac{D\sigma}{\sqrt{T}}
        \\
        &=
        \frac{\bar D_1^2}{\sqrt{2}T}
        \left(\frac{\bar D_2}{\widehat D}\right)^a H_1
        +
        \sqrt{\frac{\bar D_1^2}{T}}
        \frac{2\bar D_2}{\widehat D}\sigma
        +
        \frac{D\sigma}{\sqrt{T}}
        +
        \lim_{\delta\to0}
        \left\{
        \sqrt{\frac{\bar D_1^2}{T}}
        \frac{2\sqrt{2}\,\bar D_2}{\widehat D}
        \sqrt{\delta H_1}
        \right\}
        \\
        &=
        \frac{\bar D_1^2}{\sqrt{2}T}
        \left(\frac{\bar D_2}{\widehat D}\right)^a H_1
        +
        \frac{\sigma}{\sqrt{T}}
        \left(
        \frac{2\bar D_1\bar D_2}{\widehat D}+D
        \right).
    \end{align*}
    Hence, relation~\eqref{eps_holder_main} holds for any $\nu \in [0,1]$, with the convention $0^0 = 1$. 
\end{proof}
%
%


\section{Missing details of Section~\ref{sec_dist_ada}}\label{app_dist_ada}

There can exist settings where the underlying set is complex, for example multiple intersections of convex sets. In these settings, a closed form expression of projection might not be easily computable and one might resort to solving the EGM prox subproblem upto a desired accuracy. So, it is theoretically interesting to investigate whether the estimate of the true diameter can be done adaptively. This question is also motivated by the broader goal of extending parameter-free VI methods to problems with an unknown effective diameter, including VIs over unbounded feasible sets or the entire space. A complete treatment of such problems requires additional ideas and is beyond the scope of this work. As an initial step toward this goal, we remain in the bounded-set setting but remove the need to specify a sufficiently accurate diameter estimate in advance.

\subsection{Algorithmic Development for Distance Adaptation}\label{app_subsec_da}

We consider the modified ExtraGradient Method, aka. modified EGM, which we present again here. The method starts with the initialization $v_0 \in X$ and $H_0 = 0$, and update for all $k\geq 0$ as
\begin{align}
    &x_{k+1} = \argmin_{x \in X} \left\{ D_k \langle \widehat F(v_k, \omega_k),\, x - v_k \rangle + \frac{H_{k}}{2} \|x - v_k\|^2 \right\}, \label{mod_EGM} \\
    &v_{k+1} = \argmin_{x \in X} \left\{ D_k \langle \widehat F(x_{k+1}, \omega_{k+1}),\, x - x_{k+1} \rangle + \frac{H_{k}}{2} \|x - v_k\|^2 \right\} , \nonumber
\end{align}
where $D_k$ is the adaptive diameter and is updated for all $k \geq 0$ as
\begin{align*}
    D_{k+1} = \max \{ D_k, \|v_{k+1} - v_0\|, r_k \} \qquad \text{with} \quad D_0 = \widehat D > 0.
\end{align*}
We present a lemma from where the error balance equations will be derived.
\begin{lemma}\label{lem_gen_dist_ada}
    For the modified EGM, the following holds for any $x \in X$ and any $k \geq 0$:
    \begin{align*}
        D_k &\langle F(x_{k+1}), x_{k+1} - x \rangle \leq \widehat \Delta_k(H_{k+1}) + \frac{H_{k+1} - H_k}{2} r_k^2  + D_k \langle Z_{k+1}, h_k - x_{k+1} \rangle \nonumber\\
        & + \frac{D_k \|Z_{k+1}\|_*^2}{2 \eta} + \frac{\eta D_k}{2} [\|x - h_k\|^2 - \|x - h_{k+1}\|^2] + \frac{H_{k}}{2} [\|x - v_k\|^2 - \|x - v_{k+1}\|^2] ,
    \end{align*}
    where the sequence $\{h_k\}$ is defined in Lemma~\ref{lem_generic_reln}, $\eta > 0$ is any constant, $Z_k$ is the stochastic error (cf.\ relation~\eqref{stoc_errors}), $D_k$ is the adaptive diameter, and $\widehat \Delta_k(H) := \widehat \xi_k - \frac{H}{2} r_k^2$, where
    \begin{align*}
        &r_k^2 = \|v_k - x_{k+1}\|^2 + \|v_{k+1} - x_{k+1}\|^2 \\
        &\widehat \xi_k = D_k \langle \widehat F(x_{k+1}, \omega_{k+1}), x_{k+1} - v_{k+1} \rangle - \frac{H_k}{2} \|v_k - v_{k+1}\|^2 + \frac{H_k}{2} r_k^2 \\
        & \quad \leq D_k \la \widehat F(x_{k+1},\omega_{k+1}) - \widehat F(v_k,\omega_k), x_{k+1} - v_{k+1} \ra \equiv \widehat \xi_k' .
    \end{align*}
    Moreover,
    \begin{align}
        \widehat \Delta_k(H) \leq \widehat \xi_k' - \frac{H}{2} r_k^2 \leq \frac{D_k^{2}}{2 H} \|\widehat F(x_{k+1},\omega_{k+1}) - \widehat F(v_k,\omega_k)\|_*^2 - \frac{H}{2} \|x_{k+1} - v_k\|^2 \equiv \widetilde \Delta_k(H) . \label{tilde_delta_upbd}
    \end{align}
\end{lemma}
\begin{proof}
    Similar to the start of Lemma~\ref{lem_generic_reln}, we apply $H_k$ strong convexity of $v_{k+1}$ update and simplify terms to obtain
    \begin{align}
        D_k \langle \widehat F(x_{k+1}, \omega_{k+1}), x_{k+1} - x \rangle &\leq \widehat \xi_k - \frac{H_{k+1}}{2} r_k^2 + \frac{H_{k+1} - H_k}{2} r_k^2 \nonumber\\
        & + \frac{H_{k}}{2} [\|x - v_k\|^2 - \|x - v_{k+1}\|^2] , \label{dis_ada_main1}
    \end{align}
    where $r_k^2 = \|v_k - x_{k+1}\|^2 + \|v_{k+1} - x_{k+1}\|^2$ and
    \begin{align*}
        \widehat \xi_k = D_k \langle \widehat F(x_{k+1}, \omega_{k+1}), x_{k+1} - v_{k+1} \rangle - \frac{H_k}{2} \|v_k - v_{k+1}\|^2 + \frac{H_k}{2} r_k^2 .
    \end{align*}
    We denote 
    \begin{align*}
        \widehat \Delta_k(H_{k+1}) := \widehat \xi_k - \frac{H_{k+1}}{2} r_k^2 .
    \end{align*}
    The rest of the proof follows along lines similar to that of Lemma~\ref{lem_generic_reln}, and applying Lemma~\ref{lem_stoc_error} to simplify the left-hand side of relation~\eqref{dis_ada_main1}, we can finally obtain
    \begin{align*}
        D_k \langle F(x_{k+1}), x_{k+1} - x \rangle &\leq \widehat \Delta_k(H_{k+1}) + \frac{H_{k+1} - H_k}{2} r_k^2  + D_k \langle Z_{k+1}, h_k - x_{k+1} \rangle \nonumber\\
        & \qquad + \frac{D_k \|Z_{k+1}\|_*^2}{2 \eta} + \frac{\eta D_k}{2} [\|x - h_k\|^2 - \|x - h_{k+1}\|^2] \nonumber\\
        & \qquad + \frac{H_{k}}{2} [\|x - v_k\|^2 - \|x - v_{k+1}\|^2].
    \end{align*}
    Moreover from $H_k$ strong convexity of $x_{k+1}$ update, we can obtain
    \begin{align*}
        D_k \left\langle
        \widehat{F}(v_k,\omega_k),
        x_{k+1}-v_{k+1}
        \right\rangle
        \leq
        \frac{H_k}{2}\left\|v_k-v_{k+1}\right\|^2
        -
        \frac{H_k}{2}r_k^2 .
    \end{align*}
    Applying the preceding relation to upper bound $\widehat \Delta_k(H_{k+1})$, we obtain
    \begin{align*}
        \widehat \Delta_k(H_{k+1}) &\leq \underbrace{D_k \la \widehat F(x_{k+1},\omega_{k+1}) - \widehat F(v_k,\omega_k), x_{k+1} - v_{k+1} \ra}_{\widehat \xi_k'} - \frac{H_{k+1}}{2} r_k^2 \\
        & \leq \underbrace{\frac{D_k^{2}}{2 H_{k+1}} \|\widehat F(x_{k+1},\omega_{k+1}) - \widehat F(v_k,\omega_k)\|_*^2 - \frac{H_{k+1}}{2} \|x_{k+1} - v_k\|^2}_{\widetilde \Delta_k(H_{k+1})}  .
    \end{align*}
    This concludes the proof.
\end{proof}
Similar to Section~\ref{sec_algo}, here we also design balance equations to update the sequence $H_k$ and control the error. We present two such examples below:
\begin{align}
    \frac{H_{k+1} - H_k}{2} D_{k+1}^2 &= \left[\widehat \Delta_k(H_{k+1}) \right]_+ \qquad &\text{(Method 1)} , \nonumber\\ 
    \frac{H_{k+1} - H_k}{2} D_{k+1}^2 &= \left[\widehat \xi_k' - \frac{H_{k+1}}{2} r_k^2 \right]_+ \qquad &\text{(Method 2)} . \nonumber 
\end{align}
Since the balance equations in Section~\ref{sec_algo} are labeled as ``Rule", we use ``Method" with the corresponding number here to distinguish. Applying Lemma~\ref{lem_balance_eqns}, a unique solution of $H_{k+1}$ for both the methods can be obtained, i.e., for all $k \geq 0$,
\begin{align*}
    H_{k+1} &= H_k + \frac{2 \left[\widehat \xi_k - \frac{H_k}{2} r_k^2 \right]_+}{\widehat R_{k}^2} \\
    &= H_k + \frac{2 \left[D_k \langle \widehat F(x_{k+1}, \omega_{k+1}), x_{k+1} - v_{k+1} \rangle - \frac{H_k}{2} \|v_k - v_{k+1}\|^2 \right]_+}{\widehat R_{k}^2} \qquad &\text{(Method 1)} \\
    \text{and } H_{k+1} &= H_k + \frac{2 \left[\widehat \xi_k' - \frac{H_k}{2} r_k^2 \right]_+}{\widehat R_{k}^2} \qquad &\text{(Method 2)},
\end{align*}
where the quantity $\widehat R_{k}$ is defined as
\begin{align}
    \widehat R_{k}^2 = D_{k+1}^2 + r_k^2 = D_{k+1}^2 + \| v_k - x_{k + 1} \|^2 + \| v_{k + 1} - x_{k + 1} \|^2 . \label{rad_dist_ada}
\end{align}
We present the pseudocode of both methods in Algorithm~\ref{algo_dist_ada_vi}. Both balance equations of Methods~1 and 2 yield an upper bound on $H_k$ for all $k \geq 0$, which we present next.

%
%
%
%
\begin{algorithm}
\caption{Line search-free modified EGM using Distance Adaptation}
\label{algo_dist_ada_vi}
\begin{algorithmic}[1]
\STATE \textbf{Input:} $v_0 \in X$, $H_0 = 0 $, small estimate $D_0 = \widehat D>0$
\FOR{$k = 0, 1, 2, \dots, T-1$}
    \STATE Compute $x_{k+1}$ and $v_{k+1}$ according to modified EGM (relation~\eqref{mod_EGM})
    \STATE Compute $D_{k+1} = \max \{ D_k, \|v_{k+1} - v_0\|, r_k \}$
    \STATE Compute $\widehat R_k$ as in relation~\eqref{rad_dist_ada}
    \STATE $H_{k+1} = H_{k} + \frac{2 \left[\bar \xi_k - \frac{H_k}{2} r_k^2 \right]_+}{\widehat R_{k}^2}$ where $\bar \xi_k = \widehat \xi_k$ or $\widehat \xi_k'$ as in Lemma~\ref{lem_gen_dist_ada}
\ENDFOR
\end{algorithmic}
\end{algorithm}
%
%
%
%

Note that the effective step size of Algorithm~\ref{algo_dist_ada_vi} is $D_k / H_k$ for all $k \geq 1$. There is a beautiful analogy between Rules~1 and 2 of Algorithm~\ref{algo_params_free_vi} and Methods~1 and 2 of Algorithm~\ref{algo_dist_ada_vi}. Let us denote $H_k = M_k D_k$. Then, the balance equations of Methods~1 and 2 can be equivalently written as
\begin{align*}
    \frac{M_{k+1} D_{k+1} - M_k D_k}{2} D_{k+1}^2 &= \left[\widehat \Delta_k(M_{k+1} D_{k+1}) \right]_+ \qquad &\text{(Method 1)} , \\
    \frac{M_{k+1} D_{k+1} - M_k D_k}{2} D_{k+1}^2 &= \left[\widehat \xi_k' - \frac{M_{k+1} D_{k+1}}{2} r_k^2 \right]_+ \qquad &\text{(Method 2)} .
\end{align*}
In the special case where $D_k=\widehat D$ for all $k\geq0$, Methods~1 and~2 reduce to Rules~1 and~2, respectively, of Algorithm~\ref{algo_params_free_vi}, with $M_k$ corresponding to $H_k$. Next, we present a lemma to obtain the upper bound on $H_k$ for Algorithm~\ref{algo_dist_ada_vi}.

\begin{lemma}\label{lem_H_upbd}
    Under Assumptions~\ref{asum_set} and \ref{asum_map_growth}, and with the initialization $H_0 = 0$, the balance equations of Methods~1 and 2 (cf.\ relations~\eqref{bal_ada1} and \eqref{bal_ada2}) satisfy the following for any $k \geq 1$:
    \begin{align*}
        H_k \leq \sqrt{2} L(\delta) D + 4 \sqrt{L(\delta) \delta k} + 2 \sqrt{\sum_{i = 0}^{k-1} \|Z_{i+1} - Y_i\|_*^2} .
    \end{align*}
\end{lemma}
\begin{proof}
    Applying relation~\eqref{tilde_delta_upbd} of Lemma~\ref{lem_gen_dist_ada}, for both Methods~1 and~2, we have
    \begin{align*}
        \frac{H_{k+1} - H_k}{2} D_{k+1}^2 &\leq \left[ \frac{D_k^{2}}{2 H_{k+1}} \|\widehat F(x_{k+1},\omega_{k+1}) - \widehat F(v_k,\omega_k)\|_*^2 - \frac{H_{k+1}}{2} \|x_{k+1} - v_k\|^2 \right]_+ \\
        & \leq \left[ \frac{D_k^{2}}{H_{k+1}} \|F(x_{k+1}) - F(v_k)\|_*^2 - \frac{H_{k+1}}{2} \|x_{k+1} - v_k\|^2 \right]_+ \\
        & \hspace{6cm} + \frac{D_k^{2}}{H_{k+1}} \|Z_{k+1} - Y_k\|_*^2 ,
    \end{align*}
    where the quantity $Z_{k+1}$ and $Y_k$ are the stochastic errors and are defined in relation~\eqref{stoc_errors}. Now application of Assumption~\ref{asum_map_growth} on the preceding relation yields
    \begin{align*}
        \frac{H_{k+1} - H_k}{2} D_{k+1}^2 \leq \left[ \frac{D_k^{2}}{H_{k+1}} L^2 - \frac{H_{k+1}}{2} \right]_+ \|x_{k+1} - v_k\|^2 + \frac{4 L \delta D_k^{2}}{H_{k+1}} + \frac{D_k^{2}}{H_{k+1}} \|Z_{k+1} - Y_k\|_*^2 ,
    \end{align*}
    where for notational simplicity, we denoted $L(\delta) \equiv L$ for any $\delta > 0$. We multiply both sides of the preceding relation by $4 H_{k+1}$ and obtain
    \begin{align}
        2 H_{k+1} (H_{k+1} \!-\! H_k) D_{k+1}^2 &\leq 2 \left[ 2D_k^{2} L^2 - H_{k+1}^2 \right]_+ \|x_{k+1} \!-\! v_k\|^2 + 16 L \delta D_k^{2} + 4 D_k^{2} \|Z_{k+1} - Y_k\|_*^2 \nonumber\\
        & \leq 2 \left[ 2D^{2} L^2 - H_{k+1}^2 \right]_+ D_{k+1}^2 + 16 L \delta D_k^{2} + 4 D_k^{2} \|Z_{k+1} - Y_k\|_*^2 \label{H_upbd1} ,
    \end{align}
    where the last relation is obtained using Assumption~\ref{asum_set}. Note that for any $k \geq 0$,
    \begin{align*}
        2 H_{k+1} (H_{k+1} \!-\! H_k) = (H_{k+1} + H_{k+1}) (H_{k+1} \!-\! H_k) \geq (H_{k+1} + H_{k}) (H_{k+1} \!-\! H_k) = H_{k+1}^2 \!-\! H_k^2 .
    \end{align*}
    Applying the preceding relation to lower bound the left-hand side of relation~\eqref{H_upbd1}, we obtain
    \begin{align*}
        (H_{k+1}^2 - H_k^2) D_{k+1}^2 \leq 2 \left[ 2D^{2} L^2 - H_{k+1}^2 \right]_+ D_{k+1}^2 + 16 L \delta D_k^{2} + 4 D_k^{2} \|Z_{k+1} - Y_k\|_*^2 .
    \end{align*}
    We divide both sides of the relation by $D_{k+1}^2$ and using the fact that $D_k^2/D_{k+1}^2 \leq 1$, we obtain
    \begin{align*}
        H_{k+1}^2 - H_k^2 \leq 2 \left[ 2 D^2 L^2 - H_{k+1}^2 \right]_+ + 16 L \delta + 4 \|Z_{k+1} - Y_k\|_*^2.
    \end{align*}
    Applying Lemma~\ref{lem_aux_seq} with $\gamma_k = H_k^2$, $a = 2$, $V = 2 L^2 D^2$, $e_k = 16 L \delta + 4 \|Z_{k+1} - Y_k\|_*^2$, and noting that $\gamma_0 = H_0^2 = 0$ (based on  the initialization $H_0 = 0$), we have for all $k \geq 1$,
\begin{align*}
    H_k^2 \leq 2 L^2 D^2 + 16 L \delta k + 4 \sum_{i = 0}^{k-1} \|Z_{i+1} - Y_i\|_*^2 .
\end{align*}
Taking square root on both sides of the preceding relation and using sub-additivity of the square root function yield the final relation of the lemma.
\end{proof}
Lemma~\ref{lem_H_upbd} serves as a counterpart to Lemma~\ref{lem_algo1_2_gen}. However, their analyses differ. Lemma~\ref{lem_algo1_2_gen} also accommodates the classical AdaGrad-norm update (cf.\ Algorithm~\ref{algo_adagd_norm}) within the same universal analysis by introducing the proxy constant $\bar H$. However, in the distance-adaptive setting, we no longer consider the classical AdaGrad-norm update, which allows us to improve the constants and considerably simplify the analysis. A similar analysis can also be applied to Algorithm~\ref{algo_params_free_vi}.

\subsection{Convergence analysis for Distance Adaptive VI algorithm}

In this subsection, we present theorems that establish the convergence rate of Algorithm~\ref{algo_dist_ada_vi} under Assumption~\ref{asum_map_growth} and also study the rate for H\"older continuous mapping as a special case.

\begin{theorem}\label{thm_distada}
    Let Assumptions~\ref{asum_set}, \ref{asum_map_growth} and \ref{asum_stochastic} hold, and let $\rho = \max \{1, D/\widehat D \}$. Then, for any $T \geq 1$, the iterates generated by Algorithm~\ref{algo_dist_ada_vi} converge as
    \begin{align*}
        \EXP{\epsilon_{\tau(T)}}
        \leq \left( \rho \right)^{1/T} \log\left( e \rho \right)
        \left[\frac{3\sqrt{2}L(\delta) D^2}{T}
        +
        \frac{
        12D\sqrt{L(\delta) \delta}
        +(6\sqrt{2}+5)\sigma D
        }{\sqrt{T}}
        \right],
    \end{align*}
    where ${\tau(T)} := \argmin_{1 \leq t \leq T} \frac{D_t}{\sum_{k=0}^{t-1} D_k}$ and $\epsilon_{\tau(T)}$ is the corresponding residual function (cf.\ relation~\eqref{residual_val_eq}).
\end{theorem}
\begin{proof}
    We start from the first relation of Lemma~\ref{lem_gen_dist_ada}, and see that since $\widehat \Delta_k(H) \leq \widehat \xi_k' - \frac{H}{2} r_k^2$ from relation~\eqref{tilde_delta_upbd}, hence using the balance equations of both Methods~1 and 2 can yield
    \begin{align*}
        D_k \langle F(x_{k+1}), x_{k+1} - x \rangle &\leq \frac{H_{k+1} - H_k}{2} (r_k^2 + D_{k+1}^2)  + D_k \langle Z_{k+1}, h_k - x_{k+1} \rangle \nonumber\\
        & \qquad + \frac{D_k \|Z_{k+1}\|_*^2}{2 \eta} + \frac{\eta D_k}{2} [\|x - h_k\|^2 - \|x - h_{k+1}\|^2] \nonumber\\
        & \qquad + \frac{H_{k}}{2} [\|x - v_k\|^2 - \|x - v_{k+1}\|^2].
    \end{align*}
    We apply the fact that $r_k \leq D_{k+1}$ in order to simplify the first quantity on the right-hand side. Then we sum both sides of the preceding relation from $k = 0$ to $T-1$ for any $T \geq 1$, yielding
    \begin{align}
        \sum_{k=0}^{T-1} D_k \langle F(x_{k+1}), x_{k+1} - x \rangle \leq & (H_T - H_0) D_T^2 + \sum_{k=0}^{T-1} D_k \langle Z_{k+1}, h_k - x_{k+1} \rangle \nonumber\\
        & + \frac{\sum_{k=0}^{T-1} D_k \|Z_{k+1}\|_*^2}{2 \eta} + \frac{\eta D_T}{2} [\|x - h_0\|^2 - \|x - h_{T}\|^2] \nonumber\\
        & + \sum_{k=0}^{T-1} \frac{H_{k}}{2} [\|x - v_k\|^2 - \|x - v_{k+1}\|^2] , \label{adadist_main1}
    \end{align}
    where in some of the terms, we use the fact that $D_k \leq D_T$ for any $0 \leq k \leq T$ and any $T \geq 1$. We initialized $H_0 = 0$. We also denote $d_k = \|x-v_k\|$. Then, the last term on the right-hand side of relation~\eqref{adadist_main1} can be simplified as
    \begin{align*}
        \frac{1}{2} \sum_{k=0}^{T-1} H_{k}[d_k^2 - d_{k+1}^2] &= \frac{1}{2} \left[H_0 d_0^2 - H_{T-1} d_T^2 + \sum_{k=1}^{T-1} d_k^2 (H_k - H_{k-1}) \right] \nonumber\\
        & = \frac{1}{2} \left[H_0 d_0^2 - H_{T-1} d_T^2 + \bar d_T^2 (H_{T-1} - H_{0}) \right] \nonumber\\
        & = \frac{H_{T-1}}{2} (\bar d_T^2 - d_T^2) \leq H_T \bar d_T (\bar d_T - d_T) ,
    \end{align*}
    where $\bar d_T = \max_{1 \leq k \leq T-1} d_k$. Using Assumption~\ref{asum_set}, we obtain $\bar d_T \leq D$. Let us denote the index $\iota = \argmax_{1 \leq k \leq T-1} d_k$. Then, 
    \begin{align*}
        \bar d_T - d_T = d_\iota - d_T \leq \|v_T - v_\iota\| \leq \|v_T - v_0\| + \|v_\iota - v_0\| \leq 2 D_{T}.
    \end{align*}
    Hence, we obtain
    \begin{align*}
        \sum_{k=0}^{T-1} \frac{H_{k}}{2} [\|x - v_k\|^2 - \|x - v_{k+1}\|^2] \leq 2 D D_T H_T .
    \end{align*}
    Moreover using Assumption~\ref{asum_set}, the first quantity on the right-hand side of relation~\eqref{adadist_main1} can be upper bounded as
    \begin{align*}
        (H_T - H_0) D_T^2 \leq D D_T H_T .
    \end{align*}
    Making use of the preceding two relationships, relation~\eqref{adadist_main1} simplifies as
    \begin{align*}
        \sum_{k=0}^{T-1} D_k \langle F(x_{k+1}), x_{k+1} - x \rangle &\leq 3 D D_T H_T + \sum_{k=0}^{T-1} D_k \langle Z_{k+1}, h_k - x_{k+1} \rangle \nonumber\\
        & \qquad + \frac{D_T \sum_{k=0}^{T-1} \|Z_{k+1}\|_*^2}{2 \eta} + \frac{\eta D_T}{2} \|x - h_0\|^2 .
    \end{align*}
    Taking maximum with respect to $x \in X$ and dividing by $\sum_{k=0}^{T-1} D_k$ on both sides and using Assumption~\ref{asum_set}, we obtain the upper bound on the residual function $\epsilon_T = \epsilon_T((x_t, D_t)_{t=1}^T)$ as
    \begin{align}
        \epsilon_T &\leq \frac{D_T}{\sum_{k=0}^{T-1} D_k} \left[ 3 D H_T + \frac{\sum_{k=0}^{T-1} \|Z_{k+1}\|_*^2}{2 \eta} + \frac{\eta D^2}{2} \right] + \frac{\sum_{k=0}^{T-1} D_k \langle Z_{k+1}, h_k - x_{k+1} \rangle}{\sum_{k=0}^{T-1} D_k} . \nonumber
    \end{align}
    Next, we substitute the value of $H_T$ from Lemma~\ref{lem_H_upbd} in the preceding relation and change the variable name from $T$ to $t$ to obtain
    \begin{align}
        \epsilon_t \leq \frac{D_t}{\sum_{k=0}^{t-1} D_k} &\left[ 3\sqrt{2} L D^2 + 12 D \sqrt{L \delta t} + 6 D \sqrt{\sum_{i = 0}^{t-1} \|Z_{i+1} - Y_i\|_*^2} + \frac{\sum_{k=0}^{t-1} \|Z_{k+1}\|_*^2}{2 \eta} + \frac{\eta D^2}{2} \right] \nonumber\\
        & \hspace{6cm} +\frac{\sum_{k=0}^{t-1} D_k \langle Z_{k+1}, h_k - x_{k+1} \rangle}{\sum_{k=0}^{t-1} D_k}  , \label{adadist_main2}
    \end{align}
    where for simplicity, we considered $L(\delta) \equiv L$ for any $\delta > 0$. Let us denote
    \begin{align*}
        \tau \equiv \tau(T) := \argmin_{1 \leq t \leq T} \frac{D_t}{\sum_{k=0}^{t-1} D_k}  .
    \end{align*}
    For notational simplicity, we will write $\tau(T)$ as $\tau$ throughout the rest of the proof. Hence, we see
    \begin{align*}
        \min_{1 \leq t \leq T} \frac{D_t}{\sum_{k=0}^{t-1} D_k} = \frac{D_\tau}{\sum_{k=0}^{\tau-1} D_k} .
    \end{align*}
    Note that $1 \leq \tau \leq T$. Then, relation~\eqref{adadist_main2} can be simplified as 
    \begin{align}
        \epsilon_\tau &\leq \left[ \frac{D_\tau}{\sum_{k=0}^{\tau-1} D_k} \right] \Bigg[ 3\sqrt{2} L D^2 + 12 D \sqrt{L \delta T} + 6 D \sqrt{\sum_{i = 0}^{T-1} \|Z_{i+1} - Y_i\|_*^2} \nonumber\\
        & \quad\qquad  + \frac{\sum_{k=0}^{T-1} \|Z_{k+1}\|_*^2}{2 \eta} + \frac{\eta D^2}{2} \Bigg] + \left[\frac{D_\tau}{\sum_{k=0}^{\tau-1} D_k} \right] \left| \frac{\sum_{k=0}^{\tau-1} D_k \langle Z_{k+1}, h_k - x_{k+1} \rangle}{D_\tau} \right|. \label{adadist_main3}
    \end{align}
    For notational simplicity, let us denote 
    \begin{align}
        q_\tau = \frac{D_\tau}{\sum_{k=0}^{\tau-1} D_k}, \quad \rho = \max \left\{1, \frac{D}{\widehat D} \right\}, \quad b_k = \langle Z_{k+1}, h_k - x_{k+1} \rangle , \quad B_t = \sum_{i=0}^{t-1} b_i \text{ with } B_0 = 0, \label{adadist_main4}
    \end{align}
    where the definition of $\rho$ captures both the cases $\widehat D < D$ and $\widehat D \geq D$, and we note that $D_0 = \widehat D \leq D_T$. 
    By Lemma~\ref{lem_adadist_bound} and Assumption~\ref{asum_set}, we obtain
    \begin{align}
        q_\tau \leq \frac{\left( \frac{D_T}{D_0} \right)^{\frac{1}{T}} \log \frac{eD_T}{D_0}}{T} \leq \frac{\left( \rho \right)^{\frac{1}{T}} \log (e \rho)}{T} , \label{q_upbd}
    \end{align}
    Moreover, using the notation in~\eqref{adadist_main4}, for the last summation term on the right-hand side of relation~\eqref{adadist_main3} we have
    \begin{align*}
        \sum_{k=0}^{\tau-1} D_k b_k
        =
        \sum_{k=0}^{\tau-1} D_k \left(B_{k+1}-B_k\right) =
        D_{\tau-1}B_\tau
        -
        \sum_{k=0}^{\tau-2}
        \left(D_{k+1}-D_k\right)B_{k+1} ,
    \end{align*}
    which holds for any $\tau \geq 2$, and $\sum_{k=0}^{\tau-1} D_k b_k = D_0 b_0$
    for $\tau = 1$. Now, using the triangle inequality, and the fact that $(D_k)_{k\geq 0}$ is positive and non-decreasing, we obtain for any $\tau \geq 2$,
    \begin{align}
        \left|
        \frac{\sum_{k=0}^{\tau-1}D_k b_k}{D_\tau}
        \right| &= \frac{D_{\tau - 1}}{D_{\tau}} \left| B_\tau \right| + \frac{\sum_{k=0}^{\tau - 2} (D_{k+1} - D_k) \left|B_{k+1} \right|}{D_\tau} \nonumber\\
        &\leq
        \frac{
        D_{\tau-1}
        +\sum_{k=0}^{\tau-2}(D_{k+1}-D_k)
        }{D_\tau}
        \max_{1\leq t\leq \tau}|B_t| \nonumber\\
        &=
        \frac{2D_{\tau-1}-D_0}{D_\tau}
        \max_{1\leq t\leq \tau}|B_t| \nonumber\\
        &\leq
        2\max_{1\leq t\leq T}|B_t|. \label{adadist_main5}
    \end{align}
    Note that the right-hand side of the preceding relation also gives an upper bound for $\tau = 1$. Hence, substituting relations~\eqref{q_upbd} and \eqref{adadist_main5} back to \eqref{adadist_main3}, we obtain
    \begin{align*}
        \epsilon_\tau \leq \frac{\left( \rho \right)^{\frac{1}{T}} \log (e \rho)}{T} &\Bigg[ 3\sqrt{2} L D^2 + 12 D \sqrt{L \delta T} + 6 D \sqrt{\sum_{i = 0}^{T-1} \|Z_{i+1} - Y_i\|_*^2} \\
        & \hspace{3cm} + \frac{\sum_{k=0}^{T-1} \|Z_{k+1}\|_*^2}{2 \eta} + \frac{\eta D^2}{2} + 2 \max_{1 \leq t \leq T} |B_t| \Bigg] .
    \end{align*}
    Taking total expectation on both sides of the preceding relation, we obtain for any $T \geq 1$,
    \begin{align}
        \EXP{\epsilon_\tau} \leq \frac{\left( \rho \right)^{\frac{1}{T}} \log (e \rho)}{T} &\Bigg[ 3\sqrt{2} L D^2 + 12 D \sqrt{L \delta T} + 6 D \EXP{\sqrt{\sum_{i = 0}^{T-1} \|Z_{i+1} - Y_i\|_*^2} } \nonumber\\
        & \hspace{2cm} + \frac{\EXP{\sum_{k=0}^{T-1} \|Z_{k+1}\|_*^2}}{2 \eta} + \frac{\eta D^2}{2} + 2 \EXP{\max_{1 \leq t \leq T} |B_t|} \Bigg] . \label{adadist_main6}
    \end{align}
    Note that the third term on the right-hand side can be upper bounded using Jensen's inequality, relation~\eqref{reln_stoc_split}, and Assumption~\ref{asum_stochastic} to yield
    \begin{align}
        6 D \EXP{\sqrt{\sum_{i = 0}^{T-1} \|Z_{i+1} - Y_i\|_*^2} } \leq 6 D \sqrt{ \EXP{\sum_{i = 0}^{T-1} \|Z_{i+1} - Y_i\|_*^2}} \leq 6\sqrt{2} \sigma D \sqrt{T} . \label{adadist_main7}
    \end{align}
    The fourth term on the right-hand side of relation~\eqref{adadist_main6} can be upper bounded using Assumption~\ref{asum_stochastic} yielding
    \begin{align}
        \frac{\EXP{\sum_{k=0}^{T-1} \|Z_{k+1}\|_*^2}}{2 \eta} \leq \frac{\sigma^2 T}{2 \eta} . \label{adadist_main8}
    \end{align}
    Next, we focus on the sixth term on the right-hand side of relation~\eqref{adadist_main6}. We have
    \begin{align*}
        \EXP{b_k \mid \cF_k} = \langle \EXP{Z_{k+1} \mid \cF_k}, h_k - x_{k+1} \rangle = 0 .
    \end{align*}
    Next, we show that the sequence $(B_t)_{t = 1}^T$ is a martingale. Using the definition of $B_{t+1}$ from relation~\eqref{adadist_main4}, we have
    \begin{align*}
        \EXP{B_{t+1} \mid \cF_t} = \EXP{B_t+b_t\mid\cF_t} = B_t+\EXP{b_t\mid\cF_t} = B_t.
    \end{align*}
    Hence, the sequence $(B_t)_{t = 1}^T$ is a martingale. Hence, by Lemma~\ref{lem_doob}, we obtain
    \begin{align*}
        \EXP{\left(\max_{1 \leq t \leq T} |B_t| \right)^2} = \EXP{\max_{1 \leq t \leq T} |B_t|^2} \leq 4 \EXP{|B_T|^2} .
    \end{align*}
    Next, for $\EXP{|B_T|^2}$ we have
    \begin{align*}
        \EXP{|B_T|^2} &= \EXP{\left| \sum_{k=0}^{T-1} \la Z_{k+1}, h_k - x_{k+1} \ra \right|^2} \\
        & = \sum_{k=0}^{T-1} \EXP{ \left(\langle Z_{k+1},h_k-x_{k+1} \rangle \right)^2}
        \\
        &\quad + 2\sum_{0\leq i<j\leq T-1} \EXP{\left[ \left\langle Z_{i+1},h_i-x_{i+1} \right\rangle \left\langle Z_{j+1},h_j-x_{j+1} \right\rangle \right] } \\
        & \leq \sum_{k=0}^{T-1} \EXP{\|Z_{k+1}\|^2 \underbrace{\|h_k-x_{k+1} \|^2}_{\leq D^2}} \\
        &\quad + 2\sum_{0\leq i<j\leq T-1} \EXP{\left[ \left\langle Z_{i+1},h_i-x_{i+1} \right\rangle \underbrace{\EXP{\left\langle Z_{j+1},h_j-x_{j+1} \right\rangle \mid \cF_j}}_{=0} \right] } \\
        & \leq \sigma^2 D^2 T ,
    \end{align*}
    where we used Assumptions~\ref{asum_set} and \ref{asum_stochastic}. Hence, using Jensen's inequality, we obtain
    \begin{align}
        \EXP{\max_{1 \leq t \leq T} |B_t|} = \EXP{\sqrt{\left( \max_{1 \leq t \leq T} |B_t| \right)^2}} \leq \sqrt{\EXP{\left( \max_{1 \leq t \leq T} |B_t| \right)^2}} \leq 2 \sigma D \sqrt{T} . \label{adadist_main9}
    \end{align}
    Now substituting relations~\eqref{adadist_main7}, \eqref{adadist_main8}, and \eqref{adadist_main9} back in relation~\eqref{adadist_main6}, we further obtain
    \begin{align}
        \EXP{\epsilon_\tau} \leq \frac{\left( \rho \right)^{\frac{1}{T}} \log (e \rho)}{T} \Bigg[ 3\sqrt{2} L D^2 + 12 D \sqrt{L \delta T} + 6\sqrt{2} \sigma D \sqrt{T} + \frac{\sigma^2 T}{2 \eta} + \frac{\eta D^2}{2} + 4 \sigma D \sqrt{T} \Bigg] . \label{adadist_main10}
    \end{align}
    We minimize the preceding expression with respect to $\eta$ and obtain the optimal selection
    \begin{align*}
        \eta^* = \frac{\sigma\sqrt{T}}{D} .
    \end{align*}
    Substituting $\eta = \eta^*$ in relation~\eqref{adadist_main10}, we obtain the final relation of the theorem.
\end{proof}

It is evident from the results of Theorem~\ref{thm_distada} that Algorithm~\ref{algo_dist_ada_vi} is universal. Moreover, compared to Theorem~\ref{thm_universal}, we are able to remove the constant $\frac{D}{\widehat D}$ directly, which can get large when the estimate is poor $\widehat D << D$. We just have the term $\left( \rho \right)^{1/T}$ which decreases very fast to $1$ as $T$ increases, where $\rho = \max \{1, D/ \widehat D \}$. Moreover, when $T \geq \log \rho$, then $\left( \rho \right)^{1/T} \leq e$ which is a modest constant. Next, we present a corollary to show the convergence rate for H\"older continuous mapping.

\begin{corollary}
\label{cor_adadist_holder}
    Suppose that Assumptions~\ref{asum_set} and~\ref{asum_stochastic}
    hold, and $F$ is a H\"older continuous mapping for some $\nu\in[0,1]$
    and $H_\nu>0$. Then, for any $T \geq 1$, the iterates generated by
    Algorithm~\ref{algo_dist_ada_vi} satisfy
    \begin{align*}
        \EXP{\epsilon_{\tau(T)}}
        &\leq
        \left( \rho \right)^{1/T} \log\left( e \rho \right)
        \left[
        \frac{
        3 \times 2^{1-\nu/2}H_\nu D^{1+\nu}
        }{
        (1-\nu)^{(1-\nu)/2}T^{(1+\nu)/2}
        }
        +
        (6\sqrt{2}+5)\frac{\sigma D}{\sqrt{T}}
        \right],
    \end{align*}
    where $\rho = \max \{1, D/\widehat D \}$ and we use the convention $0^0=1$.
\end{corollary}
\begin{proof}
The proof follows along the similar lines as Corollary~\ref{cor_holder}. We start from the final relation of Theorem~\ref{thm_distada} and substitute the value of $L(\delta)$ from Lemma~\ref{lem_holder_to_delta_L} to the expression to obtain
\begin{align}
\EXP{\epsilon_{\tau(T)}}
&\leq
\rho^{1/T}\log(e\rho)
\Bigg[
\frac{3\sqrt{2}D^2}{T}
\left(\frac{1-\nu}{4}\right)^{\frac{1-\nu}{1+\nu}}
\nu^{\frac{2\nu}{1+\nu}}
H_\nu^{\frac{2}{1+\nu}}
\delta^{-\frac{1-\nu}{1+\nu}}
\nonumber\\
&\hspace{3cm}
+
\frac{12D}{\sqrt{T}}
\left(\frac{1-\nu}{4}\right)^{\frac{1-\nu}{2(1+\nu)}}
\nu^{\frac{\nu}{1+\nu}}
H_\nu^{\frac{1}{1+\nu}}
\delta^{\frac{\nu}{1+\nu}}
+
(6\sqrt{2}+5)\frac{\sigma D}{\sqrt{T}}
\Bigg], \nonumber
\end{align}
where $\rho = \max \left\{1, \frac{D}{\widehat D} \right\}$. For $\nu\in(0,1)$, minimizing the first two terms on the right-hand side over $\delta>0$ yields
\begin{align}
\delta^*
=
\frac{(1-\nu)^{(3+\nu)/2}}
     {2^{(5+\nu)/2}\nu}
\frac{H_\nu D^{1+\nu}}{T^{(1+\nu)/2}}, \nonumber
\end{align}
which when substituted back yields
\begin{align*}
\EXP{\epsilon_{\tau(T)}}
&\leq
\rho^{1/T}\log(e\rho)
\left[
\frac{
(3\sqrt{2})^\nu 6^{1-\nu}
}{
(1-\nu)^{(1-\nu)/2}
}
\frac{H_\nu D^{1+\nu}}{T^{(1+\nu)/2}}
+
(6\sqrt{2}+5)\frac{\sigma D}{\sqrt{T}}
\right]
\\
&=
\rho^{1/T}\log(e\rho)
\left[
\frac{
3\times 2^{1-\nu/2}H_\nu D^{1+\nu}
}{
(1-\nu)^{(1-\nu)/2}T^{(1+\nu)/2}
}
+
(6\sqrt{2}+5)\frac{\sigma D}{\sqrt{T}}
\right].
\end{align*}
For the corner cases, i.e., when $\nu = 0$ and $\nu = 1$, the proof follows similarly to that of Corollary~\ref{cor_holder}, and we obtain the desired relation for every $\nu\in[0,1]$.
\end{proof}


\section{Experimental Details and Additional Experiments}\label{app_experiments}

For all experiments, the diameter-related parameters used by our
methods and DoG are reported in Table~\ref{tab:distance_initialization}. For projected AdaGrad-Norm in the optimization experiments, we set
$b_0=1$ and $\eta=1$, following the notation of \citep{ward2020adagrad}. We also conduct an ablation study to evaluate the sensitivity of all 
experiments to the parameters $\widehat D \in \{0.5,0.7,1,1.5,2,3\}$ for EGM Rules~1--3 and $D_0 \in \{0.01, 0.1, 1\}$ for DA EGM Methods~1 and~2. Only for the bilinear saddle-point problem, we considered $\widehat D \in \{0.05,0.08,0.1,0.5, 1,2\}$ and $D_0 \in \{10^{-10}, 10^{-5}, 10^{-2}, 10^{-1}, 1\}$. Our results show that EGM Rules~1--3 are more sensitive to the choice of $\widehat D$ than DA EGM Methods~1 and~2 are to the choice of $D_0$. This observation is consistent with our theory, as the 
convergence rates for DA EGM have better distance-dependent constants.

For the optimization problems, we report the best-iterate gap $\min_{1\leq j\leq k} \bigl\{f(x_j)-f^*\bigr\} = f(x_k^*)-f^*$, as well as the gap at the corresponding averaged iterate, $f(\widehat{x}_k)-f^*$. We also report the step sizes generated by the different methods. For the stochastic min-max dual formulation of $\ell_2$-regression problem, we report the primal-dual gap at the averaged iterate, together with the corresponding step sizes. For the bilinear saddle-point problem, we present the decrease in the residual function value along with the associated step sizes. For both optimization and VI experiments, we use uniform averaging, i.e., $\widehat{x}_k = \frac{1}{k}\sum_{t=1}^k x_t$, for our EGM Rules~1--3 (cf. Algorithms~\ref{algo_params_free_vi} and \ref{algo_adagd_norm}), UGM, UGM with AdaGrad, projected
AdaGrad-Norm, grid-tuned EGM, and adaptive GEG. Whereas, for our DA EGM Methods~1 and 2 (cf. Algorithm~\ref{algo_dist_ada_vi}), and DoG, we instead use $D_t$-weighted averaging, i.e., $\widehat{x}_k = \frac{\sum_{t=1}^k D_t x_t}{\sum_{t=1}^k D_t}$. For the stochastic min-max problem, the corresponding averaged dual iterate $\widehat{y}_k$ is computed analogously.

\begin{table*}[h]
\caption{Distance parameters used in the numerical experiments, which can be chosen arbitrarily. DoG is not included in the VI experiments.}
\label{tab:distance_initialization}
\centering
\small
\setlength{\tabcolsep}{10pt}
\begin{tabular}{lccc}
\toprule
\textbf{Problem}
&
\textbf{Estimate $\widehat D$}
&
\textbf{DoG $r_\epsilon$}
&
\textbf{Algorithm~3 $D_0$}
\\
\midrule
Softmax
&
$2$
&
$0.2$
&
$0.2$
\\
Logistic regression
&
$2$
&
$0.2$
&
$0.2$
\\
H\"older-smooth
&
$2$
&
$0.2$
&
$0.2$
\\
Min-max dual of constrained $\ell_2$-regression
&
$2\sqrt{2}$
&
Not used
&
$0.2\sqrt{2}$
\\
Bilinear saddle-point
&
$0.1$
&
Not used
&
$0.1$
\\
\bottomrule
\end{tabular}
\end{table*}

Next, we present additional experiments on Softmax problem and Logistic-regression problem. Towards the end, we also provide some of the additional experimental details on H\"older-smooth polyhedral feasibility problem and the VI problems which are also studied in Section~\ref{sec_experiments}.


\begin{figure*}[htbp]
\centering


\begin{subfigure}{0.32\textwidth}
    \centering
    \includegraphics[width=\linewidth]
    {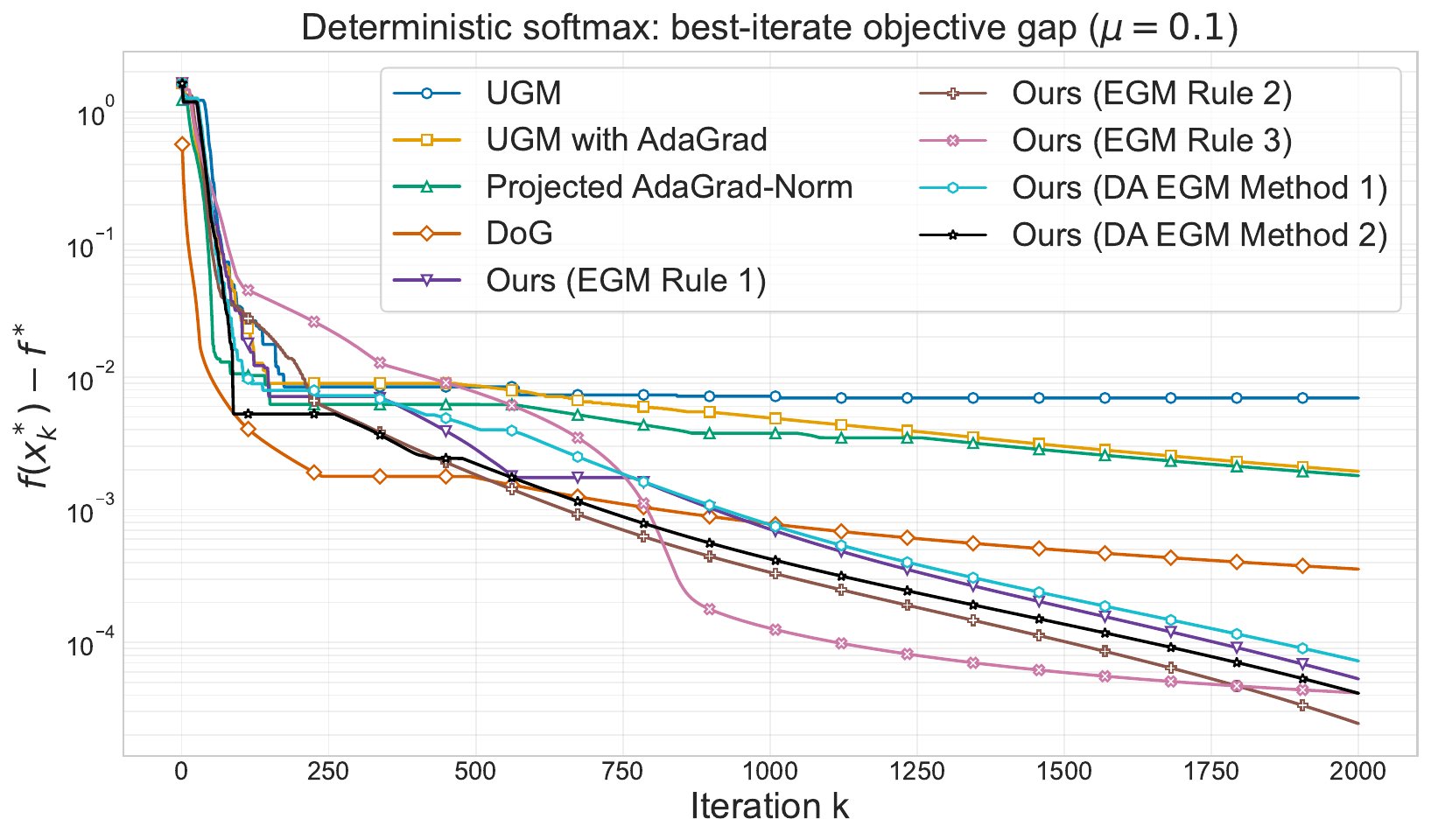}
\end{subfigure}
\hfill
\begin{subfigure}{0.32\textwidth}
    \centering
    \includegraphics[width=\linewidth]
    {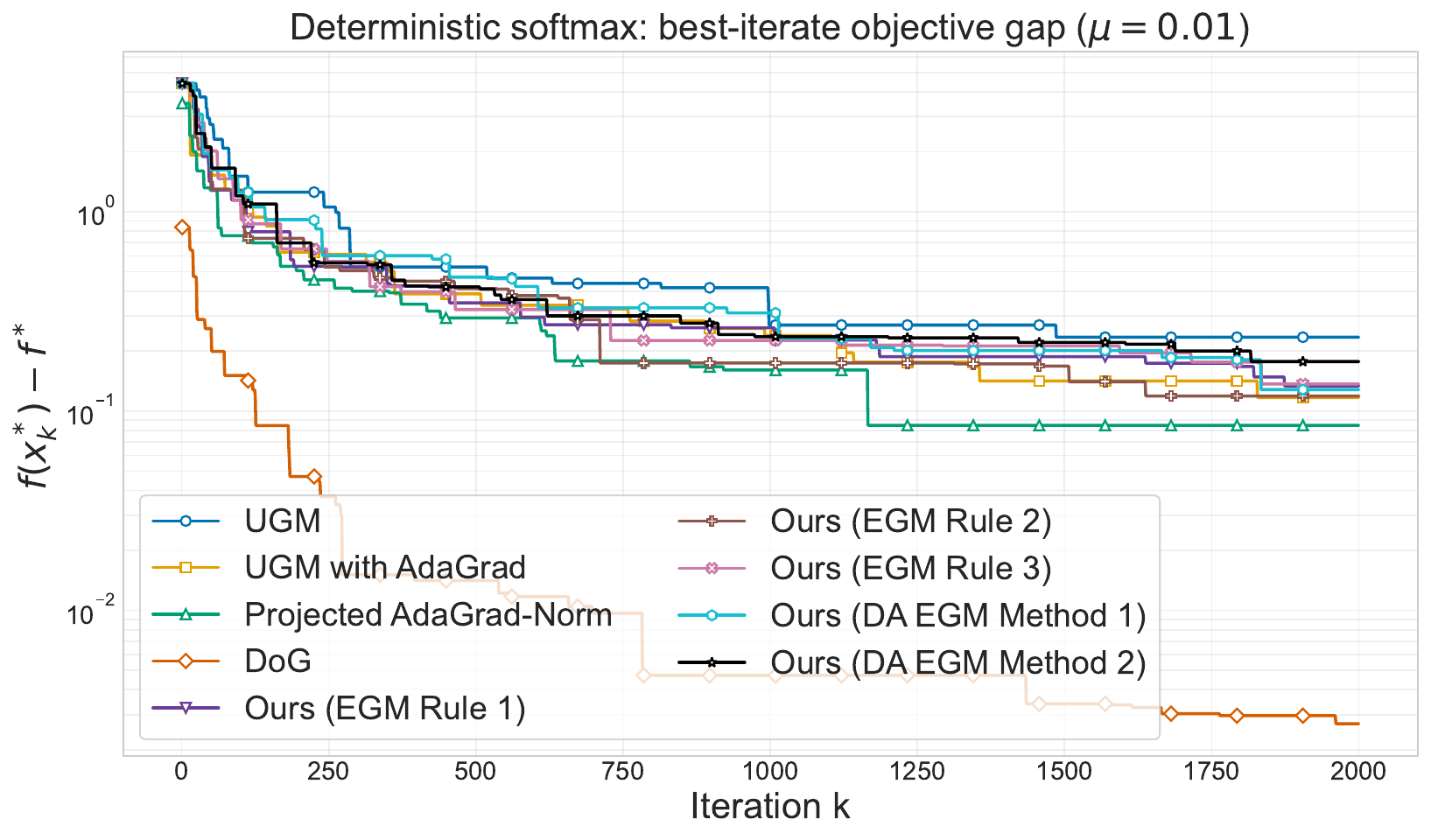}
\end{subfigure}
\hfill
\begin{subfigure}{0.32\textwidth}
    \centering
    \includegraphics[width=\linewidth]
    {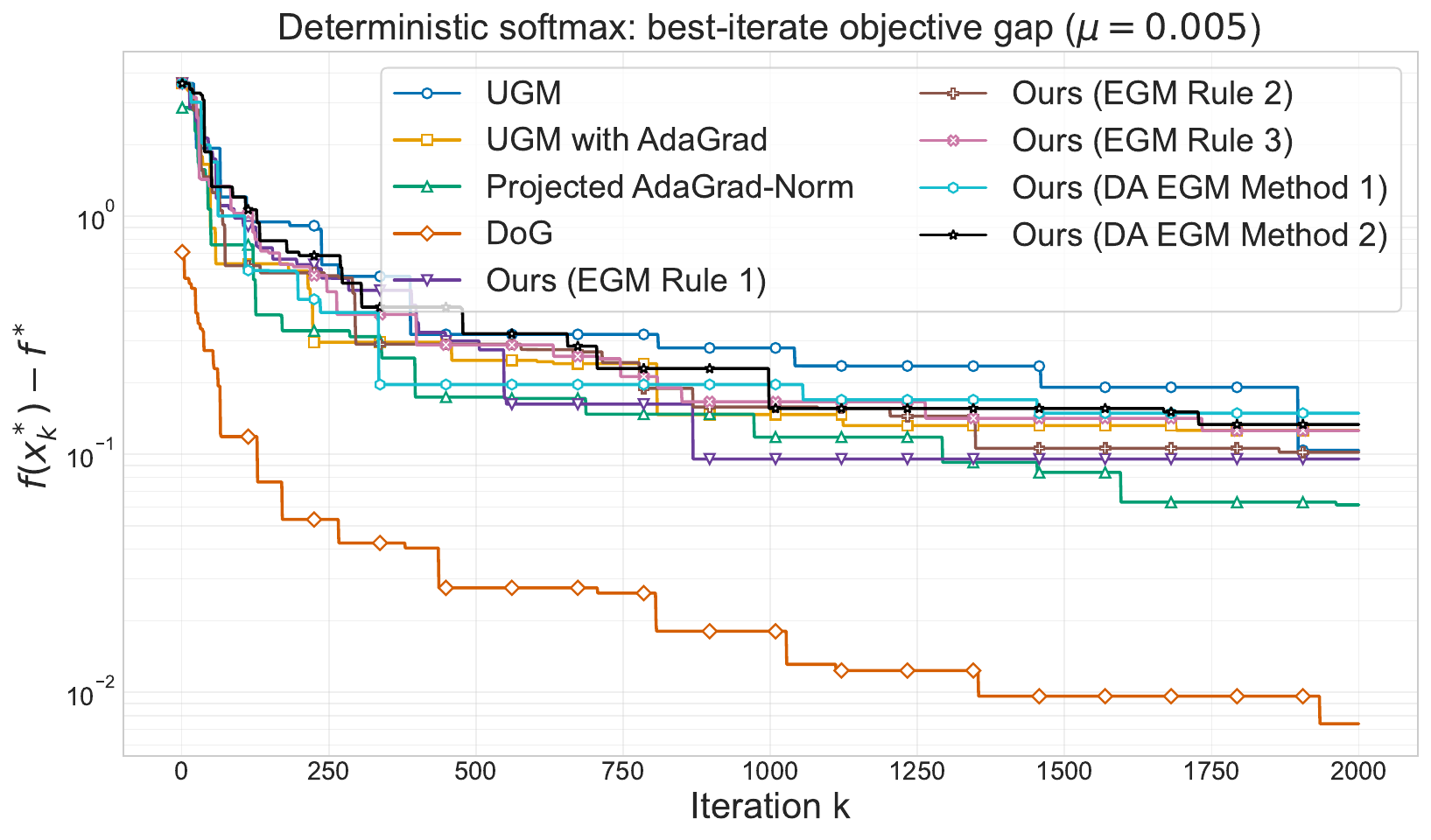}
\end{subfigure}
\par

\begin{subfigure}{0.32\textwidth}
    \centering
    \includegraphics[width=\linewidth]
    {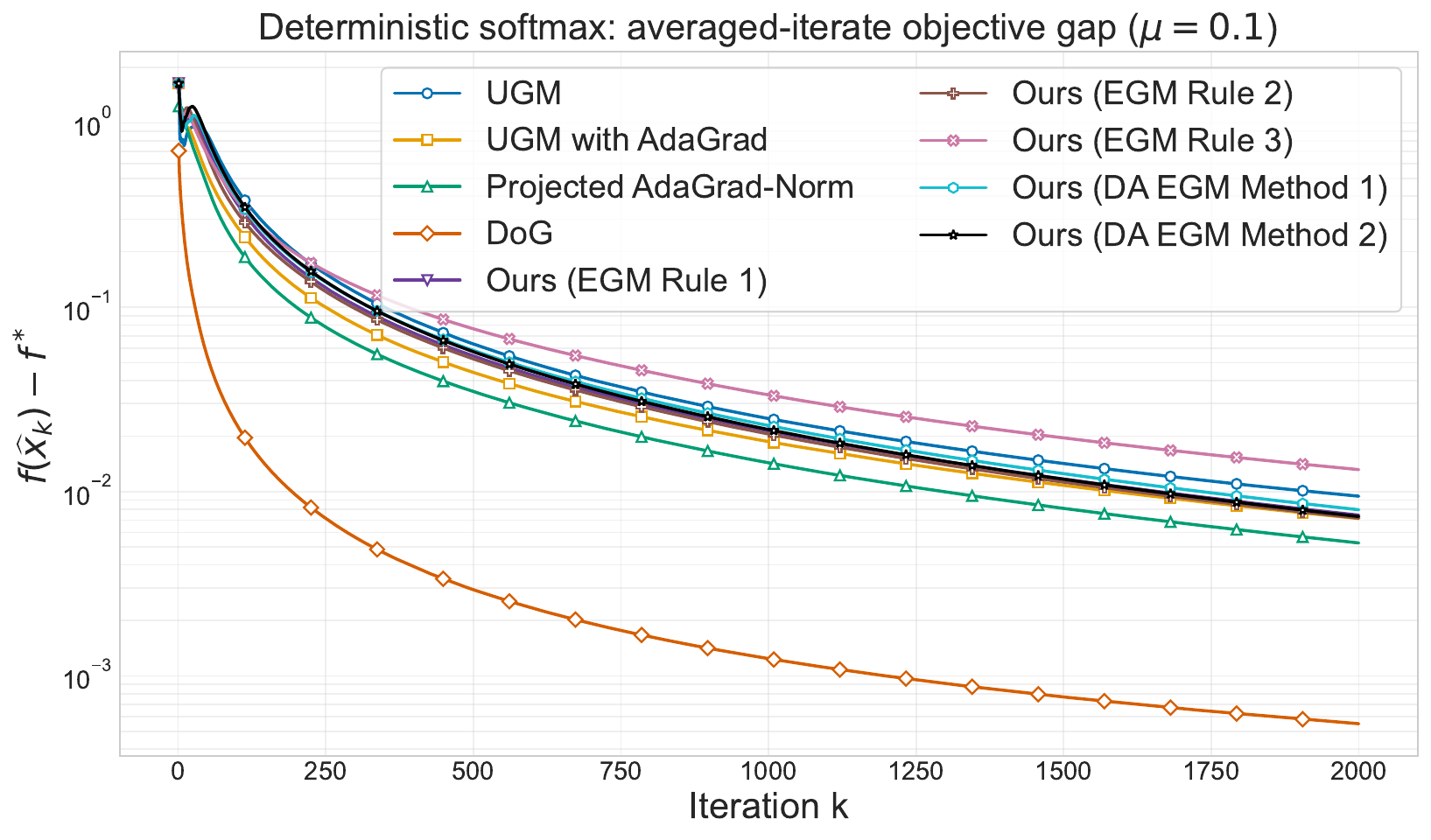}
\end{subfigure}
\hfill
\begin{subfigure}{0.32\textwidth}
    \centering
    \includegraphics[width=\linewidth]
    {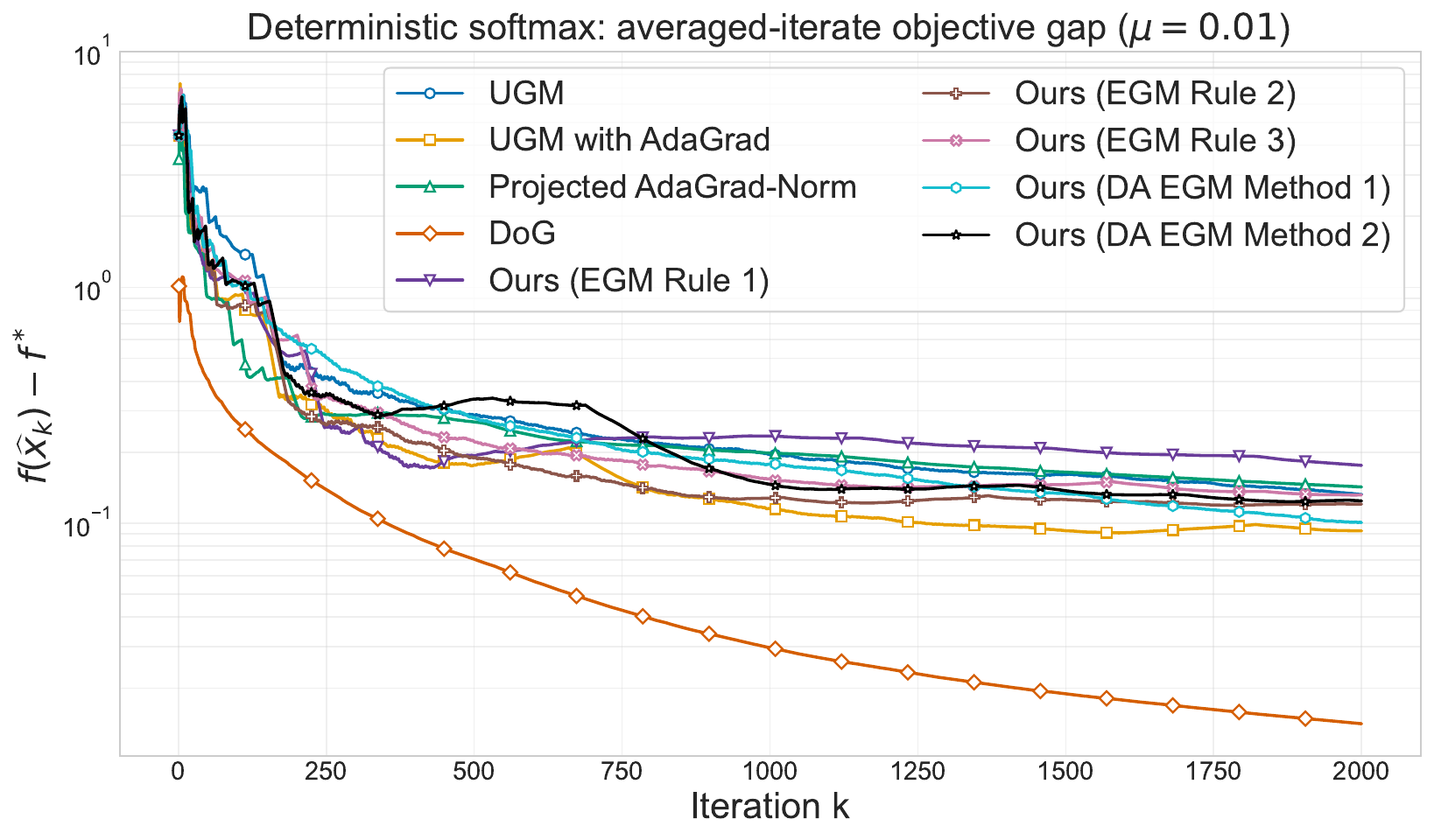}
\end{subfigure}
\hfill
\begin{subfigure}{0.32\textwidth}
    \centering
    \includegraphics[width=\linewidth]
    {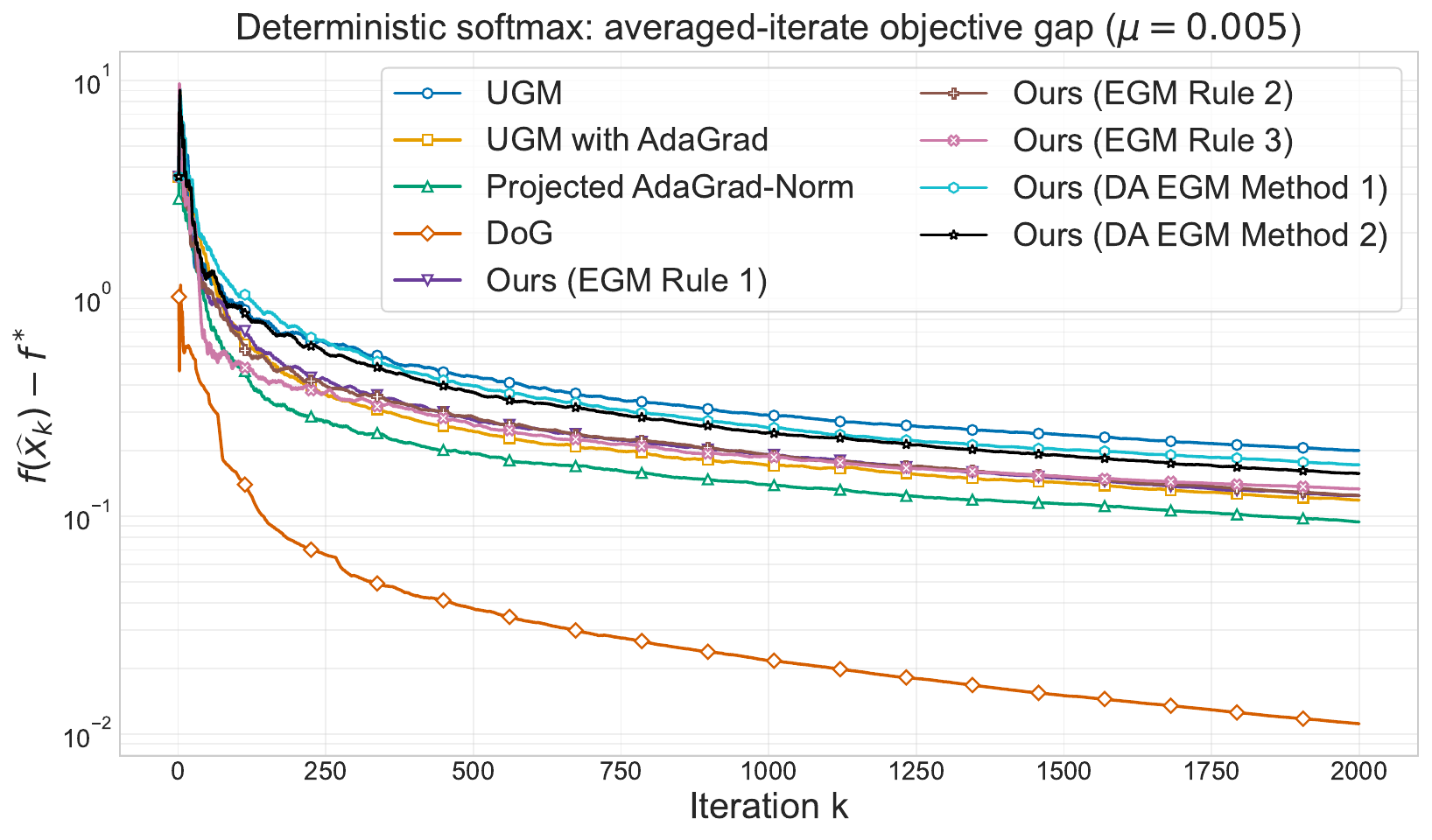}
\end{subfigure}
\par

\begin{subfigure}{0.32\textwidth}
    \centering
    \includegraphics[width=\linewidth]
    {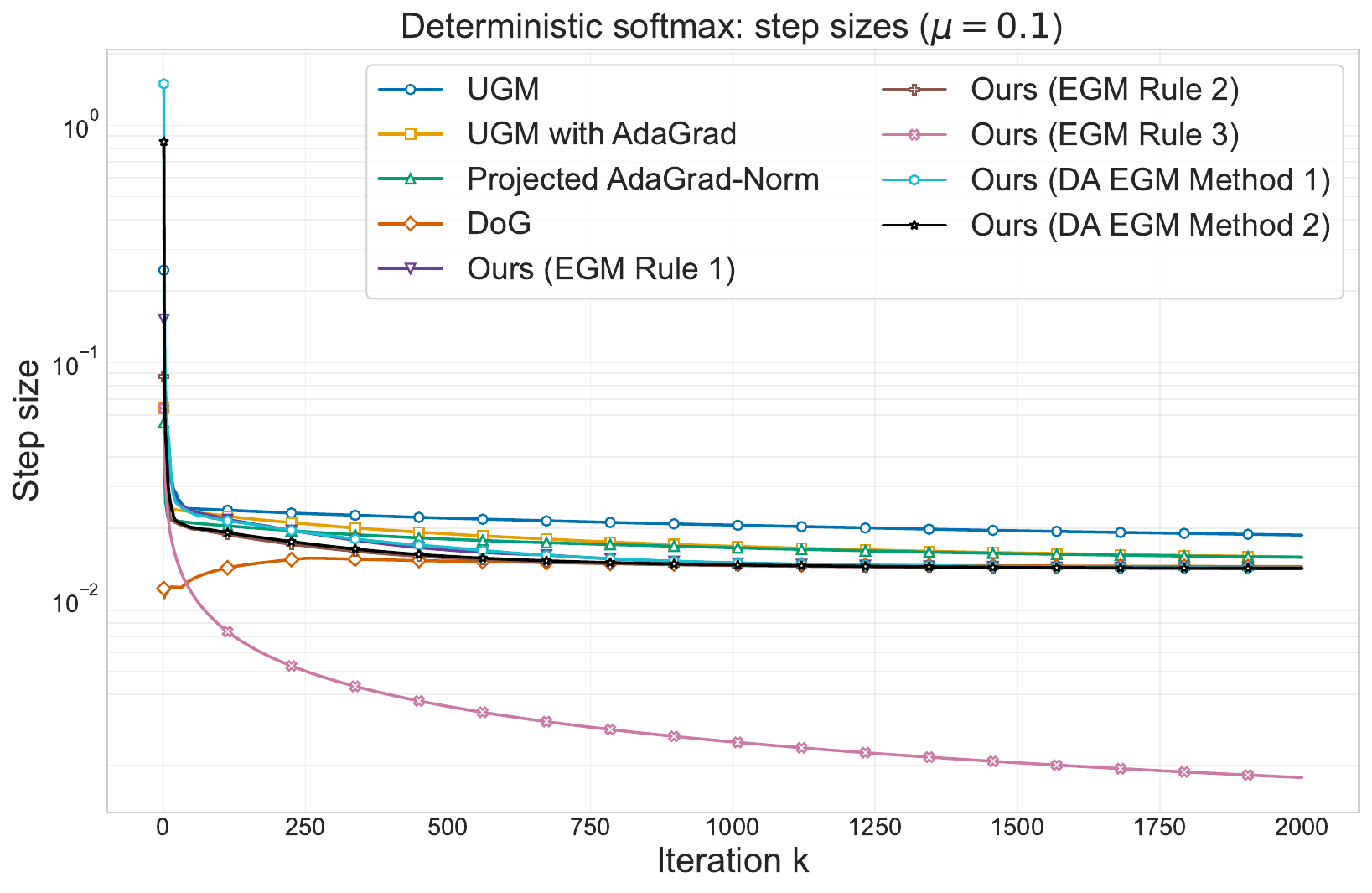}
\end{subfigure}
\hfill
\begin{subfigure}{0.32\textwidth}
    \centering
    \includegraphics[width=\linewidth]
    {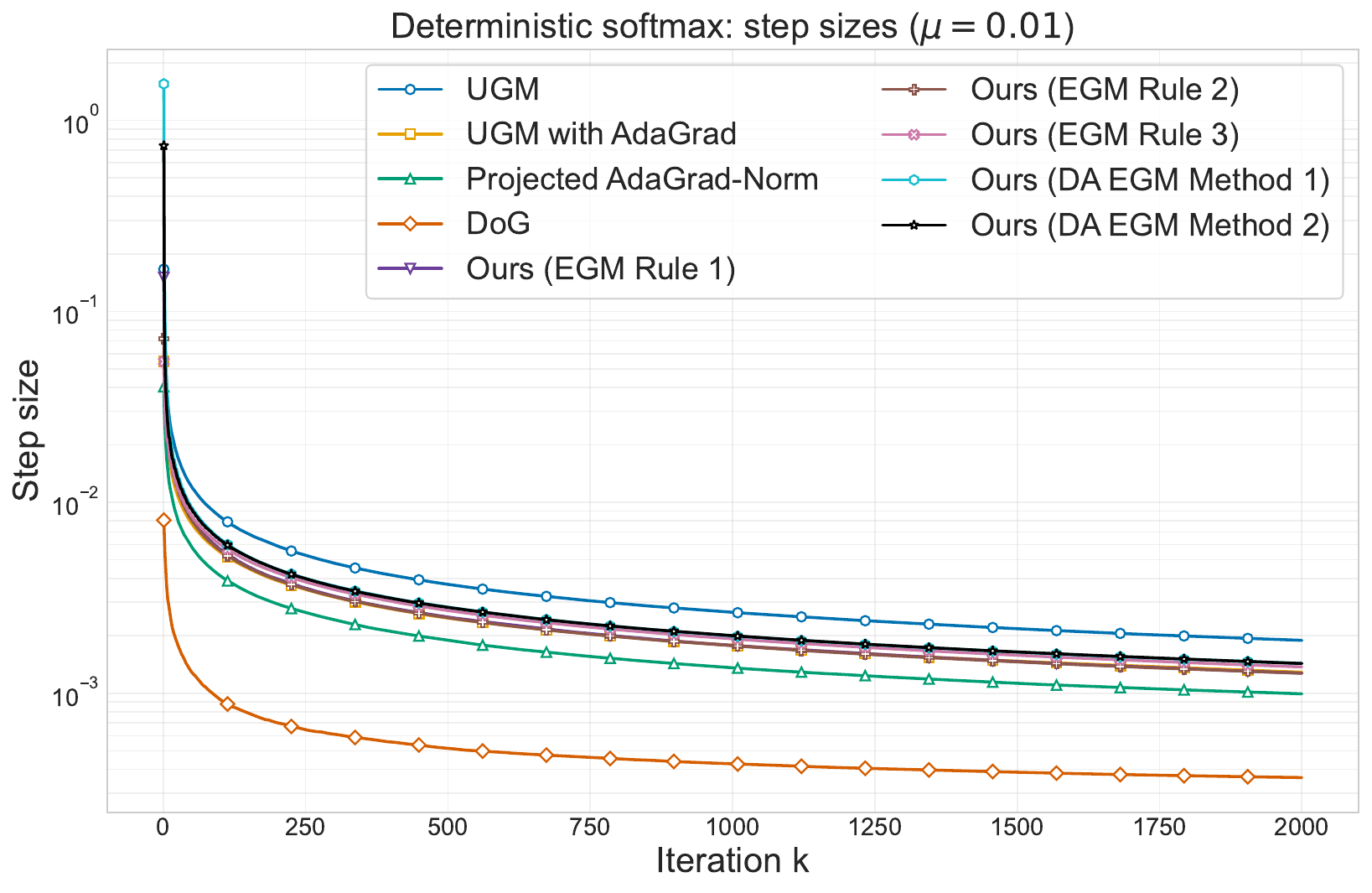}
\end{subfigure}
\hfill
\begin{subfigure}{0.32\textwidth}
    \centering
    \includegraphics[width=\linewidth]
    {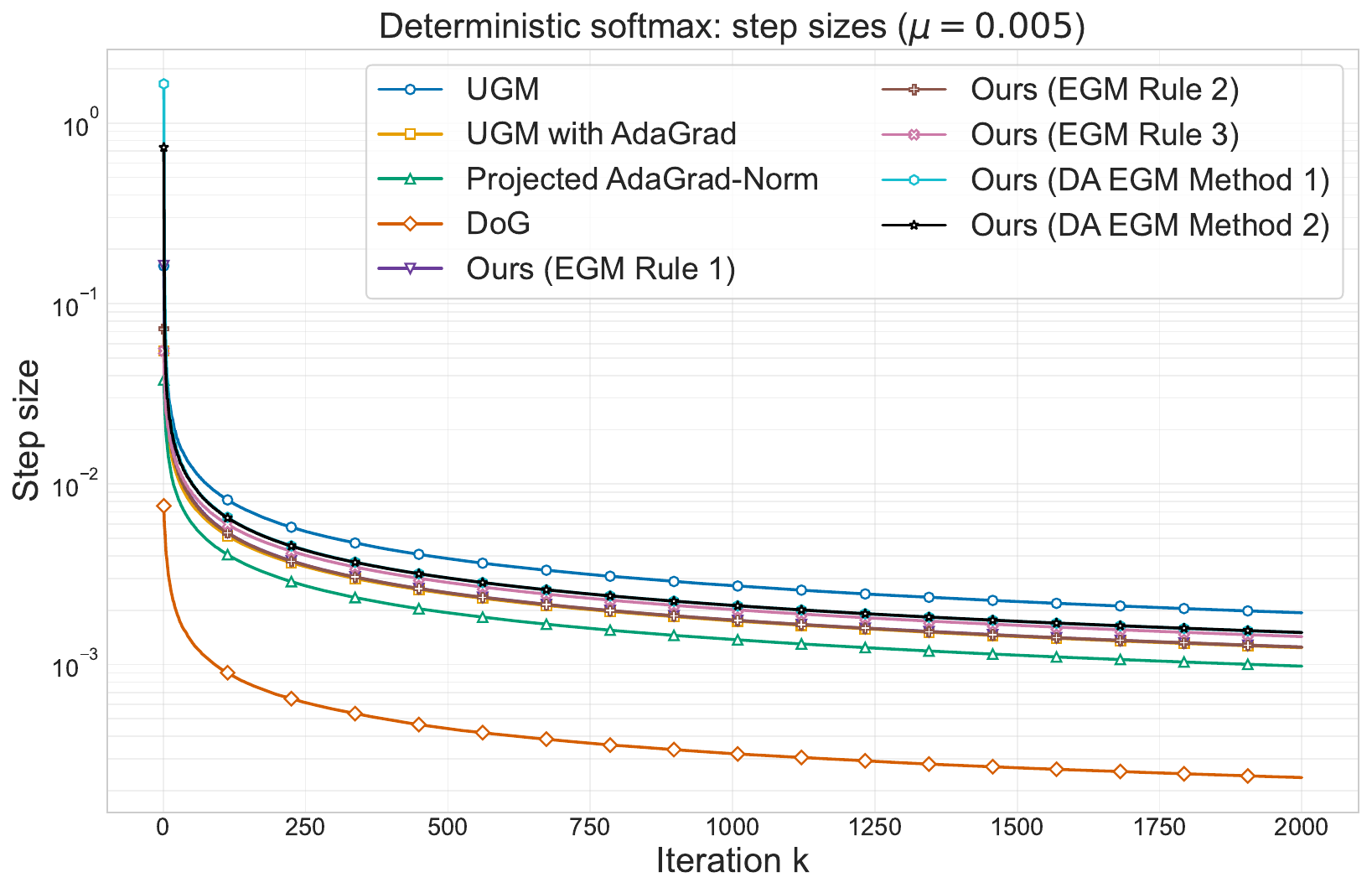}
\end{subfigure}

\caption{Results for the deterministic softmax problem. From top to
bottom, the rows show the best-iterate objective function gaps, averaged-iterate objective function gaps, and step sizes.}
\label{fig_softmax_results}
\end{figure*}




\begin{figure*}[t]
\centering

\begin{subfigure}{0.32\textwidth}
    \centering
    \includegraphics[width=\linewidth]
    {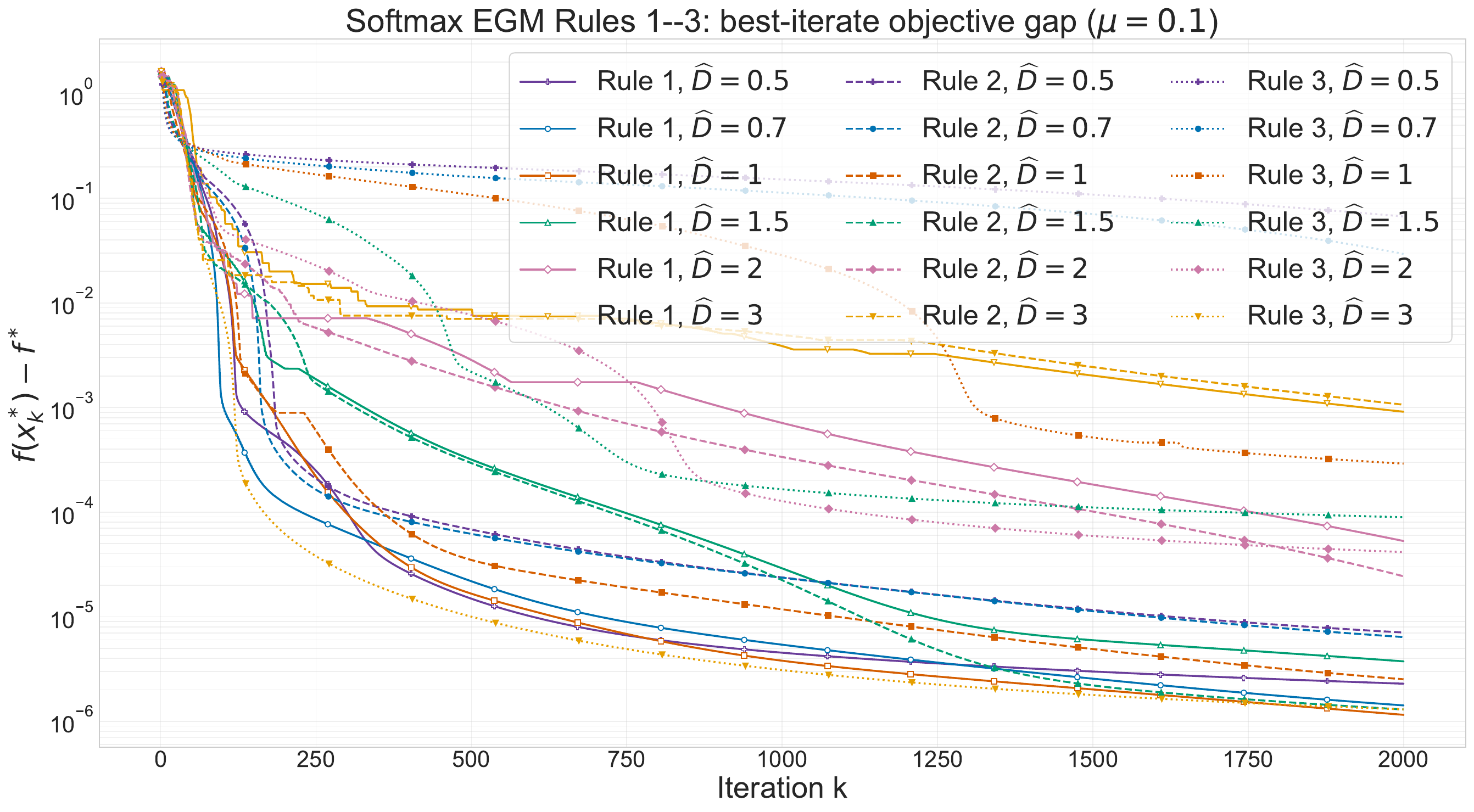}
\end{subfigure}
\hfill
\begin{subfigure}{0.32\textwidth}
    \centering
    \includegraphics[width=\linewidth]
    {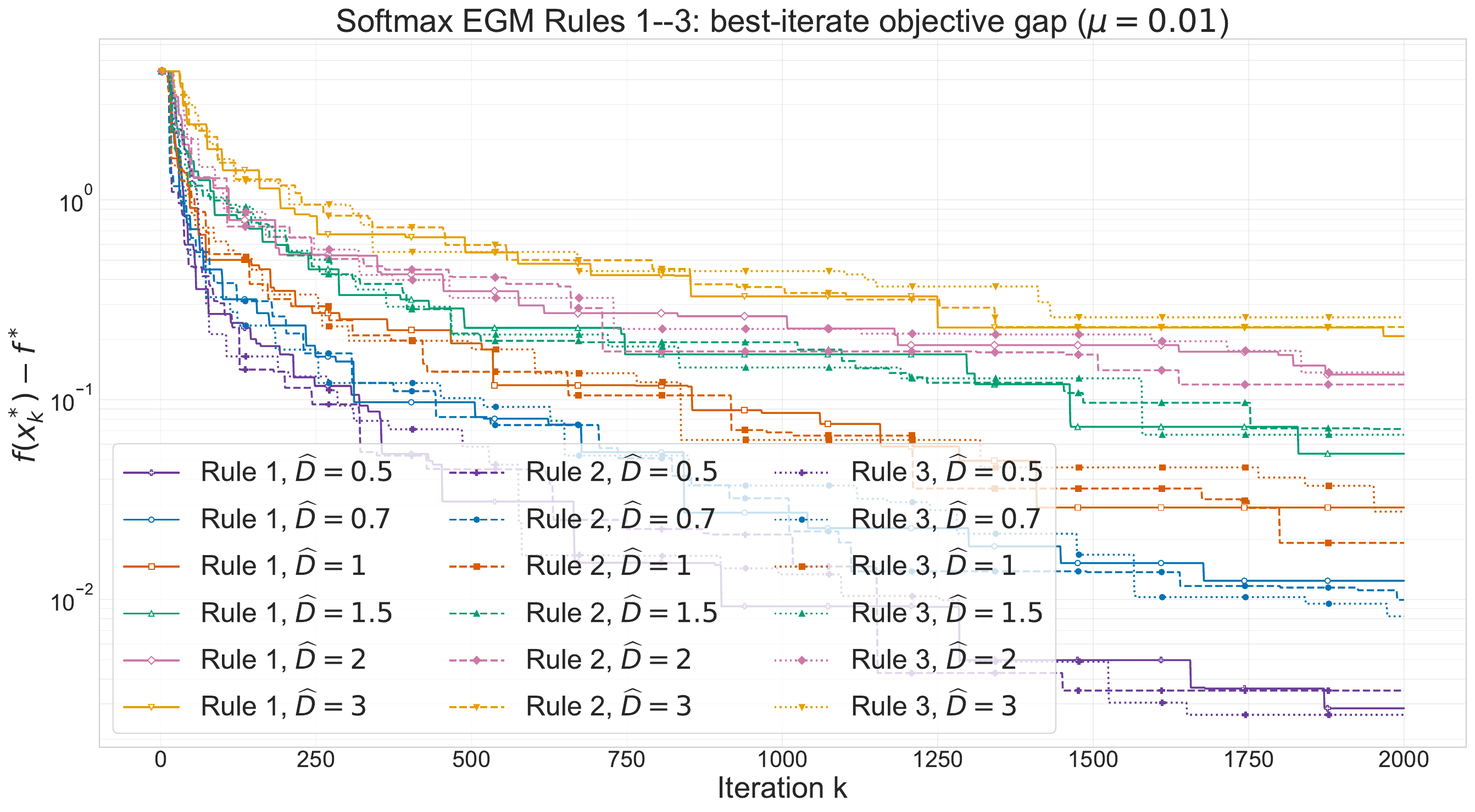}
\end{subfigure}
\hfill
\begin{subfigure}{0.32\textwidth}
    \centering
    \includegraphics[width=\linewidth]
    {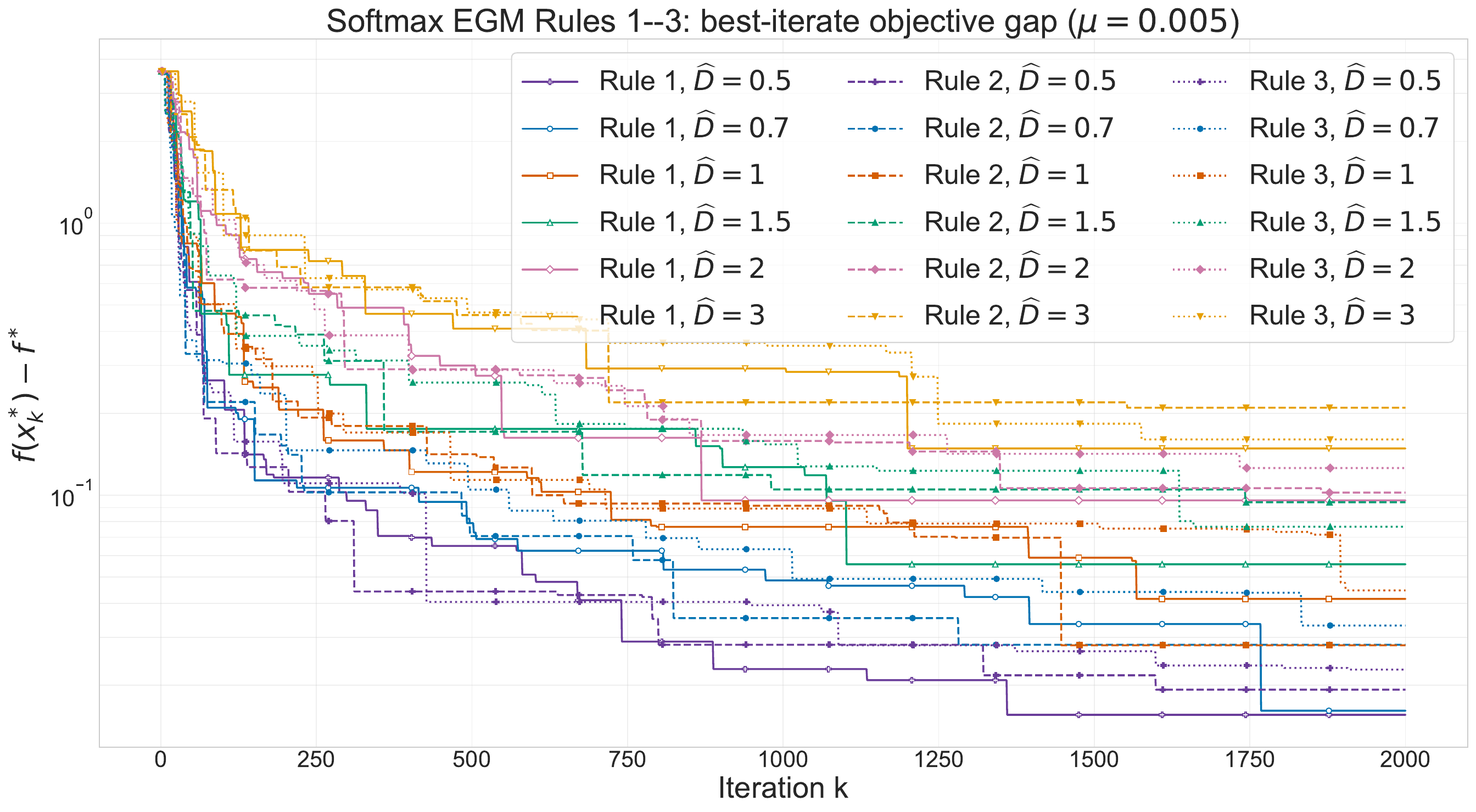}
\end{subfigure}
\par

\begin{subfigure}{0.32\textwidth}
    \centering
    \includegraphics[width=\linewidth]
    {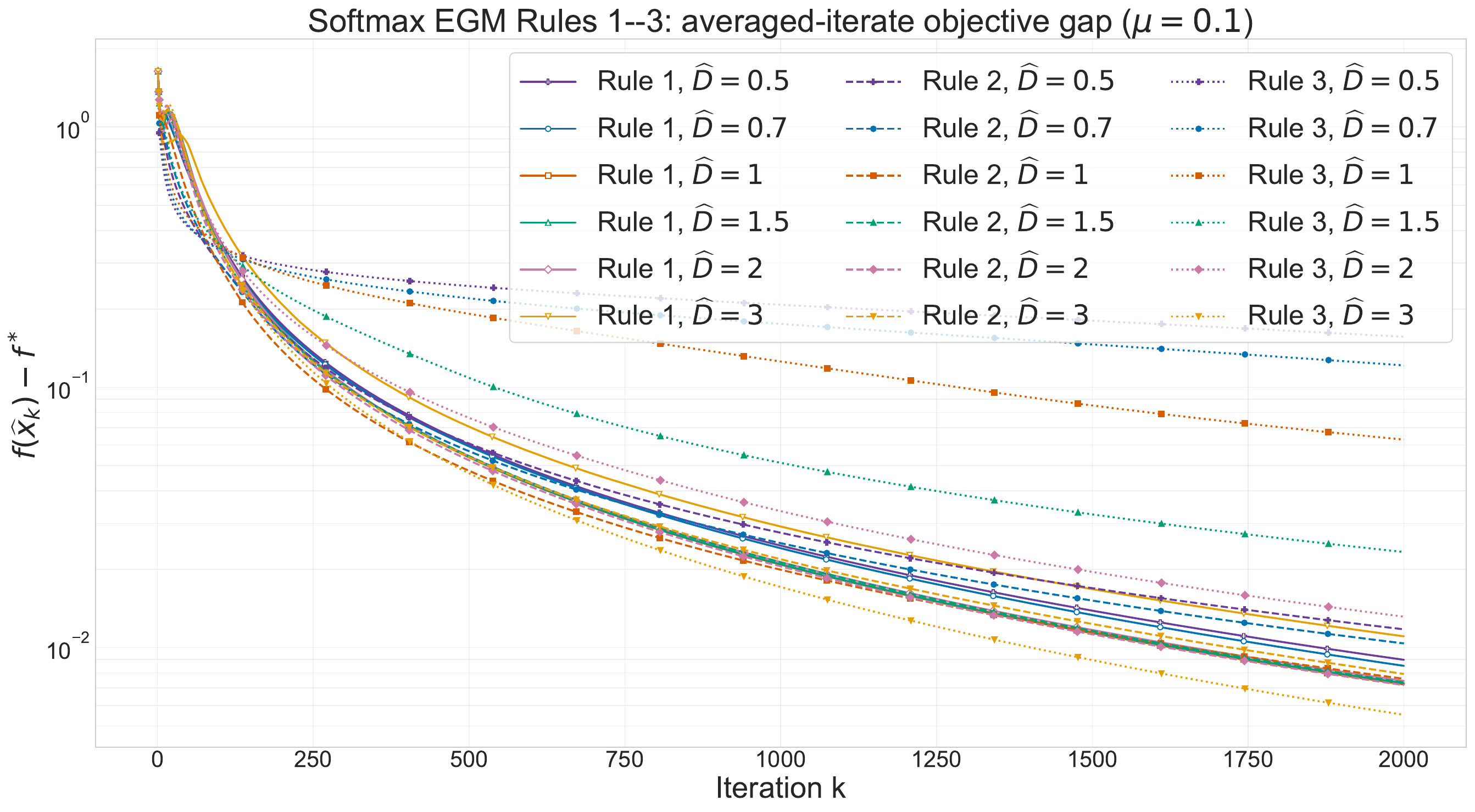}
\end{subfigure}
\hfill
\begin{subfigure}{0.32\textwidth}
    \centering
    \includegraphics[width=\linewidth]
    {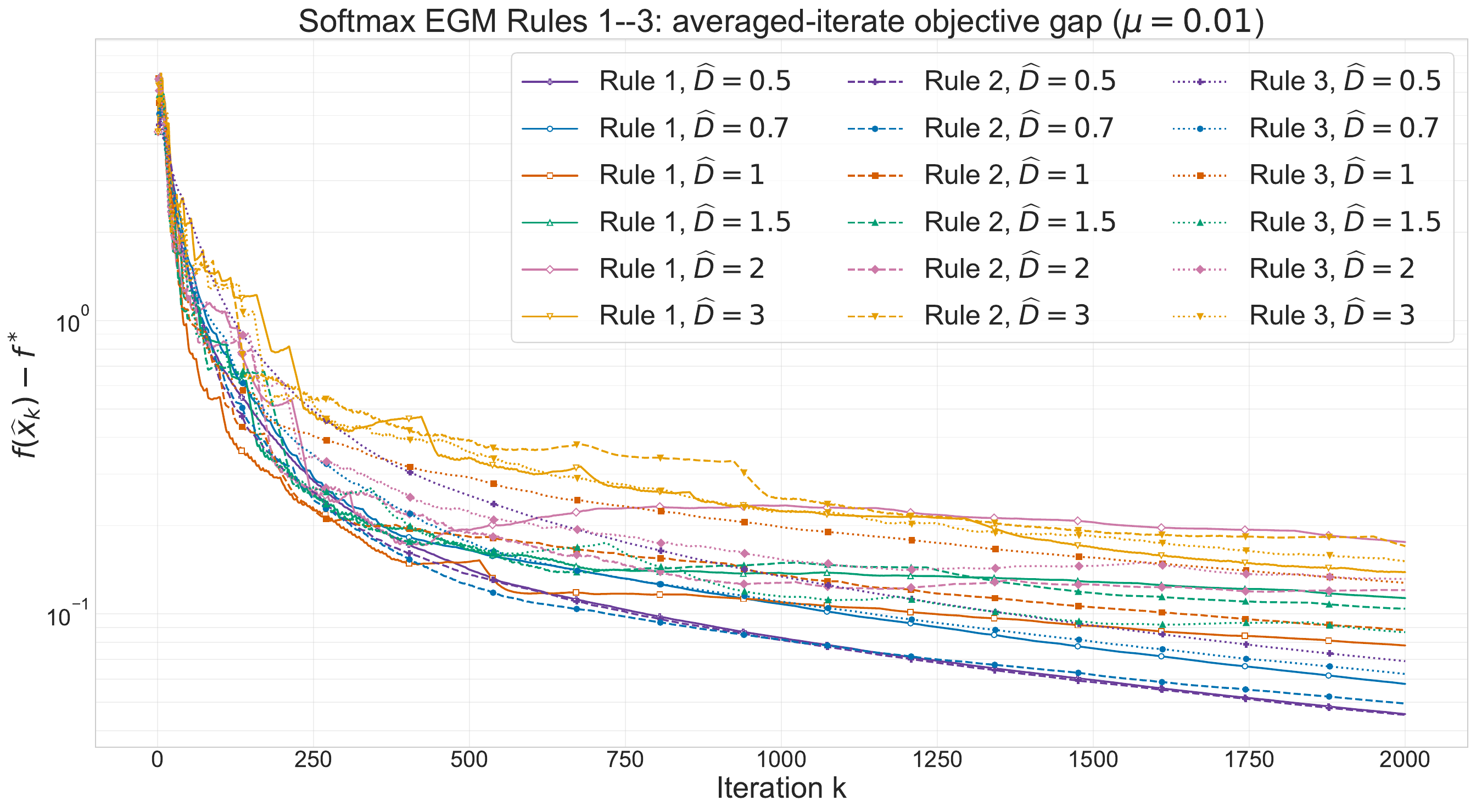}
\end{subfigure}
\hfill
\begin{subfigure}{0.32\textwidth}
    \centering
    \includegraphics[width=\linewidth]
    {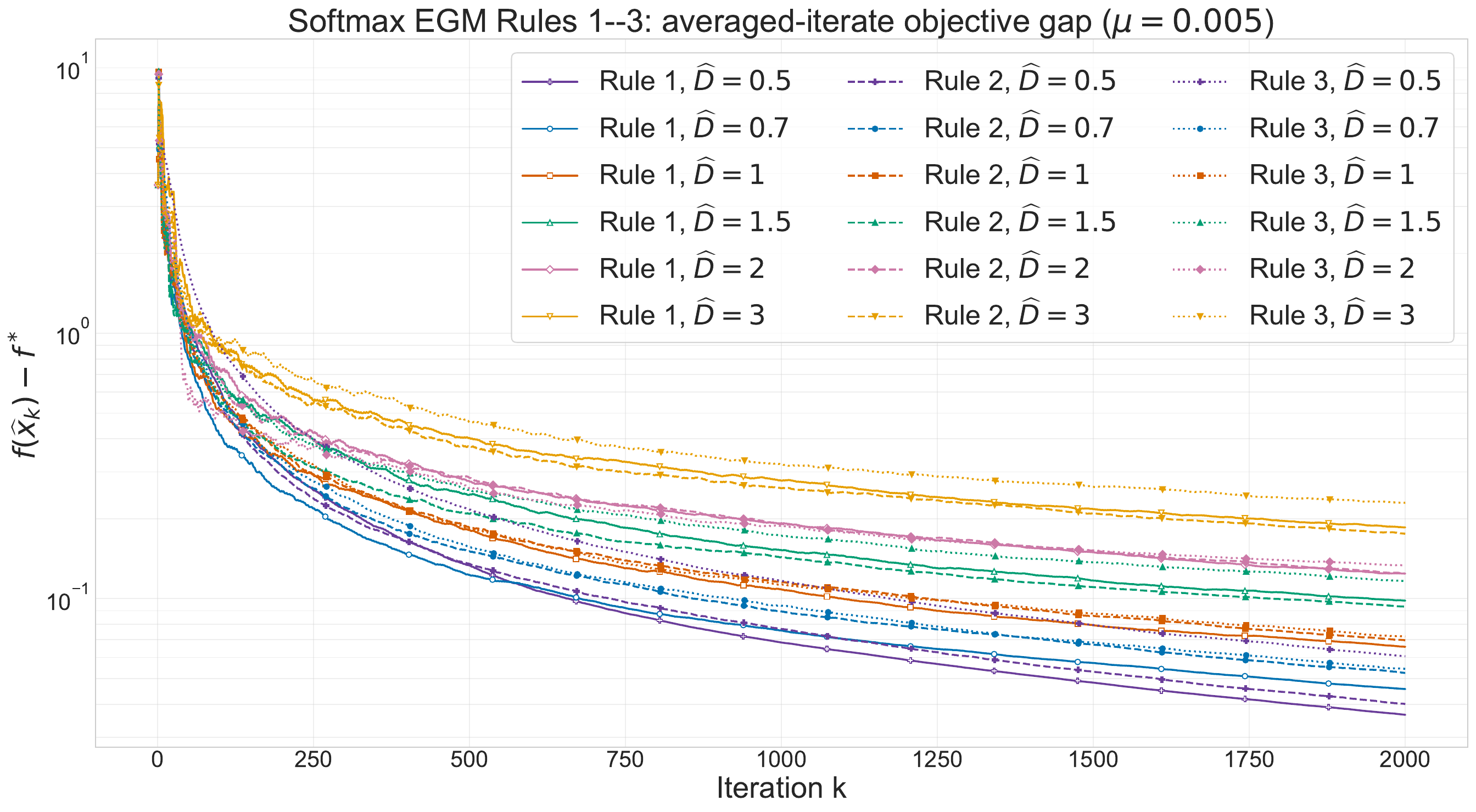}
\end{subfigure}
\par

\begin{subfigure}{0.32\textwidth}
    \centering
    \includegraphics[width=\linewidth]
    {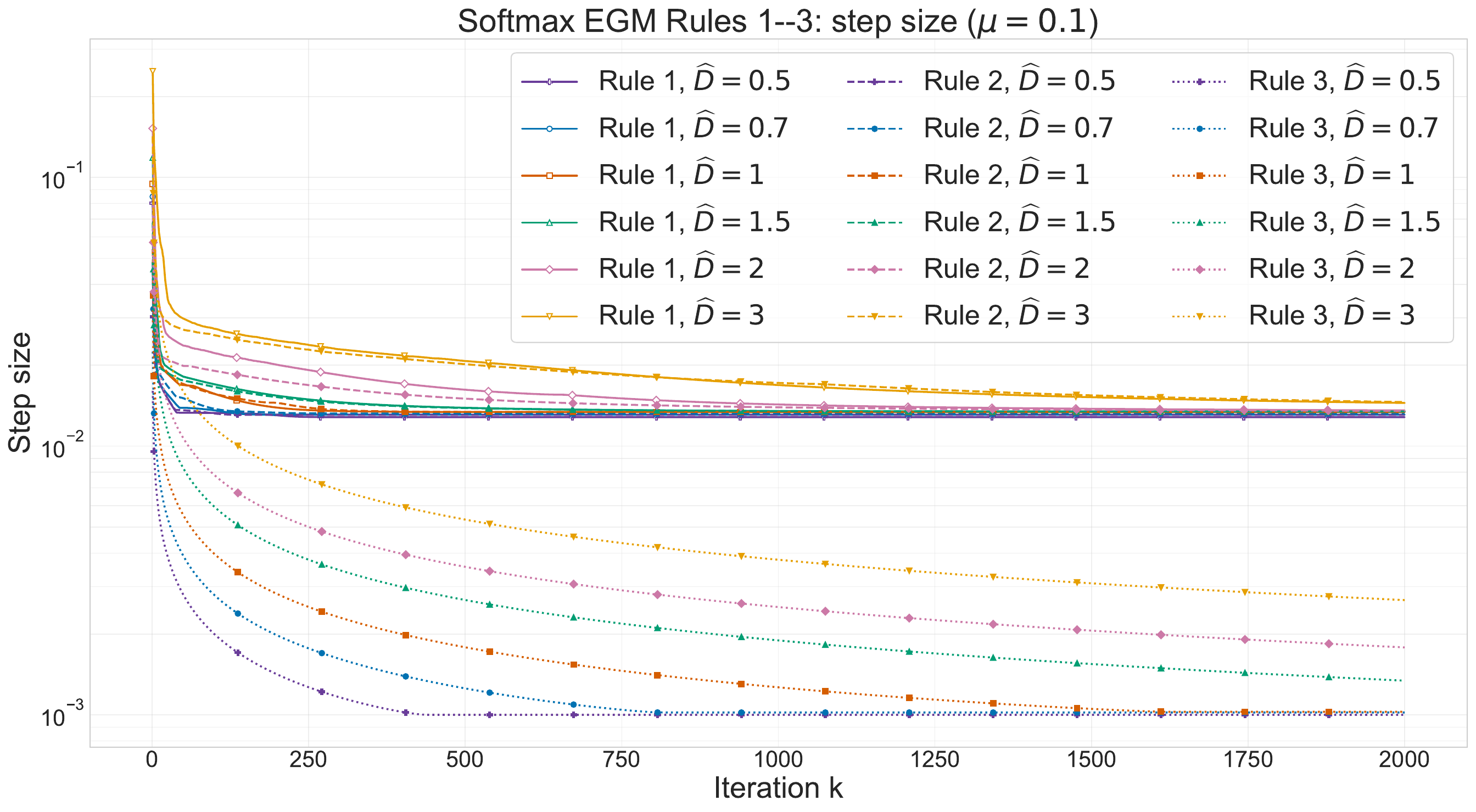}
\end{subfigure}
\hfill
\begin{subfigure}{0.32\textwidth}
    \centering
    \includegraphics[width=\linewidth]
    {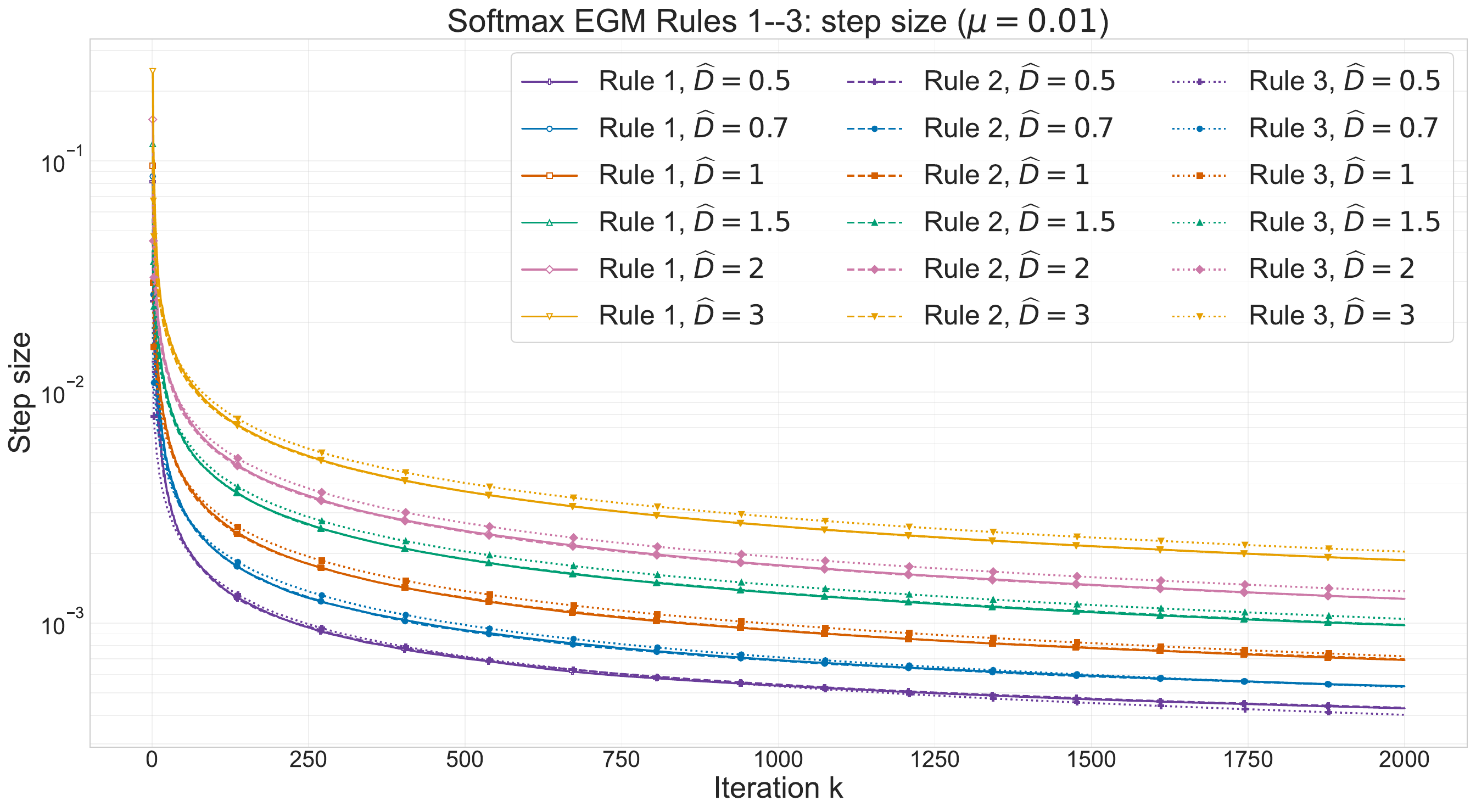}
\end{subfigure}
\hfill
\begin{subfigure}{0.32\textwidth}
    \centering
    \includegraphics[width=\linewidth]
    {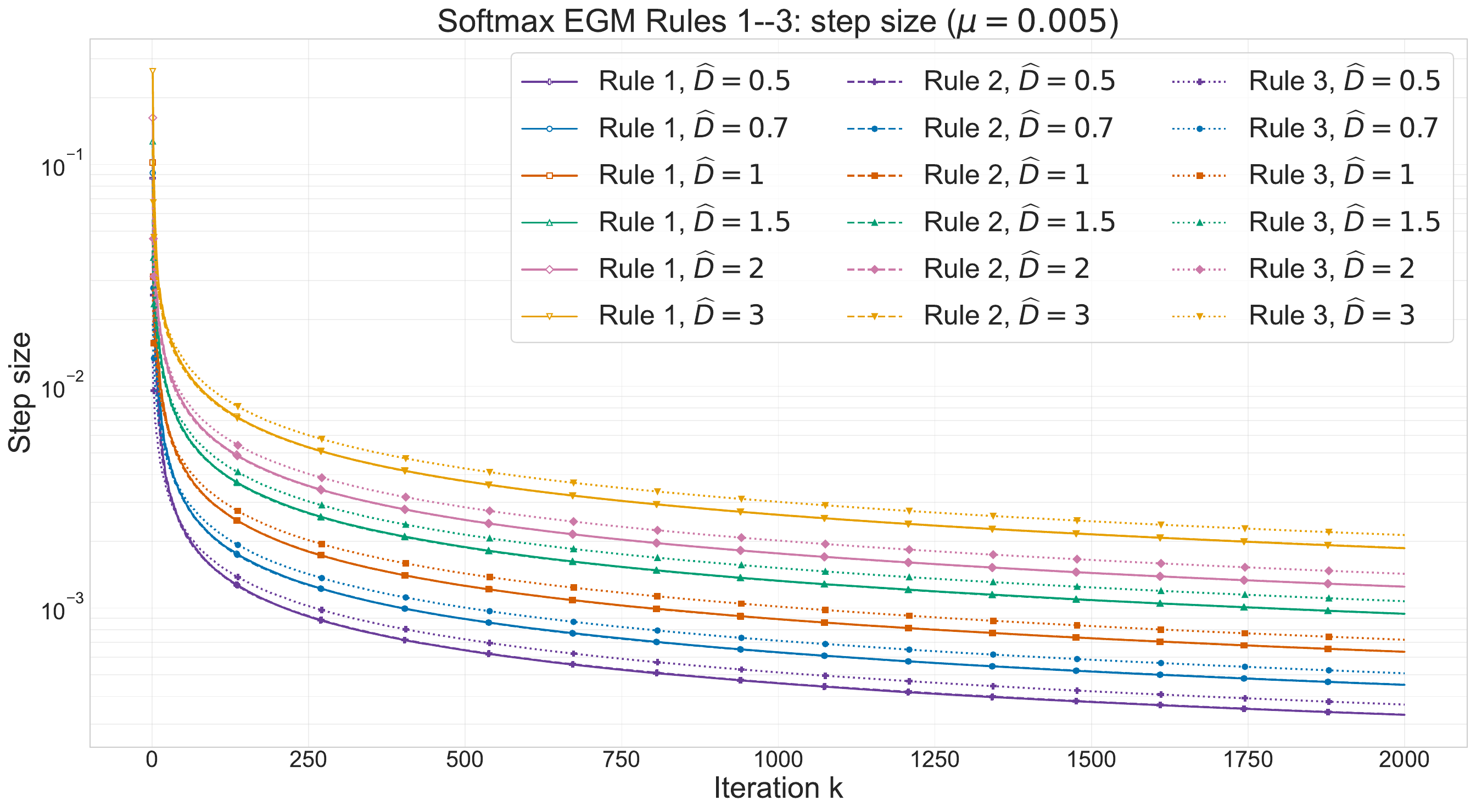}
\end{subfigure}

\caption{Sensitivity of EGM Rules 1--3 to the choice of $\widehat D$
for the deterministic softmax problem. From left to right, the columns
correspond to $\mu=0.1$, $\mu=0.01$, and $\mu=0.005$. From top to
bottom, the rows show the best-iterate objective function gaps,
averaged-iterate objective function gaps, and step sizes.}
\label{fig:softmax_D_hat_sensitivity}
\end{figure*}


\begin{figure*}[htbp]
\centering

\begin{subfigure}{0.32\textwidth}
    \centering
    \includegraphics[width=\linewidth]
    {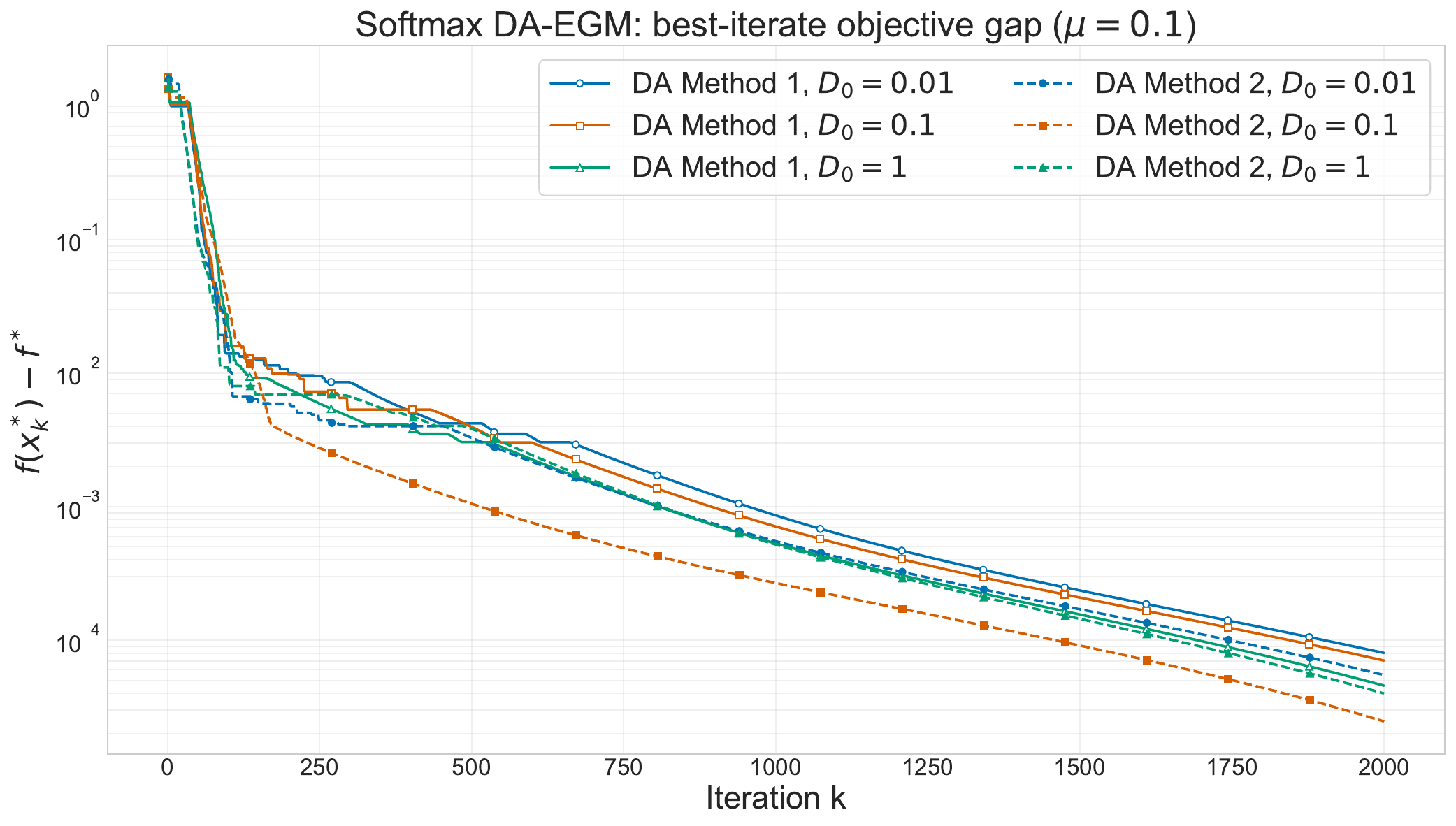}
\end{subfigure}
\hfill
\begin{subfigure}{0.32\textwidth}
    \centering
    \includegraphics[width=\linewidth]
    {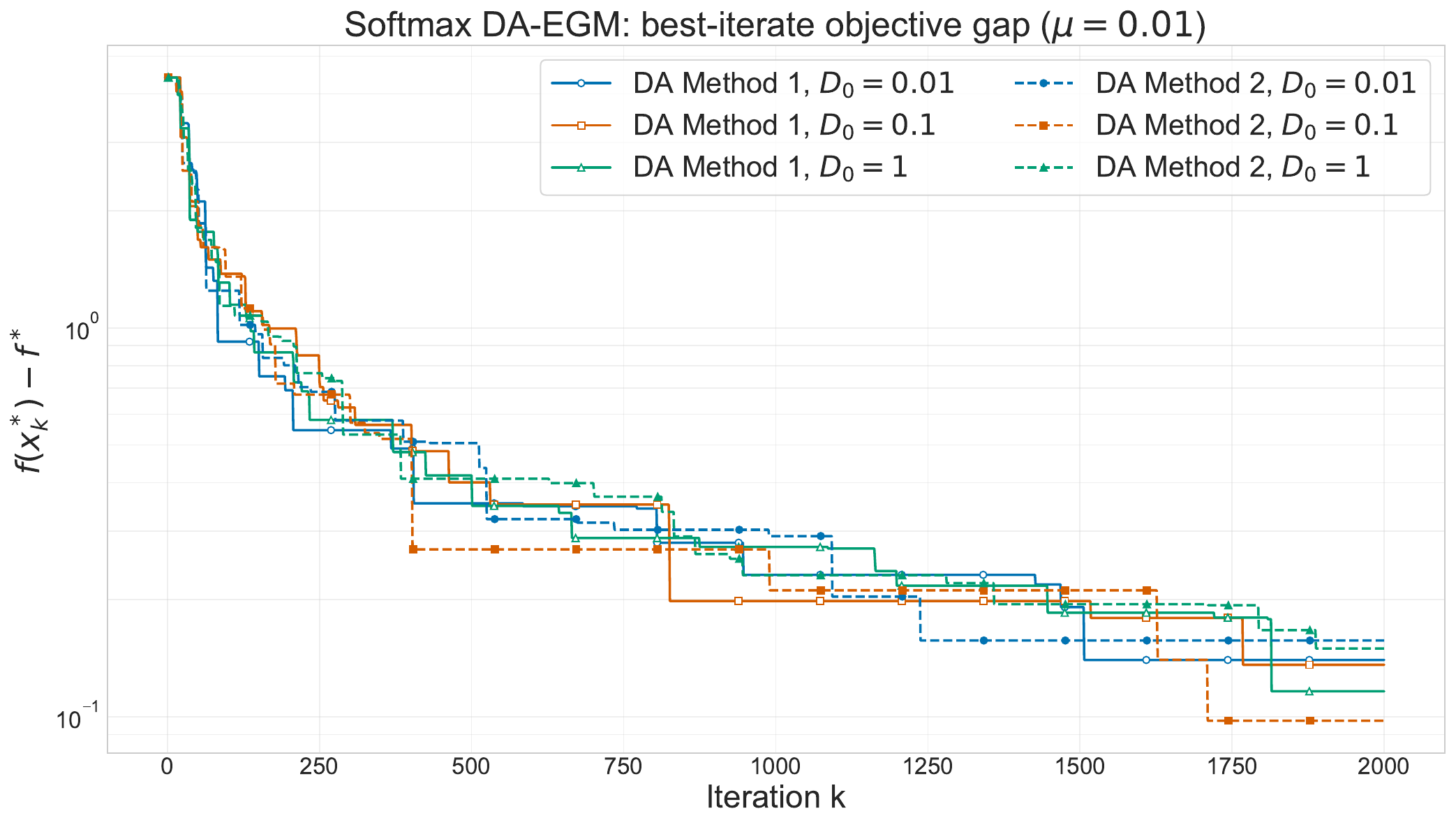}
\end{subfigure}
\hfill
\begin{subfigure}{0.32\textwidth}
    \centering
    \includegraphics[width=\linewidth]
    {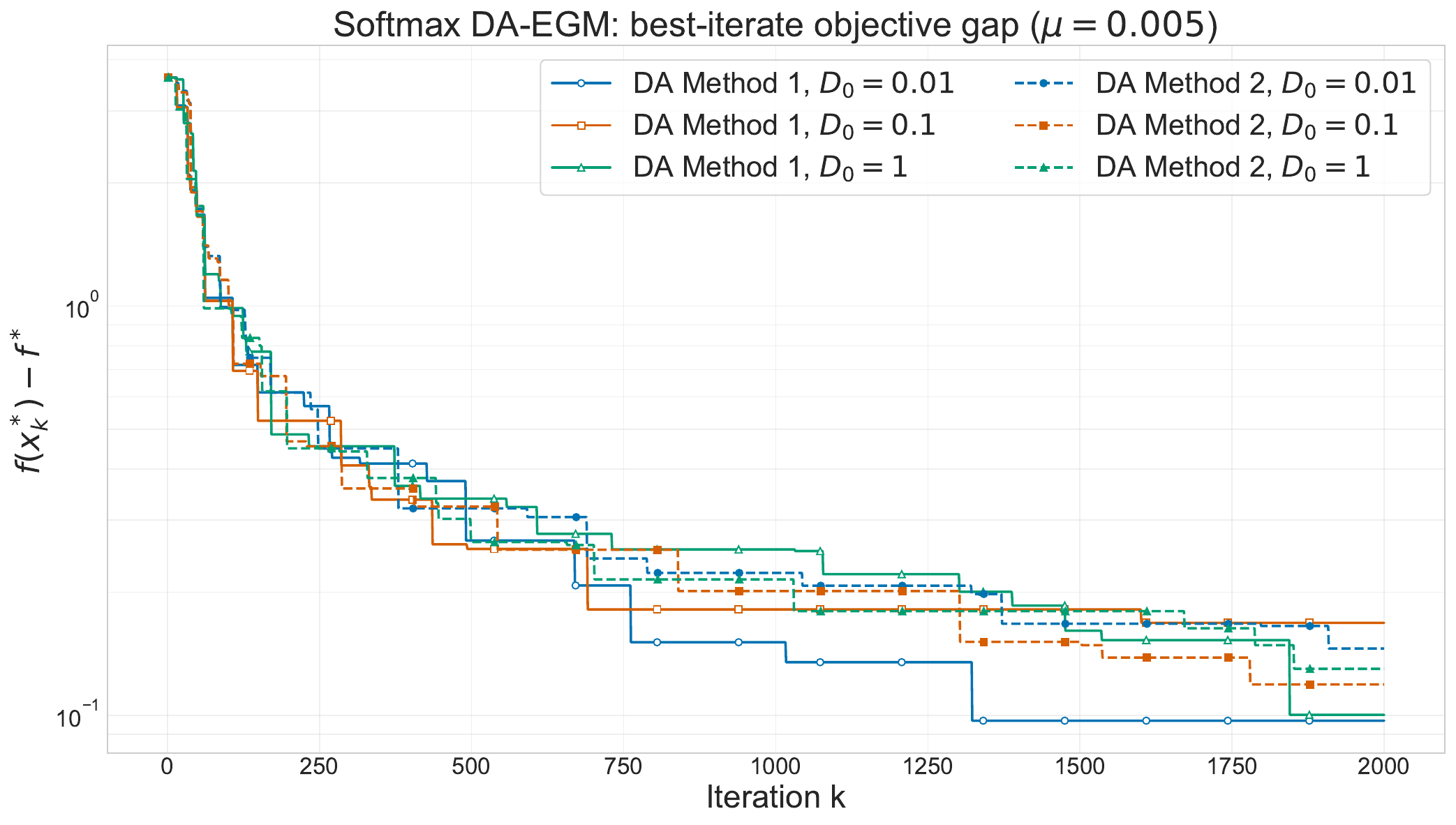}
\end{subfigure}
\par

\begin{subfigure}{0.32\textwidth}
    \centering
    \includegraphics[width=\linewidth]
    {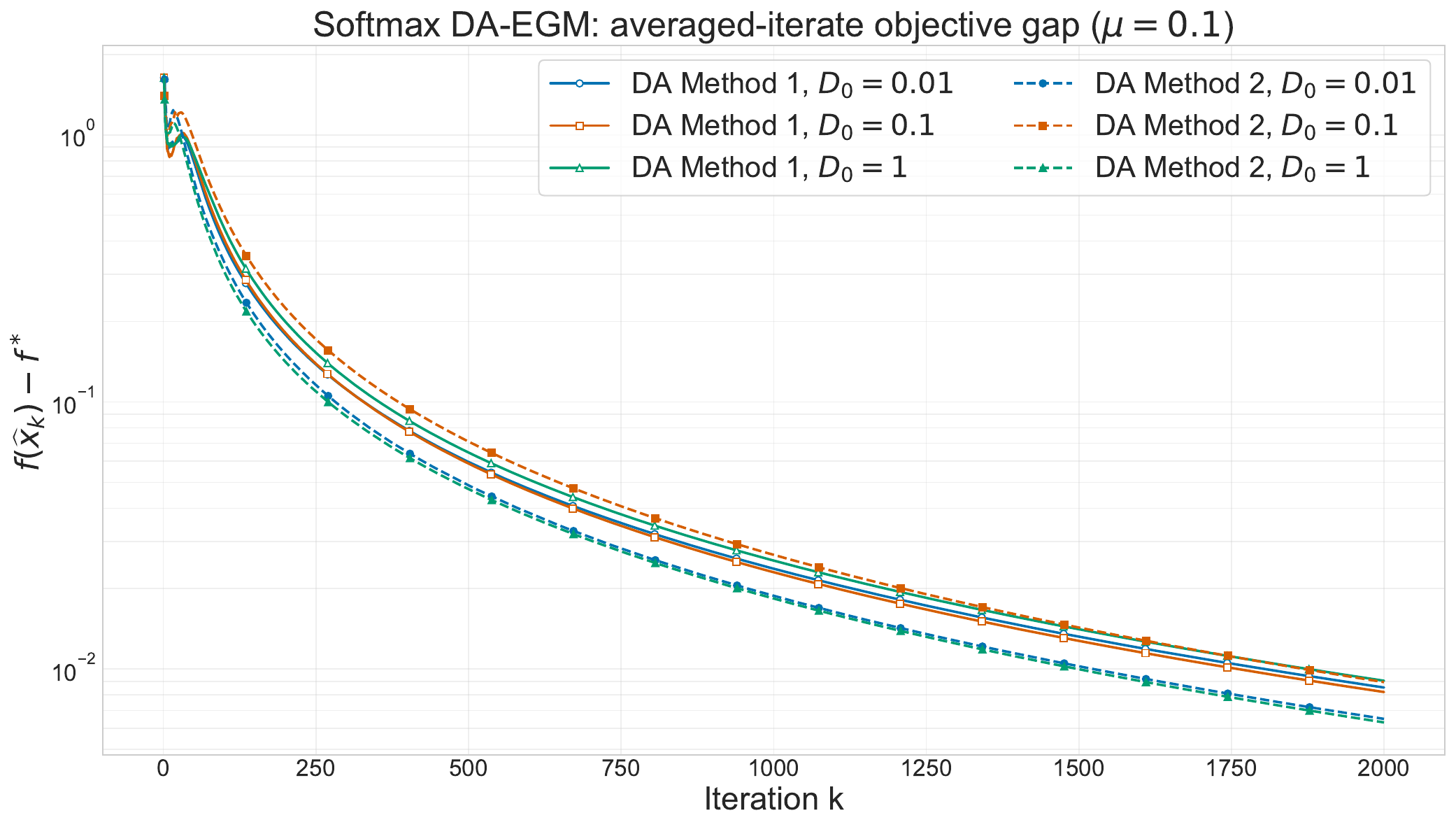}
\end{subfigure}
\hfill
\begin{subfigure}{0.32\textwidth}
    \centering
    \includegraphics[width=\linewidth]
    {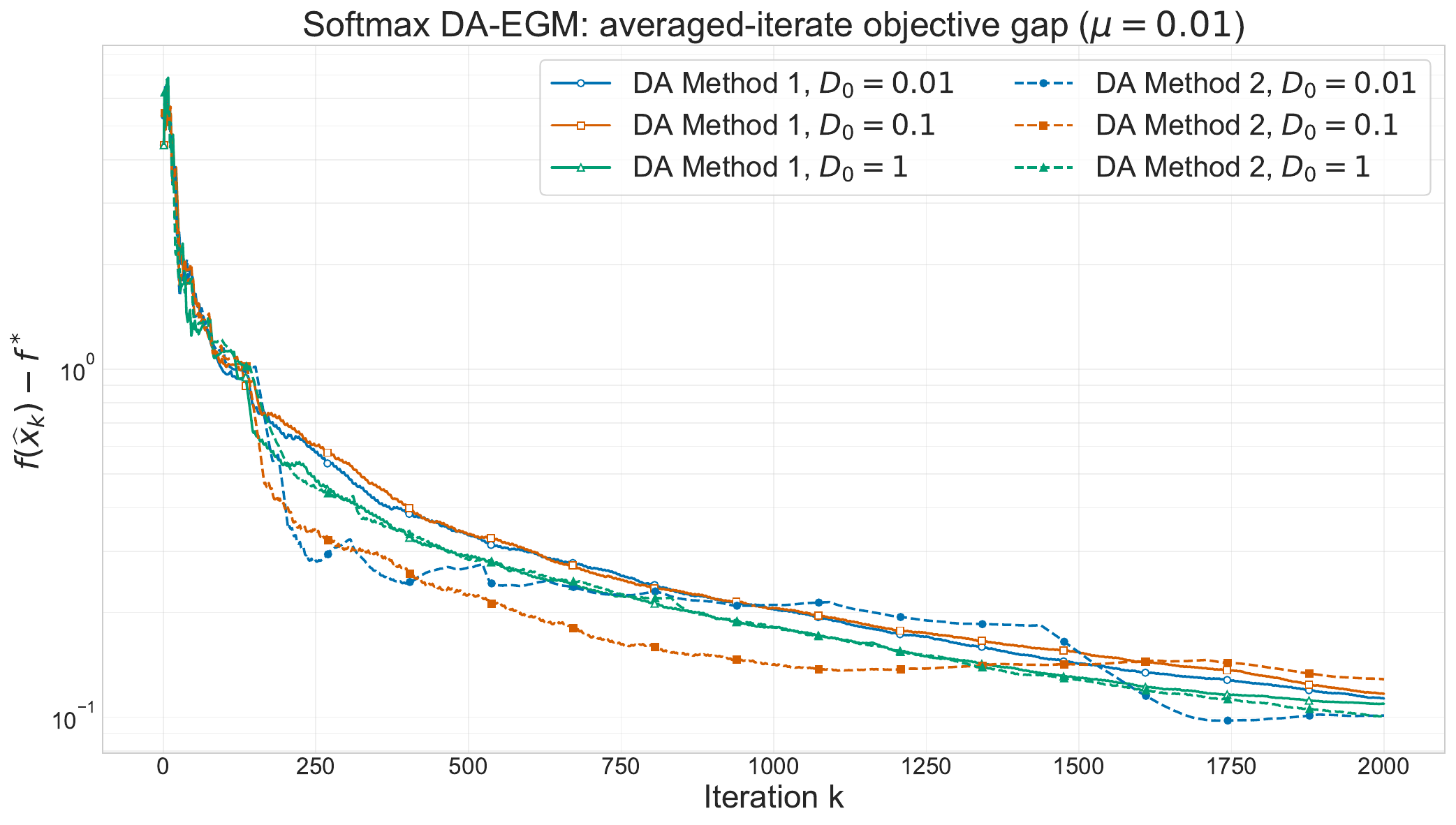}
\end{subfigure}
\hfill
\begin{subfigure}{0.32\textwidth}
    \centering
    \includegraphics[width=\linewidth]
    {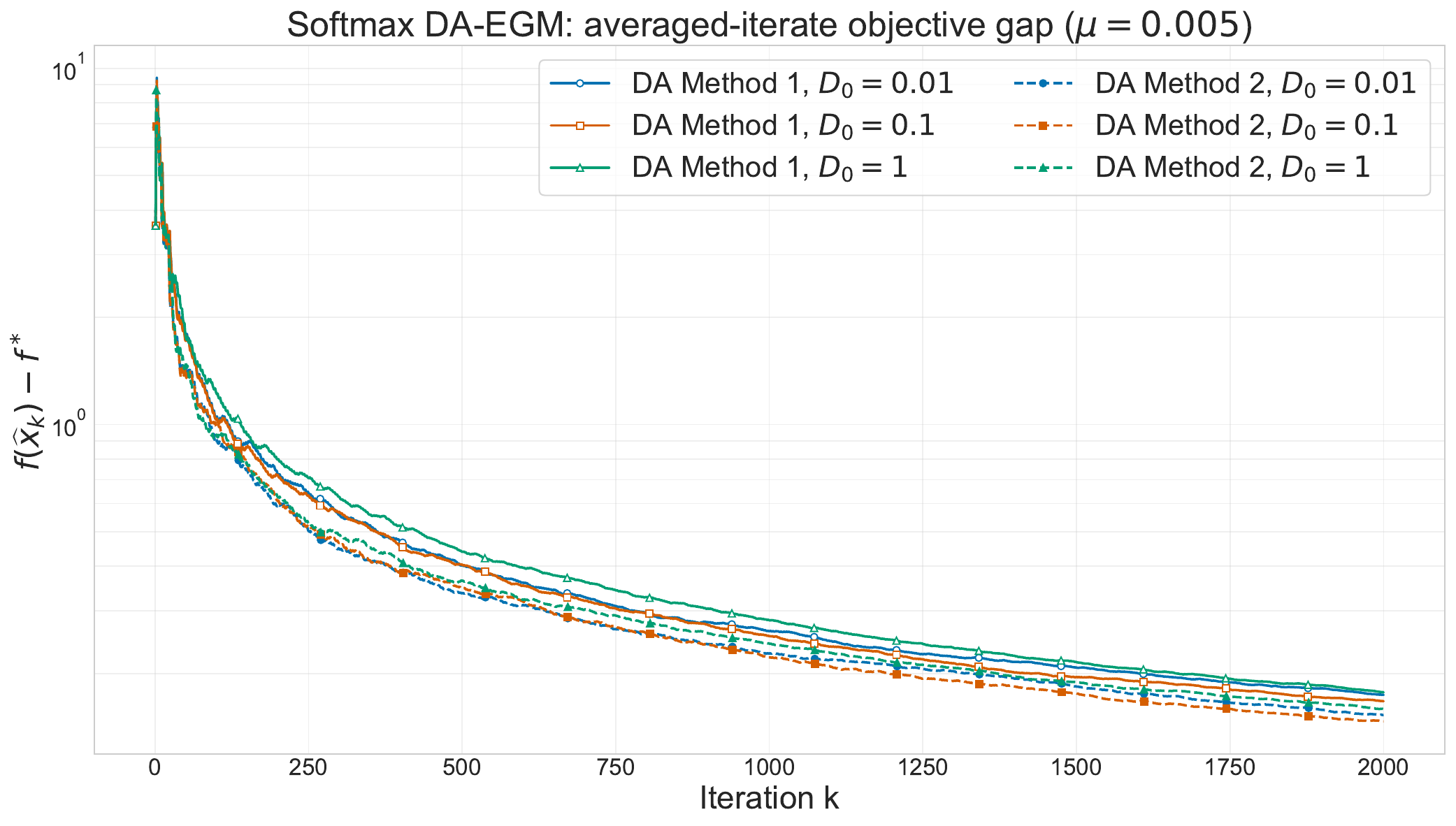}
\end{subfigure}
\par

\begin{subfigure}{0.32\textwidth}
    \centering
    \includegraphics[width=\linewidth]
    {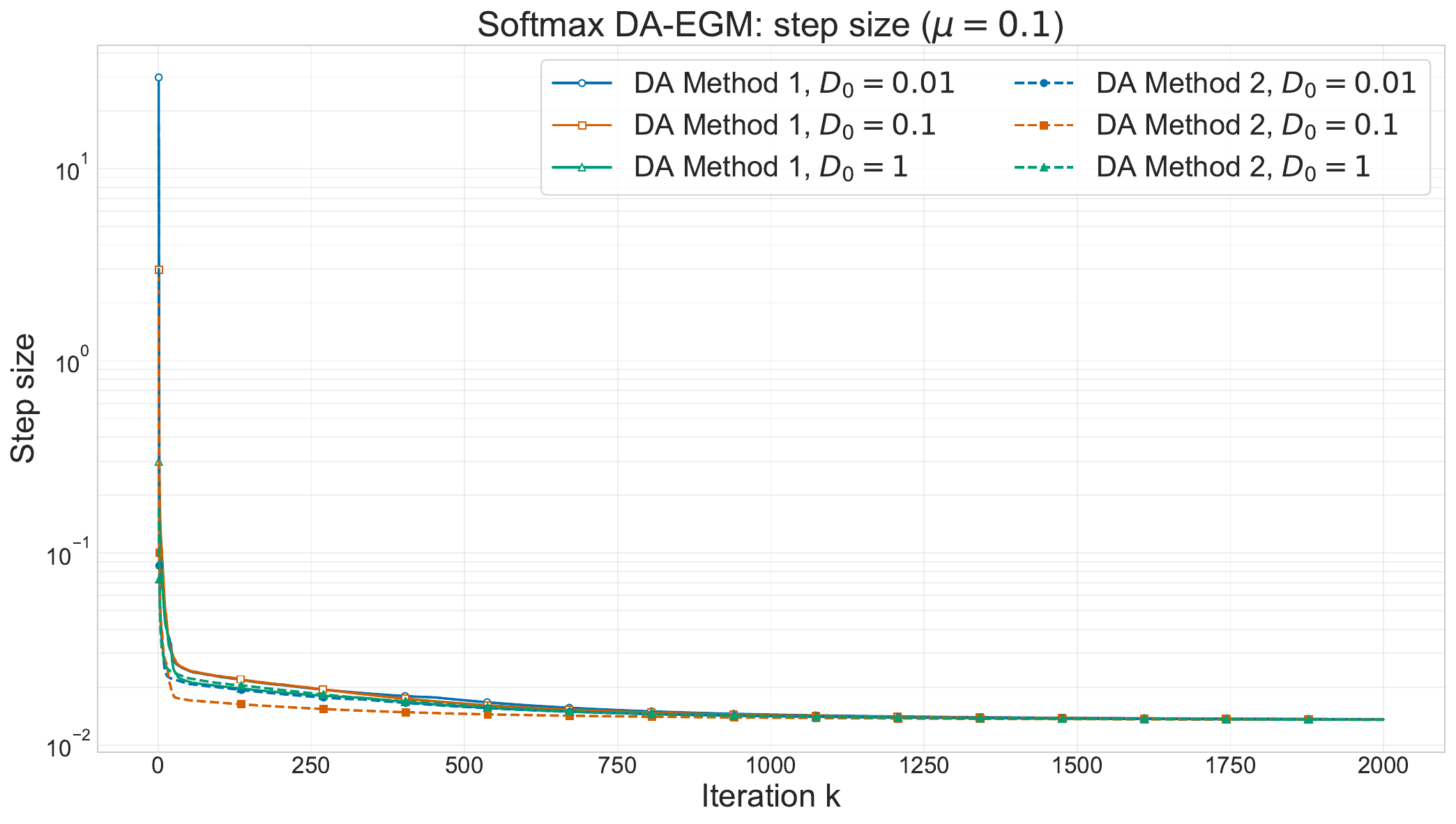}
\end{subfigure}
\hfill
\begin{subfigure}{0.32\textwidth}
    \centering
    \includegraphics[width=\linewidth]
    {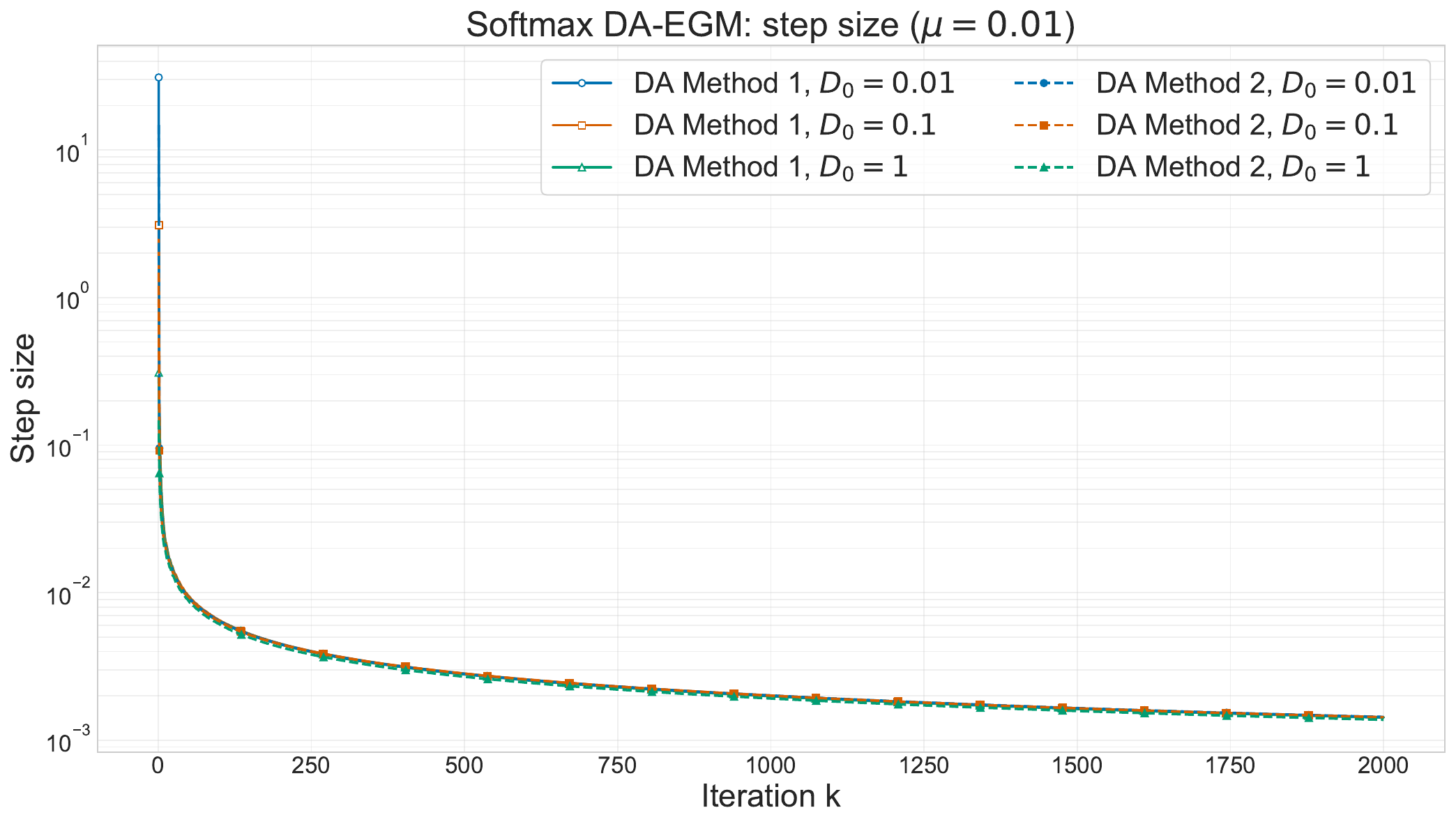}
\end{subfigure}
\hfill
\begin{subfigure}{0.32\textwidth}
    \centering
    \includegraphics[width=\linewidth]
    {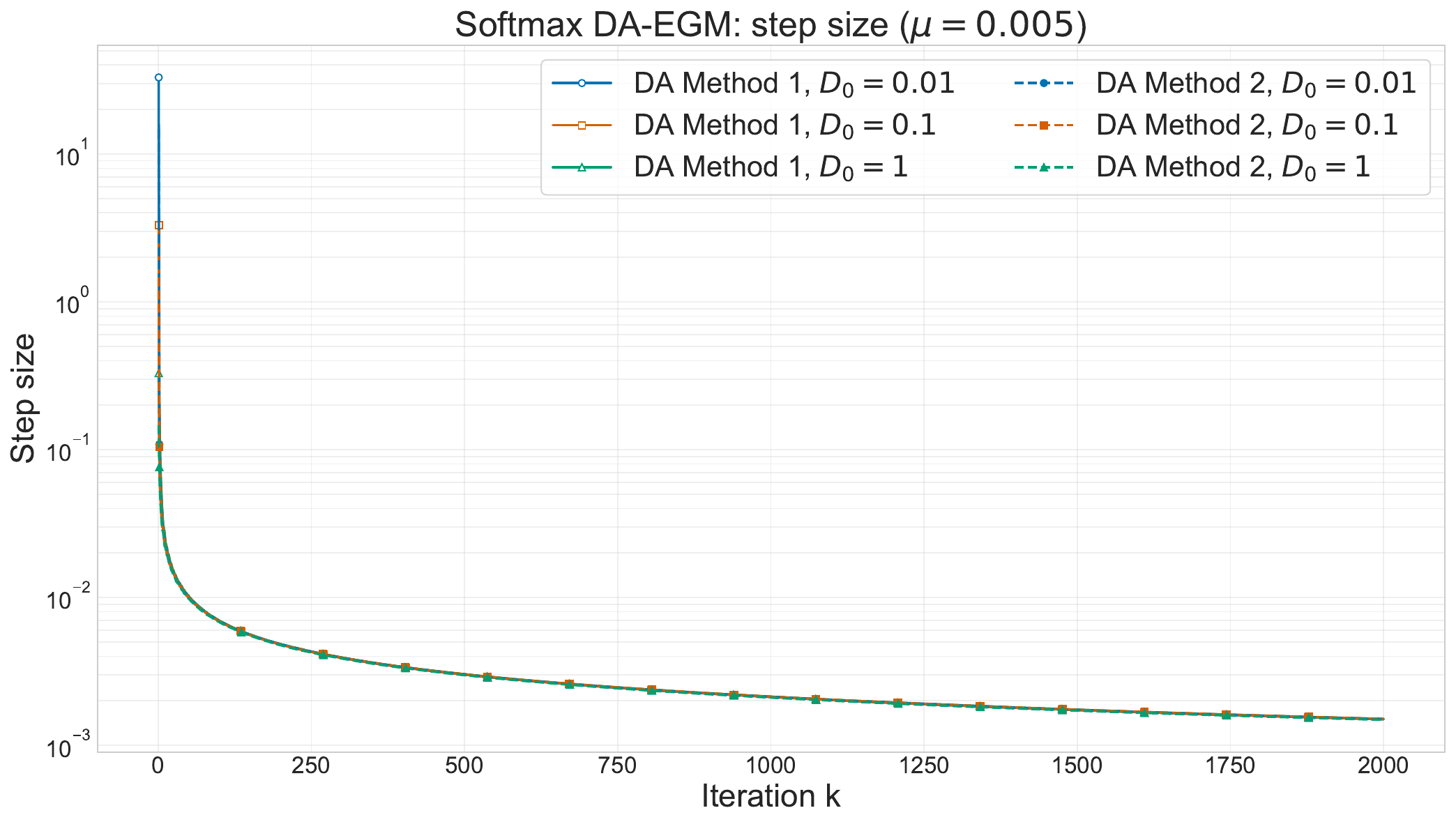}
\end{subfigure}

\caption{Sensitivity of DA-EGM Methods 1 and 2 to the initial distance
estimate $D_0$ for the deterministic softmax problem. From left to
right, the columns correspond to $\mu=0.1$, $\mu=0.01$, and
$\mu=0.005$. From top to bottom, the rows show the best-iterate
objective function gaps, averaged-iterate objective function gaps,
and step sizes.}
\label{fig:softmax_D0_sensitivity}
\end{figure*}



\paragraph{Softmax problem.}
We consider the optimization problem
\begin{equation}
    \min_{\|x\|_2\le 1}
    f(x) \quad \text{where} \quad f(x) :=
    \mu\log\left(
        \sum_{i=1}^{n}
        \exp\left(
            \frac{\langle a_i,x\rangle-b_i}{\mu}
        \right)
    \right) , \nonumber
\end{equation}
where $(a_i, b_i)_{i = 1}^n \in \R^d \times \R$ are the data points and $\mu > 0$ is a parameter. The problem is a smooth approximation of $\max_{1 \leq i \leq n} [\la a_i , x \ra - b_i]$~\citep{nesterov2005smooth}. The data are generated following the procedure in the paper~\citep[Appendix E]{moshtaghifar2026dada}. The data generation technique ensures that the optimal function value $f^* = 0$. We use \(n=1000\), \(d=2000\), and
\(\mu\in\{0.1,0.01,0.005\}\). The experiment is deterministic, and
the full batch is used at every iteration. As \(\mu\) decreases, the
softmax function becomes a sharper approximation of the non-smooth
maximum function, making the optimization problem more difficult.

The convergence results in Figure~\ref{fig_softmax_results} show that, in terms of the best-iterate objective function gap $f(x_k^*)-f^*$, our methods perform competitively on the smoother instance with $\mu=0.1$. For the less smooth instances with $\mu=0.01$ and $\mu=0.005$, DoG attains smaller objective function gaps over the considered iteration horizon, while the other methods perform similarly. In terms of the averaged-iterate objective function gap $f(\widehat{x}_k)-f^*$, DoG performs better for all tested values of $\mu$.

The purpose of these experiments is not to establish uniform empirical dominance over the existing methods. Rather, they demonstrate that our methods remain stable and consistently reduce the objective function gap across different smoothness regimes, in agreement with their theoretical universality guarantees. Although their empirical performance may be improved by tuning $\widehat D$, we use the principled and uniform choice $\widehat D=D$ in our EGM Rules~1--3, without instance-specific tuning.

We study the sensitivity of EGM and DA EGM to the choices of $\widehat D$ and $D_0$ in Figures~\ref{fig:softmax_D_hat_sensitivity} 
and~\ref{fig:softmax_D0_sensitivity}, respectively. 
Figure~\ref{fig:softmax_D_hat_sensitivity} shows that, for $\mu=0.1$, Rule~3 with $\widehat D=3$ performs well for both the best and averaged iterates. When $\mu=0.01$, Rules~1 and~2 with $\widehat D=0.5$ perform well in both cases, while Rule~3 with $\widehat D=0.5$ also performs well for the best iterate. When $\mu=0.005$, corresponding to the least smooth setting, Rule~1 with $\widehat D=0.5$ performs best. Thus, for this problem, smaller values of $\widehat D$ yield better objective-value reductions as the problem becomes 
less smooth. Figure~\ref{fig:softmax_D0_sensitivity} shows that most DA EGM variants perform similarly across the tested settings. This observation is consistent with our theory, which predicts substantially reduced sensitivity of DA EGM to the initialization $D_0$. Moreover, the best-tuned configurations in Figure~\ref{fig:softmax_D_hat_sensitivity} outperform the DoG results reported in Figure~\ref{fig_softmax_results}.


\begin{figure*}[htbp]
\centering


\begin{subfigure}{0.24\textwidth}
    \centering
    \includegraphics[width=\linewidth]{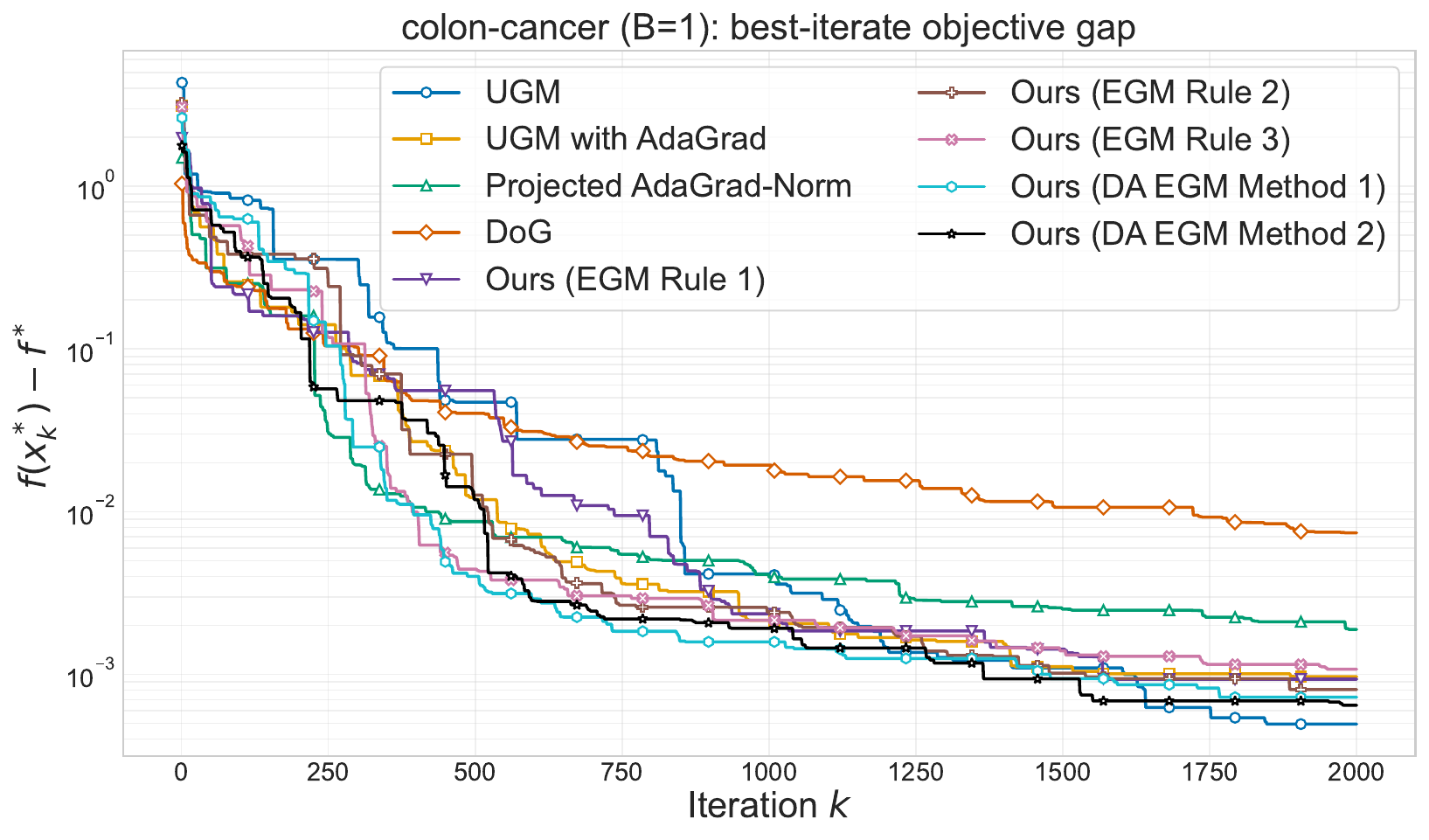}
\end{subfigure}
\hfill
\begin{subfigure}{0.24\textwidth}
    \centering
    \includegraphics[width=\linewidth]{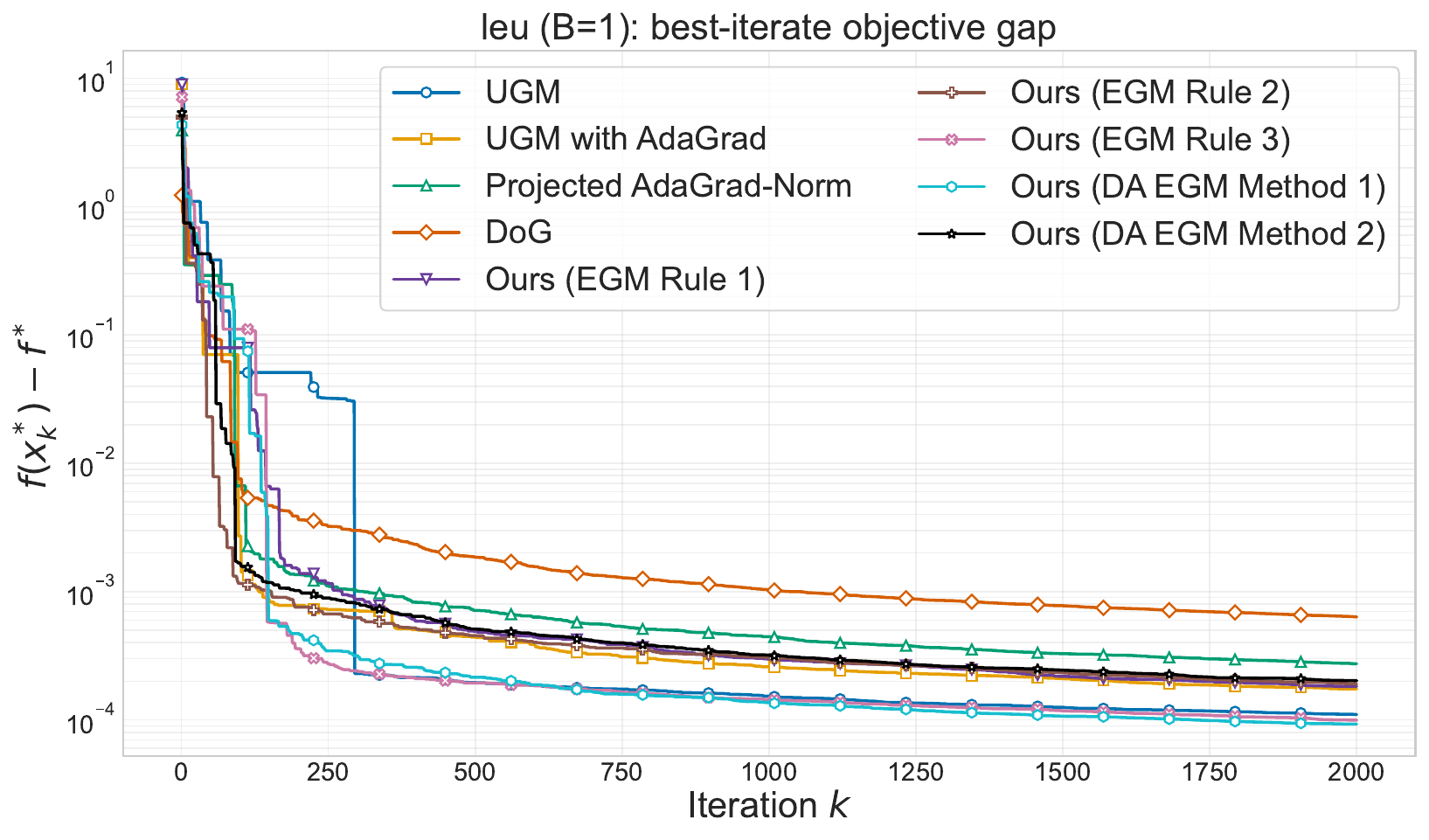}
\end{subfigure}
\hfill
\begin{subfigure}{0.24\textwidth}
    \centering
    \includegraphics[width=\linewidth]{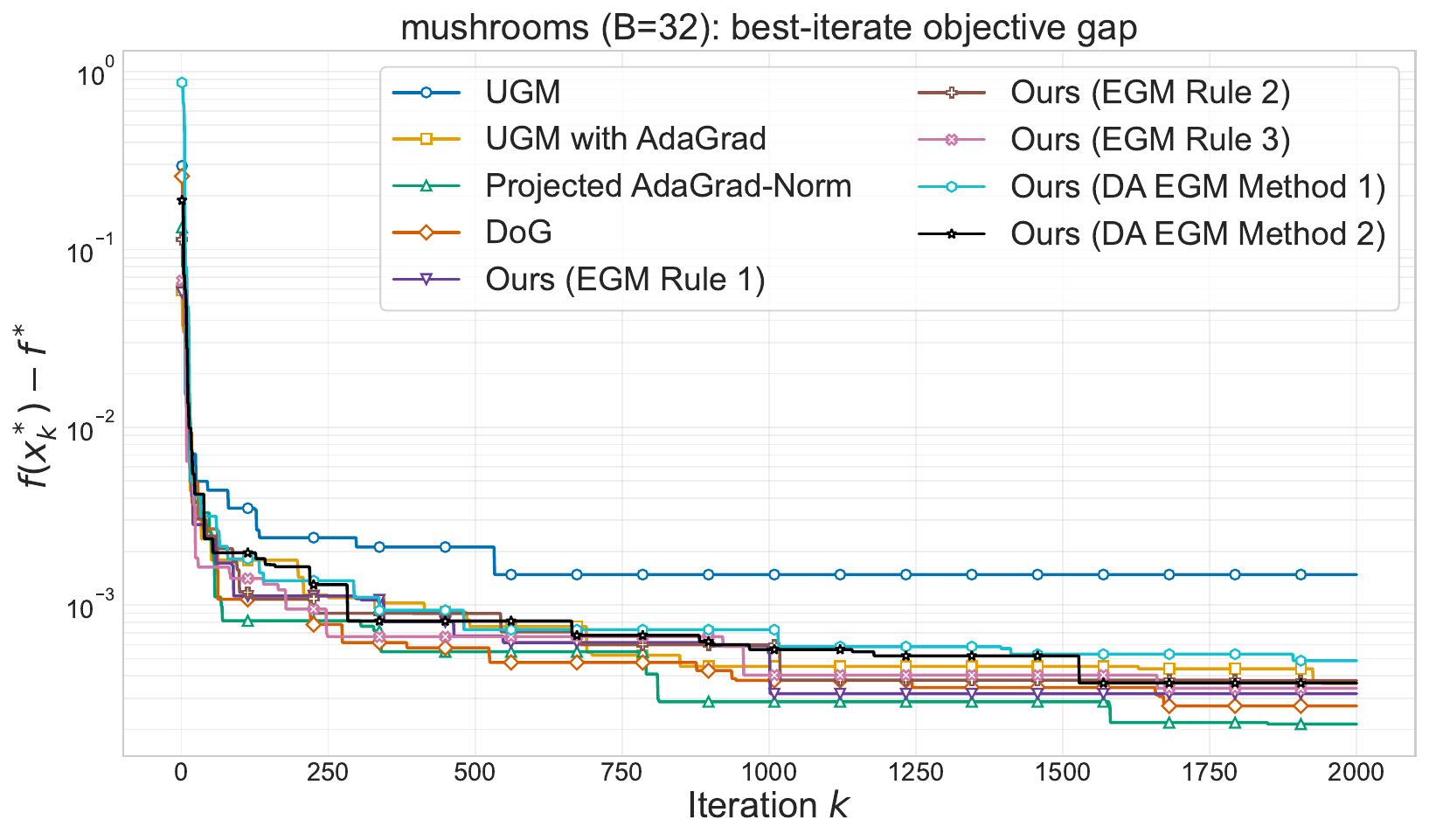}
\end{subfigure}
\hfill
\begin{subfigure}{0.24\textwidth}
    \centering
    \includegraphics[width=\linewidth]{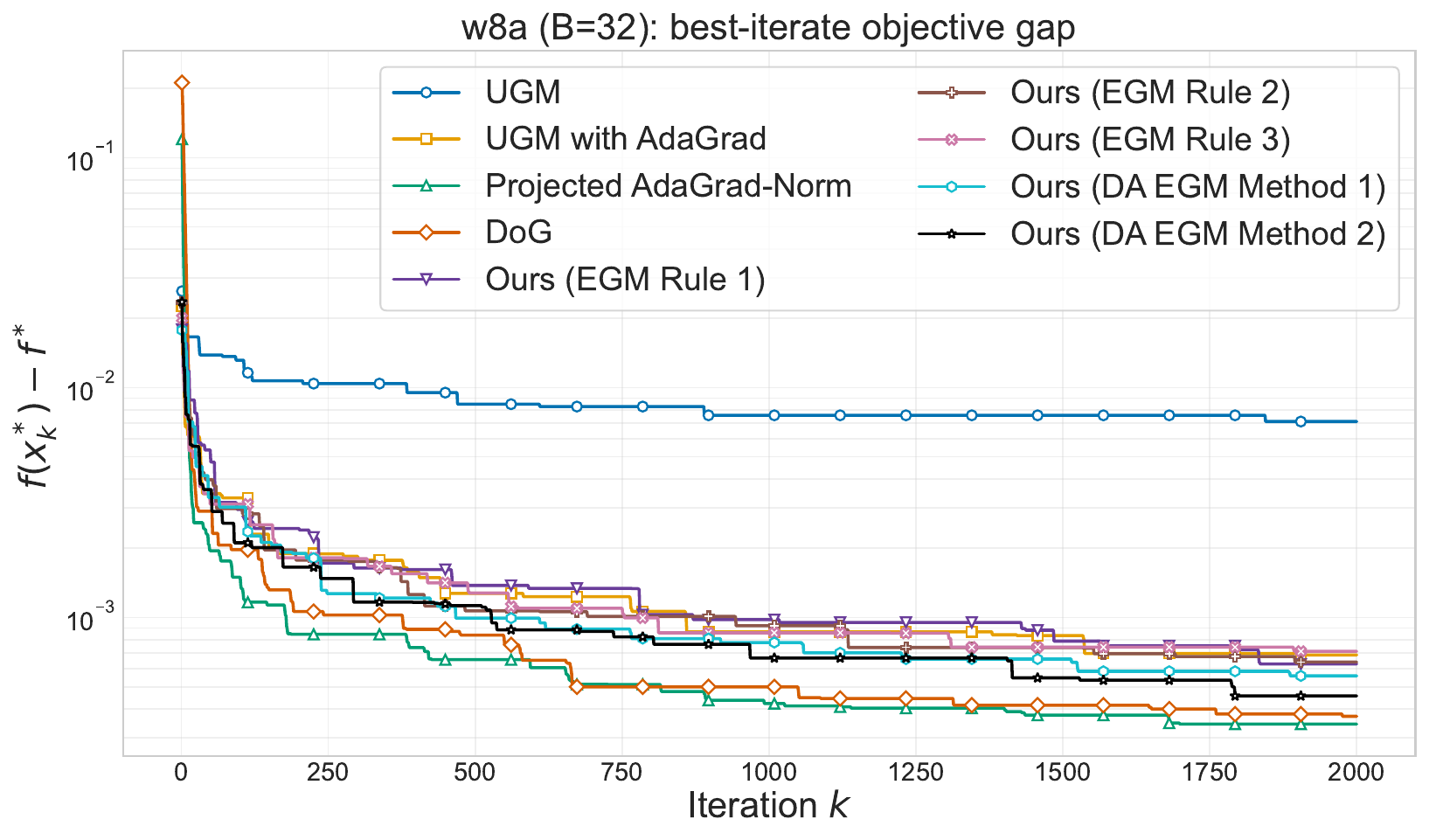}
\end{subfigure}
\par

\begin{subfigure}{0.24\textwidth}
    \centering
    \includegraphics[width=\linewidth]{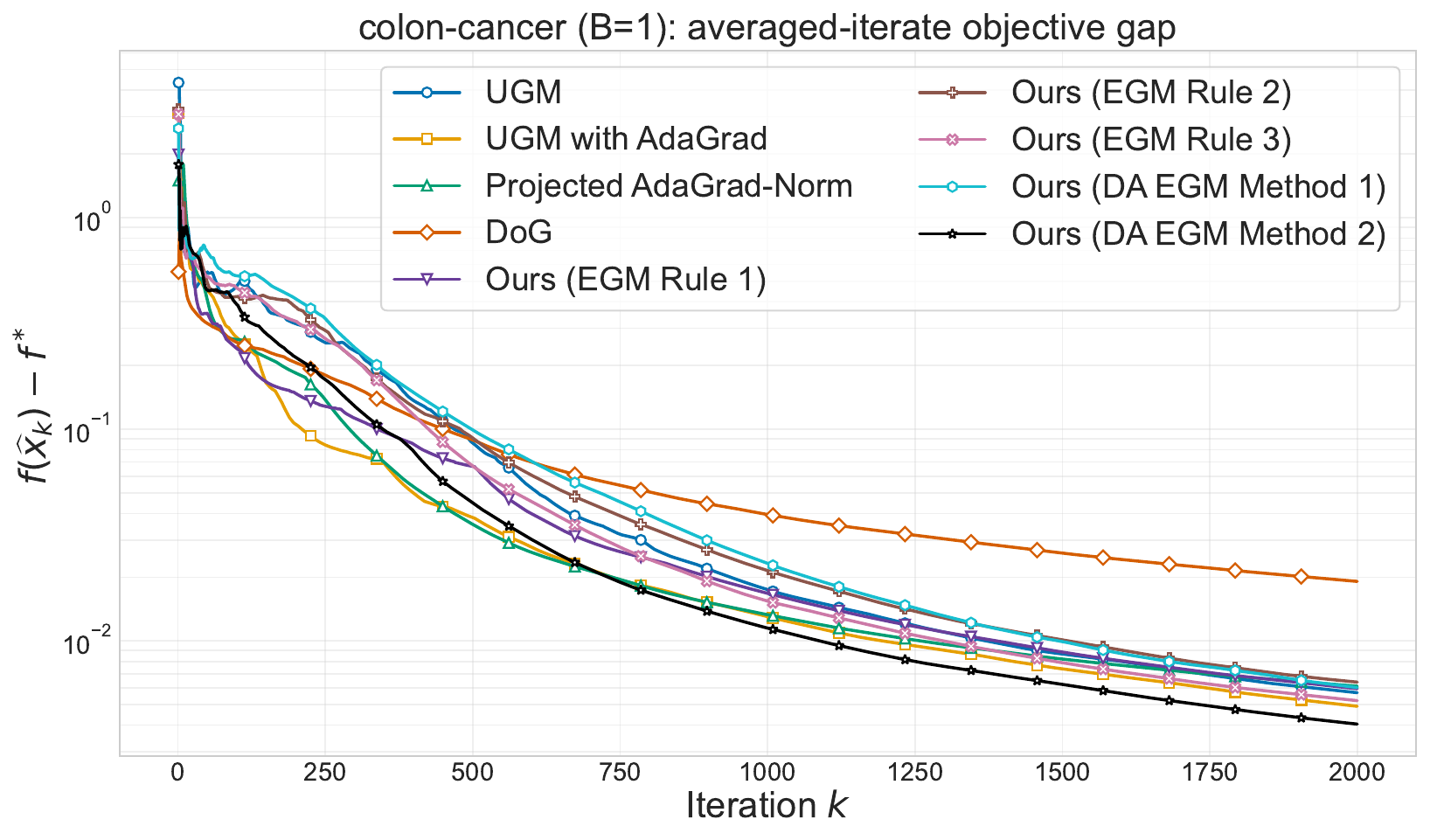}
\end{subfigure}
\hfill
\begin{subfigure}{0.24\textwidth}
    \centering
    \includegraphics[width=\linewidth]{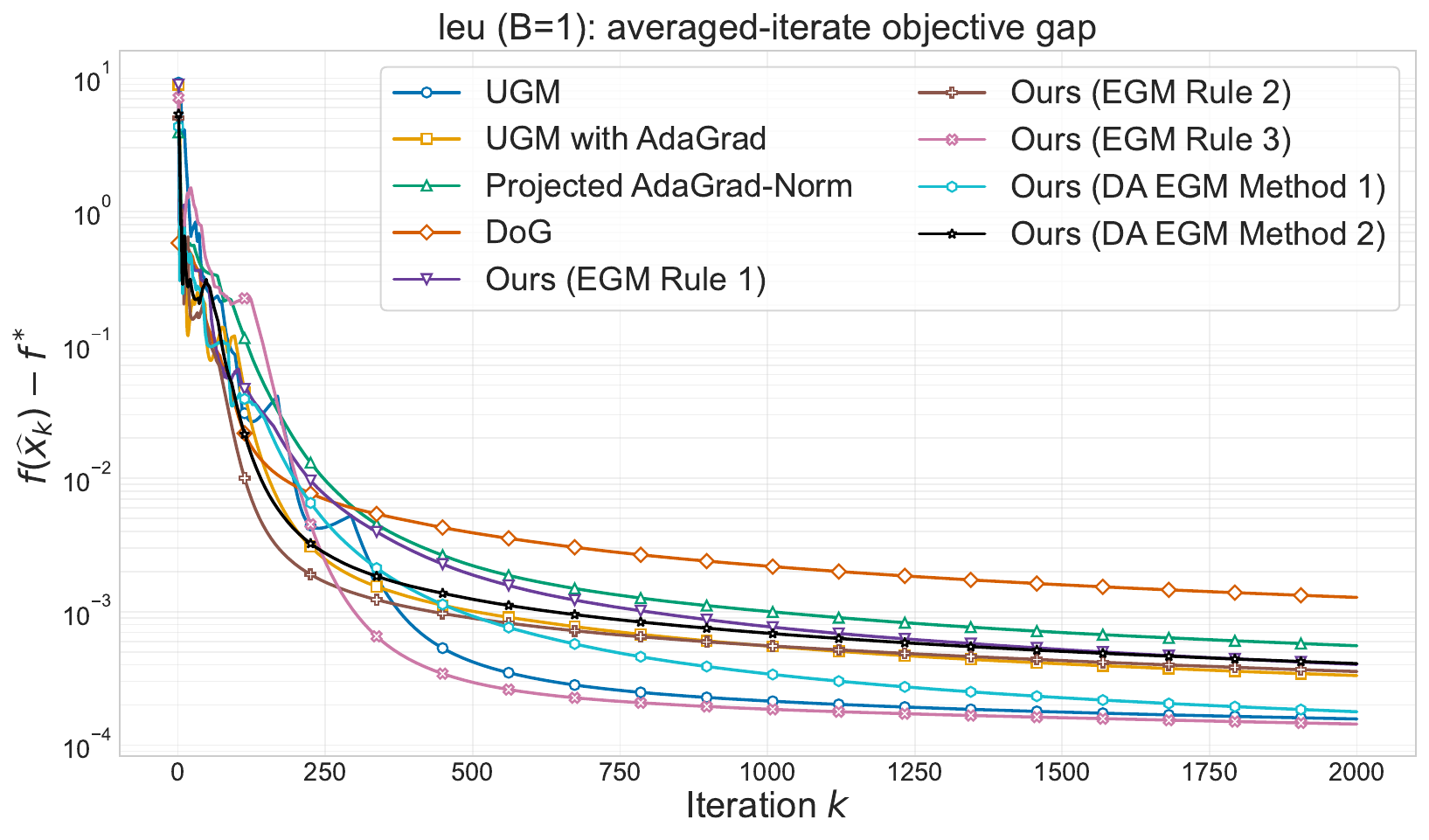}
\end{subfigure}
\hfill
\begin{subfigure}{0.24\textwidth}
    \centering
    \includegraphics[width=\linewidth]{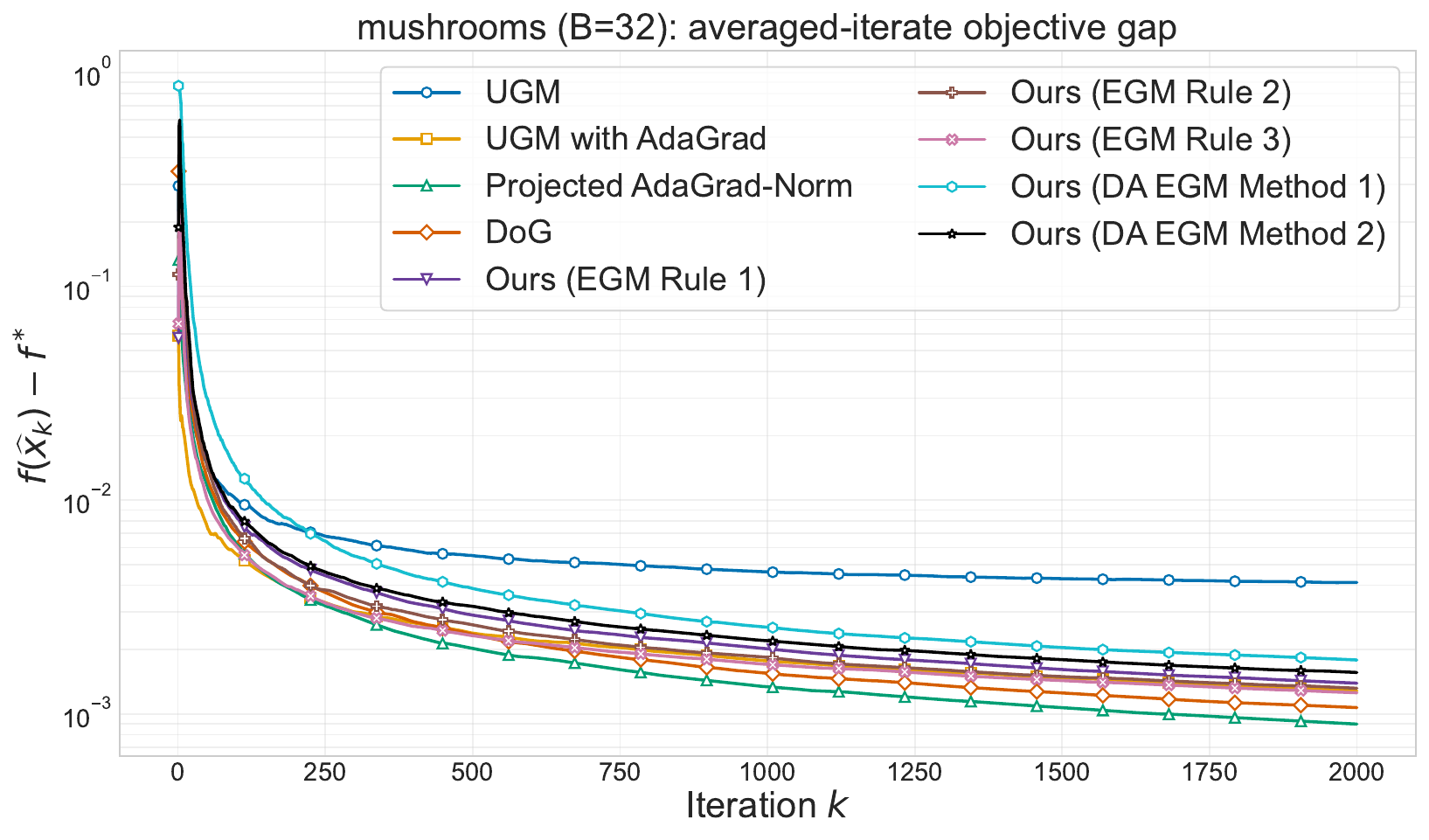}
\end{subfigure}
\hfill
\begin{subfigure}{0.24\textwidth}
    \centering
    \includegraphics[width=\linewidth]{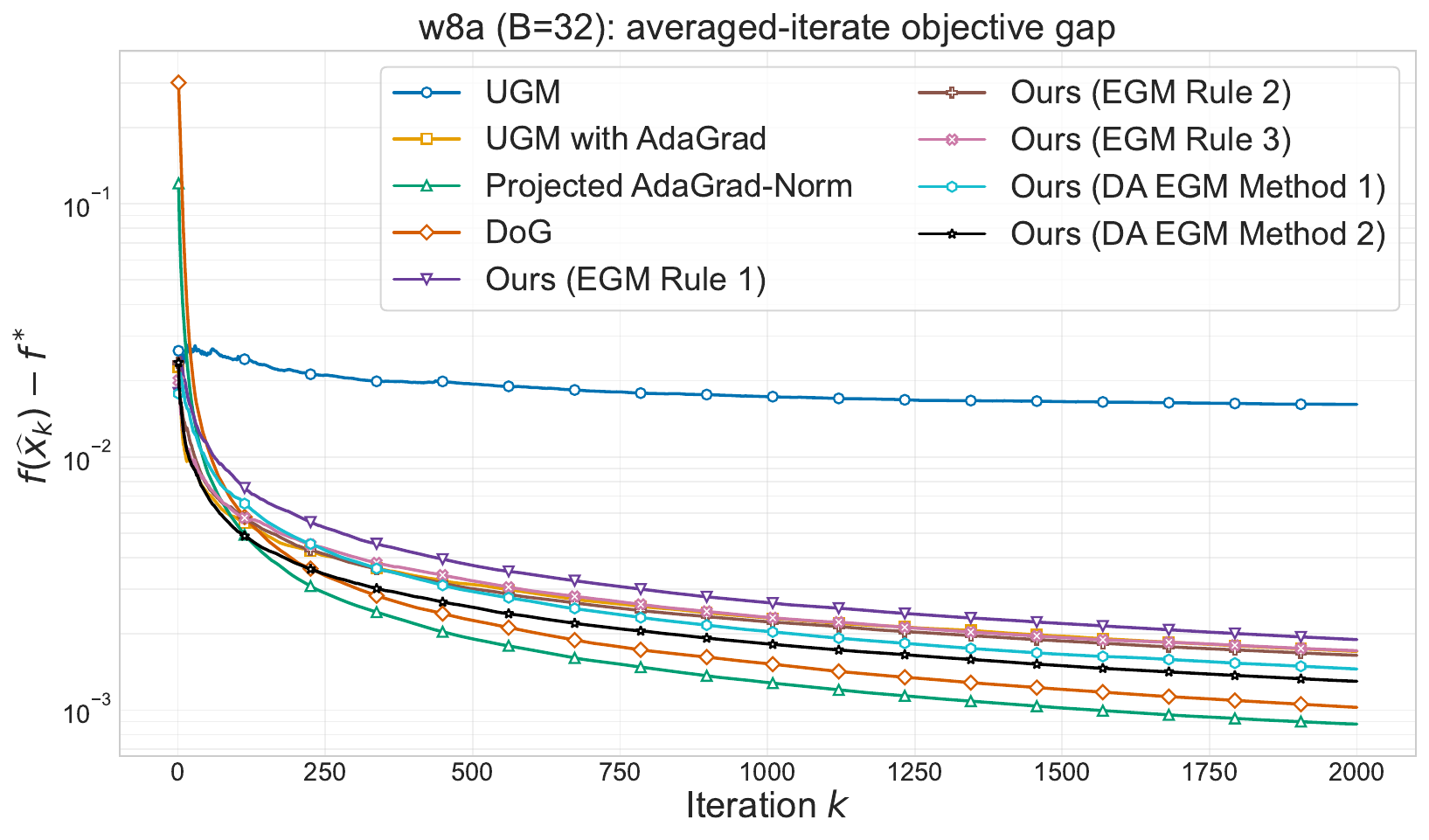}
\end{subfigure}
\par

\begin{subfigure}{0.24\textwidth}
    \centering
    \includegraphics[width=\linewidth]{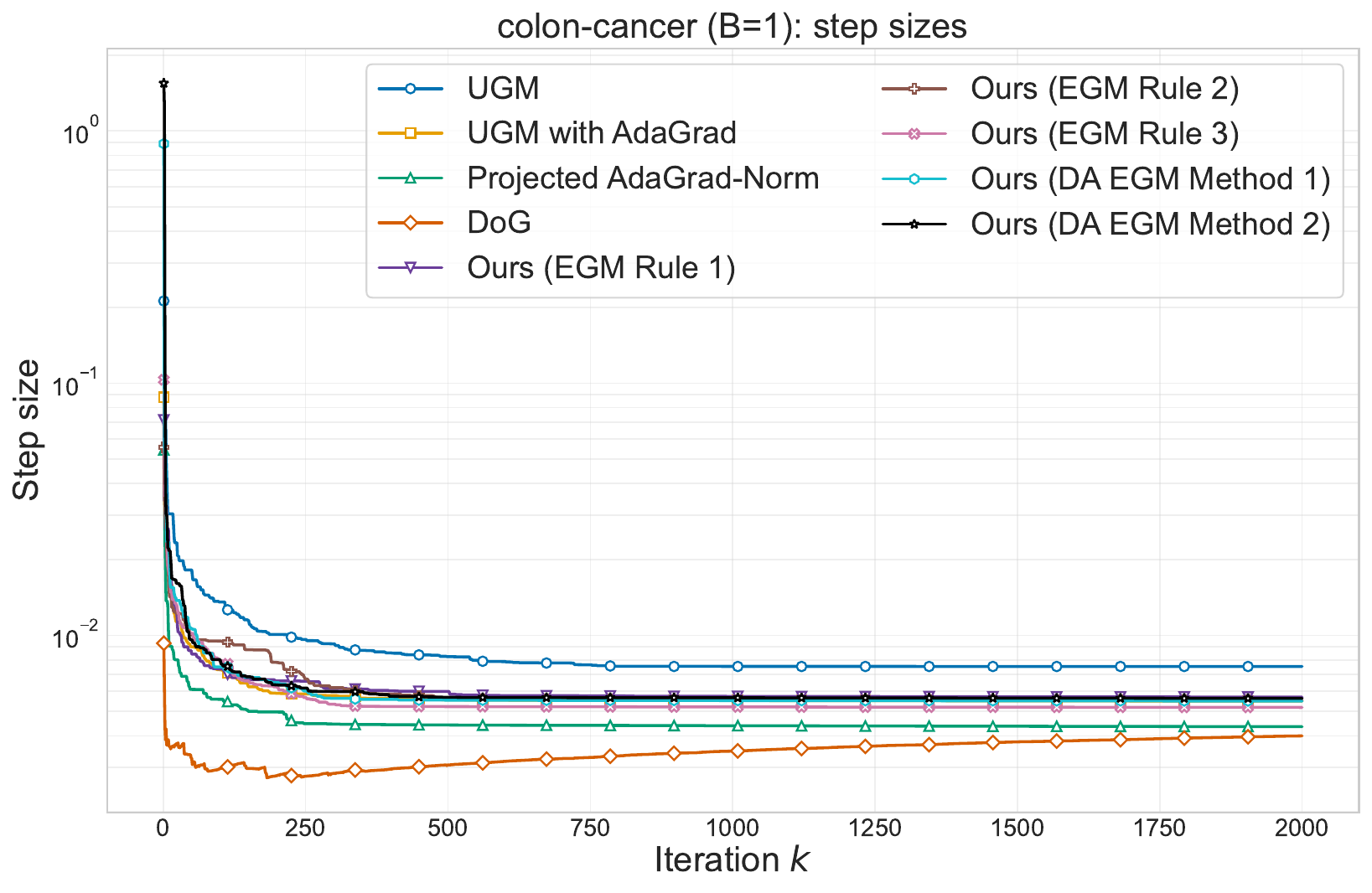}
\end{subfigure}
\hfill
\begin{subfigure}{0.24\textwidth}
    \centering
    \includegraphics[width=\linewidth]{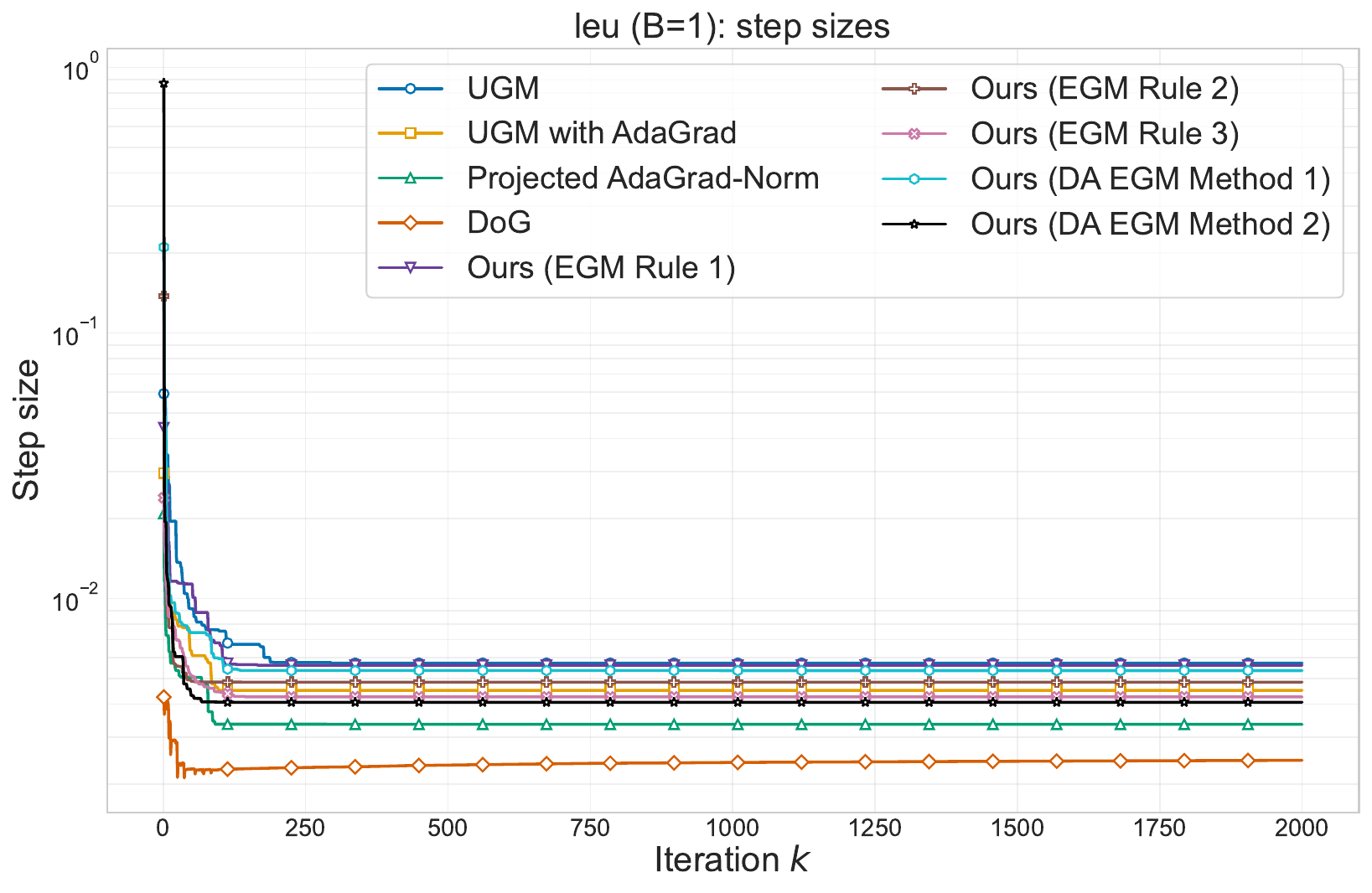}
\end{subfigure}
\hfill
\begin{subfigure}{0.24\textwidth}
    \centering
    \includegraphics[width=\linewidth]{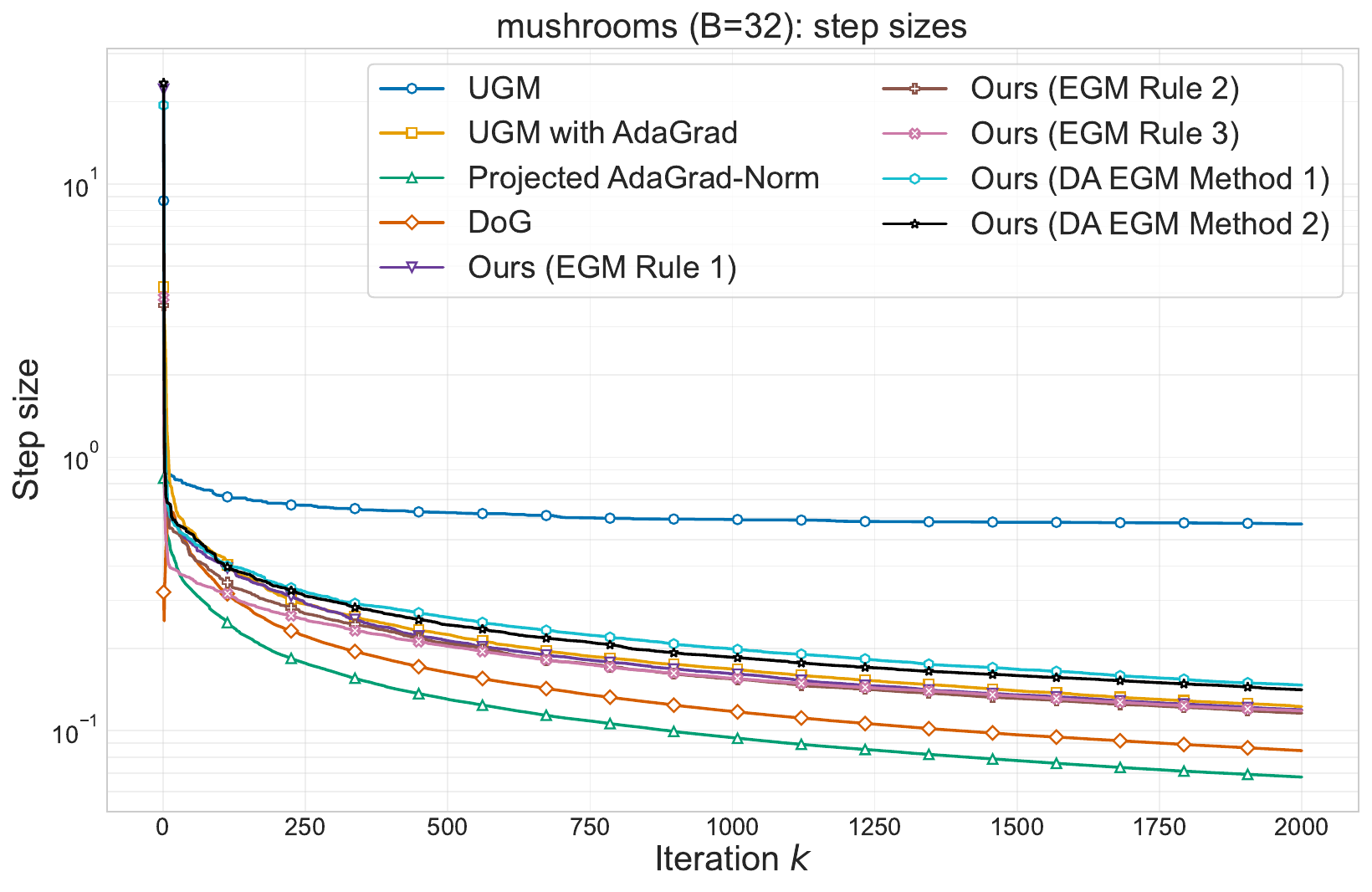}
\end{subfigure}
\hfill
\begin{subfigure}{0.24\textwidth}
    \centering
    \includegraphics[width=\linewidth]{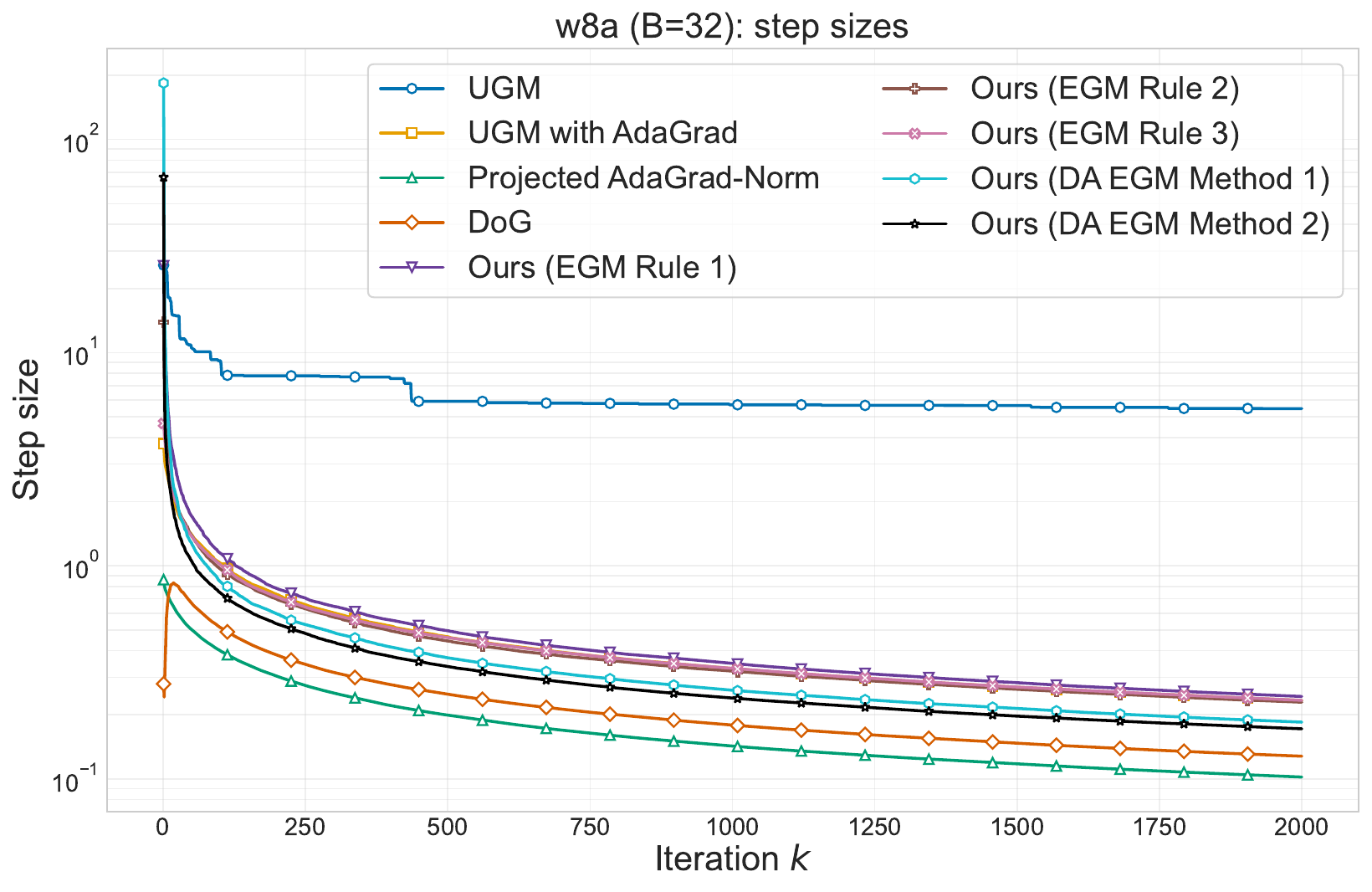}
\end{subfigure}

\caption{Results for logistic regression problem (stochastic). From top to bottom, the rows
show the best-iterate objective function gaps,
averaged-iterate objective function gaps, and step sizes.}
\label{fig_logistic}
\end{figure*}




\begin{figure*}[htbp]
\centering

\begin{subfigure}{0.24\textwidth}
    \centering
    \includegraphics[width=\linewidth]
    {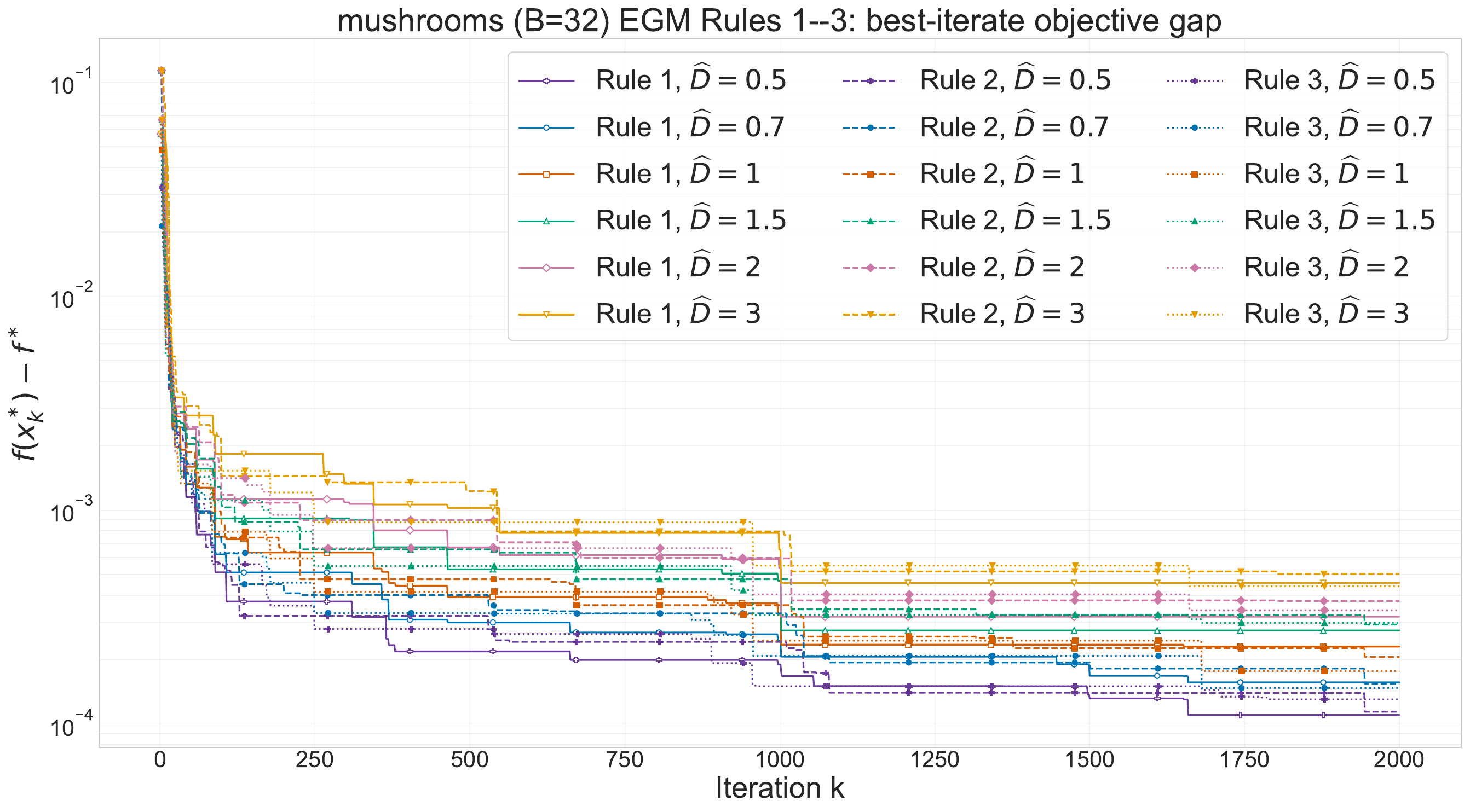}
\end{subfigure}
\hfill
\begin{subfigure}{0.24\textwidth}
    \centering
    \includegraphics[width=\linewidth]
    {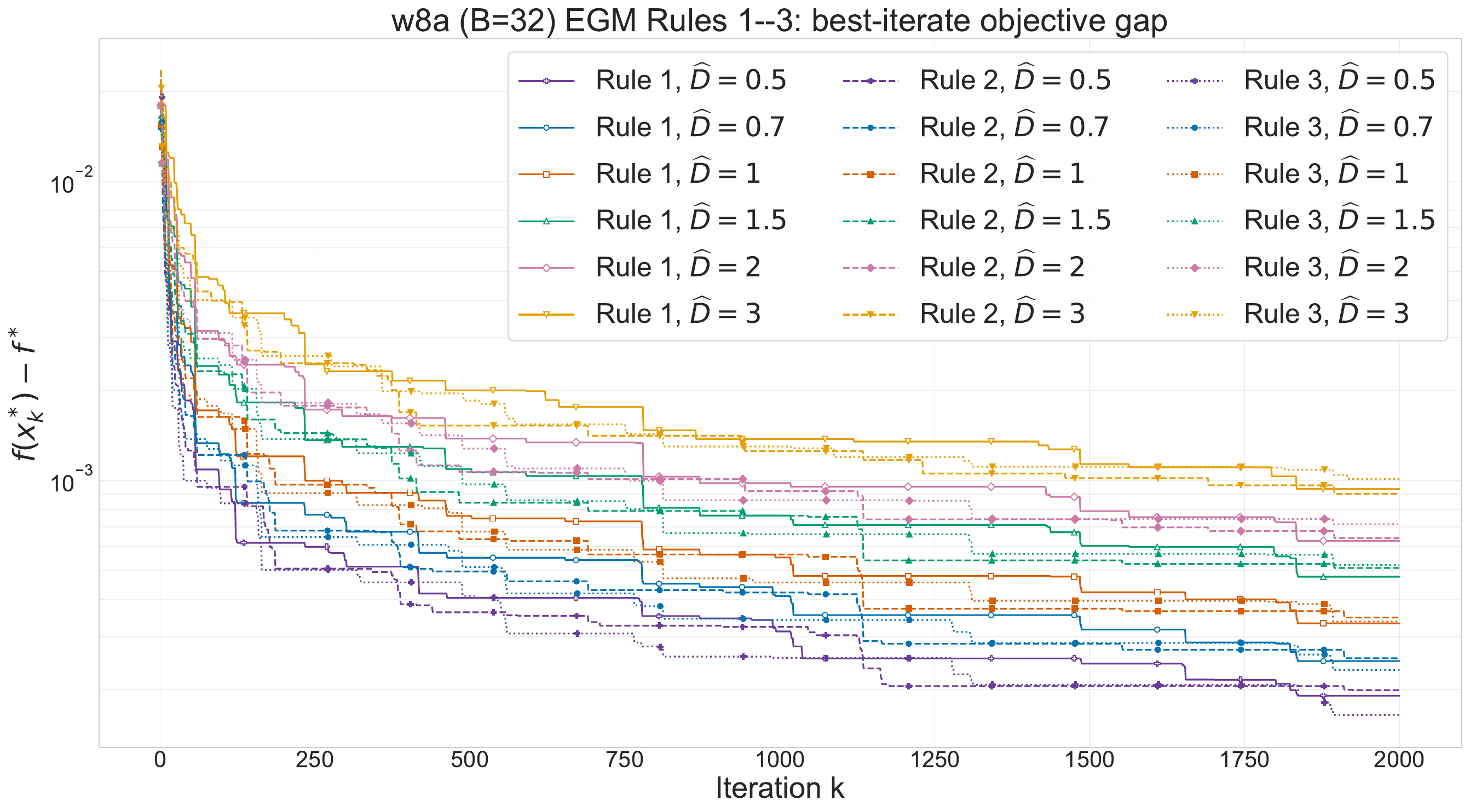}
\end{subfigure}
\hfill
\begin{subfigure}{0.24\textwidth}
    \centering
    \includegraphics[width=\linewidth]
    {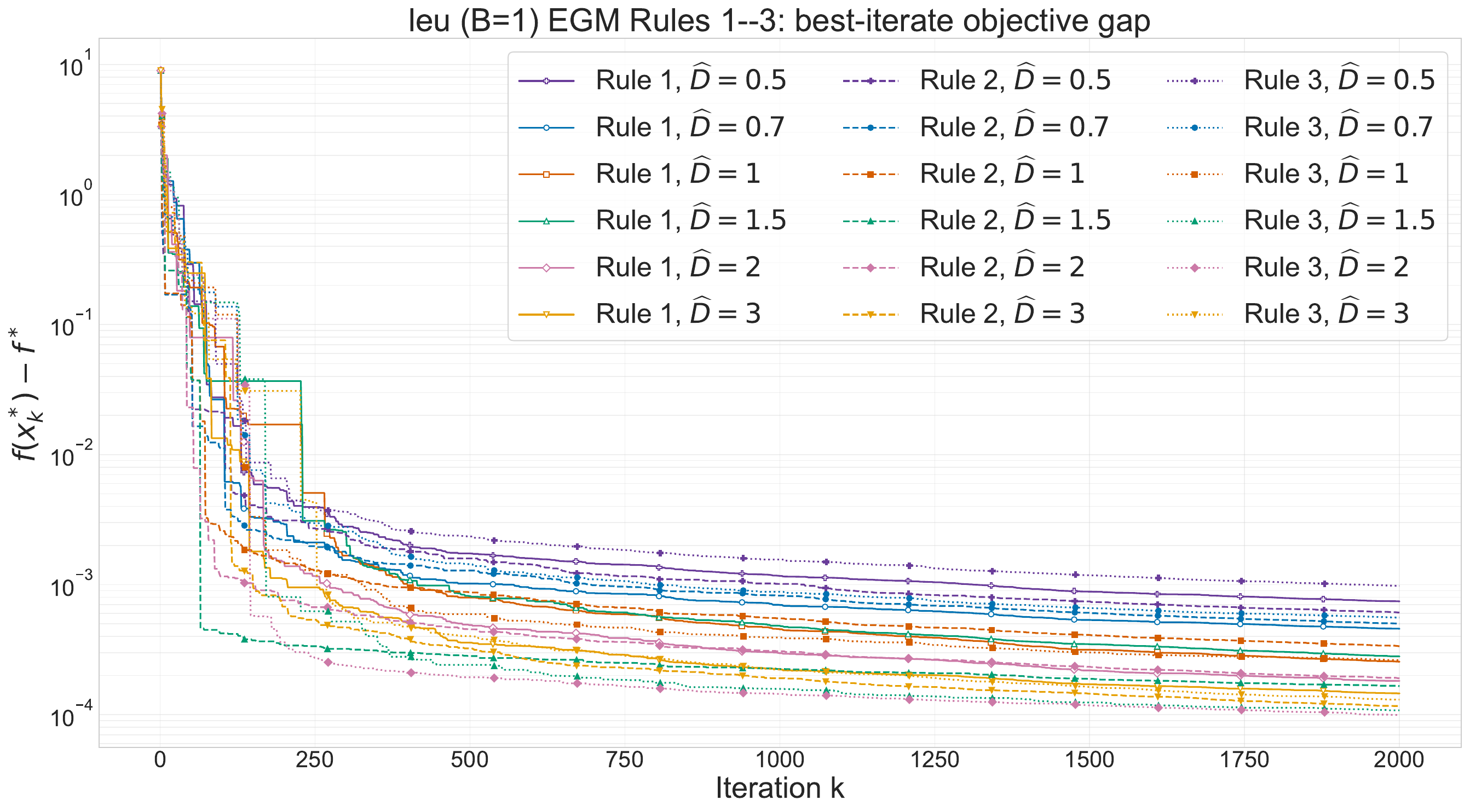}
\end{subfigure}
\hfill
\begin{subfigure}{0.24\textwidth}
    \centering
    \includegraphics[width=\linewidth]
    {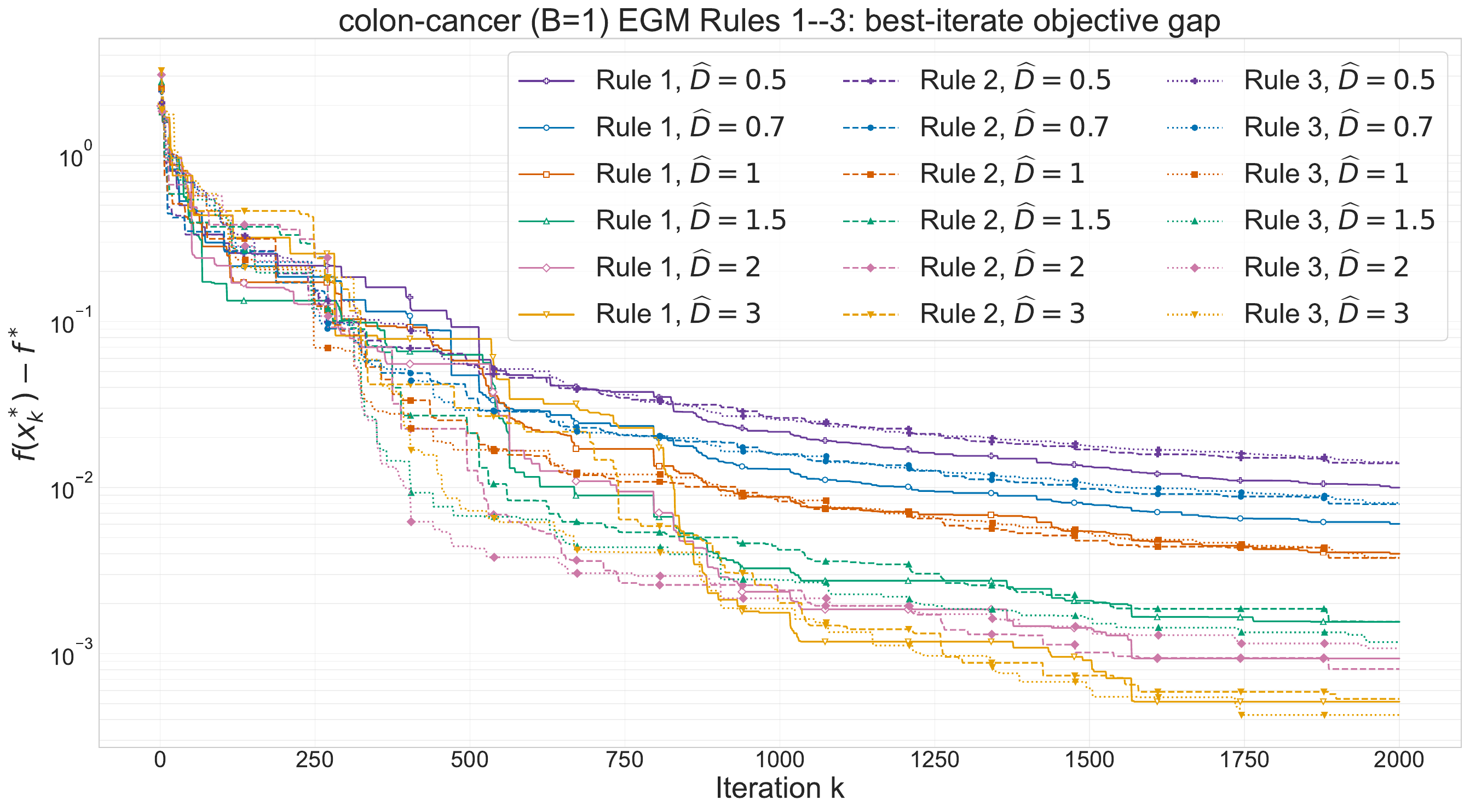}
\end{subfigure}
\par

\begin{subfigure}{0.24\textwidth}
    \centering
    \includegraphics[width=\linewidth]
    {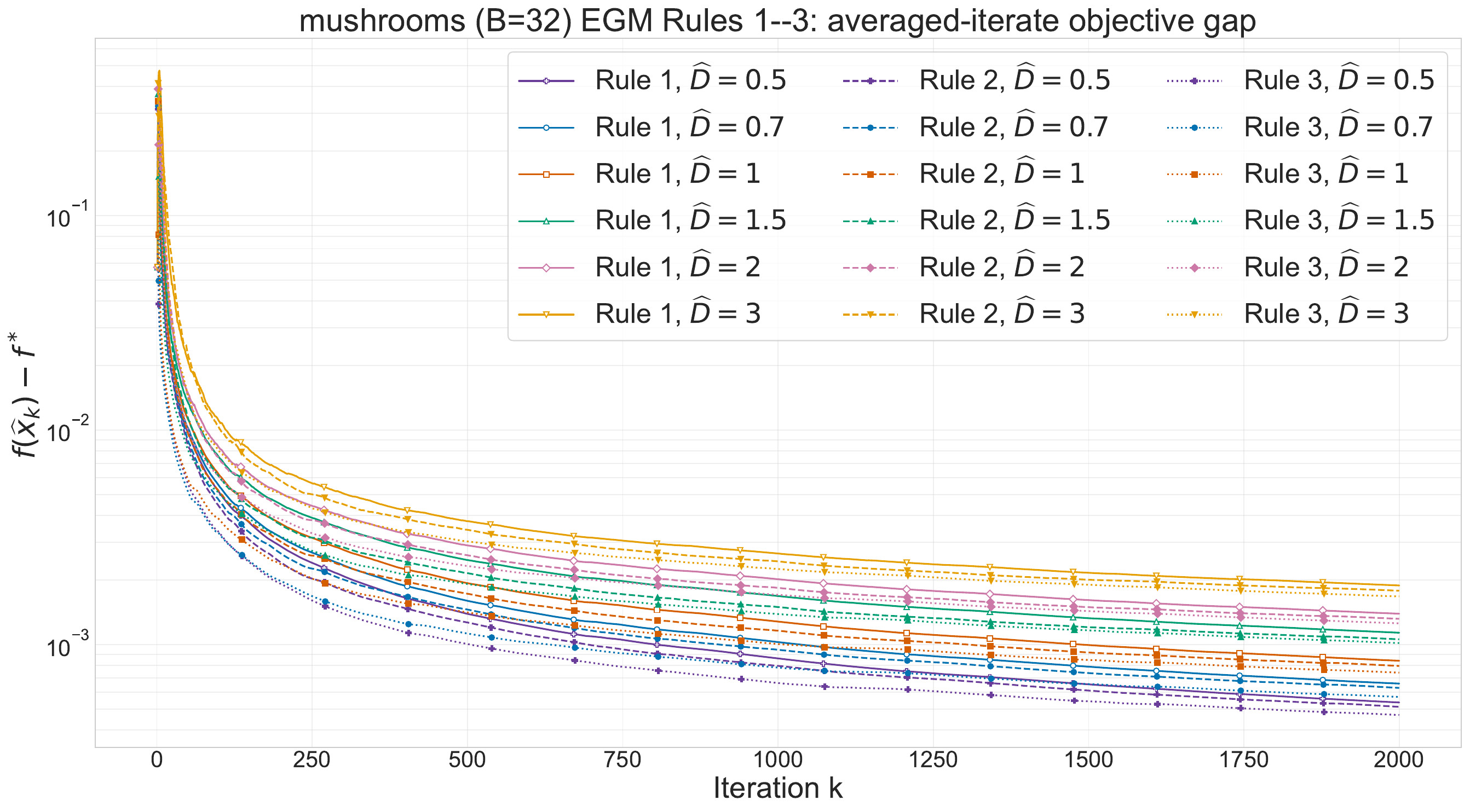}
\end{subfigure}
\hfill
\begin{subfigure}{0.24\textwidth}
    \centering
    \includegraphics[width=\linewidth]
    {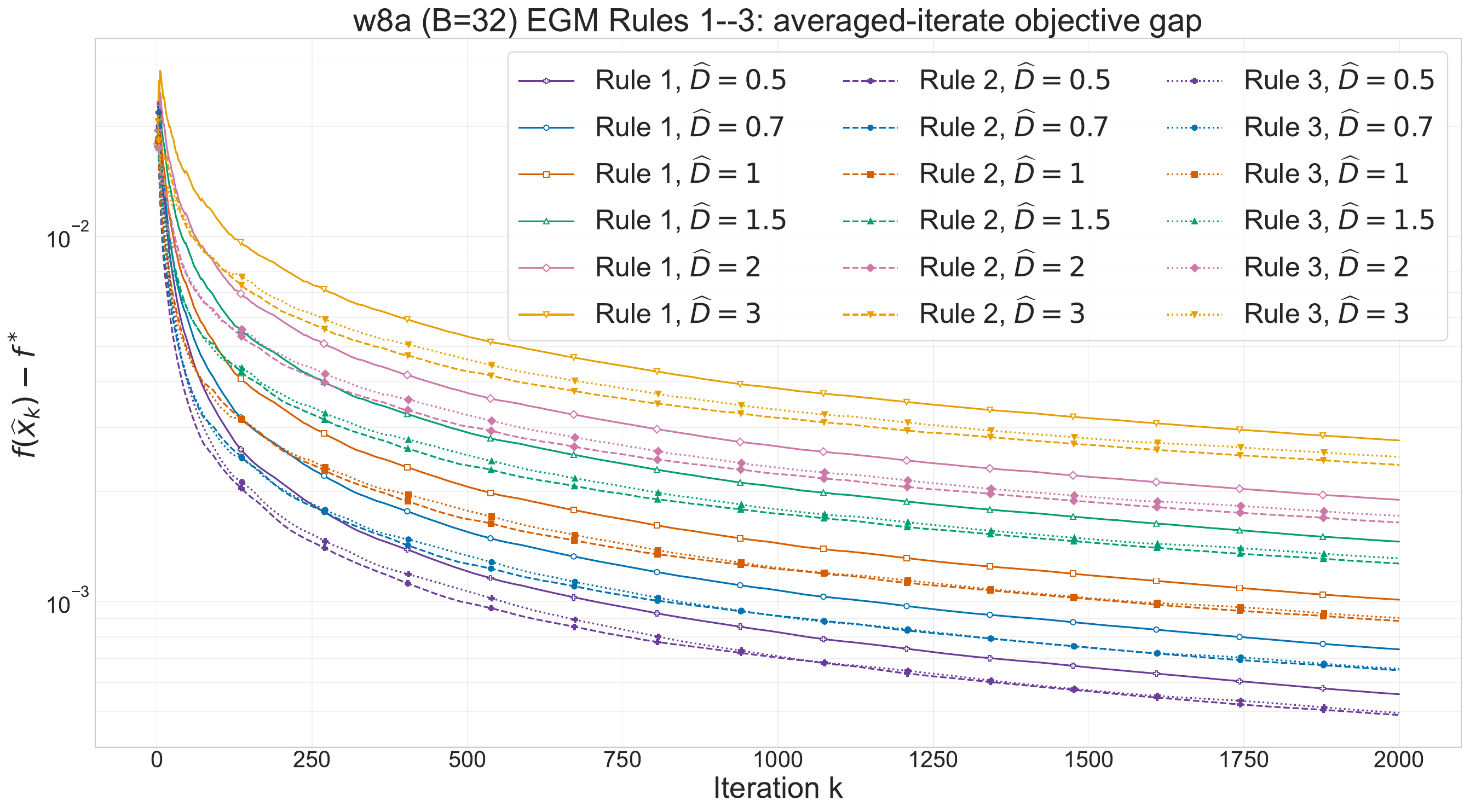}
\end{subfigure}
\hfill
\begin{subfigure}{0.24\textwidth}
    \centering
    \includegraphics[width=\linewidth]
    {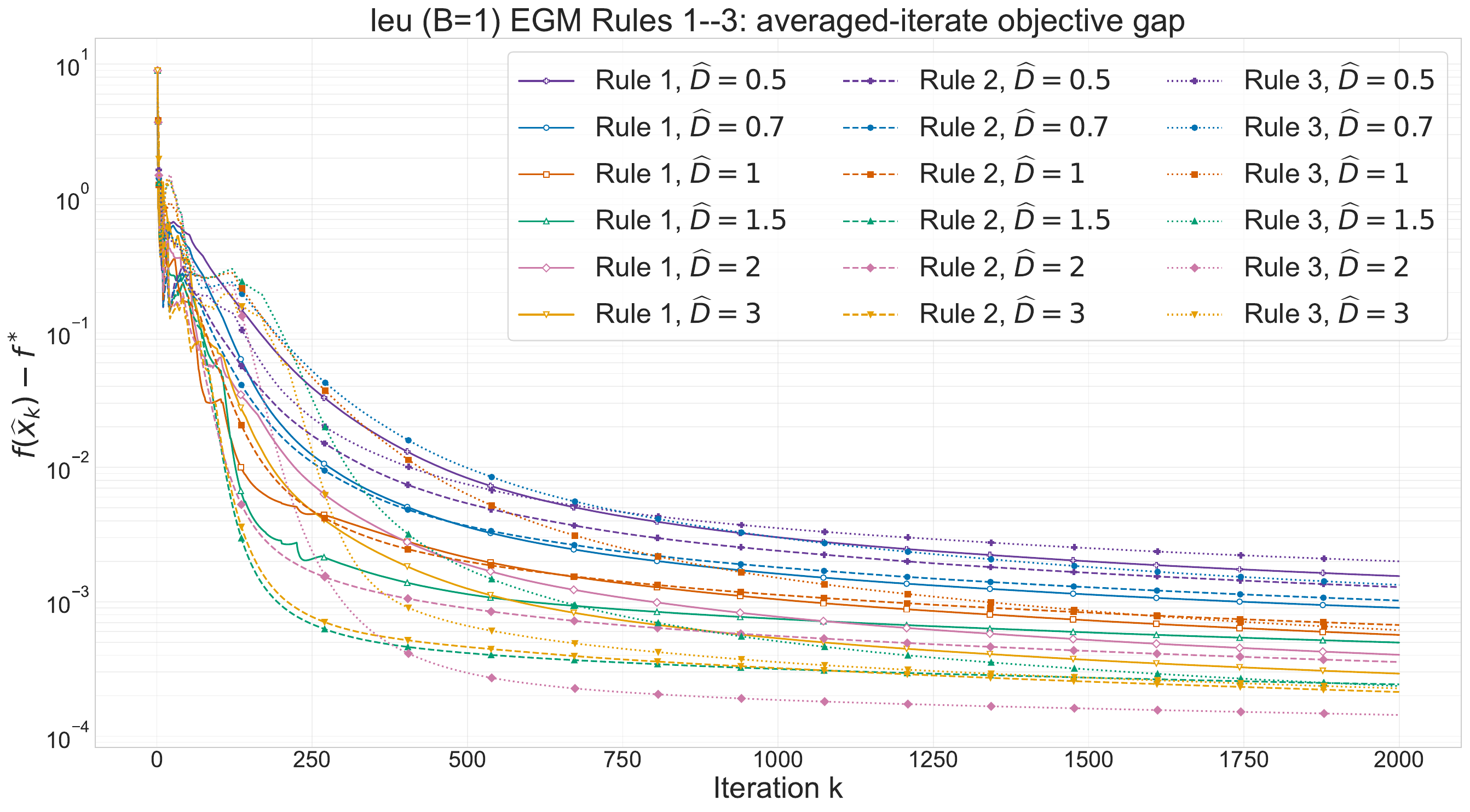}
\end{subfigure}
\hfill
\begin{subfigure}{0.24\textwidth}
    \centering
    \includegraphics[width=\linewidth]
    {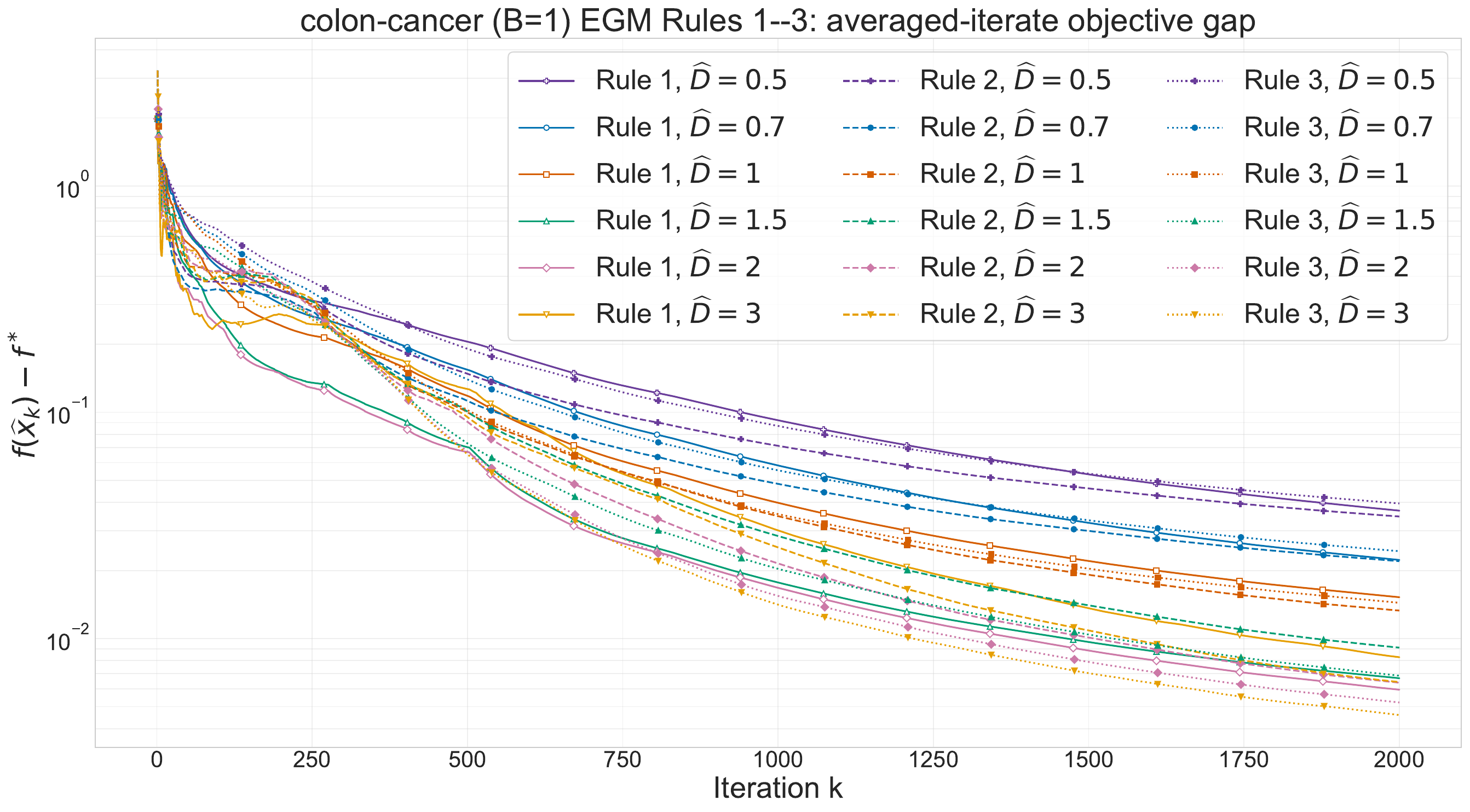}
\end{subfigure}
\par

\begin{subfigure}{0.24\textwidth}
    \centering
    \includegraphics[width=\linewidth]
    {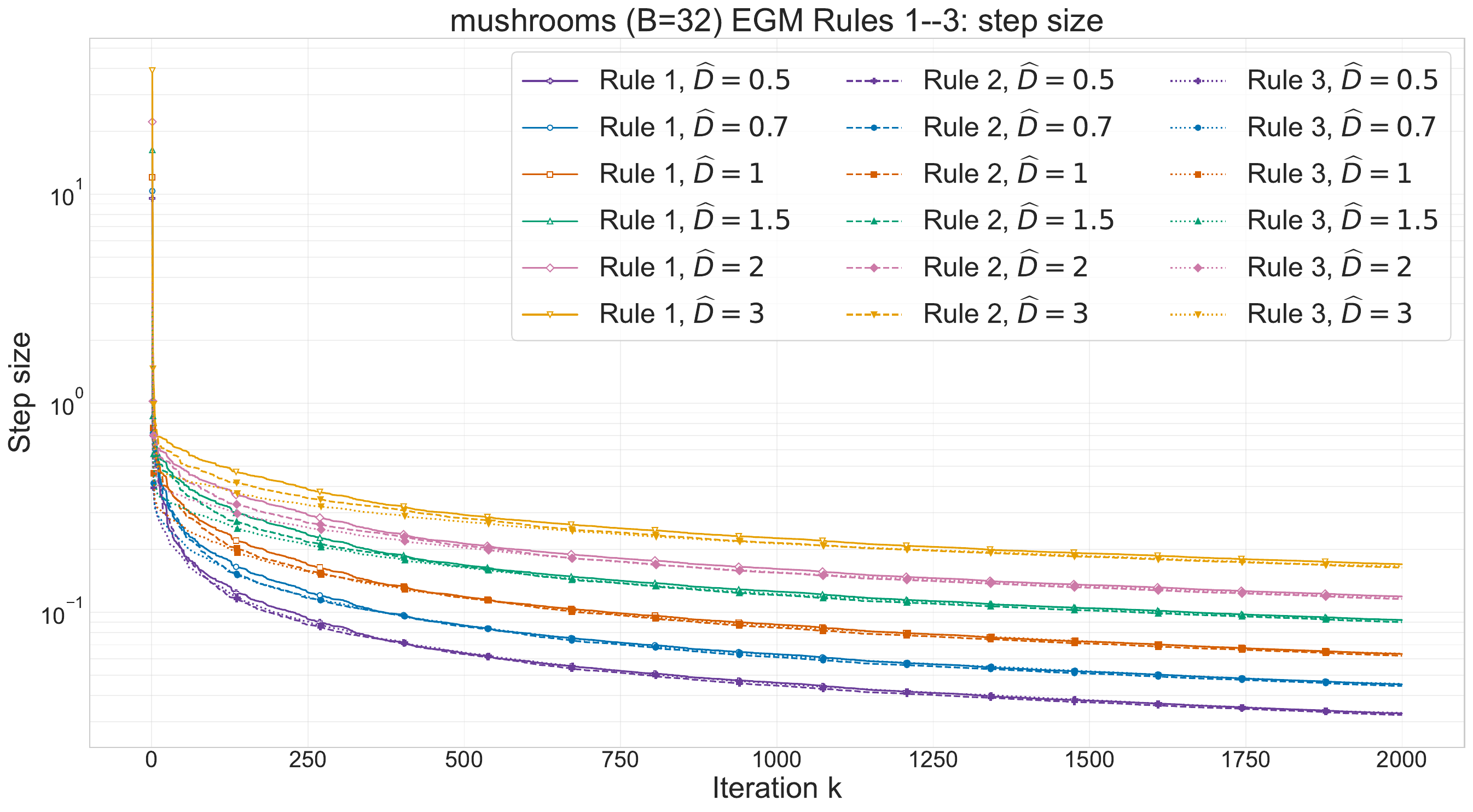}
\end{subfigure}
\hfill
\begin{subfigure}{0.24\textwidth}
    \centering
    \includegraphics[width=\linewidth]
    {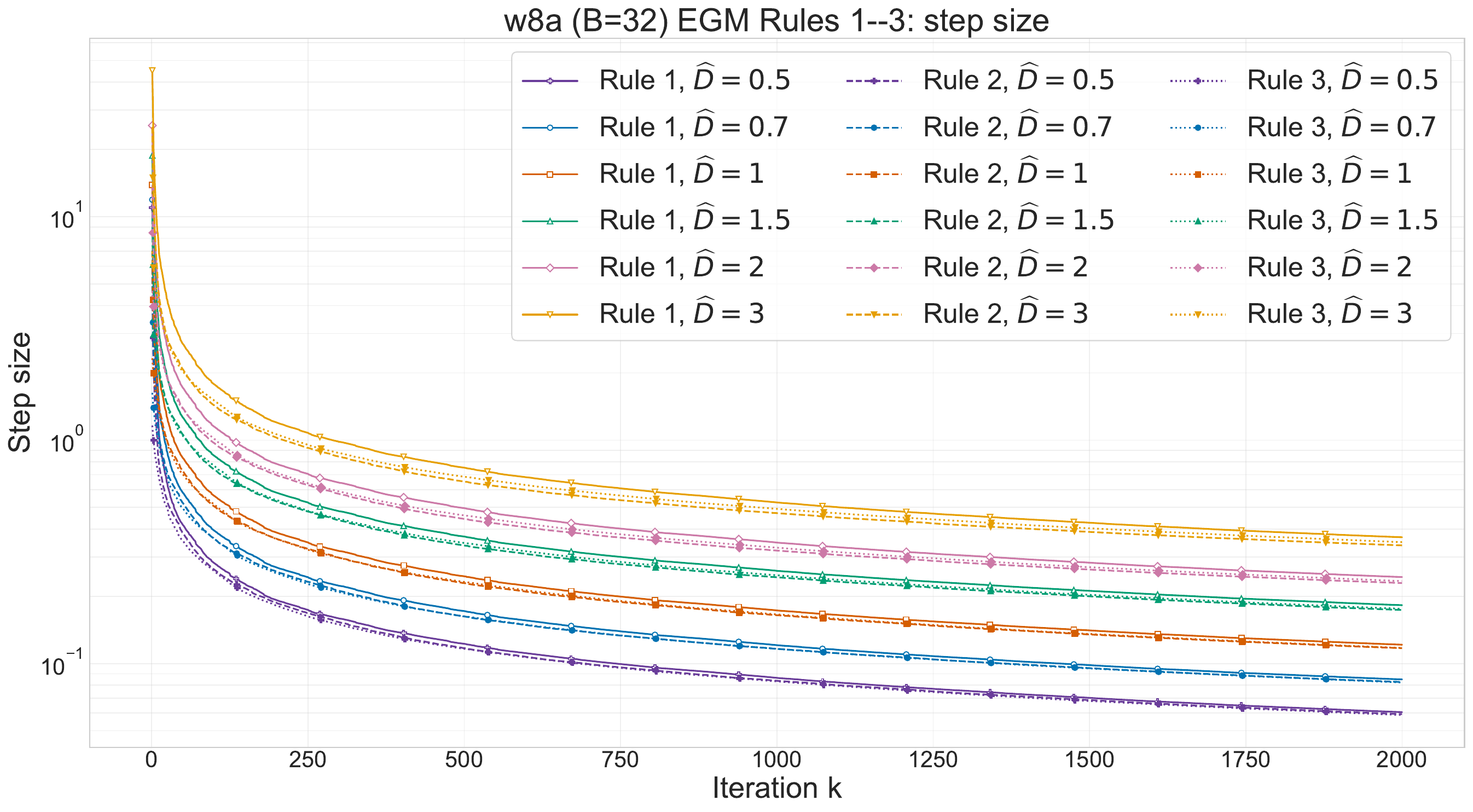}
\end{subfigure}
\hfill
\begin{subfigure}{0.24\textwidth}
    \centering
    \includegraphics[width=\linewidth]
    {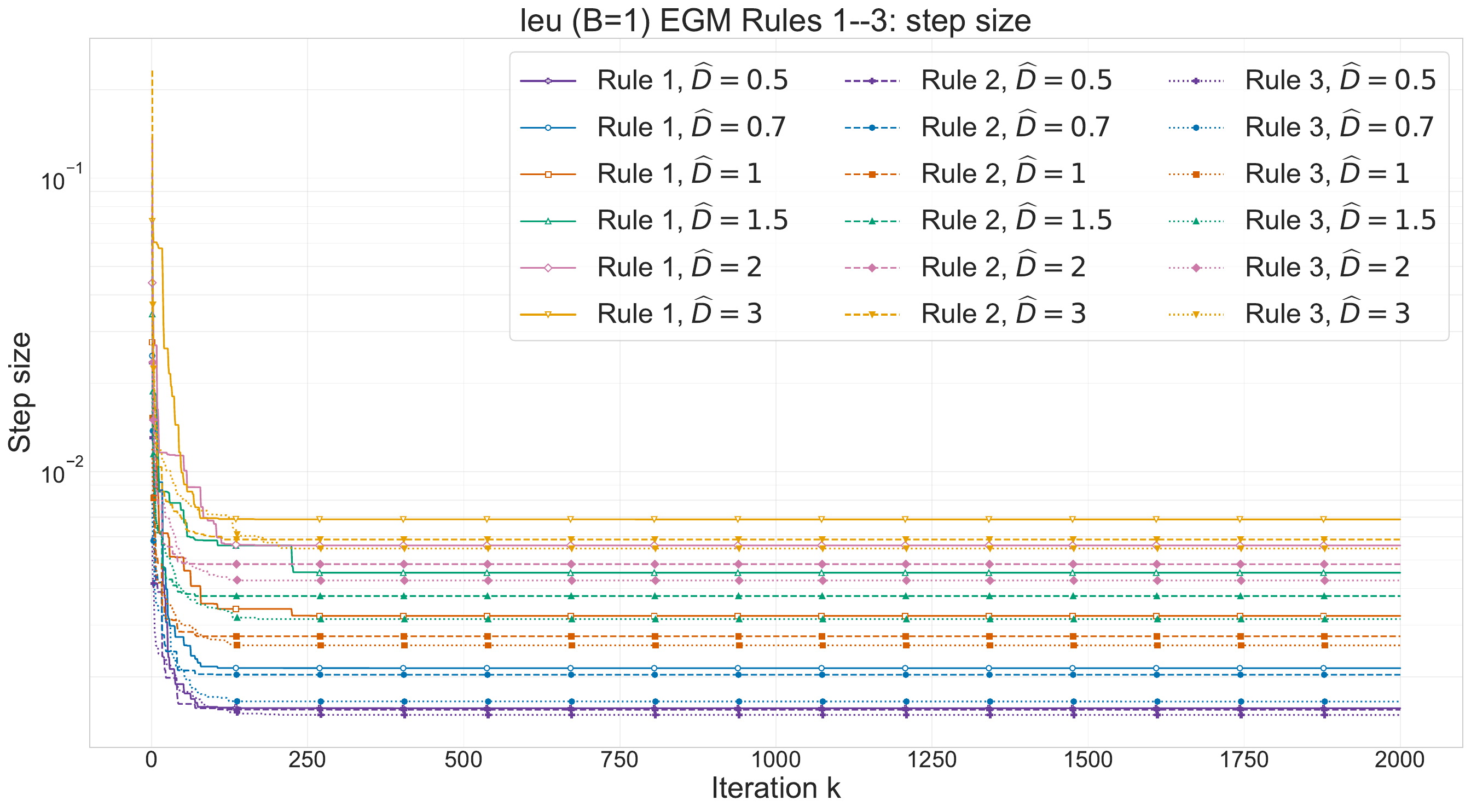}
\end{subfigure}
\hfill
\begin{subfigure}{0.24\textwidth}
    \centering
    \includegraphics[width=\linewidth]
    {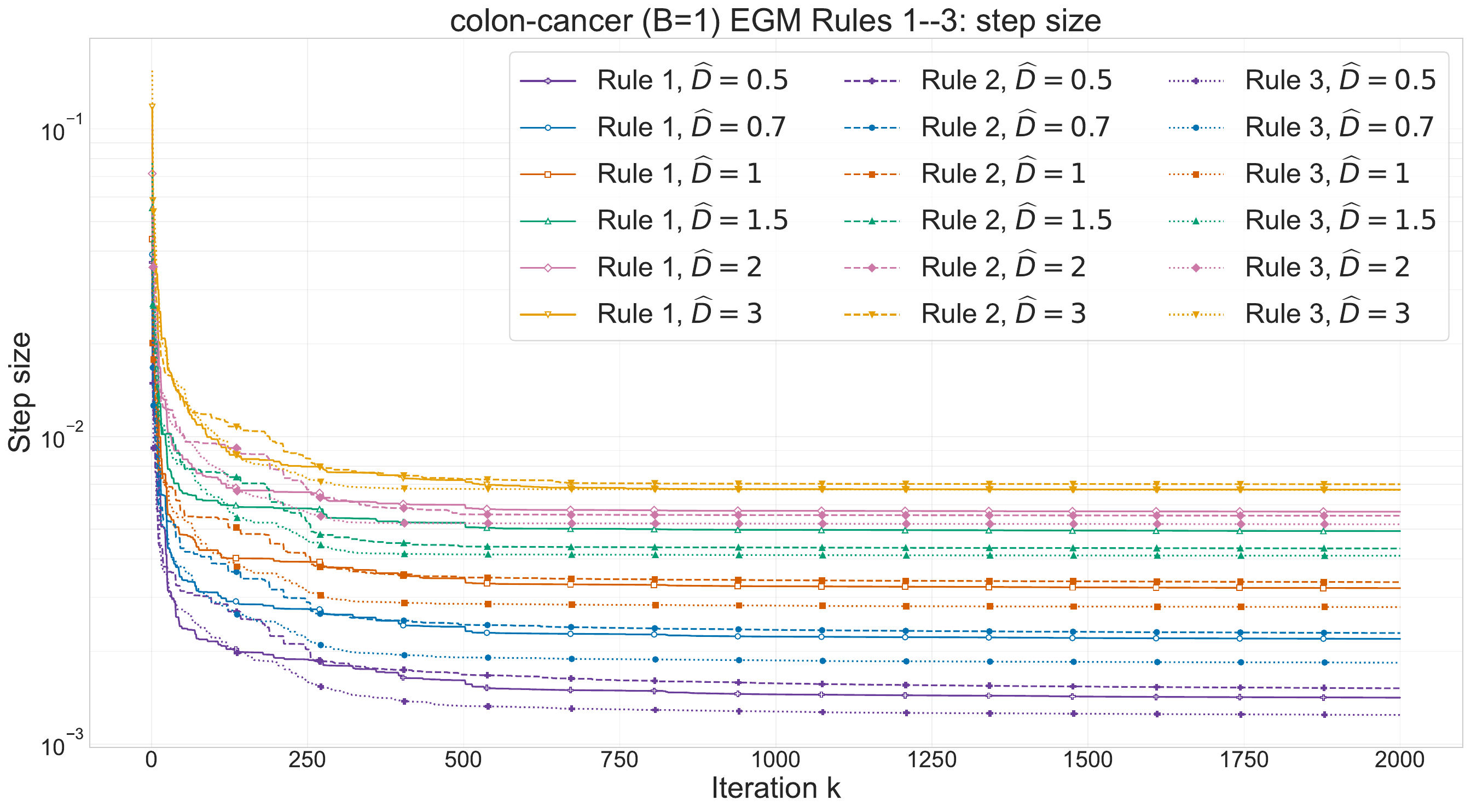}
\end{subfigure}

\caption{Sensitivity of EGM Rules 1--3 to $\widehat D$ for stochastic
logistic regression. From left to right, the columns correspond to
Mushrooms, W8a, Leu, and Colon Cancer. From top to bottom, the rows
show the best-iterate function values, averaged-iterate function
values, and step sizes.}
\label{fig:logistic_Dhat_sensitivity}
\end{figure*}


\begin{figure*}[htbp]
\centering

\begin{subfigure}{0.24\textwidth}
    \centering
    \includegraphics[width=\linewidth]
    {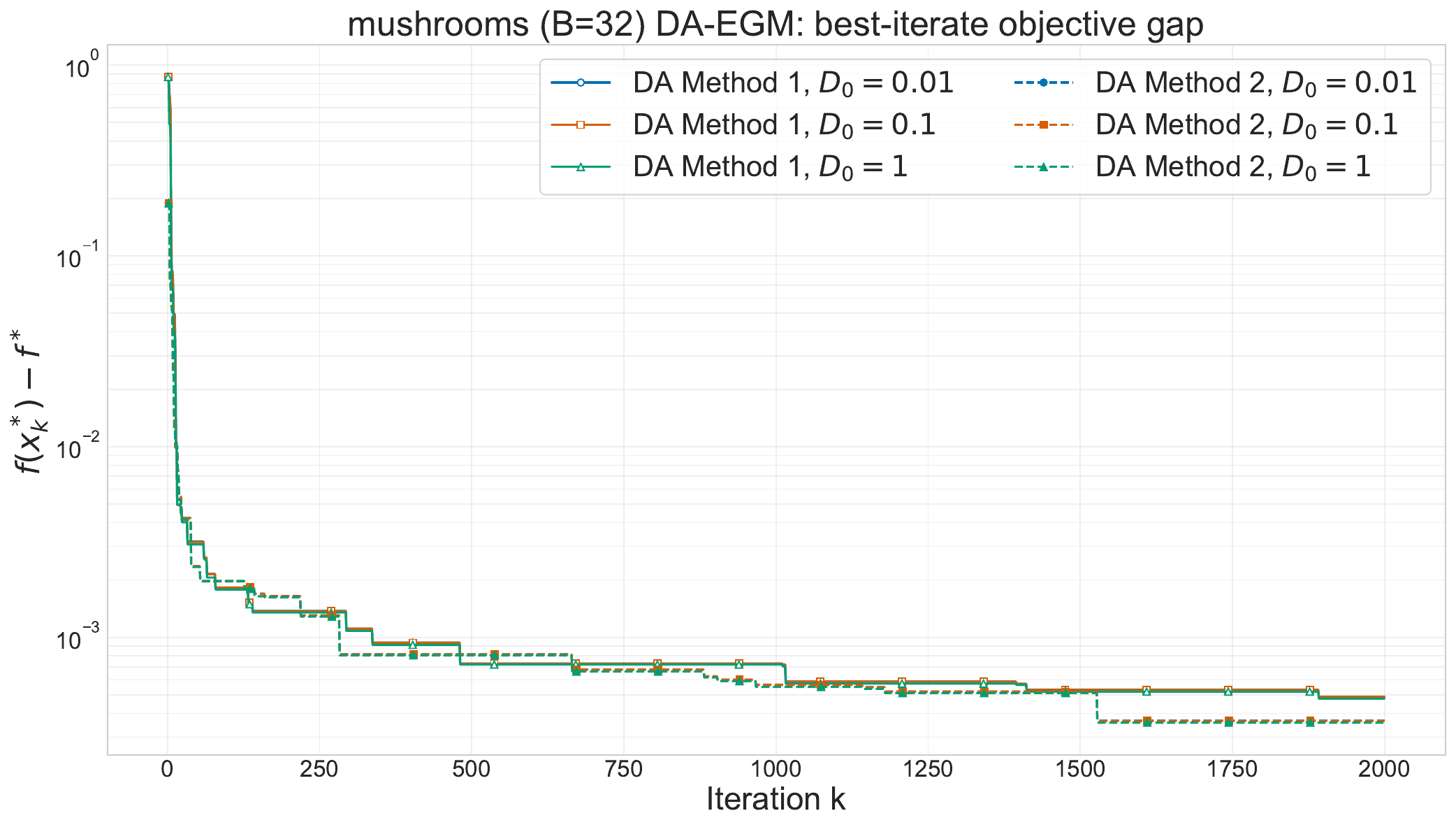}
\end{subfigure}
\hfill
\begin{subfigure}{0.24\textwidth}
    \centering
    \includegraphics[width=\linewidth]
    {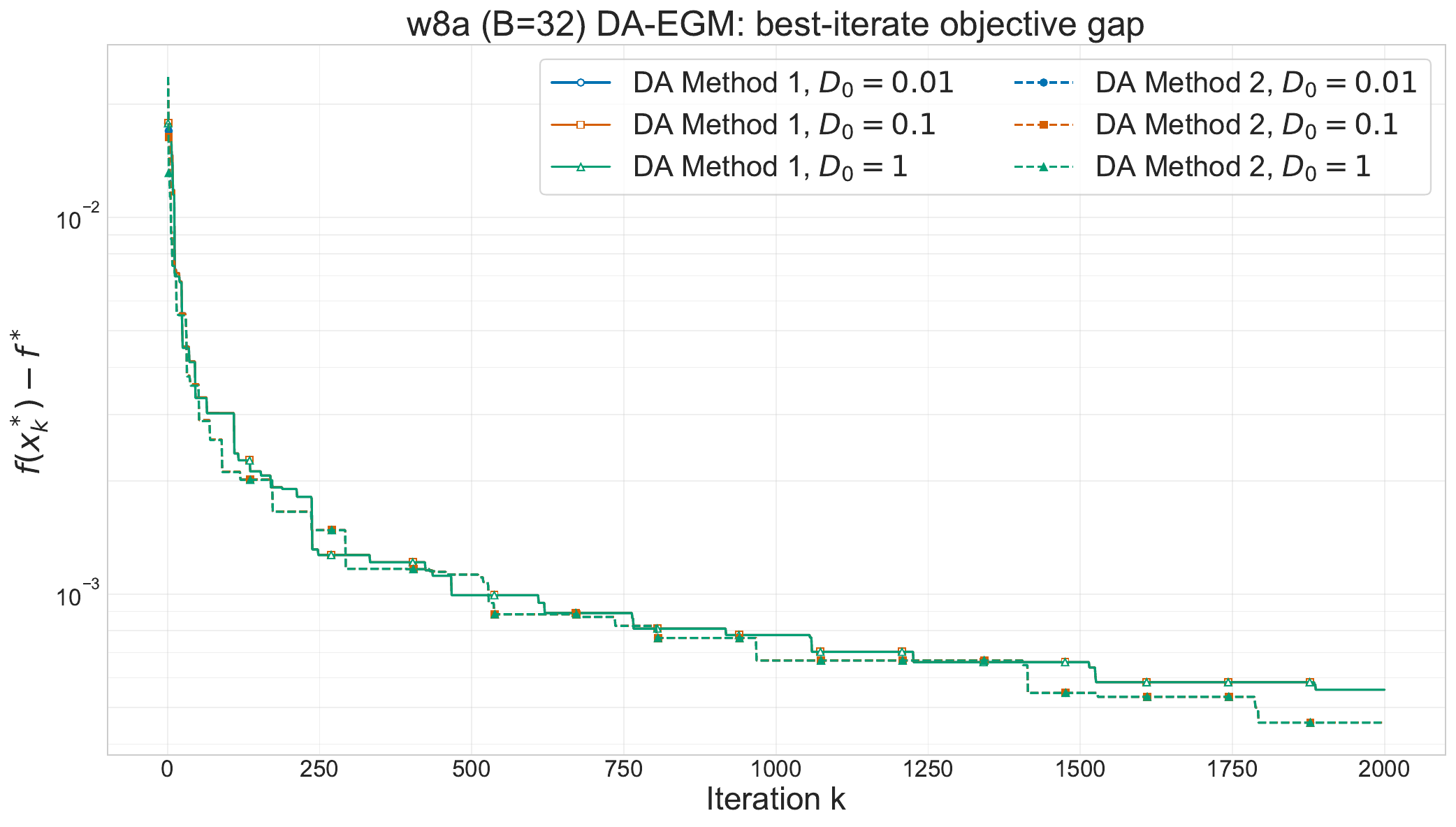}
\end{subfigure}
\hfill
\begin{subfigure}{0.24\textwidth}
    \centering
    \includegraphics[width=\linewidth]
    {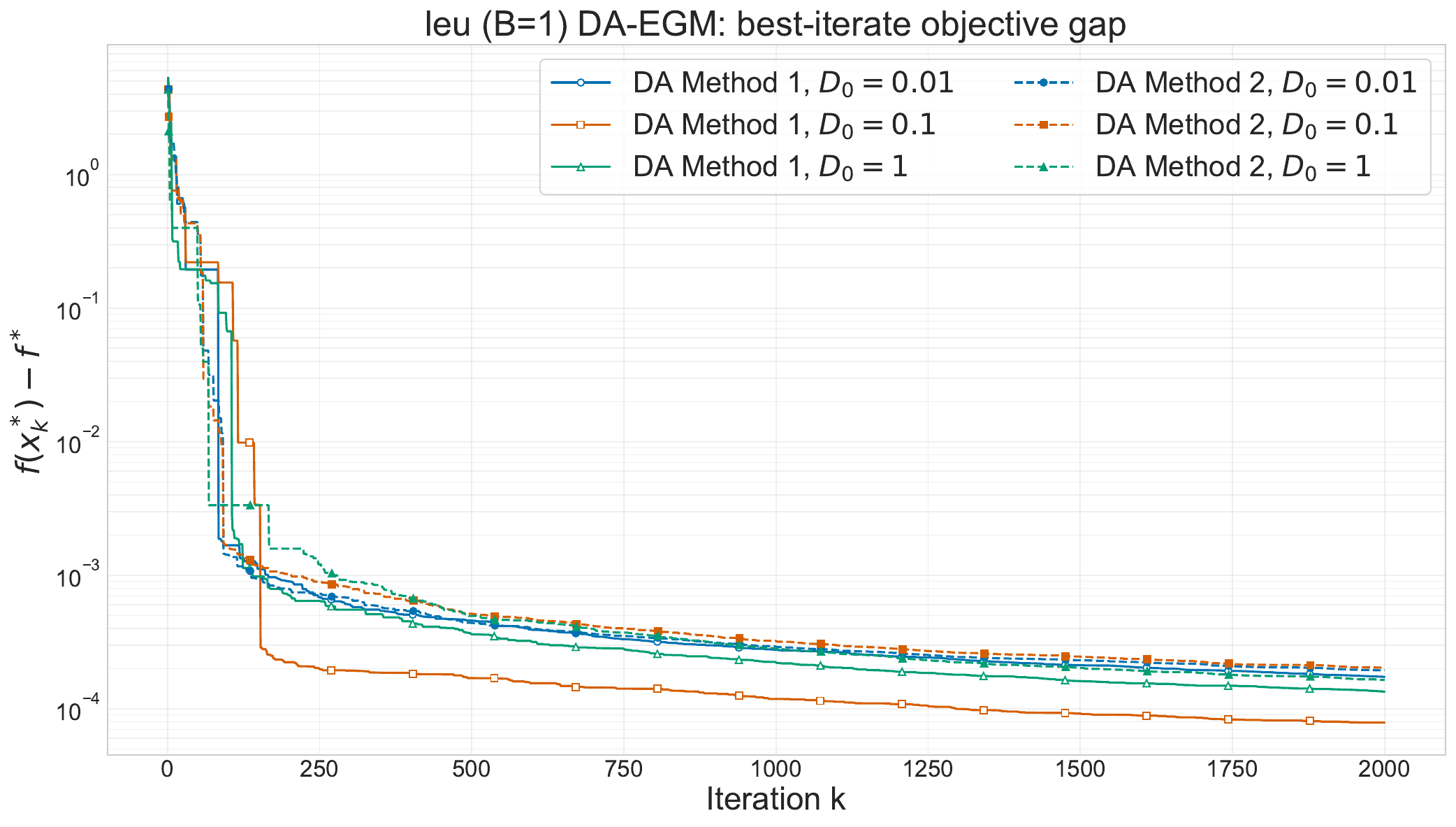}
\end{subfigure}
\hfill
\begin{subfigure}{0.24\textwidth}
    \centering
    \includegraphics[width=\linewidth]
    {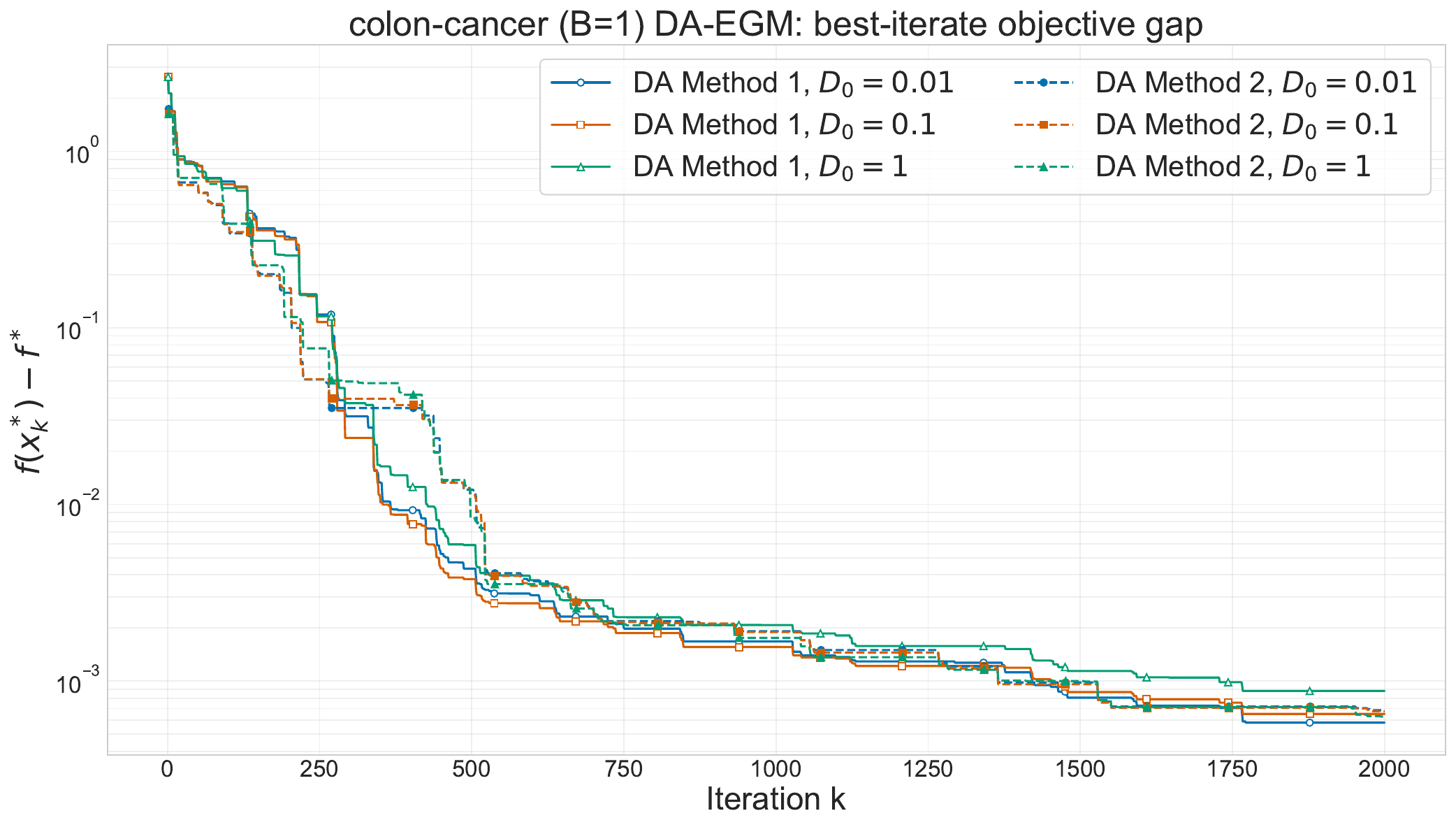}
\end{subfigure}
\par

\begin{subfigure}{0.24\textwidth}
    \centering
    \includegraphics[width=\linewidth]
    {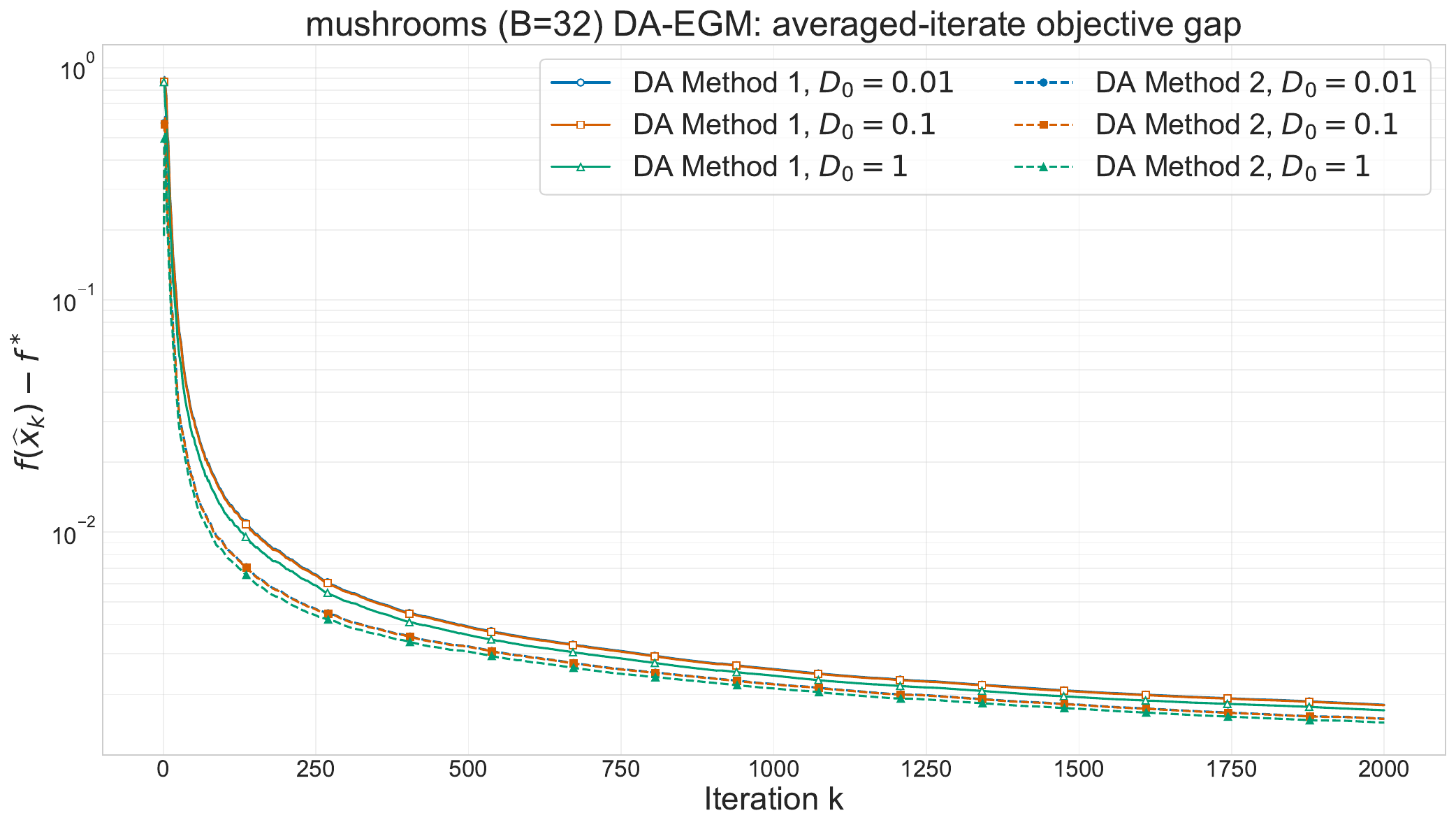}
\end{subfigure}
\hfill
\begin{subfigure}{0.24\textwidth}
    \centering
    \includegraphics[width=\linewidth]
    {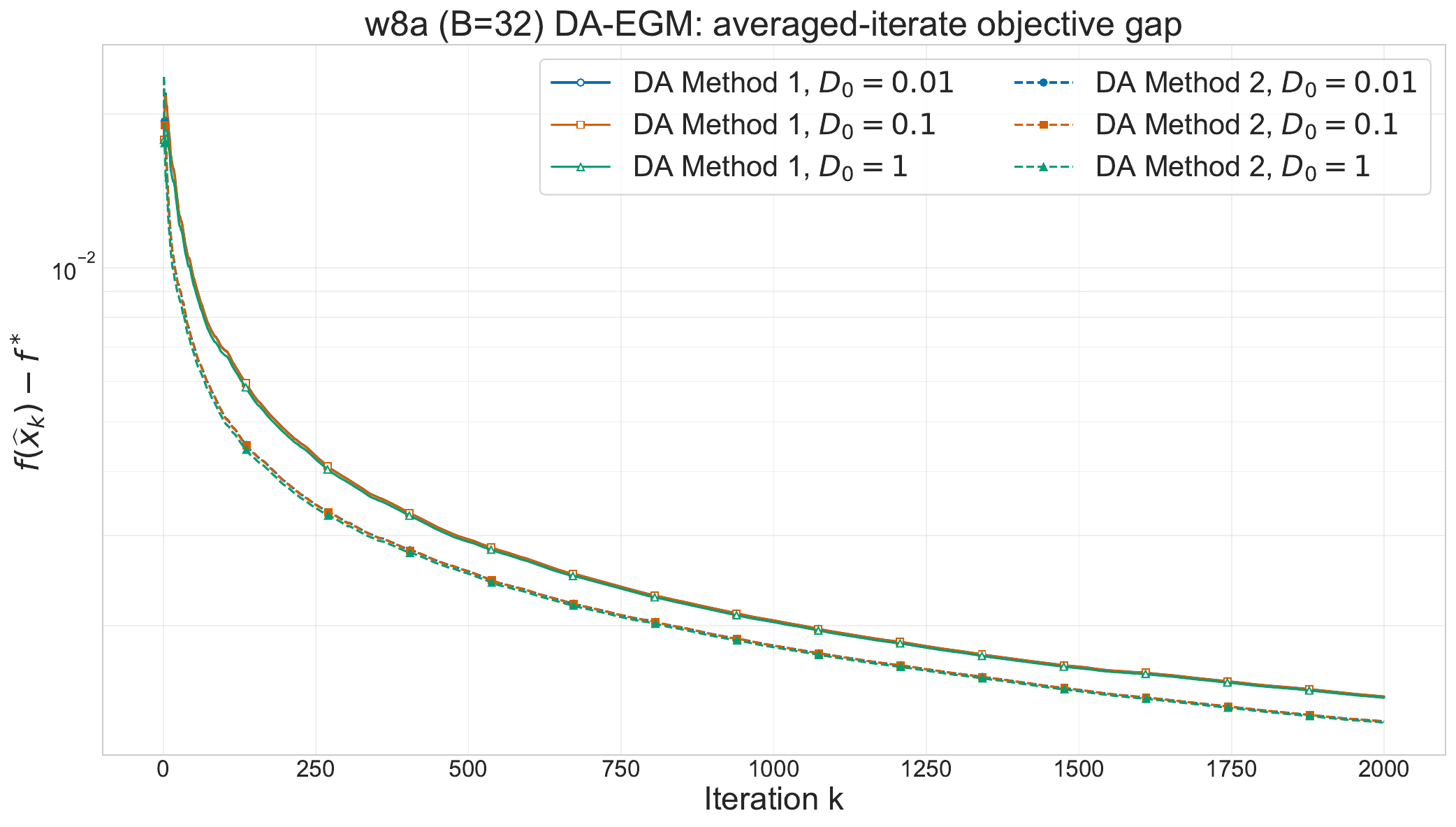}
\end{subfigure}
\hfill
\begin{subfigure}{0.24\textwidth}
    \centering
    \includegraphics[width=\linewidth]
    {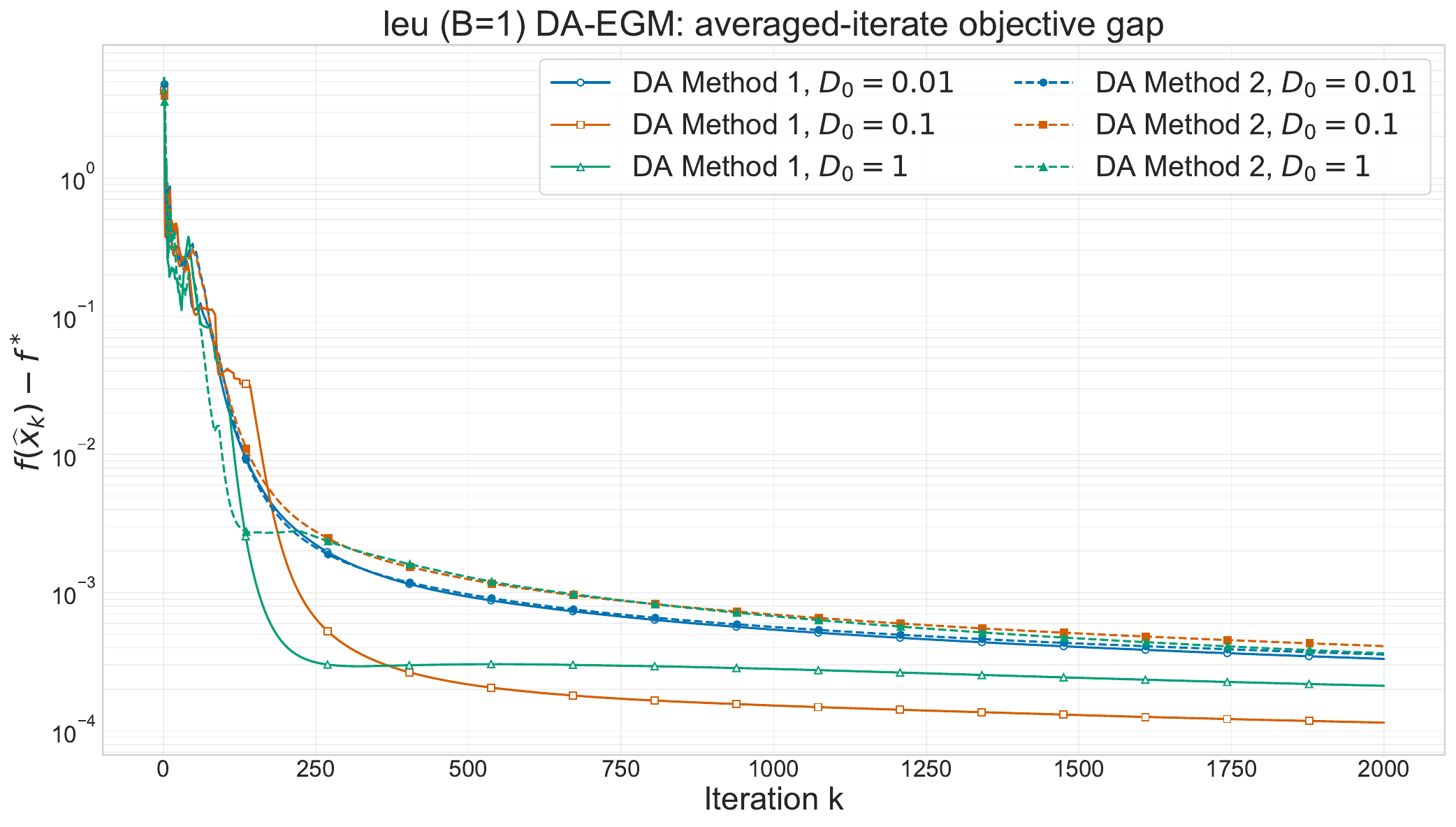}
\end{subfigure}
\hfill
\begin{subfigure}{0.24\textwidth}
    \centering
    \includegraphics[width=\linewidth]
    {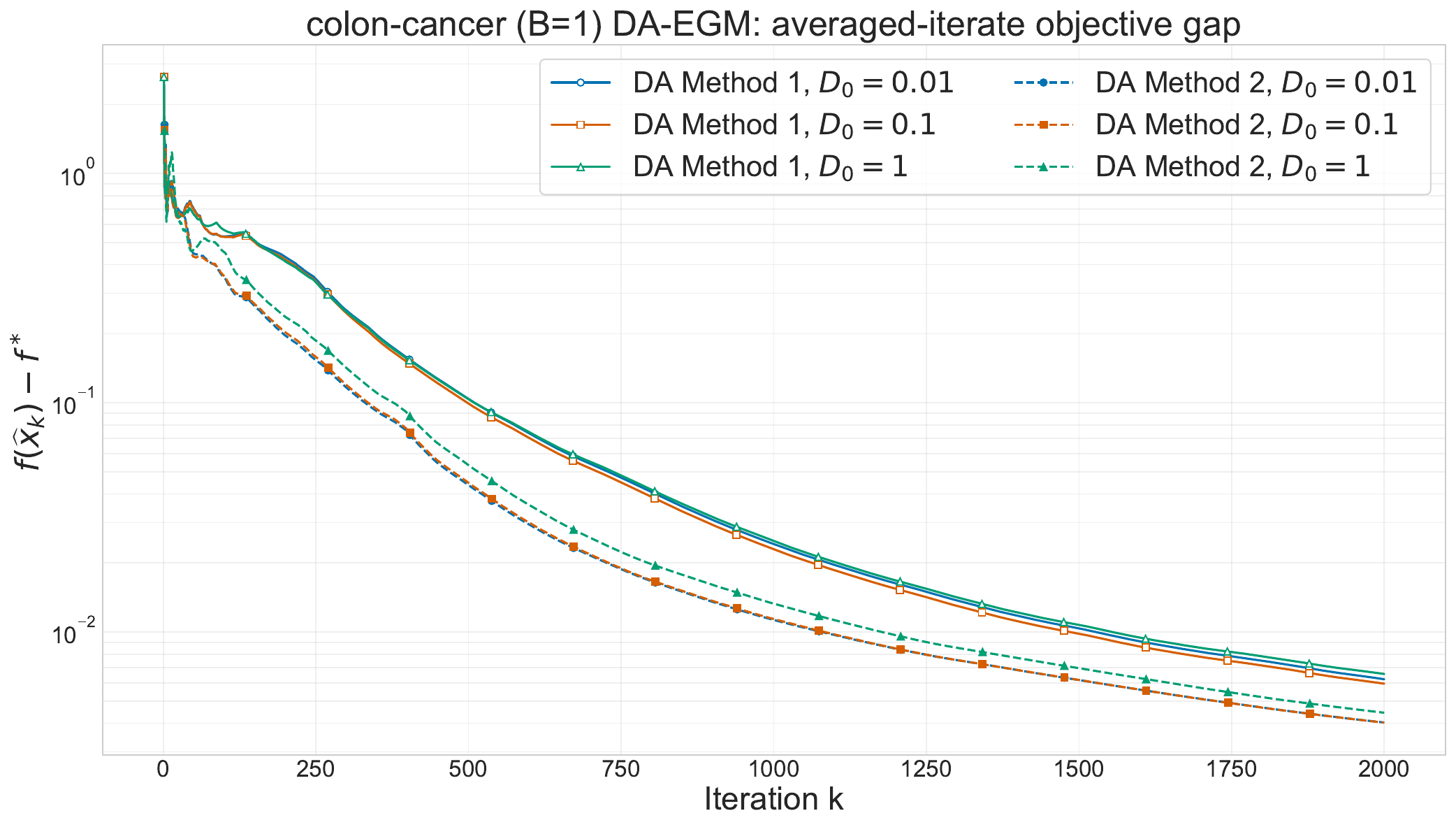}
\end{subfigure}
\par

\begin{subfigure}{0.24\textwidth}
    \centering
    \includegraphics[width=\linewidth]
    {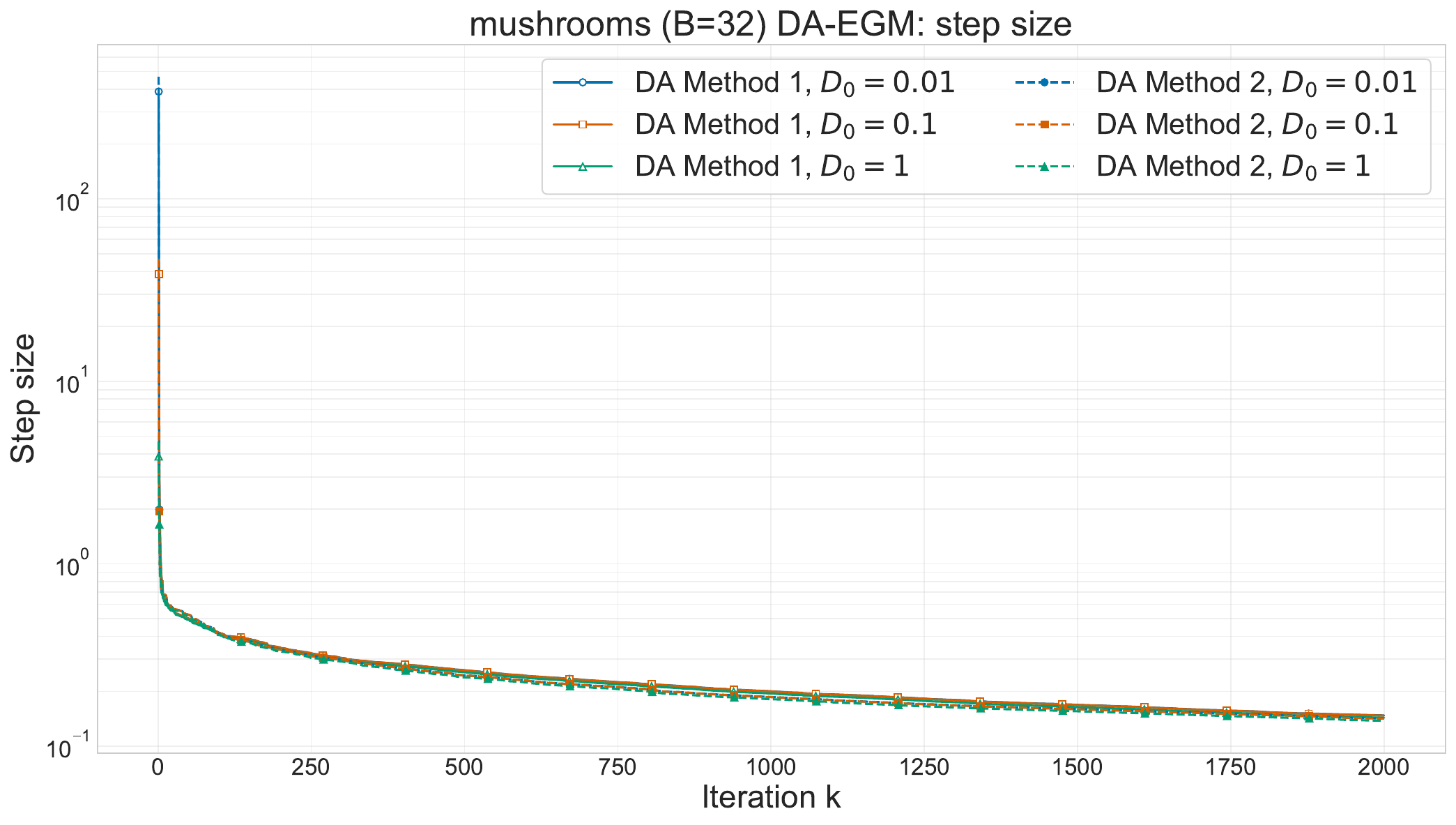}
\end{subfigure}
\hfill
\begin{subfigure}{0.24\textwidth}
    \centering
    \includegraphics[width=\linewidth]
    {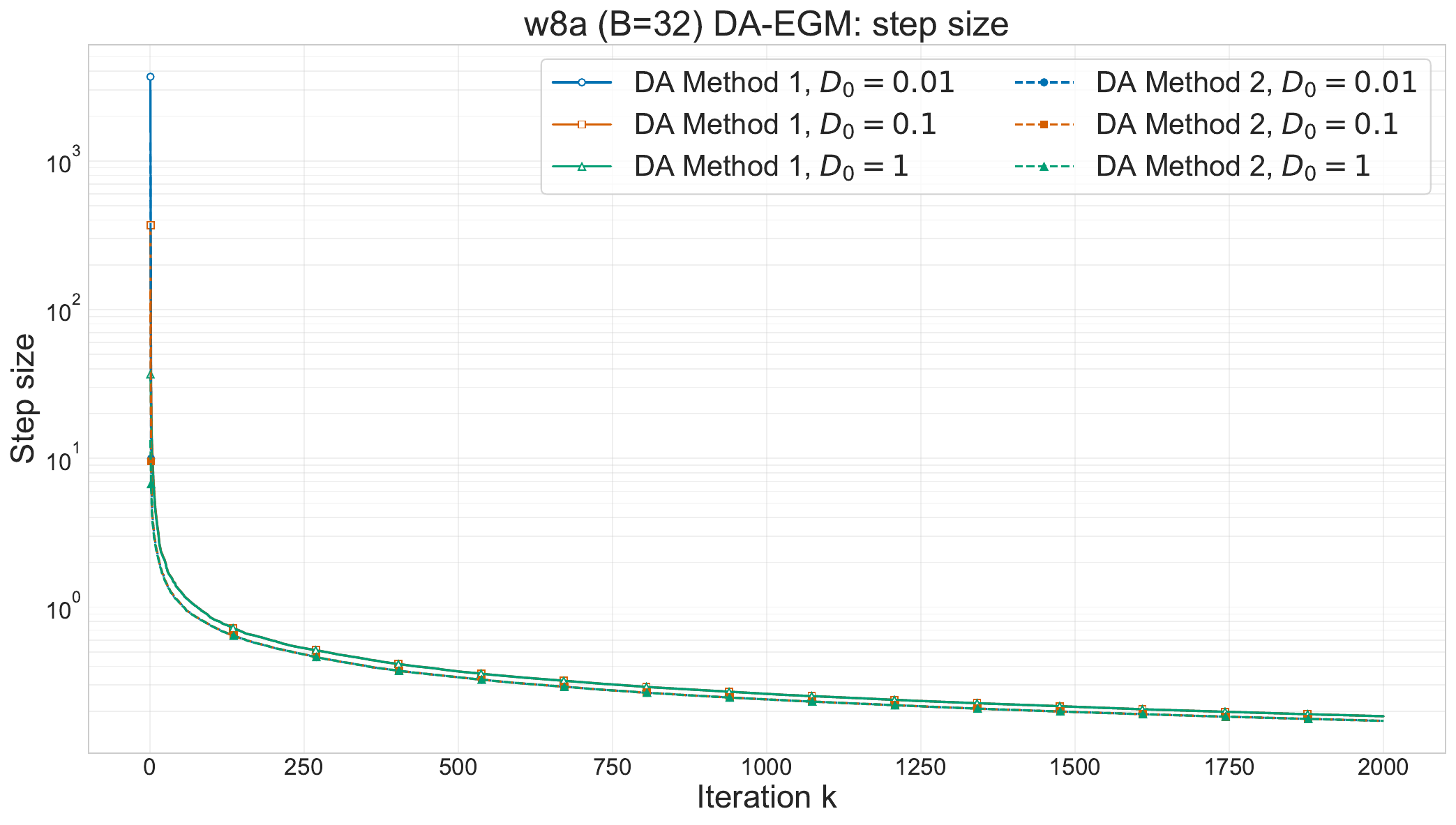}
\end{subfigure}
\hfill
\begin{subfigure}{0.24\textwidth}
    \centering
    \includegraphics[width=\linewidth]
    {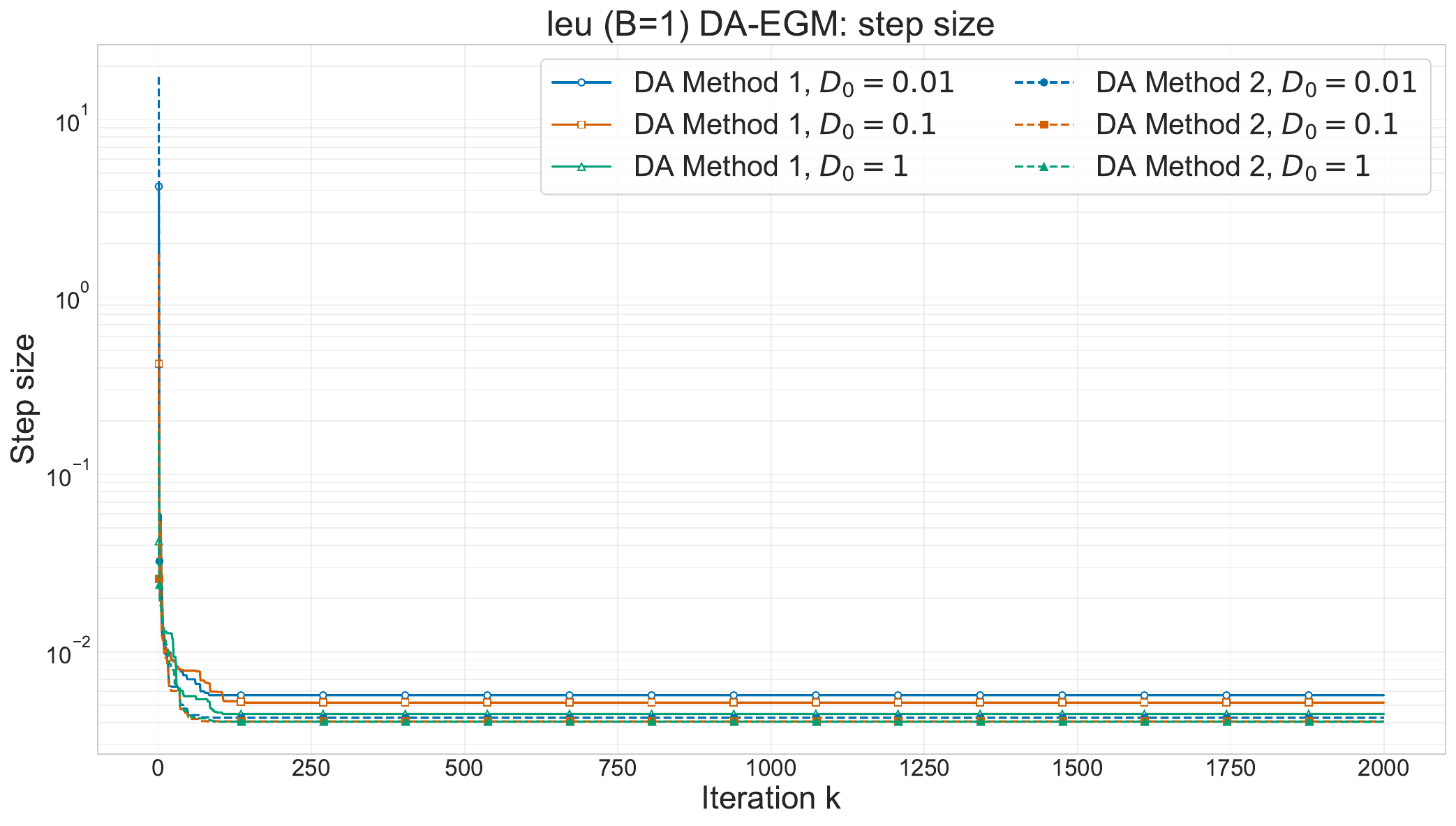}
\end{subfigure}
\hfill
\begin{subfigure}{0.24\textwidth}
    \centering
    \includegraphics[width=\linewidth]
    {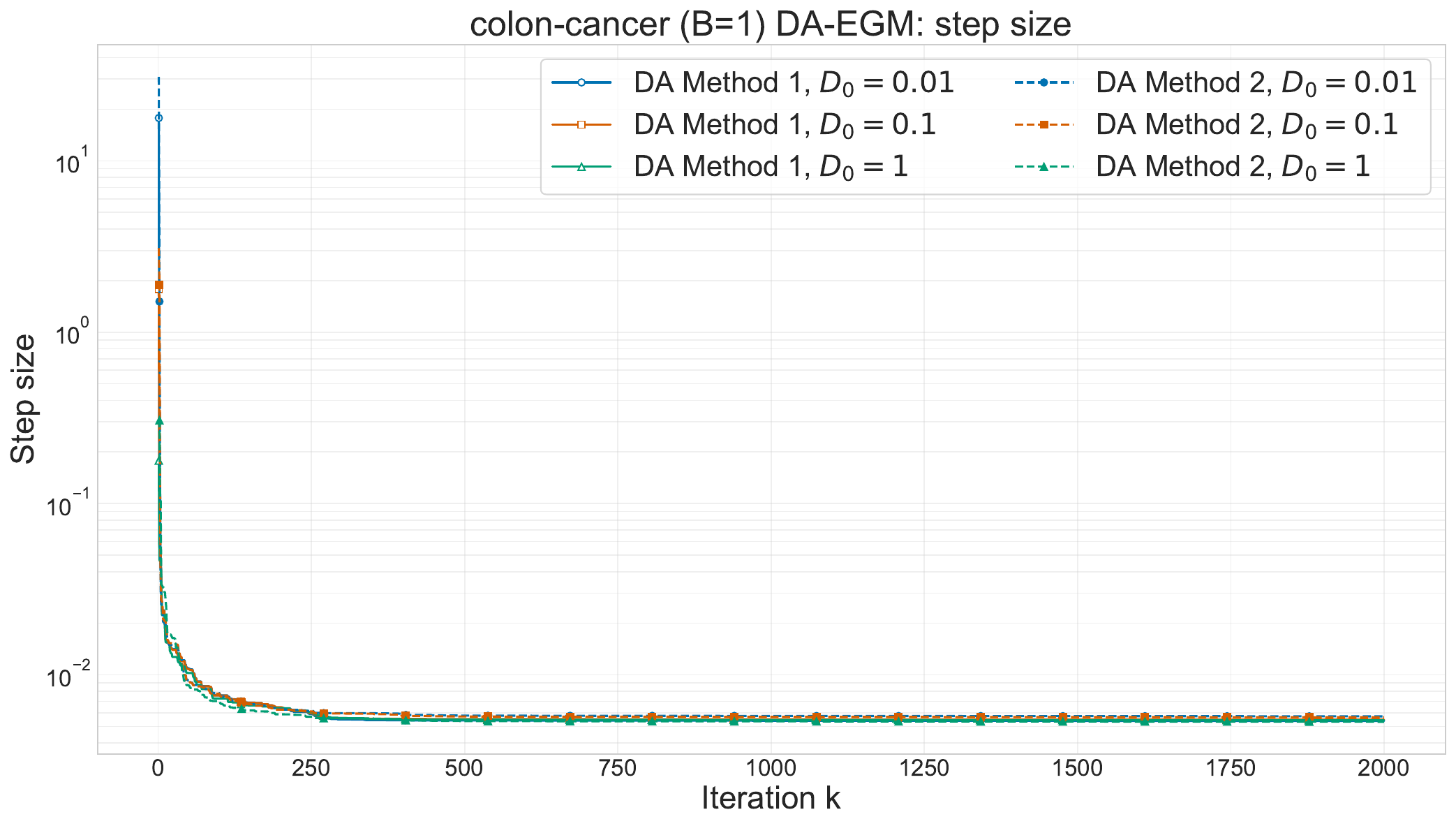}
\end{subfigure}

\caption{Sensitivity of DA-EGM Methods 1 and 2 to $D_0$ for stochastic
logistic regression. From left to right, the columns correspond to
Mushrooms, W8a, Leu, and Colon Cancer. From top to bottom, the rows
show the best-iterate function values, averaged-iterate function
values, and step sizes.}
\label{fig:logistic_D0_sensitivity}
\end{figure*}


\paragraph{Logistic-regression problem.}

We consider the stochastic constrained logistic-regression problem
\begin{equation}
    \min_{\|x\|_2\le 1} f(x) \quad \text{where} \quad f(x) 
    :=
    \frac{1}{n}
    \sum_{i=1}^{n}
    \log\left(
        1+\exp\left(-b_i\langle a_i,x\rangle\right)
    \right), \nonumber
\end{equation}
where \(a_i\in\mathbb{R}^d\) and \(b_i\in\{-1,+1\}\). We use the following datasets: (i) colon-cancer, (ii) leu, (iii) mushrooms, and (iv) w8a. We select mini-batches of sizes \(32\) for mushrooms and w8a,  and \(1\) for leu and colon-cancer. Since the exact optimal objective value $f^\star$ is unknown, for each dataset
we use the objective value returned by L-BFGS~\citep{liu1989limited} as an approximation to $f^\star$. We then use this approximation to compute the objective-function
gap.

The objective function gap values evaluated at both the best-iterate and the averaged-iterate and also the step sizes are presented in Figure~\ref{fig_logistic}. For the averaged-iterate plot, our algorithms perform better on the colon-cancer and leukemia datasets, whereas Projected AdaGrad-Norm performs better on the mushrooms and w8a datasets.

Figure~\ref{fig:logistic_Dhat_sensitivity} illustrates the sensitivity of EGM Rules~1--3 to the choice of $\widehat D$ across the different datasets. For the mushrooms and w8a datasets, all three EGM rules perform well with $\widehat D=0.5$. For the leu and colon-cancer datasets, Rule~3 performs well with $\widehat D=2$ and $\widehat D=3$, respectively. For the DA EGM methods, Figure~\ref{fig:logistic_D0_sensitivity} shows that all tested values of $D_0$ yield relatively similar performance across the datasets, although DA EGM Method~1 with $D_0=0.1$ performs slightly better on the leu dataset.




\begin{figure*}[t]
\centering


\begin{subfigure}{0.32\textwidth}
    \centering
    \includegraphics[width=\linewidth]{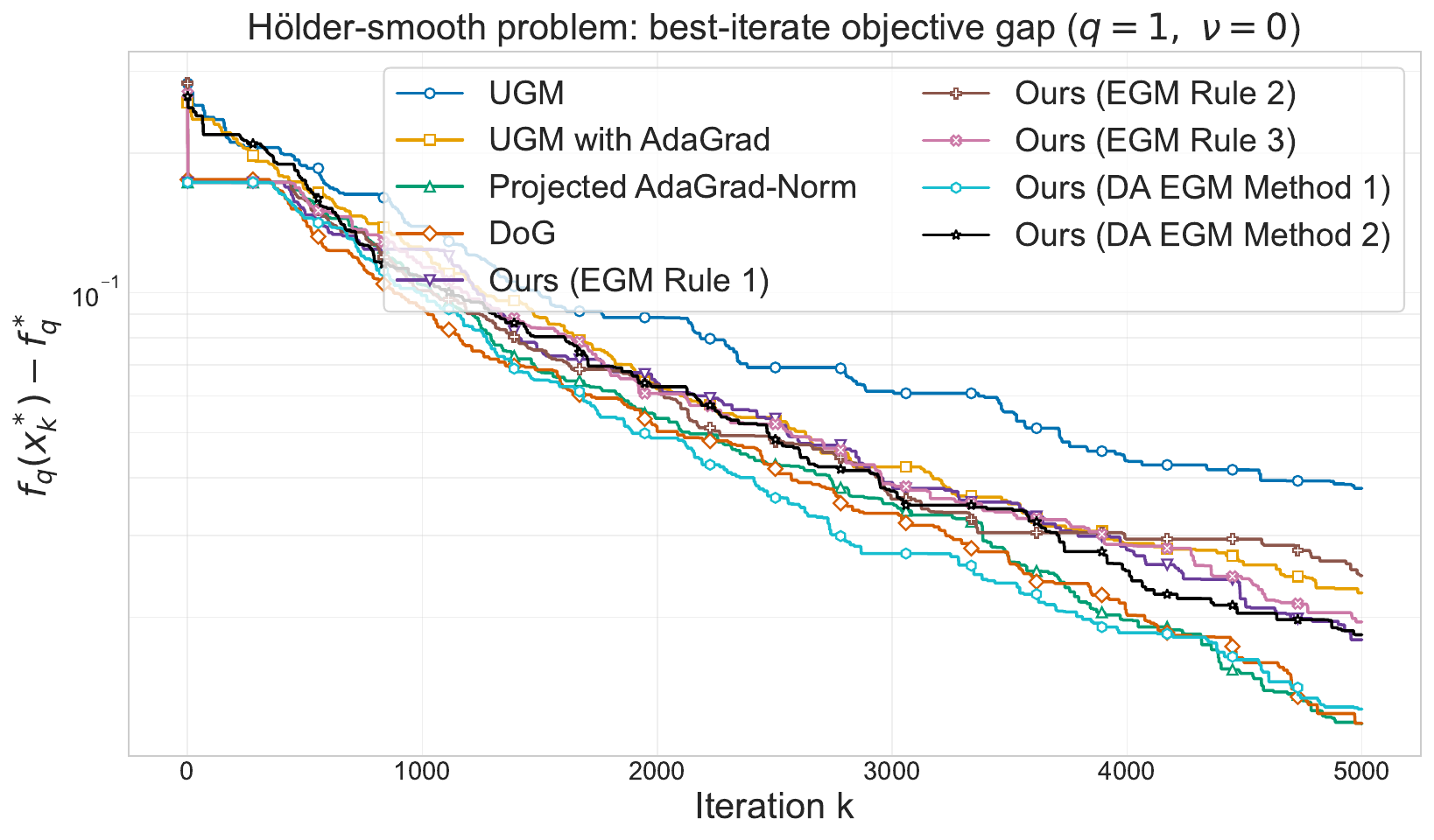}
\end{subfigure}
\hfill
\begin{subfigure}{0.32\textwidth}
    \centering
    \includegraphics[width=\linewidth]{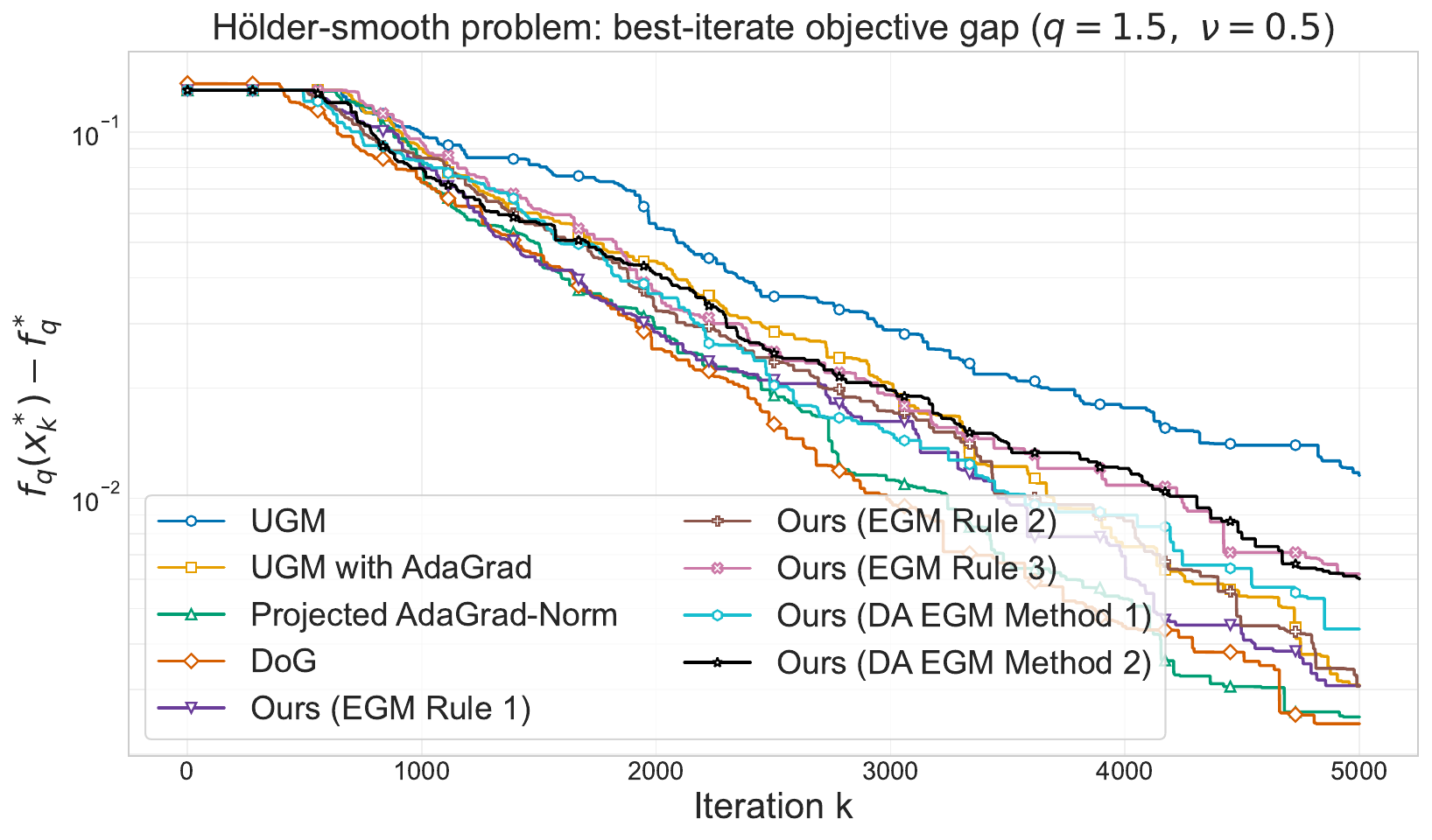}
\end{subfigure}
\hfill
\begin{subfigure}{0.32\textwidth}
    \centering
    \includegraphics[width=\linewidth]{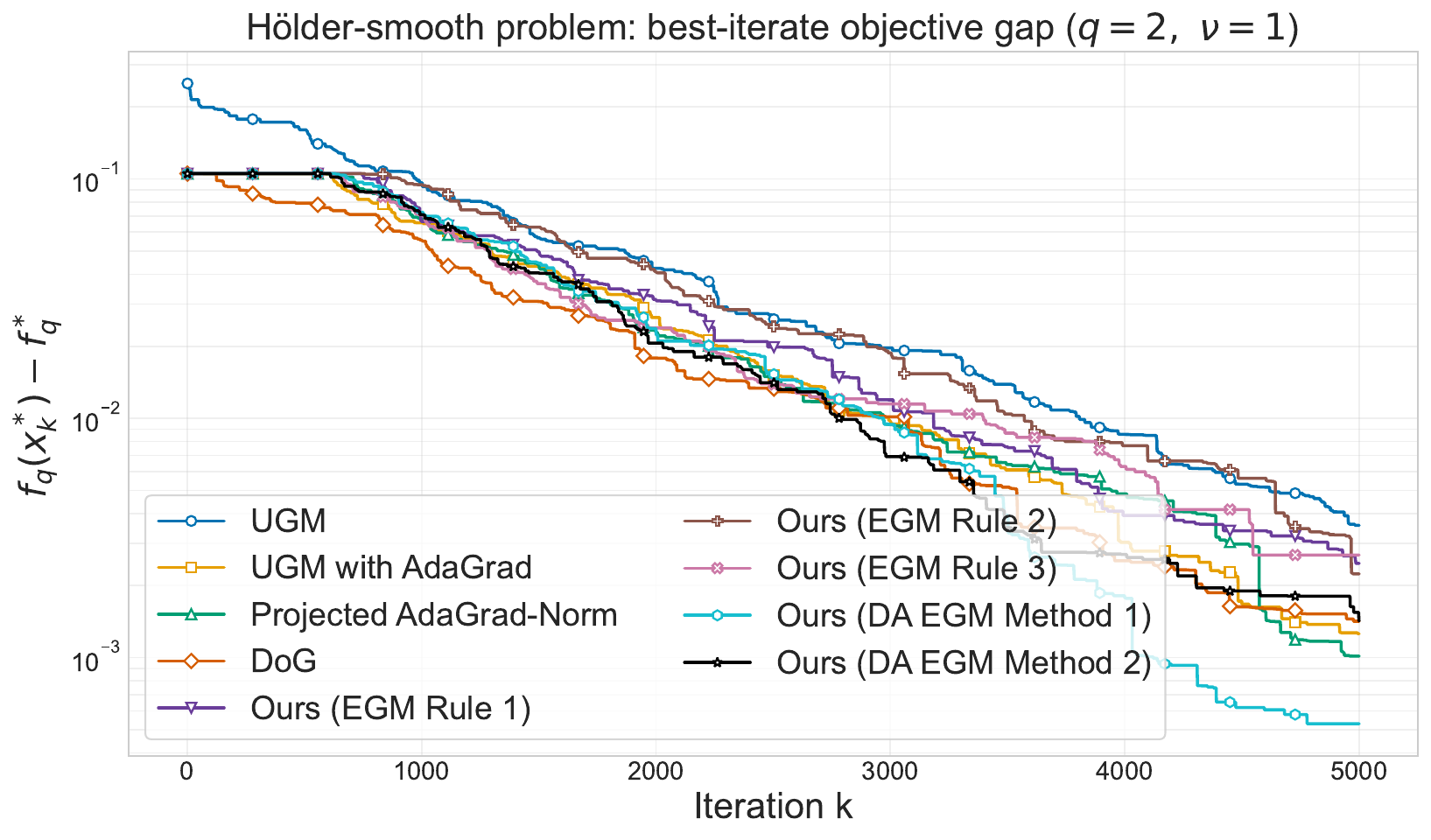}
\end{subfigure}
\par

\begin{subfigure}{0.32\textwidth}
    \centering
    \includegraphics[width=\linewidth]{pics/stoc_holder/stochastic_holder_avg_q_1p0.pdf}
\end{subfigure}
\hfill
\begin{subfigure}{0.32\textwidth}
    \centering
    \includegraphics[width=\linewidth]{pics/stoc_holder/stochastic_holder_avg_q_1p5.pdf}
\end{subfigure}
\hfill
\begin{subfigure}{0.32\textwidth}
    \centering
    \includegraphics[width=\linewidth]{pics/stoc_holder/stochastic_holder_avg_q_2p0.pdf}
\end{subfigure}
\par

\begin{subfigure}{0.32\textwidth}
    \centering
    \includegraphics[width=\linewidth]{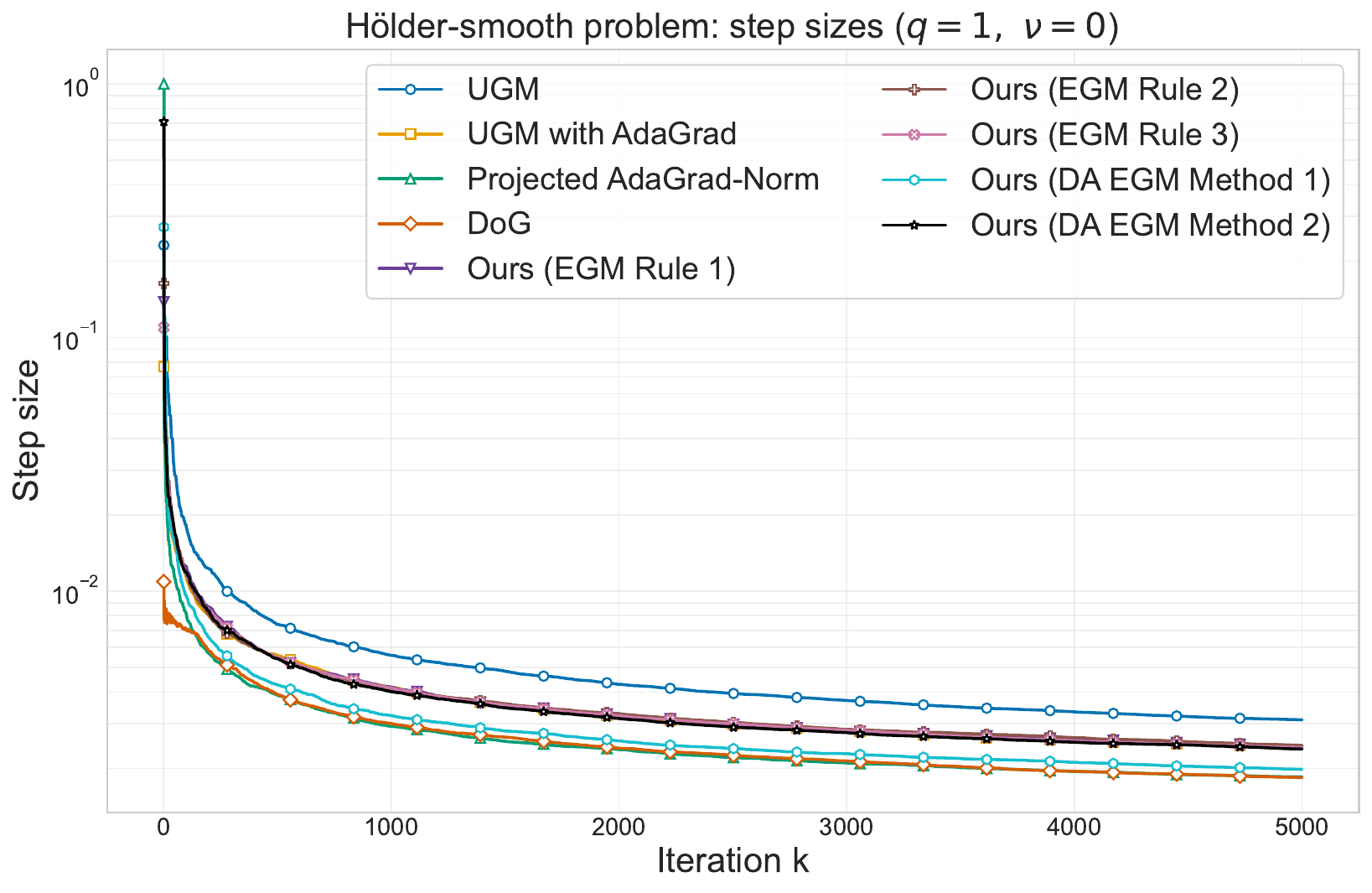}
\end{subfigure}
\hfill
\begin{subfigure}{0.32\textwidth}
    \centering
    \includegraphics[width=\linewidth]{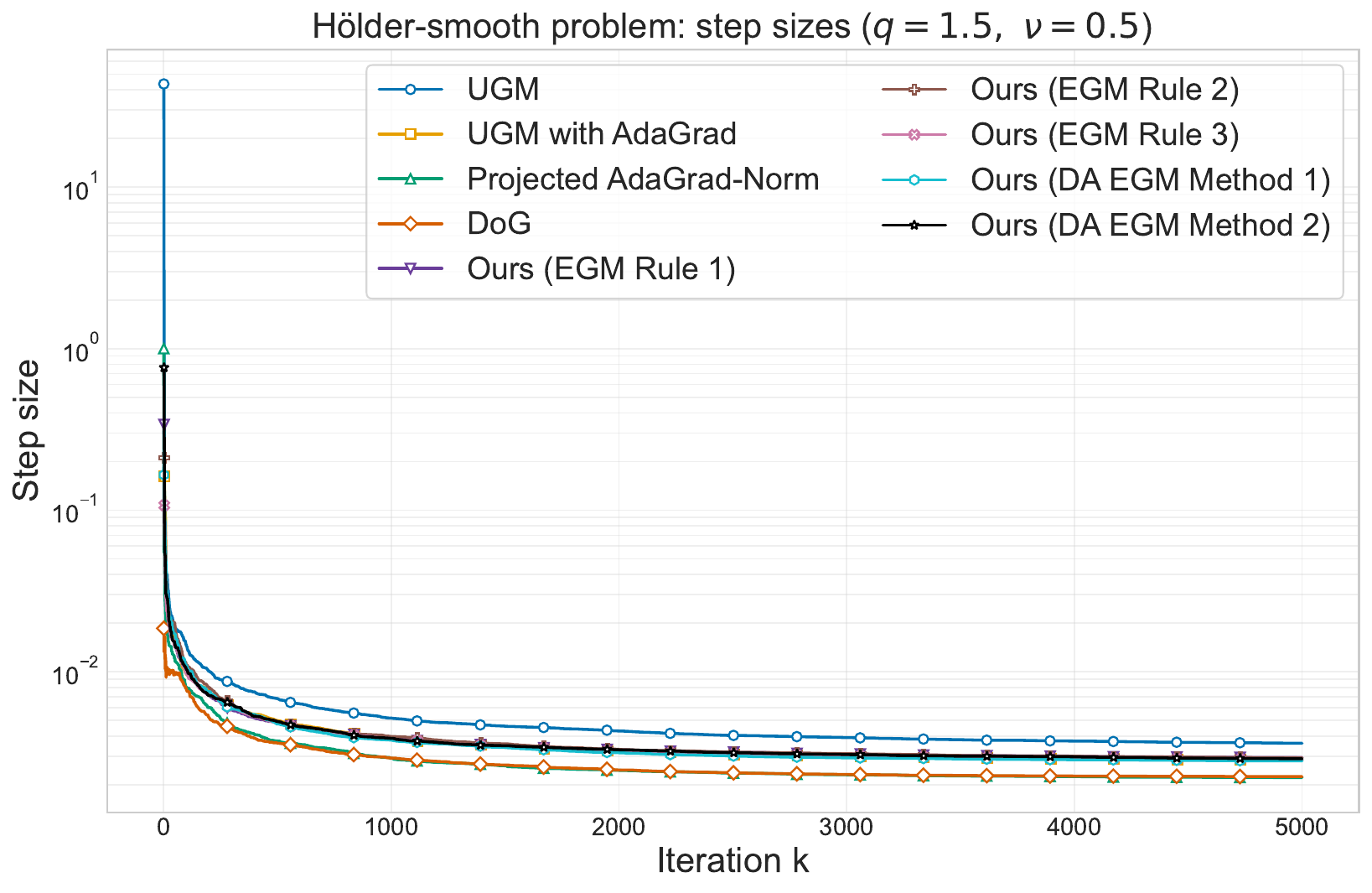}
\end{subfigure}
\hfill
\begin{subfigure}{0.32\textwidth}
    \centering
    \includegraphics[width=\linewidth]{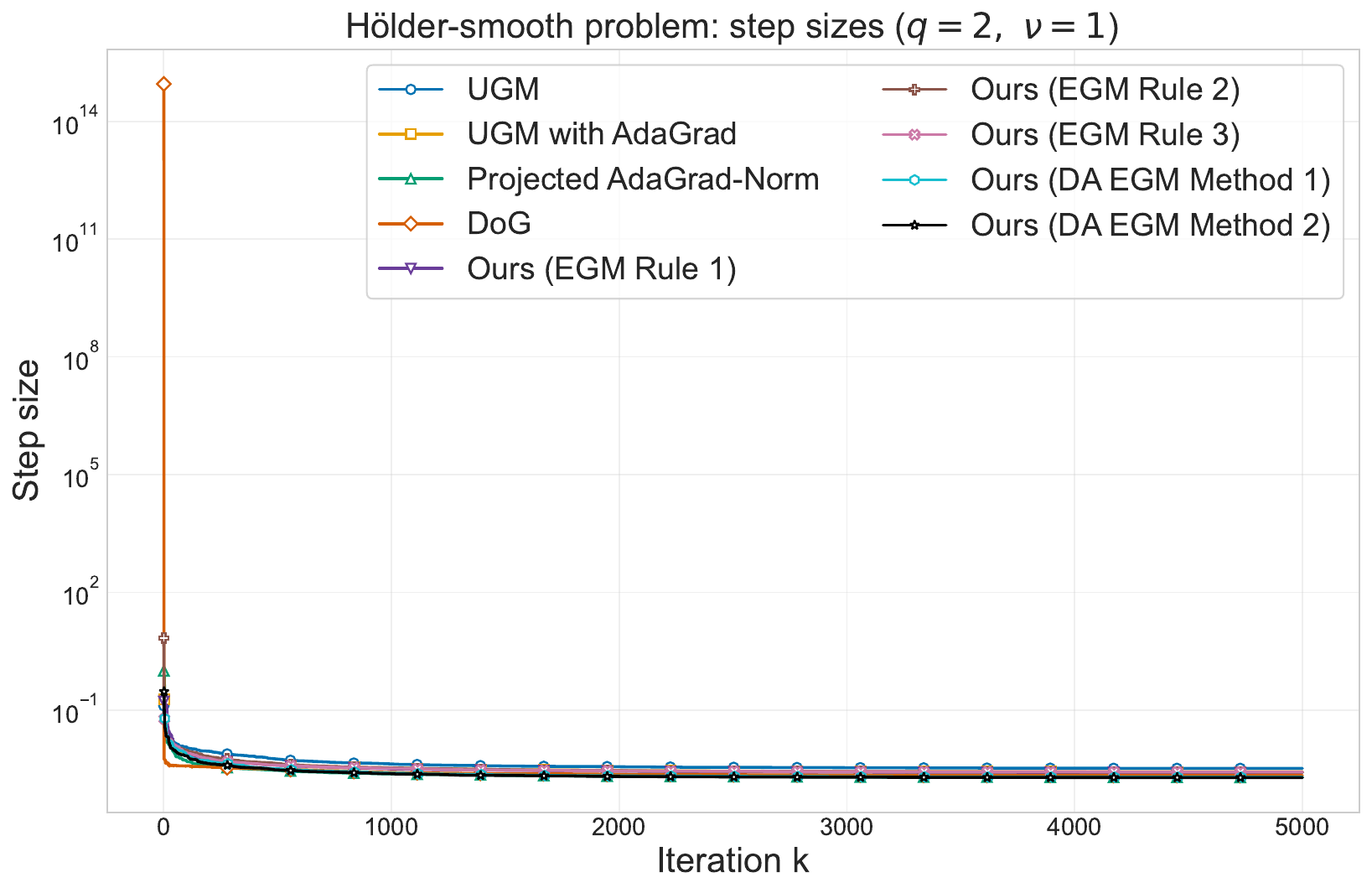}
\end{subfigure}

\caption{Results for the H\"older-smooth polyhedral feasibility problem (stochastic). From top to bottom,
the rows show the best-iterate objective function gaps, averaged-iterate objective function gaps, and step sizes. The columns correspond to $q=1$, $q=1.5$, and $q=2$, respectively.}
\label{fig_holder}
\end{figure*}



\paragraph{H\"older-smooth polyhedral feasibility problem.}

In our experiments, we use $n=1000$ inequality constraints, an optimization-variable dimension of $d=1000$, and a Euclidean-ball radius of $R=1$. We consider three values of $q$, namely, $q\in\{1,1.5,2\}$. The same realization of the generated data is used for all three values of $q$, so only the smoothness of the objective changes across the experiments.

Figure~\ref{fig_holder} shows the objective function gaps evaluated at the best and averaged iterates, along with the corresponding step sizes for all methods. Our algorithms perform comparably to the benchmark methods and, in most cases, achieve slightly smaller objective function gaps. The only exception is the averaged-iterate result for $q=2$, for which UGM with AdaGrad and DoG perform slightly better.

\begin{figure*}[htbp]
\centering

\begin{subfigure}{0.32\textwidth}
    \centering
    \includegraphics[width=\linewidth]
    {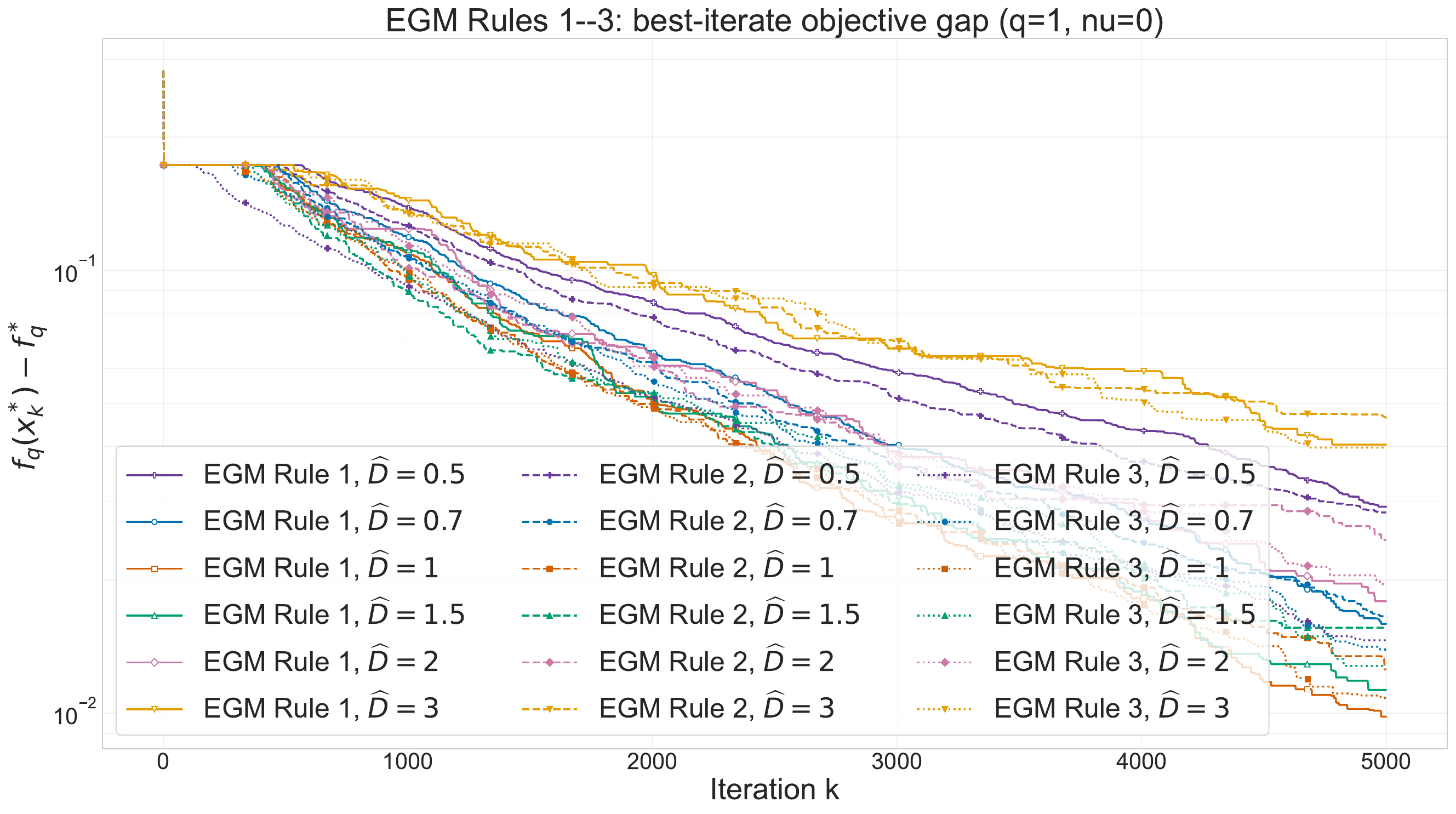}
\end{subfigure}
\hfill
\begin{subfigure}{0.32\textwidth}
    \centering
    \includegraphics[width=\linewidth]
    {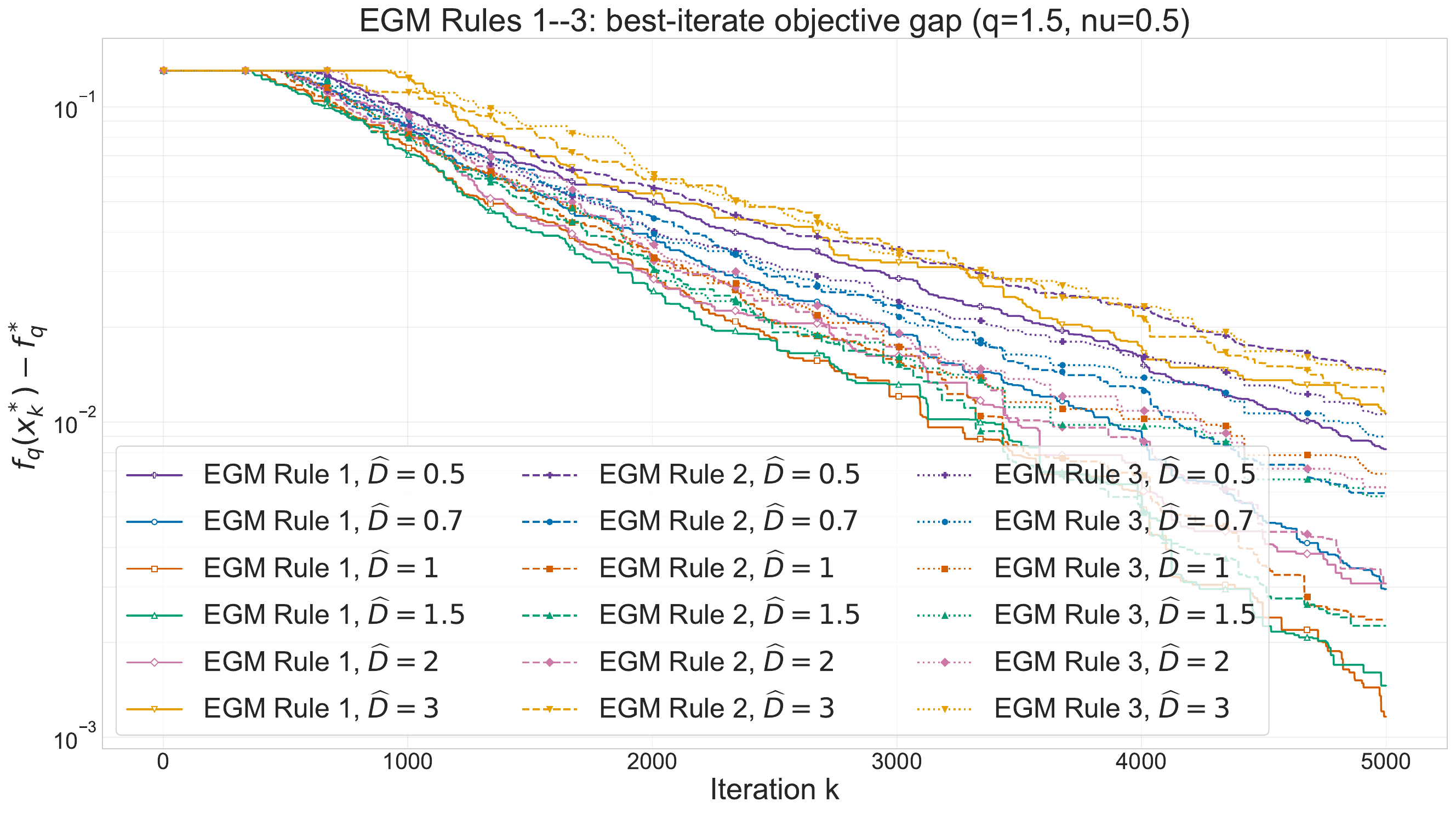}
\end{subfigure}
\hfill
\begin{subfigure}{0.32\textwidth}
    \centering
    \includegraphics[width=\linewidth]
    {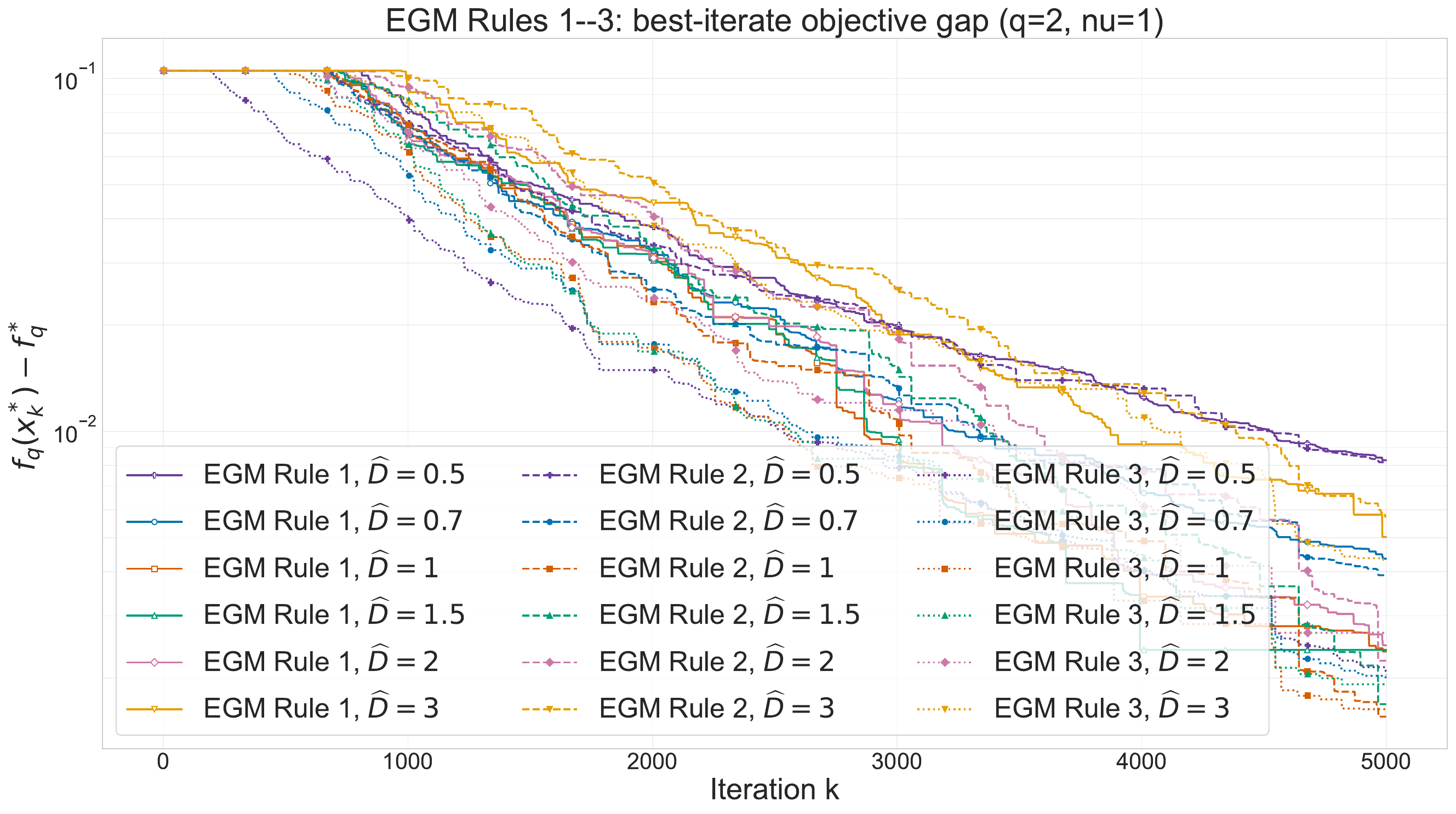}
\end{subfigure}
\par

\begin{subfigure}{0.32\textwidth}
    \centering
    \includegraphics[width=\linewidth]
    {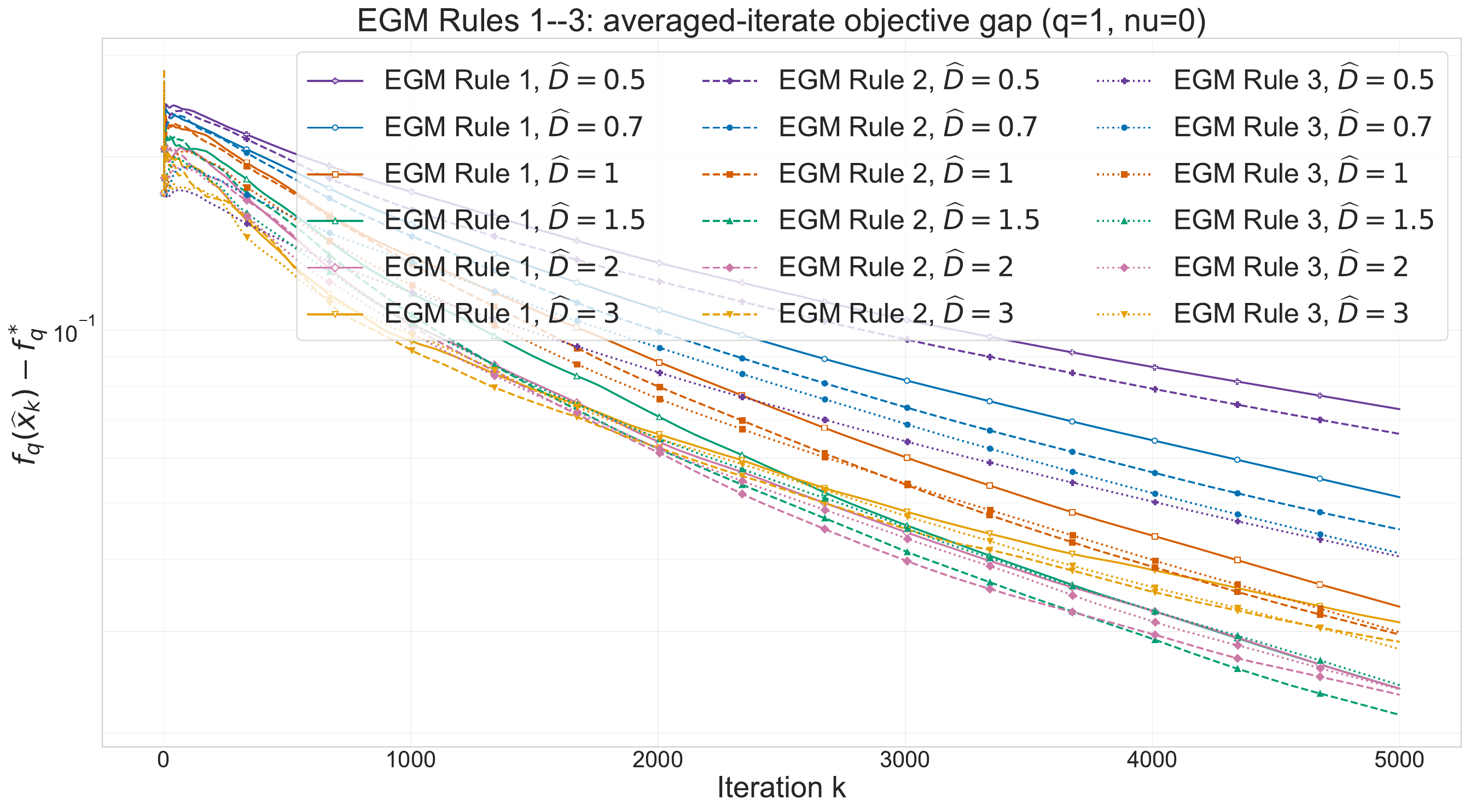}
\end{subfigure}
\hfill
\begin{subfigure}{0.32\textwidth}
    \centering
    \includegraphics[width=\linewidth]
    {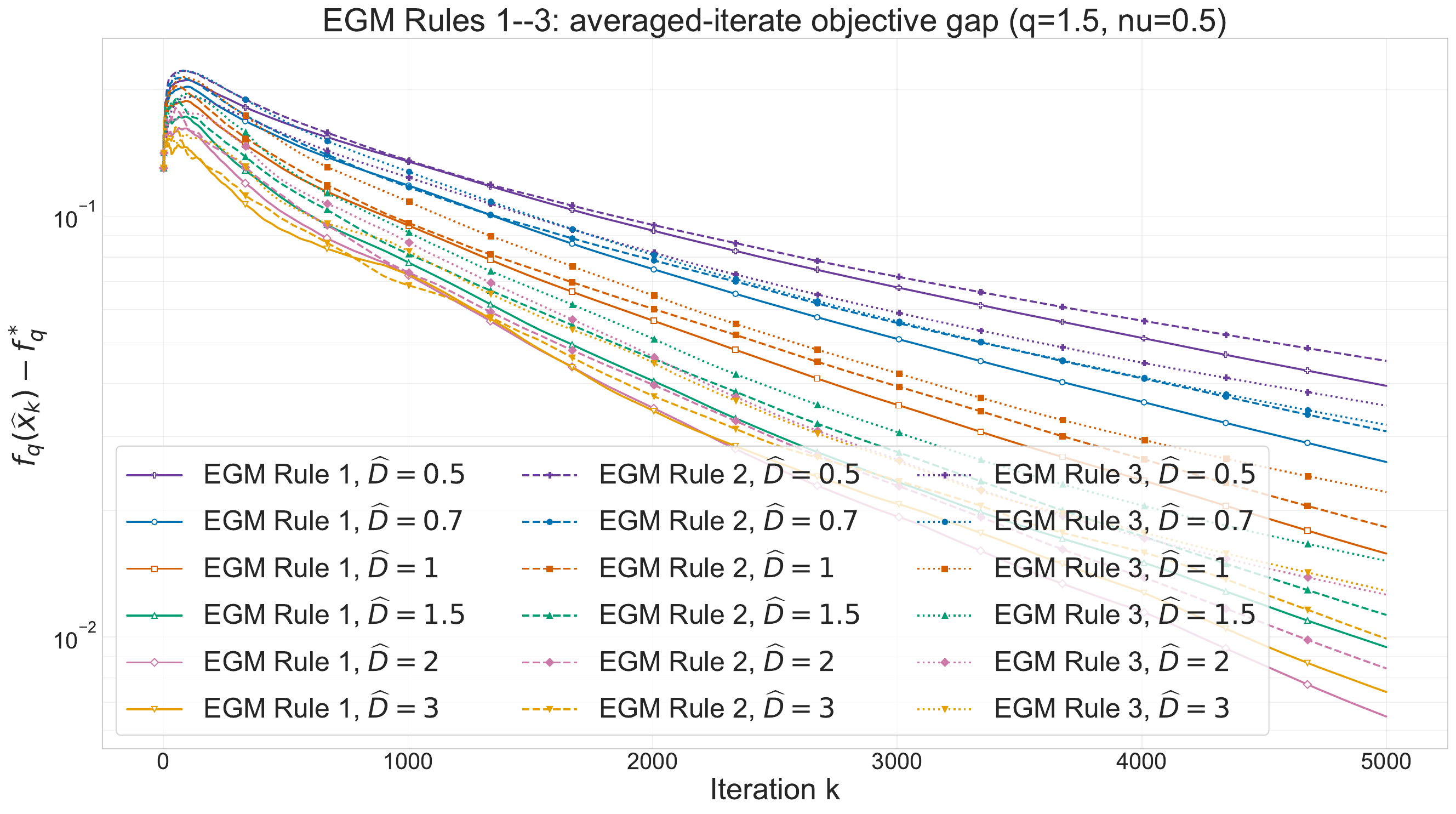}
\end{subfigure}
\hfill
\begin{subfigure}{0.32\textwidth}
    \centering
    \includegraphics[width=\linewidth]
    {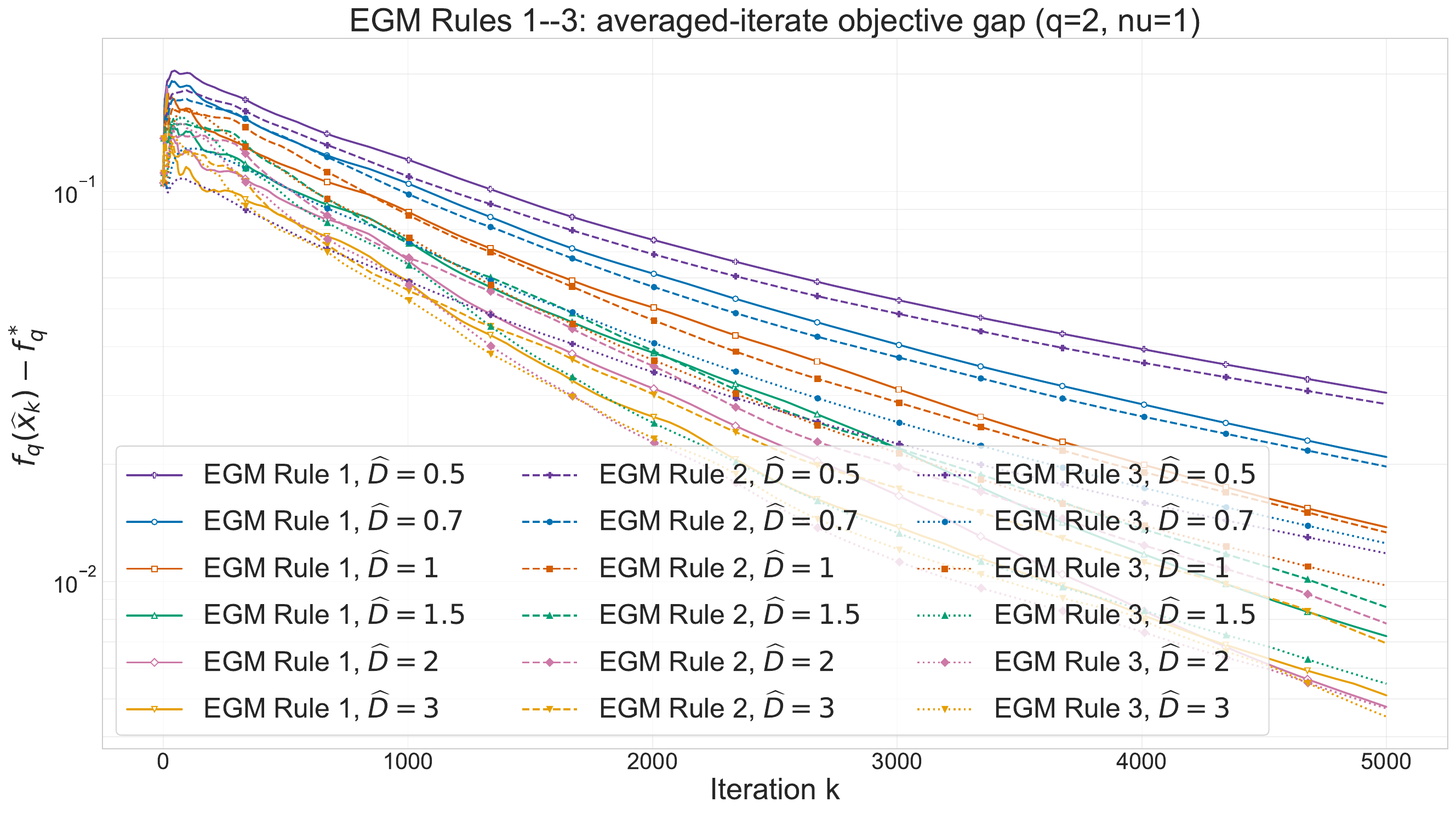}
\end{subfigure}
\par

\begin{subfigure}{0.32\textwidth}
    \centering
    \includegraphics[width=\linewidth]
    {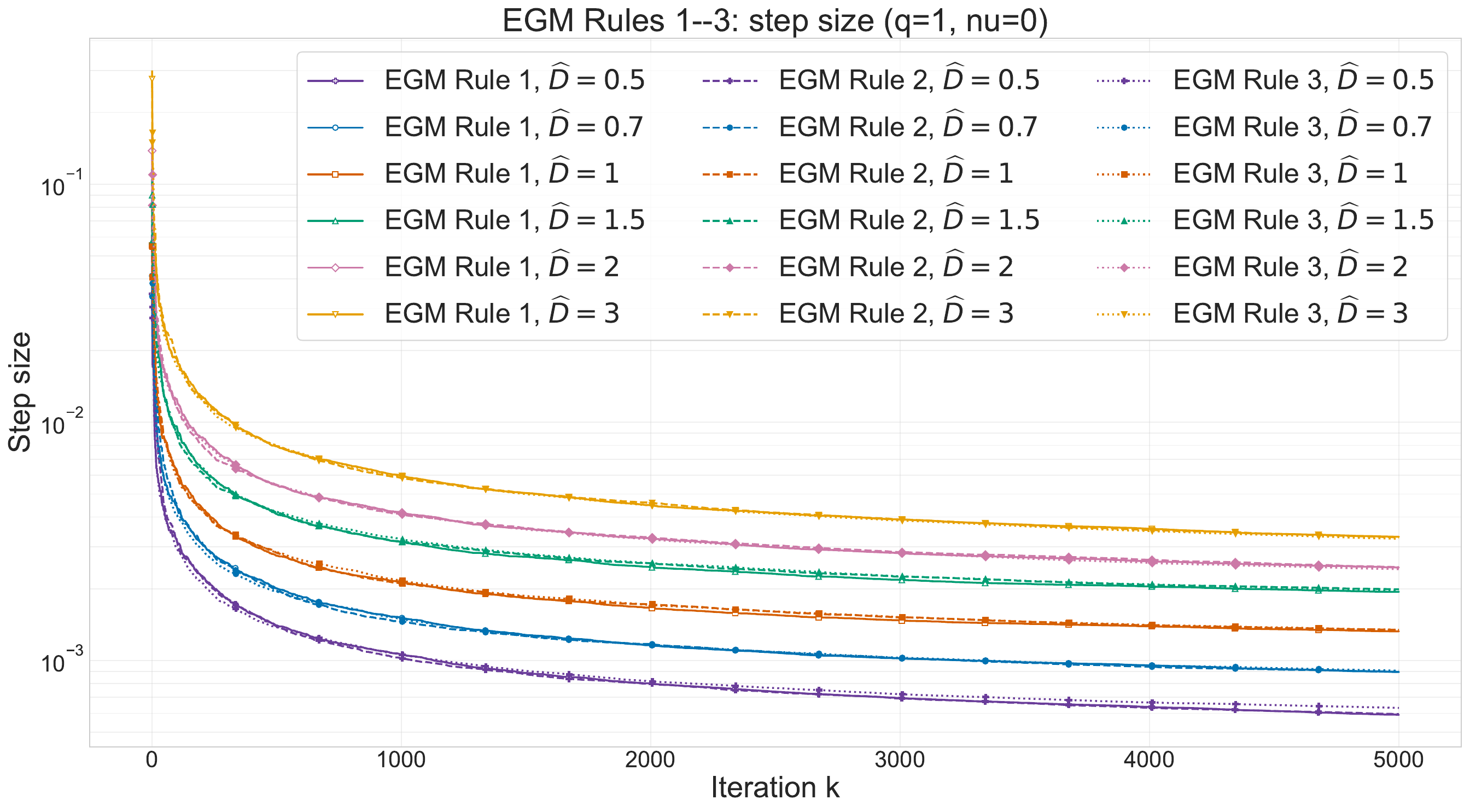}
\end{subfigure}
\hfill
\begin{subfigure}{0.32\textwidth}
    \centering
    \includegraphics[width=\linewidth]
    {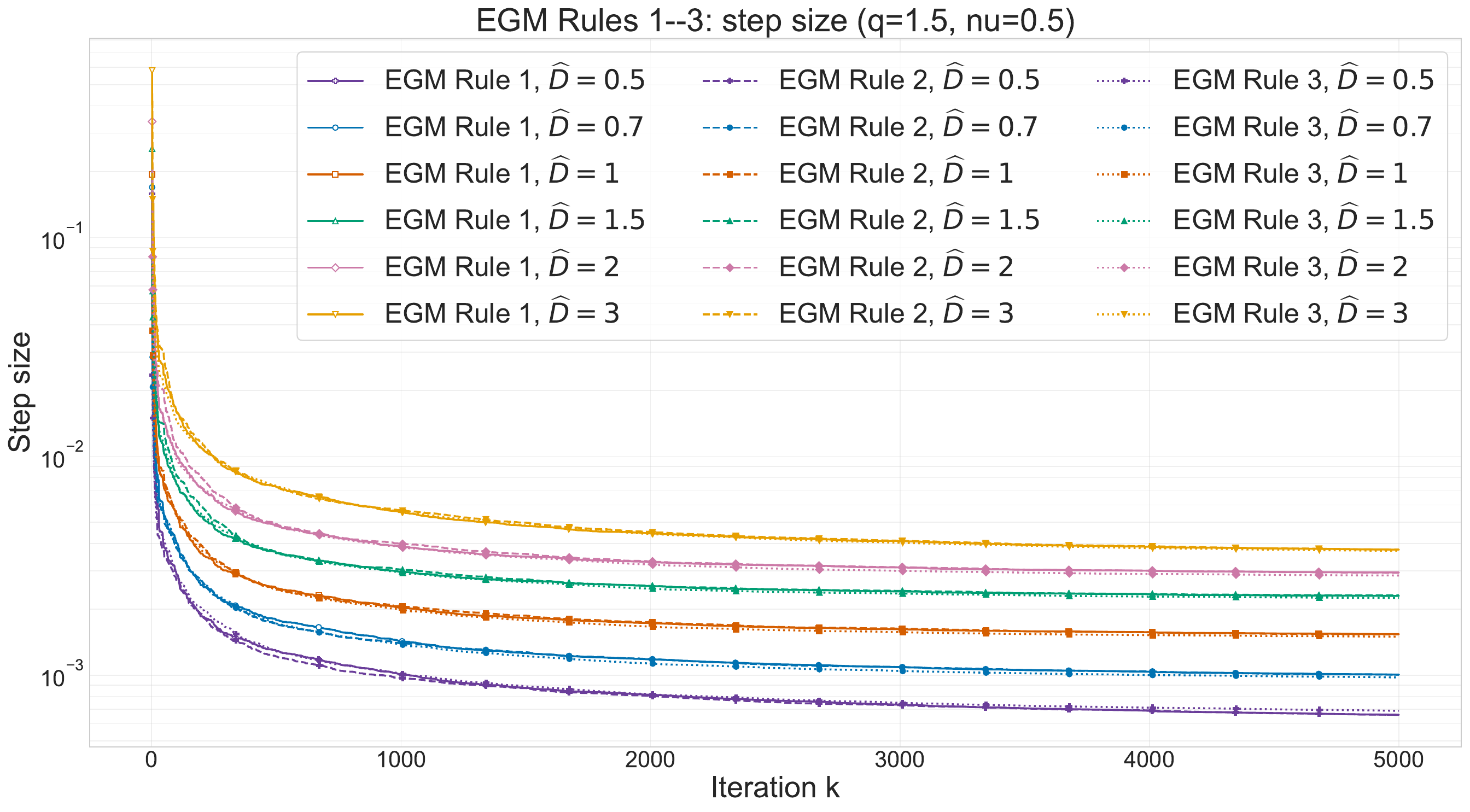}
\end{subfigure}
\hfill
\begin{subfigure}{0.32\textwidth}
    \centering
    \includegraphics[width=\linewidth]
    {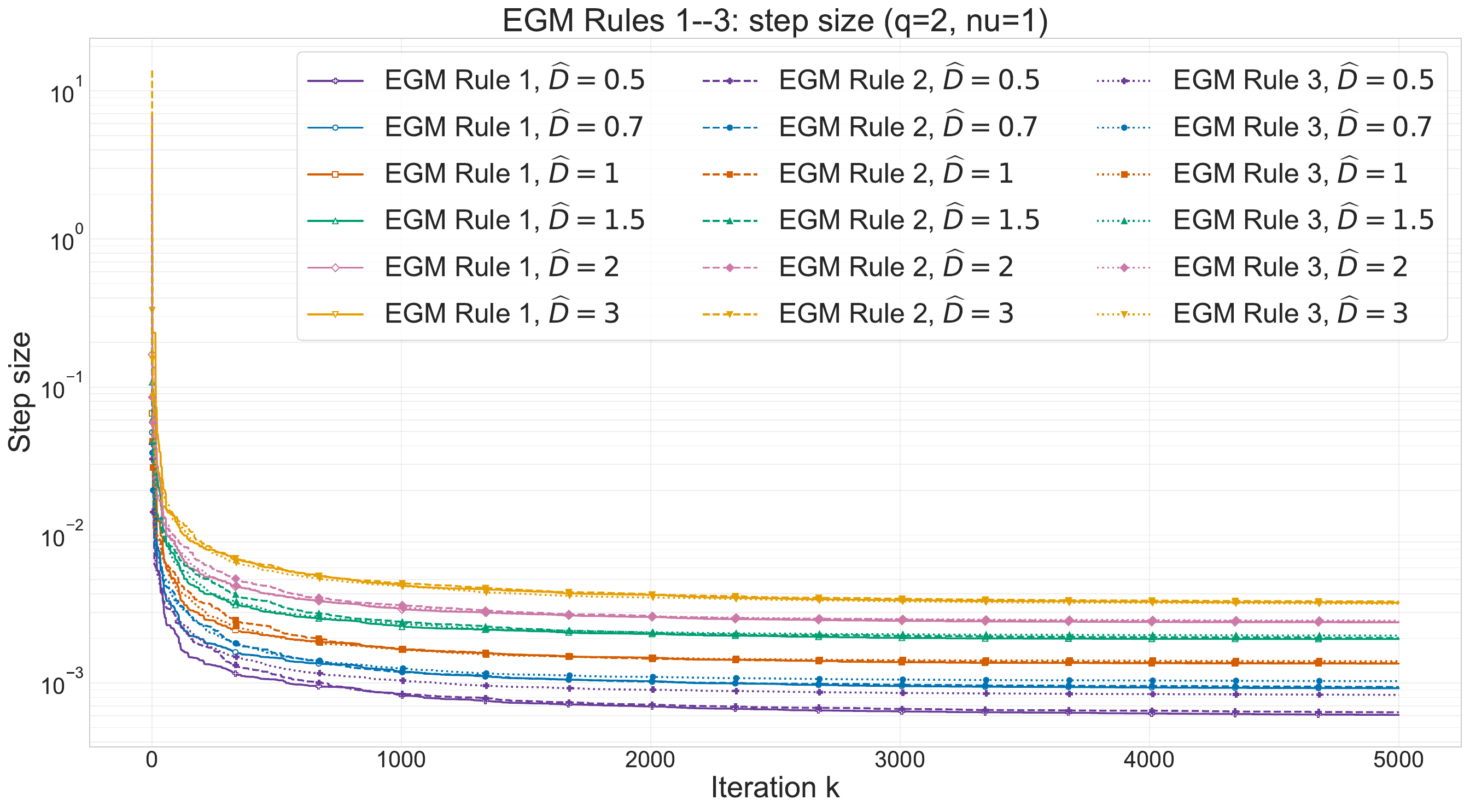}
\end{subfigure}

\caption{Sensitivity of EGM Rules 1--3 to $\widehat D$ for the
stochastic H\"older problem. From left to right, the columns correspond
to $q=1$, $q=1.5$, and $q=2$. From top to bottom, the rows show the
best-iterate function-value gaps, averaged-iterate function-value gaps,
and step sizes.}
\label{fig:holder_Dhat_sensitivity}
\end{figure*}

\begin{figure*}[htbp]
\centering

\begin{subfigure}{0.32\textwidth}
    \centering
    \includegraphics[width=\linewidth]
    {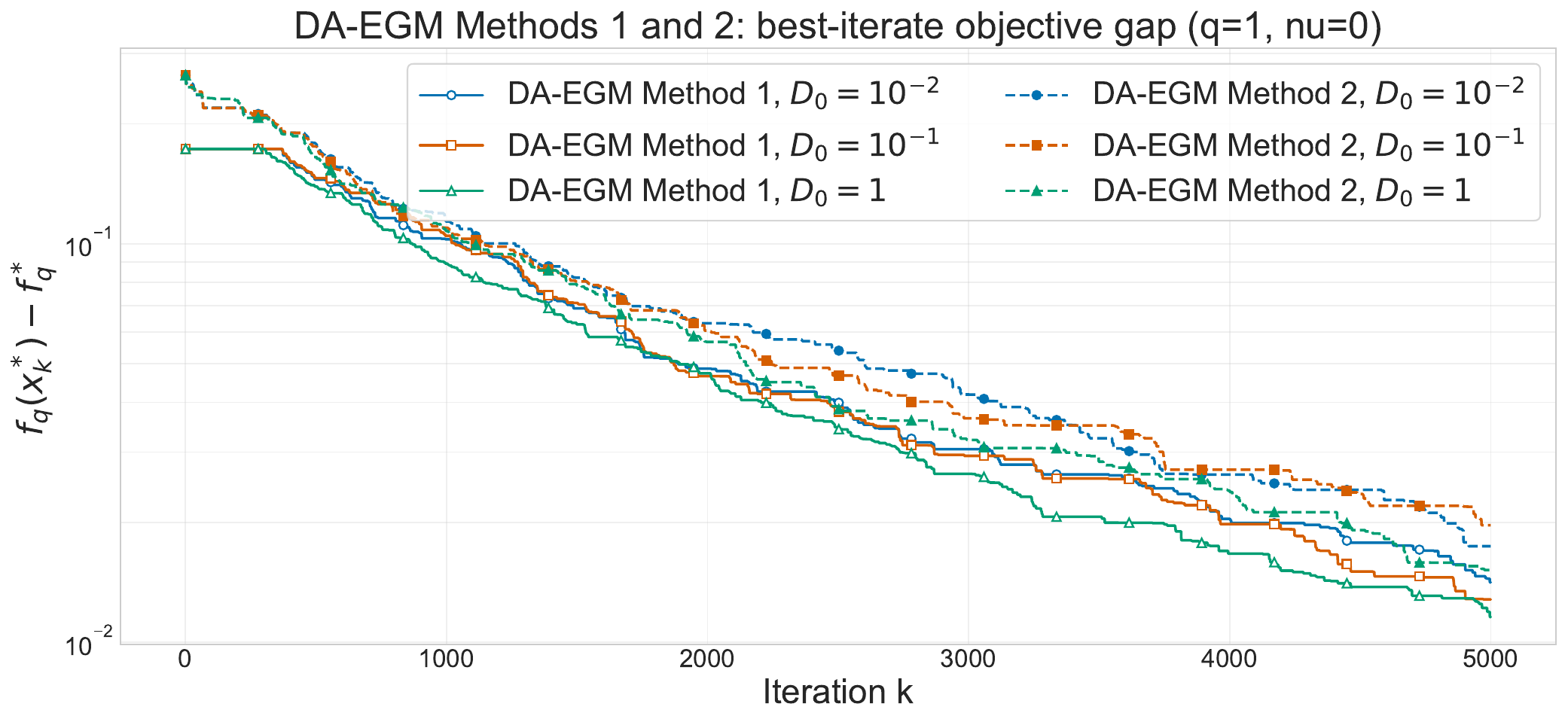}
\end{subfigure}
\hfill
\begin{subfigure}{0.32\textwidth}
    \centering
    \includegraphics[width=\linewidth]
    {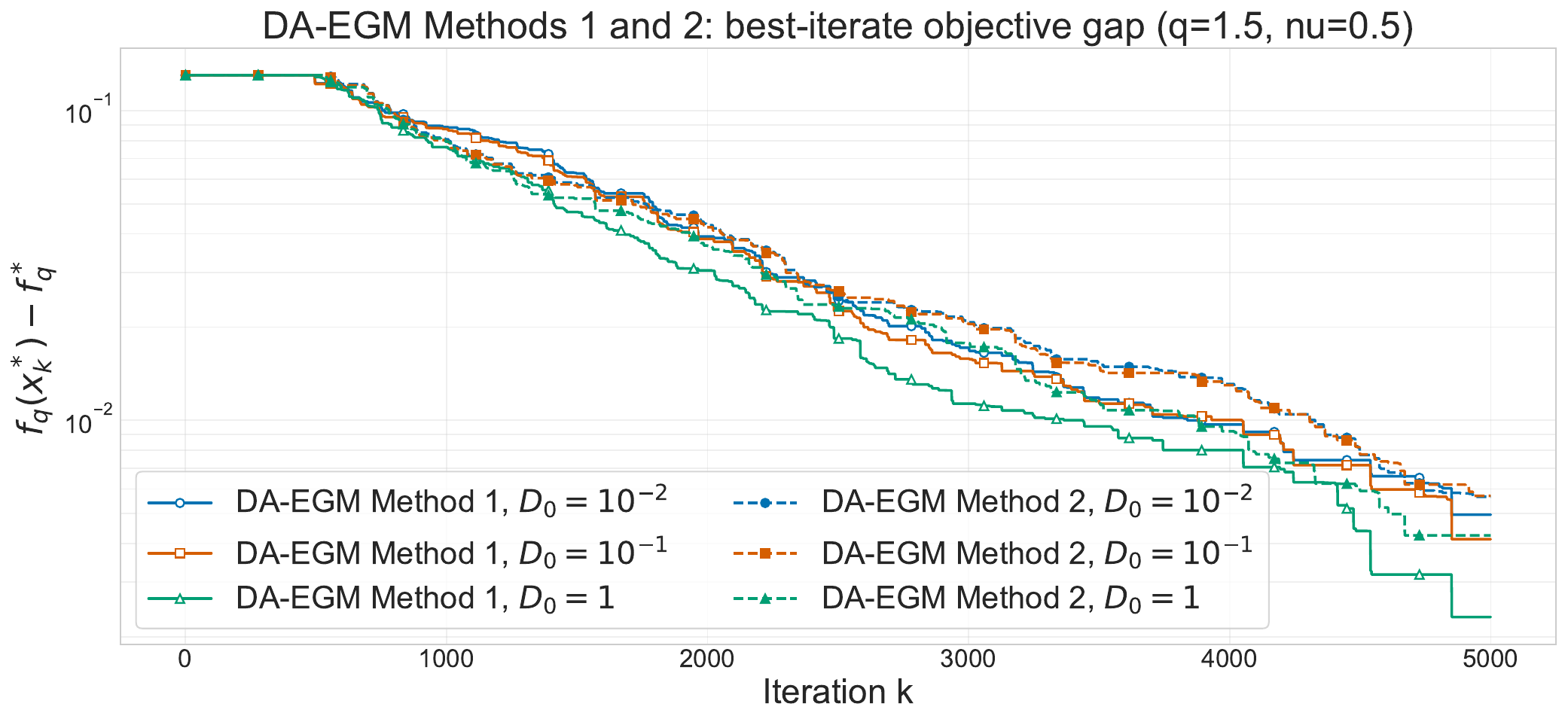}
\end{subfigure}
\hfill
\begin{subfigure}{0.32\textwidth}
    \centering
    \includegraphics[width=\linewidth]
    {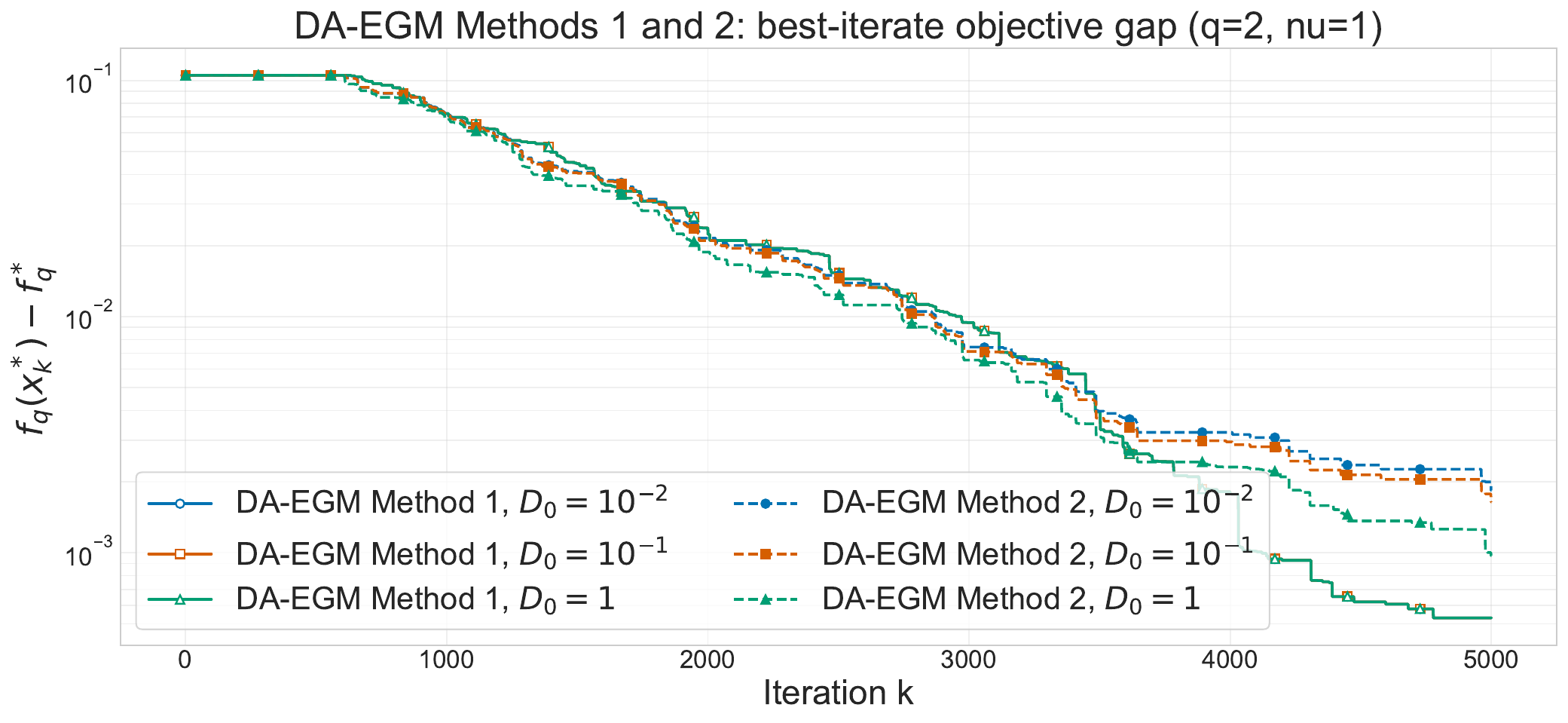}
\end{subfigure}
\par

\begin{subfigure}{0.32\textwidth}
    \centering
    \includegraphics[width=\linewidth]
    {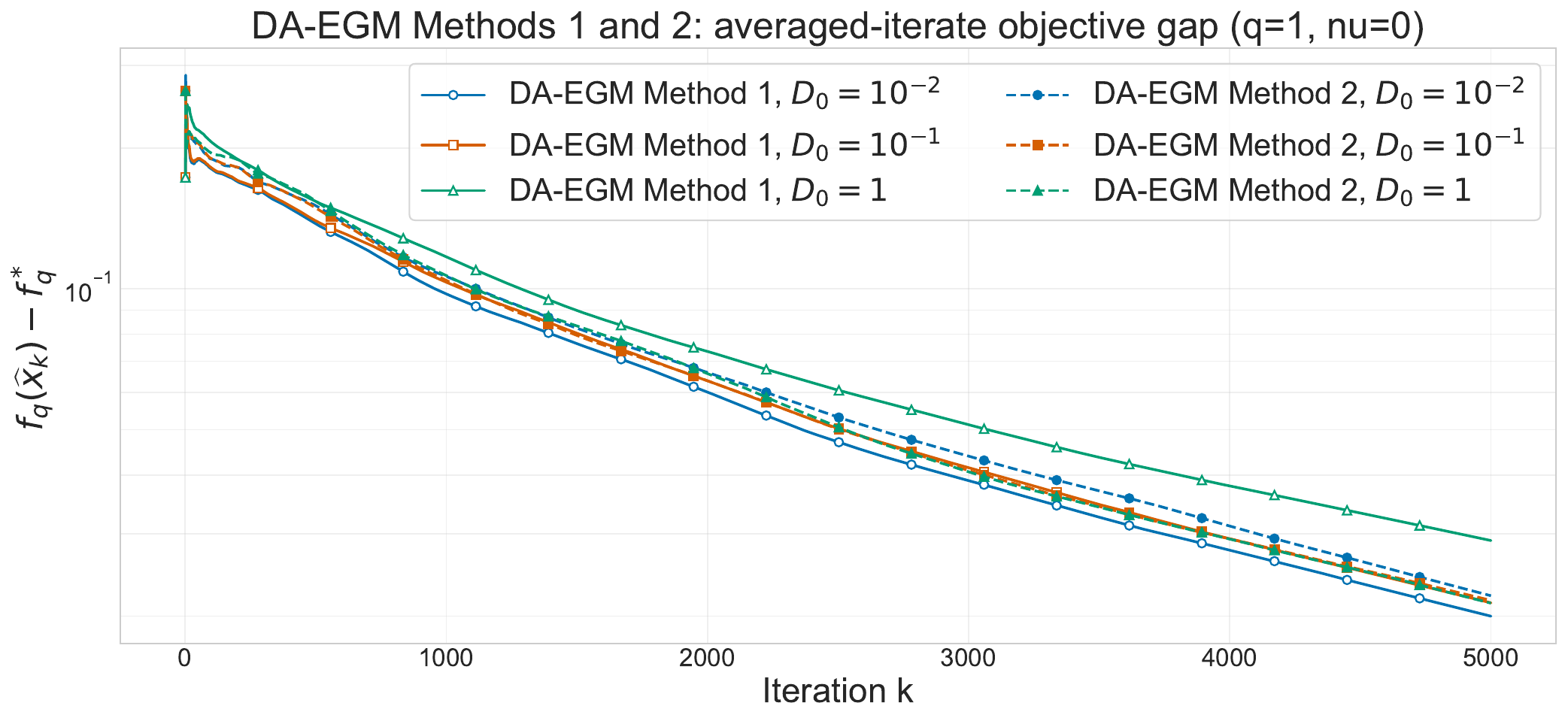}
\end{subfigure}
\hfill
\begin{subfigure}{0.32\textwidth}
    \centering
    \includegraphics[width=\linewidth]
    {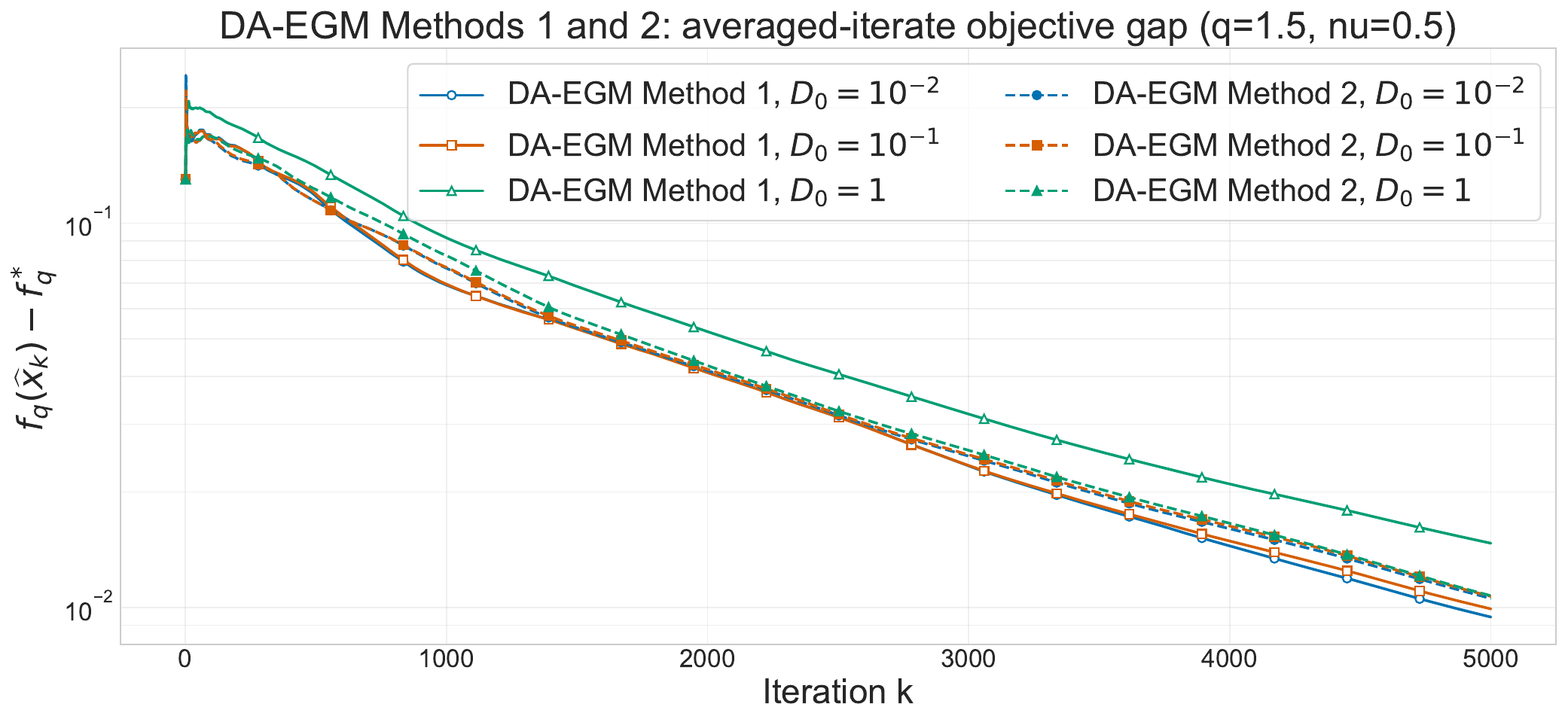}
\end{subfigure}
\hfill
\begin{subfigure}{0.32\textwidth}
    \centering
    \includegraphics[width=\linewidth]
    {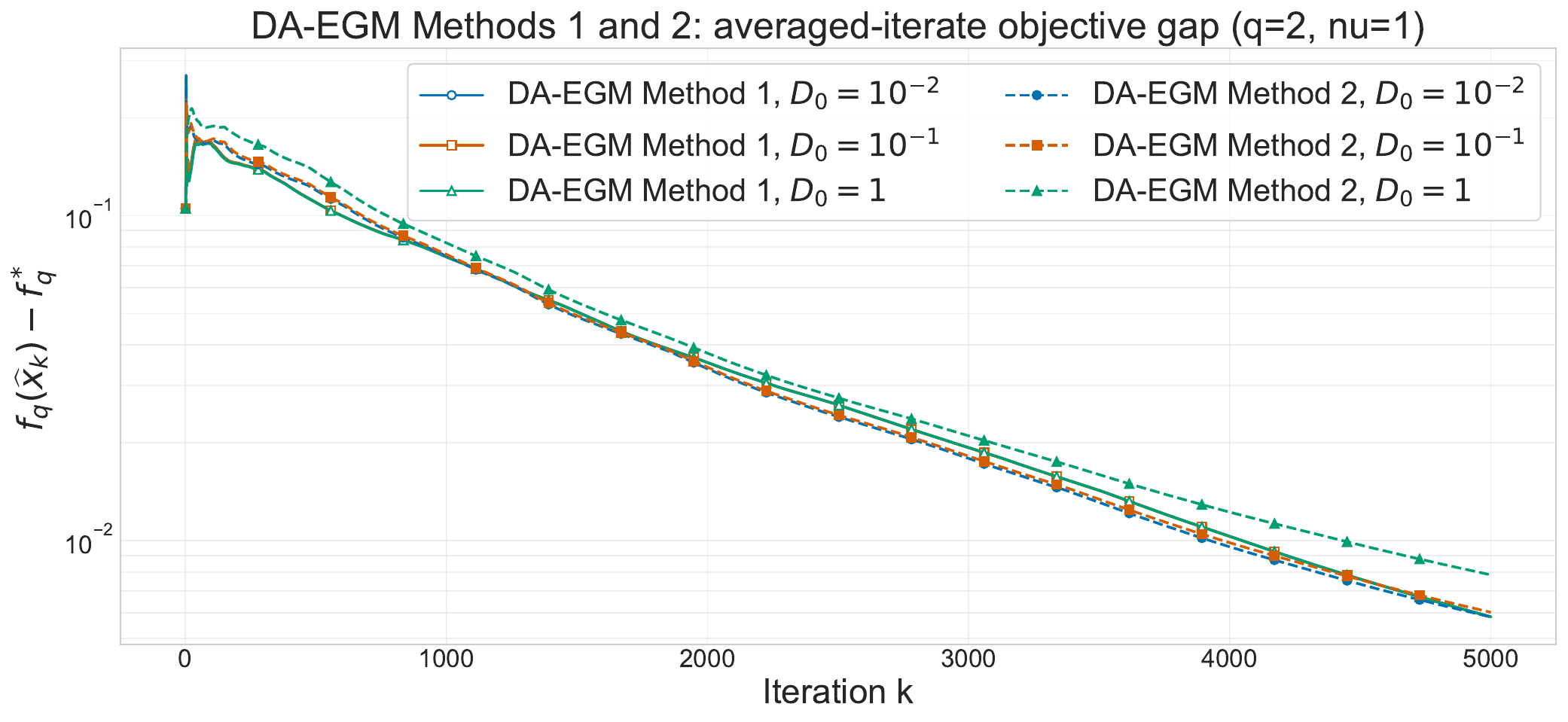}
\end{subfigure}
\par

\begin{subfigure}{0.32\textwidth}
    \centering
    \includegraphics[width=\linewidth]
    {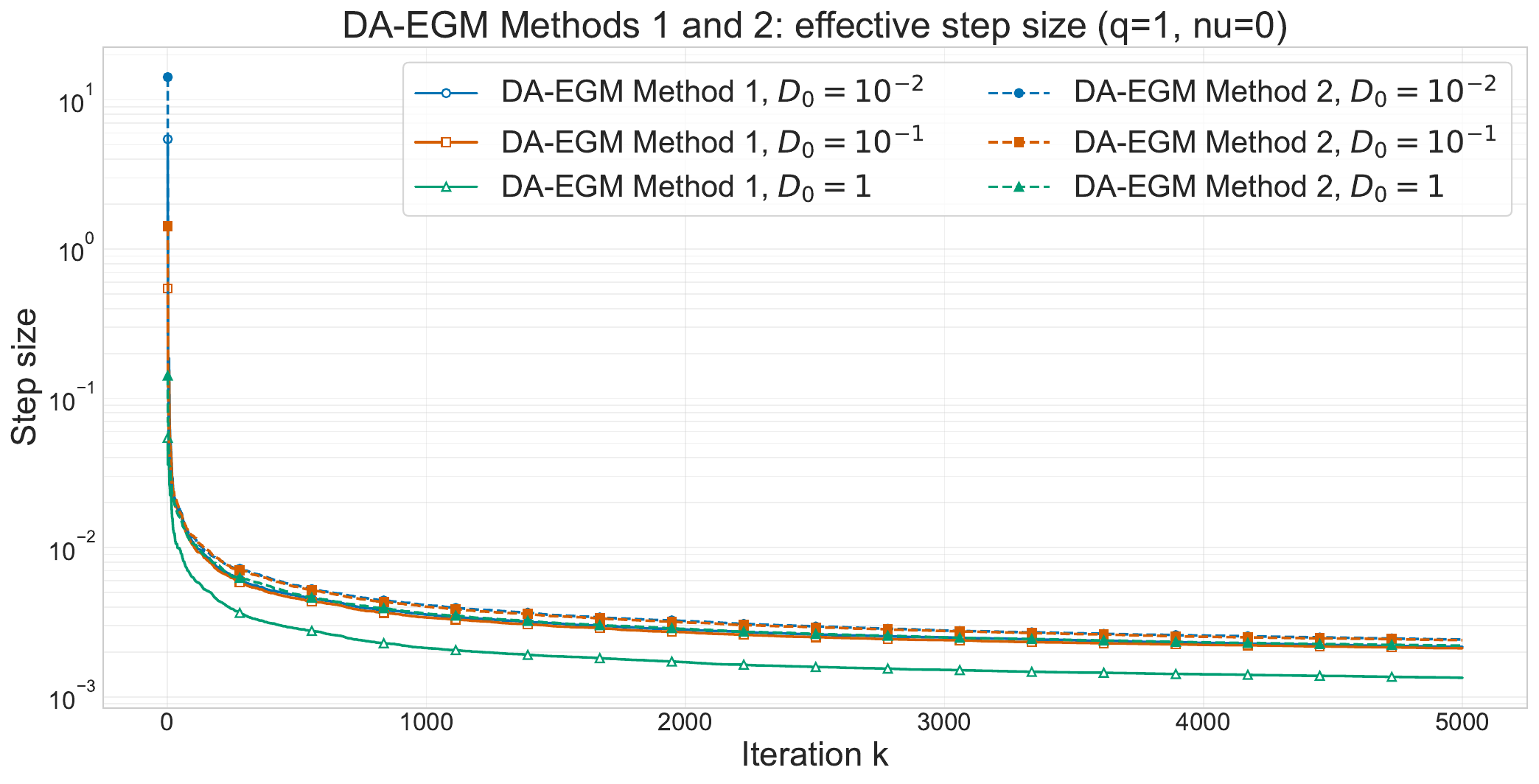}
\end{subfigure}
\hfill
\begin{subfigure}{0.32\textwidth}
    \centering
    \includegraphics[width=\linewidth]
    {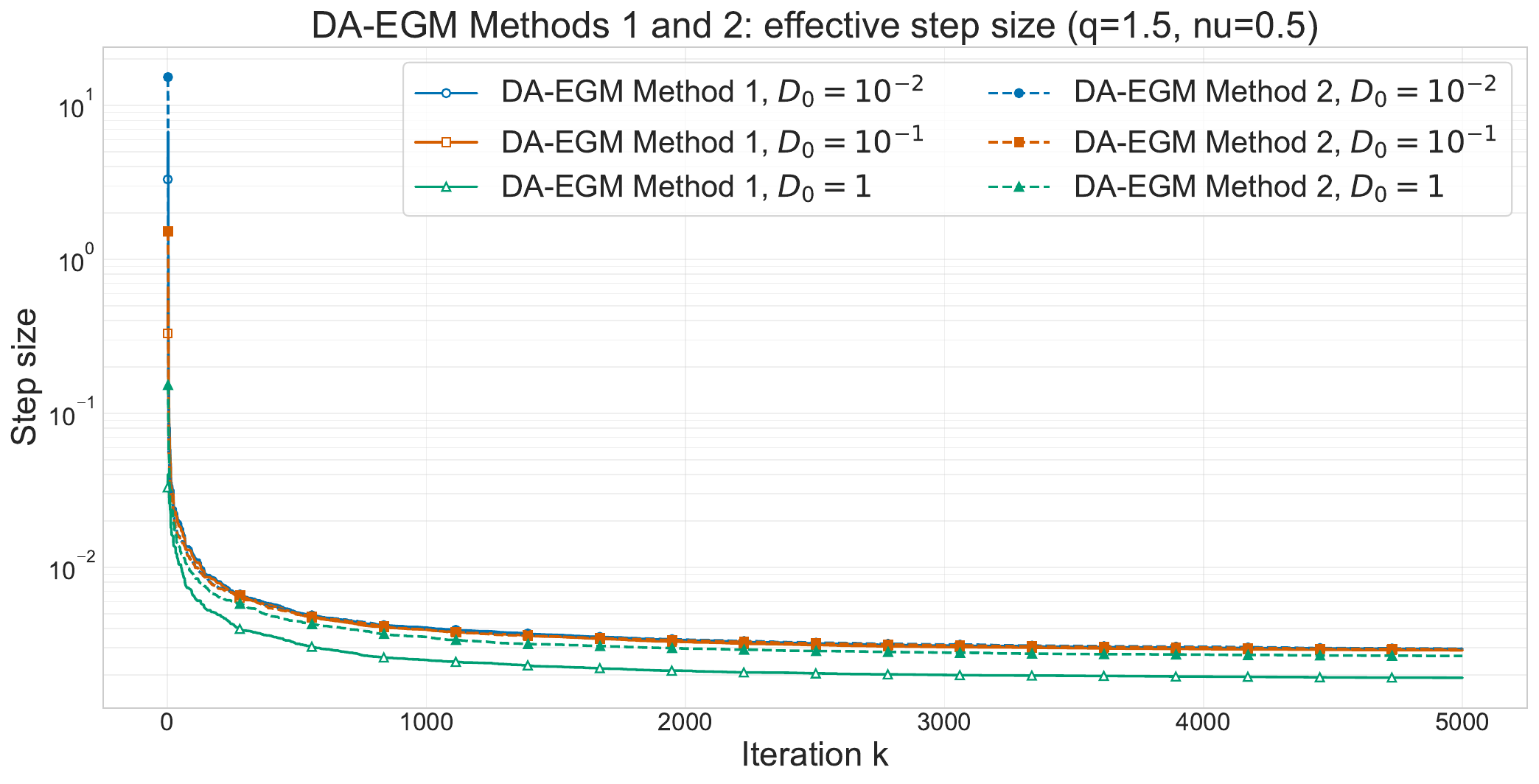}
\end{subfigure}
\hfill
\begin{subfigure}{0.32\textwidth}
    \centering
    \includegraphics[width=\linewidth]
    {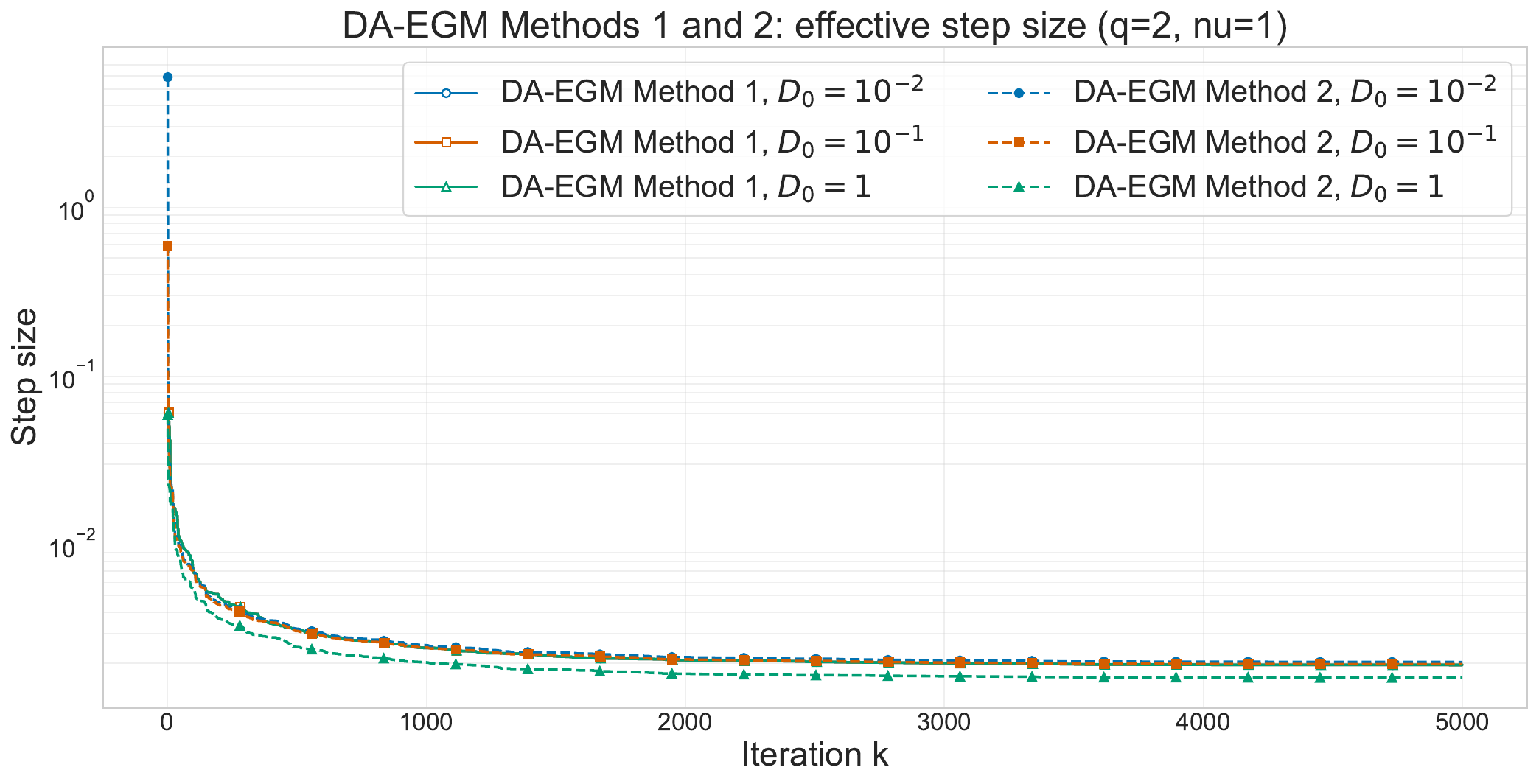}
\end{subfigure}

\caption{Sensitivity of DA-EGM Methods 1 and 2 to $D_0$ for the
stochastic H\"older problem. From left to right, the columns correspond
to $q=1$, $q=1.5$, and $q=2$. From top to bottom, the rows show the
best-iterate function-value gaps, averaged-iterate function-value gaps,
and step sizes.}
\label{fig:holder_D0_sensitivity}
\end{figure*}

Figure~\ref{fig:holder_Dhat_sensitivity} illustrates the sensitivity of EGM Rules~1--3 to the choice of $\widehat D$. For the best iterate, Rule~1 with $\widehat D=1$ performs well for the H\"older exponents $\nu=0$ and $\nu=0.5$, whereas Rules~2 and~3 with $\widehat D=1$ perform well for $\nu=1$. For the 
averaged iterate, Rule~2 with $\widehat D=1.5$ performs well for $\nu=0$, Rule~1 with $\widehat D=2$ performs well for $\nu=0.5$, and Rule~3 with $\widehat D=2$ or $3$ performs well for $\nu=1$. For the DA EGM methods, Figure~\ref{fig:holder_D0_sensitivity} shows that Method~1 with $D_0=1$ 
performs well for the best iterate across all smoothness regimes. For the averaged iterate, Method~1 with $D_0=10^{-2}$ performs slightly better across all smoothness regimes, although its performance remains comparable to that obtained with the other choices of $D_0$.


\begin{figure*}[htbp]
\centering


\begin{subfigure}{0.48\textwidth}
    \centering
    \includegraphics[width=\linewidth]
    {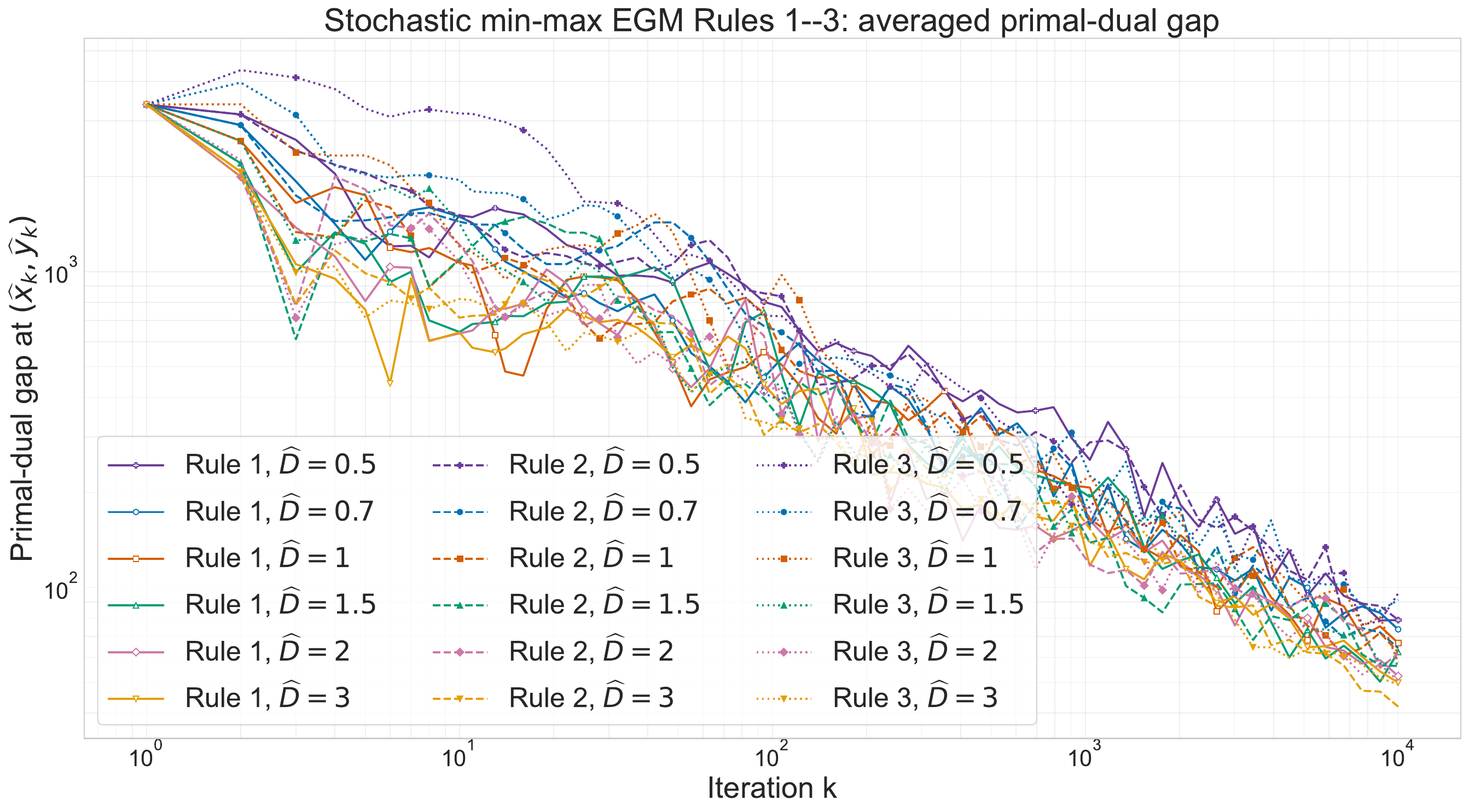}
\end{subfigure}
\hfill
\begin{subfigure}{0.48\textwidth}
    \centering
    \includegraphics[width=\linewidth]
    {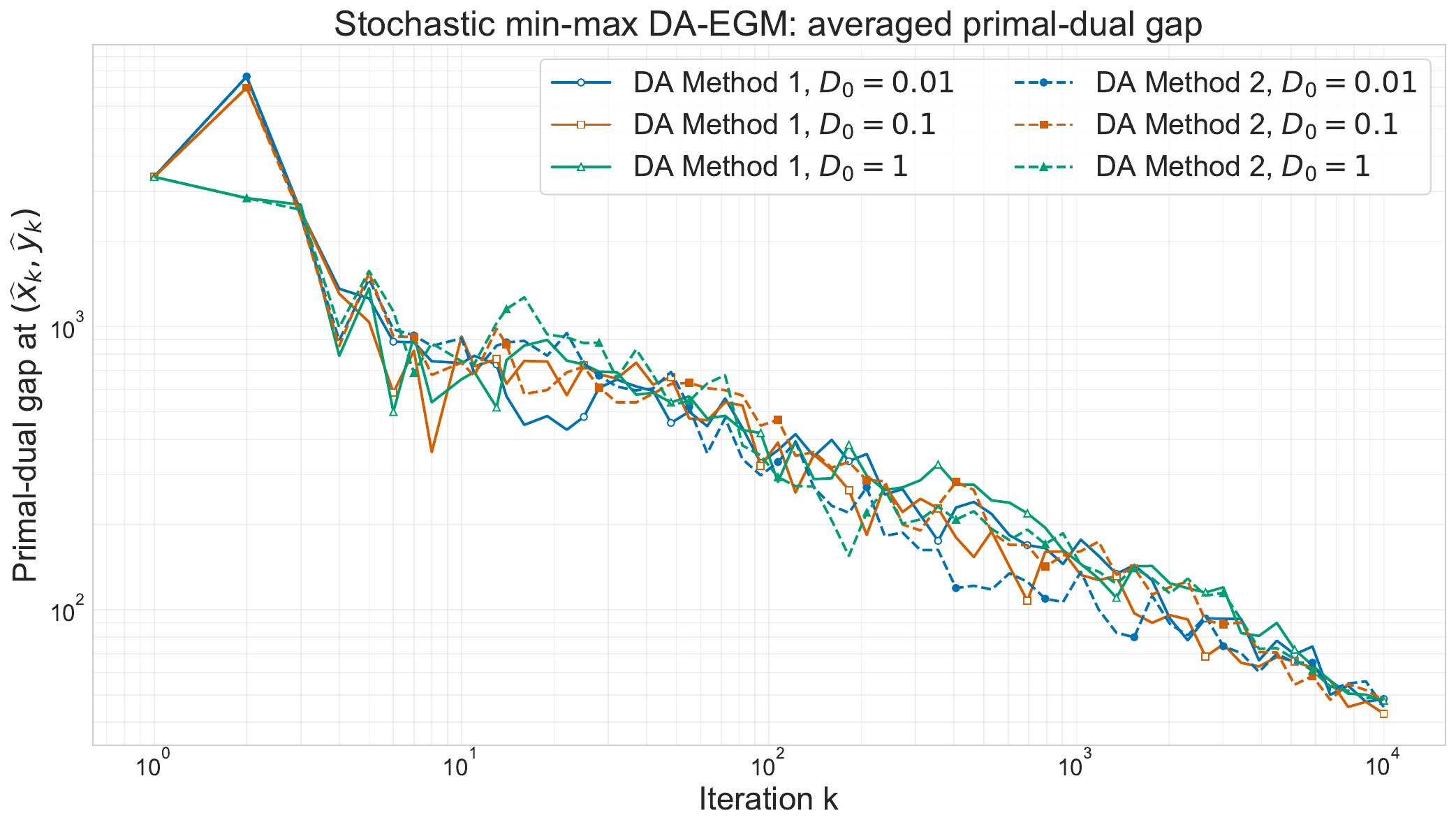}
\end{subfigure}
\par


\begin{subfigure}{0.48\textwidth}
    \centering
    \includegraphics[width=\linewidth]
    {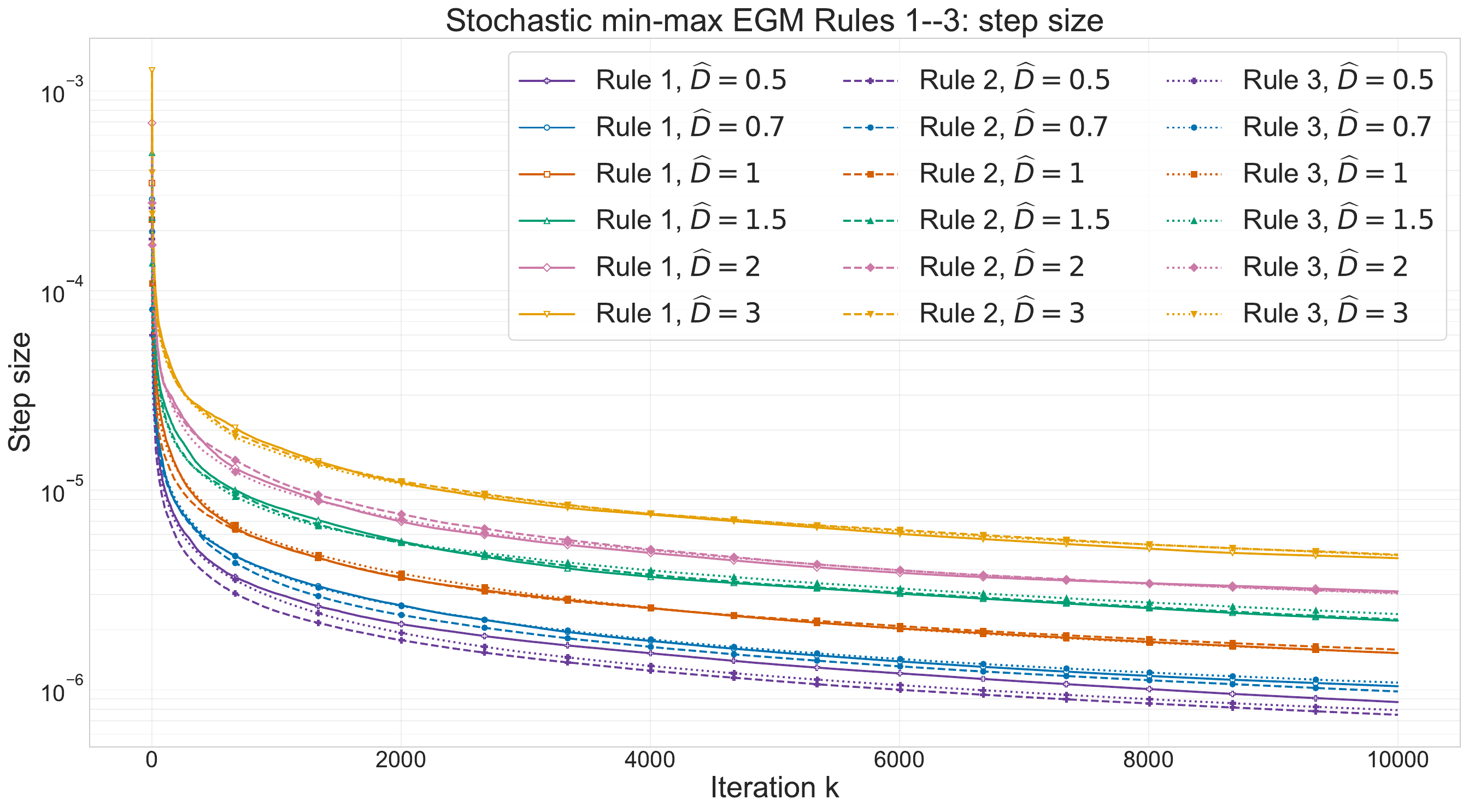}
\end{subfigure}
\hfill
\begin{subfigure}{0.48\textwidth}
    \centering
    \includegraphics[width=\linewidth]
    {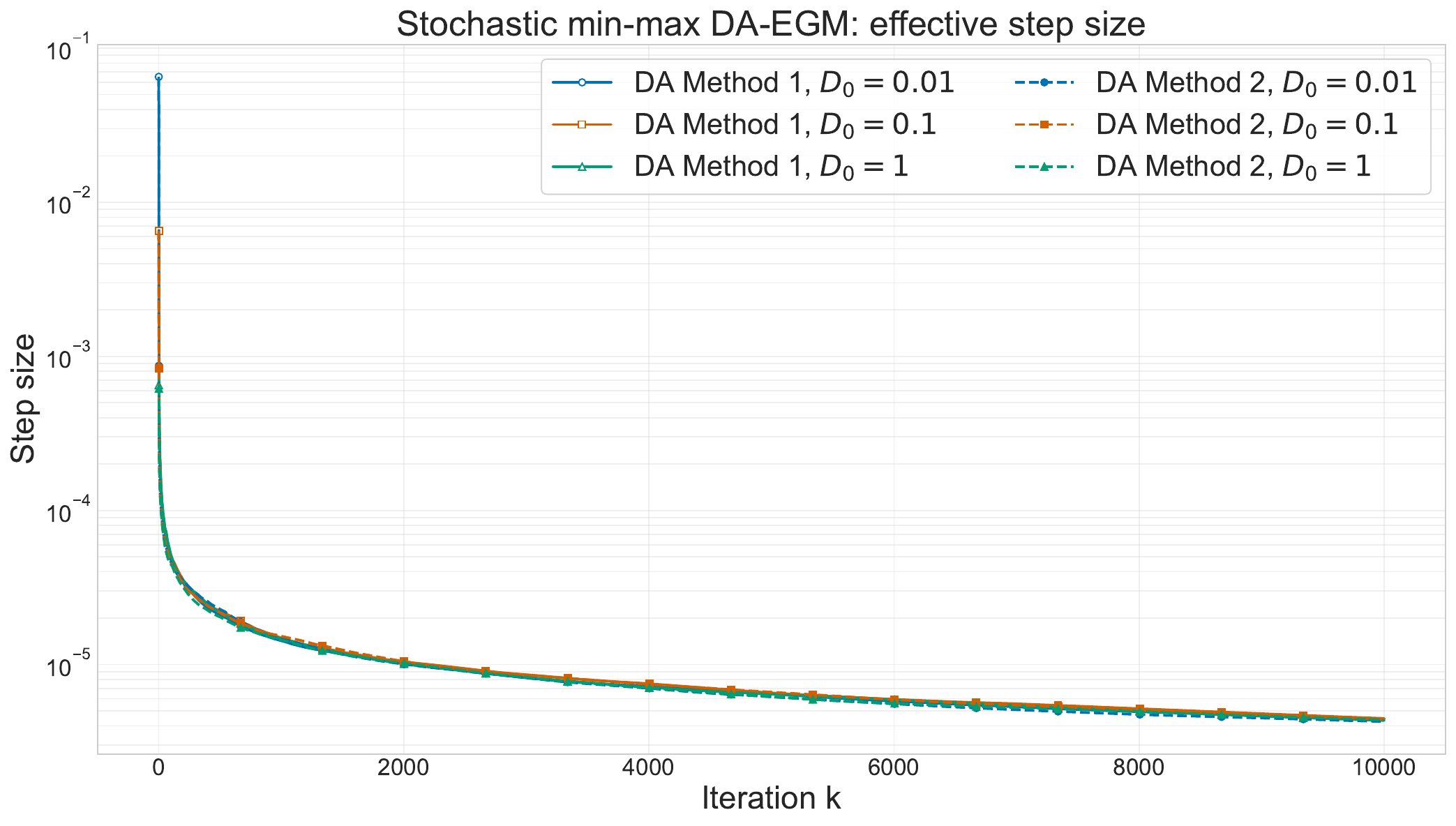}
\end{subfigure}

\caption{Sensitivity results for the stochastic least-squares
min--max problem. The left column shows EGM Rules 1--3 for varying
$\widehat D$, and the right column shows DA-EGM Methods 1 and 2 for
varying $D_0$. The first row presents the variational-inequality gaps,
and the second row presents the corresponding step sizes.}
\label{fig:minmax_sensitivity}
\end{figure*}


\paragraph{Min-max dual formulation of a constrained $\ell_2$-regression problem.}

In our experiments, we considered the matrix \(A\in\mathbb{R}^{m \times n}\) with $m = {10000}$ and $n = 8000$. The entries of the matrix $A$, denoted by $A_{ij}$, are sampled independently and identically from the uniform distribution $U[0,1)$. We generate a target vector \(x_{\mathrm{tar}}\) that satisfies
\(\|x_{\mathrm{tar}}\|_2=1.5\) and we set
\begin{align*}
    b=Ax_{\mathrm{tar}}+\xi,
    \qquad
    \|\xi\|_2=0.05\|Ax_{\mathrm{tar}}\|_2.
\end{align*}
Thus, the generating target lies outside the feasible ball, while the
additive perturbation prevents a trivial exact fit. 

Next, we discuss about the construction of the stochastic VI mapping. Following \citet[Section~3.6.6]{nesterov2013firstorder}, we construct
an unbiased stochastic estimate of the VI mapping using importance
sampling. Let $\omega=(I,J)$ denote the sampled indices. For
$x\neq 0$ and $y\neq 0$, we sample $I\in\{1,\ldots,m\}$ and
$J\in\{1,\ldots,n\}$ independently according to
\begin{align*}
    \mathbb{P}(I=i\mid y)
    &=
    \frac{|y_i|}{\|y\|_1},
    \qquad i\in\{1,\ldots,m\},\\
    \mathbb{P}(J=j\mid x)
    &=
    \frac{|x_j|}{\|x\|_1},
    \qquad j\in\{1,\ldots,n\}.
\end{align*}
We then define the stochastic matrix-vector product estimators
\begin{align*}
    \widehat{A^\top y}
    &=
    \|y\|_1\operatorname{sign}(y_I)A_{I,:}^{\top},
    \\
    \widehat{Ax}
    &=
    \|x\|_1\operatorname{sign}(x_J)A_{:,J}.
\end{align*}
When $x=0$ or $y=0$, the corresponding estimator is defined to be
zero. These estimators are conditionally unbiased because
\begin{align*}
    \mathbb{E}\!\left[\widehat{A^\top y}\mid y\right]
    &=
    \sum_{i=1}^m
    \frac{|y_i|}{\|y\|_1}
    \|y\|_1\operatorname{sign}(y_i)A_{i,:}^{\top}
    =
    \sum_{i=1}^m y_iA_{i,:}^{\top}
    =
    A^\top y,
    \\
    \mathbb{E}\!\left[\widehat{Ax}\mid x\right]
    &=
    \sum_{j=1}^n
    \frac{|x_j|}{\|x\|_1}
    \|x\|_1\operatorname{sign}(x_j)A_{:,j}
    =
    \sum_{j=1}^n x_jA_{:,j}
    =
    Ax.
\end{align*}
Consequently, we define the stochastic VI mapping by
\begin{align*}
    \widehat F([x,y],\omega)
    :=
    \begin{bmatrix}
        \widehat{A^\top y}\\
        b-\widehat{Ax}
    \end{bmatrix}.
\end{align*}
It follows that
\begin{align*}
    \mathbb{E}\!\left[
        \widehat F([x,y],\omega)\mid x,y
    \right]
    =
    \begin{bmatrix}
        A^\top y\\
        b-Ax
    \end{bmatrix}
    =
    F([x,y]).
\end{align*}
Hence, the stochastic mapping oracle is unbiased. Each oracle call
accesses only one row and one column of $A$ and has an arithmetic cost
of $\mathcal{O}(m+n)$. By contrast, evaluating the exact mapping using
the full matrix-vector products has a cost of $\mathcal{O}(mn)$ for
a dense matrix.

For a feasible averaged point
$(\widehat x_t,\widehat y_t)$ at iteration $k$ for any $k \geq 0$, we evaluate the primal-dual gap function as
\begin{align*}
    \operatorname{Gap}(\widehat x_t,\widehat y_t)
    &:=
    \max_{\|y\|_2\leq 1}
    \left\langle A\widehat x_t-b,y\right\rangle
    -
    \min_{\|x\|_2\leq 1}
    \left\langle Ax-b,\widehat y_t\right\rangle.
\end{align*}
Both optimization problems admit closed-form solutions. In particular, we obtain
\begin{align*}
    \max_{\|y\|_2\leq 1}
    \left\langle A\widehat x_t-b,y\right\rangle
    &=
    \left\|A\widehat x_t-b\right\|_2,
    \\
    \min_{\|x\|_2\leq 1}
    \left\langle Ax-b,\widehat y_t\right\rangle
    &=
    -\left\|A^\top\widehat y_t\right\|_2
    -
    \left\langle b,\widehat y_t\right\rangle.
\end{align*}
Therefore, the primal-dual gap is evaluated exactly as
\begin{align*}
    \operatorname{Gap}(\widehat x_t,\widehat y_t)
    =
    \left\|A\widehat x_t-b\right\|_2
    +
    \left\|A^\top\widehat y_t\right\|_2
    +
    \left\langle b,\widehat y_t\right\rangle.
\end{align*}
The exact full matrix-vector products in this expression are used only to evaluate the reported primal-dual gap and are not used in the
stochastic algorithmic updates.

We present an ablation study on the sensitivity of EGM and DA EGM to the initialization $\widehat D$ and $D_0$, respectively, in Figure~\ref{fig:minmax_sensitivity}. The results show that Rule~2 with $\widehat D=3$ yields a slightly smaller primal--dual gap. For DA EGM, all tested initializations eventually achieve similar performance.


\paragraph{Bilinear saddle-point problem.}

For the saddle-point problem we set
\begin{align*}
    M=10,
    \qquad
    \beta=0.5,
    \qquad \text{and} \qquad
    \sigma=0.1 .
\end{align*}
Next, we show that the mapping $F$ satisfies Assumption~\ref{asum_map_growth}. For simplicity, let us denote $z = [x,y]$. We see that the mapping $F$ consists of a 1-H\"older smooth part $F_1$ and a 0-H\"older nonsmooth part $F_0$. Hence, we can write
\begin{align*}
    F(z) = F_1(z) + F_0(z) \quad \text{where } F_1(z) = M [y, -x] \quad \text{and } F_0(z) = \beta [\text{sign}(x), \text{sign}(y)].
\end{align*}
We see that
\begin{align*}
    \|F_1(z) - F_1(z^\prime)\|_2 = M \|z - z^\prime\|_2 
\end{align*}
which shows that $F_1$ is 1-H\"older continuous and by Lemma~\ref{lem_holder_to_delta_L}, we see that $H_1 = M$ and we obtain
\begin{align*}
    L_{F_1}(\delta_1) = M .
\end{align*}
For the mapping $F_0$, we see that
\begin{align*}
    \|F_0(z) - F_0(z^\prime)\|_2 \leq 2\sqrt{2} \beta .
\end{align*}
Hence by  Lemma~\ref{lem_holder_to_delta_L}, $H_0 = 2\sqrt{2} \beta$, and we see that
\begin{align*}
    L_{F_0}(\delta_0) = \frac{(2\sqrt{2}\beta)^2}{4\delta_0}
    =
    \frac{2\beta^2}{\delta_0}.
\end{align*}
By Lemma~\ref{lem_calculus}(iii), we can take $\delta = \delta_1+\delta_0$, and obtain an admissible continuity modulus for
$F=F_1+F_0$ as
\begin{align*}
    L_F(\delta)
    =
    L_1(\delta_1)+L_0(\delta_0) = M+\frac{2\beta^2}{\delta_0}.
\end{align*}
Since, the Lipschitz component $F_1$ requires no additive
approximation error, the full $\delta$-budget can be assigned to
$F_0$, giving $\delta_0=\delta$. Thus, an admissible continuity
modulus for the complete mapping is
\begin{align*}
    L_{F}(\delta) = M+\frac{2\beta^2}{\delta}.
\end{align*}
Hence, the mapping $F$ satisfies $(\delta,L)$-continuity.


\begin{figure*}[t]
    \centering

    \begin{subfigure}{0.32\textwidth}
        \centering
        \includegraphics[width=\linewidth]
        {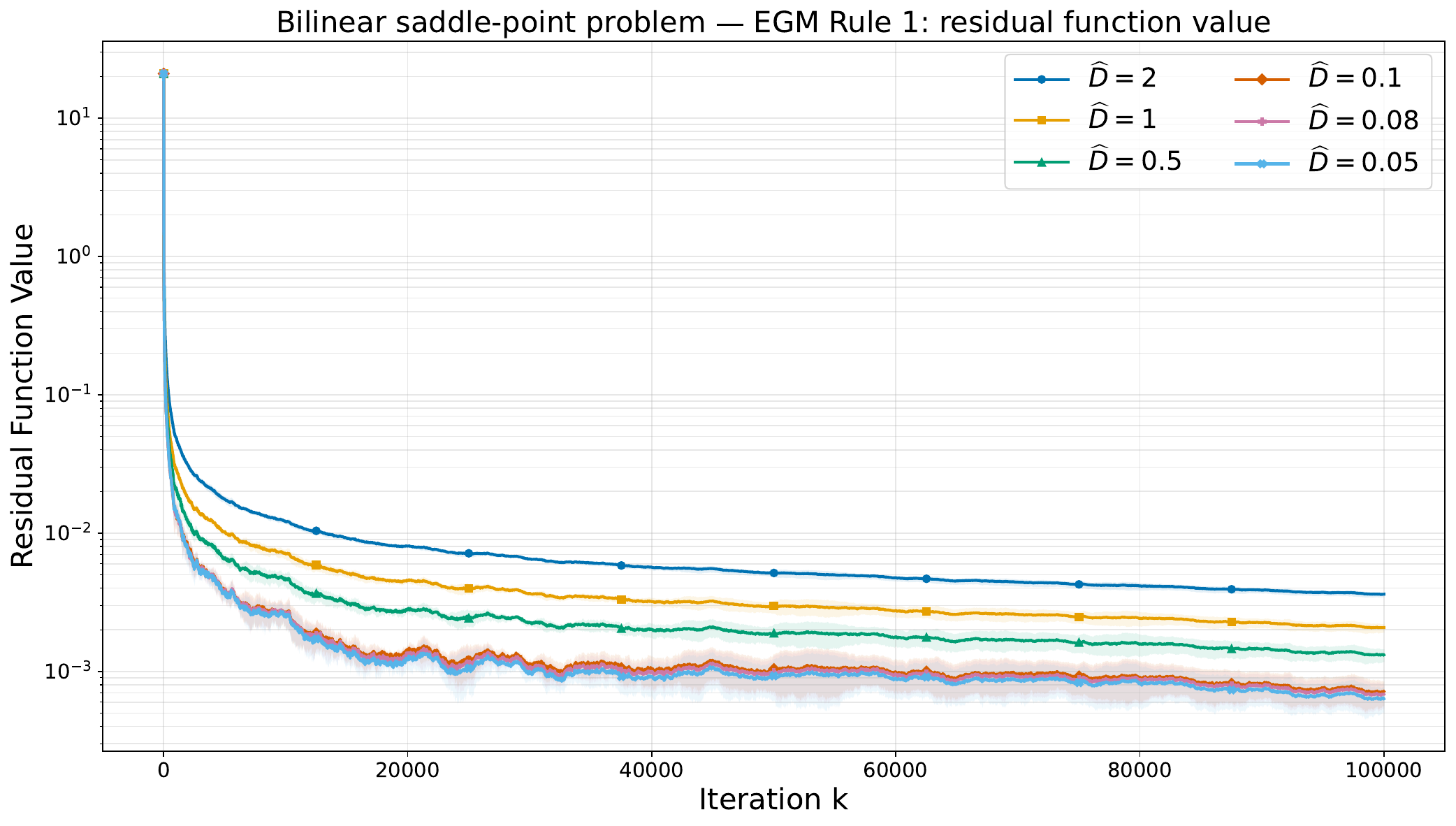}
    \end{subfigure}
    \hfill
    \begin{subfigure}{0.32\textwidth}
        \centering
        \includegraphics[width=\linewidth]
        {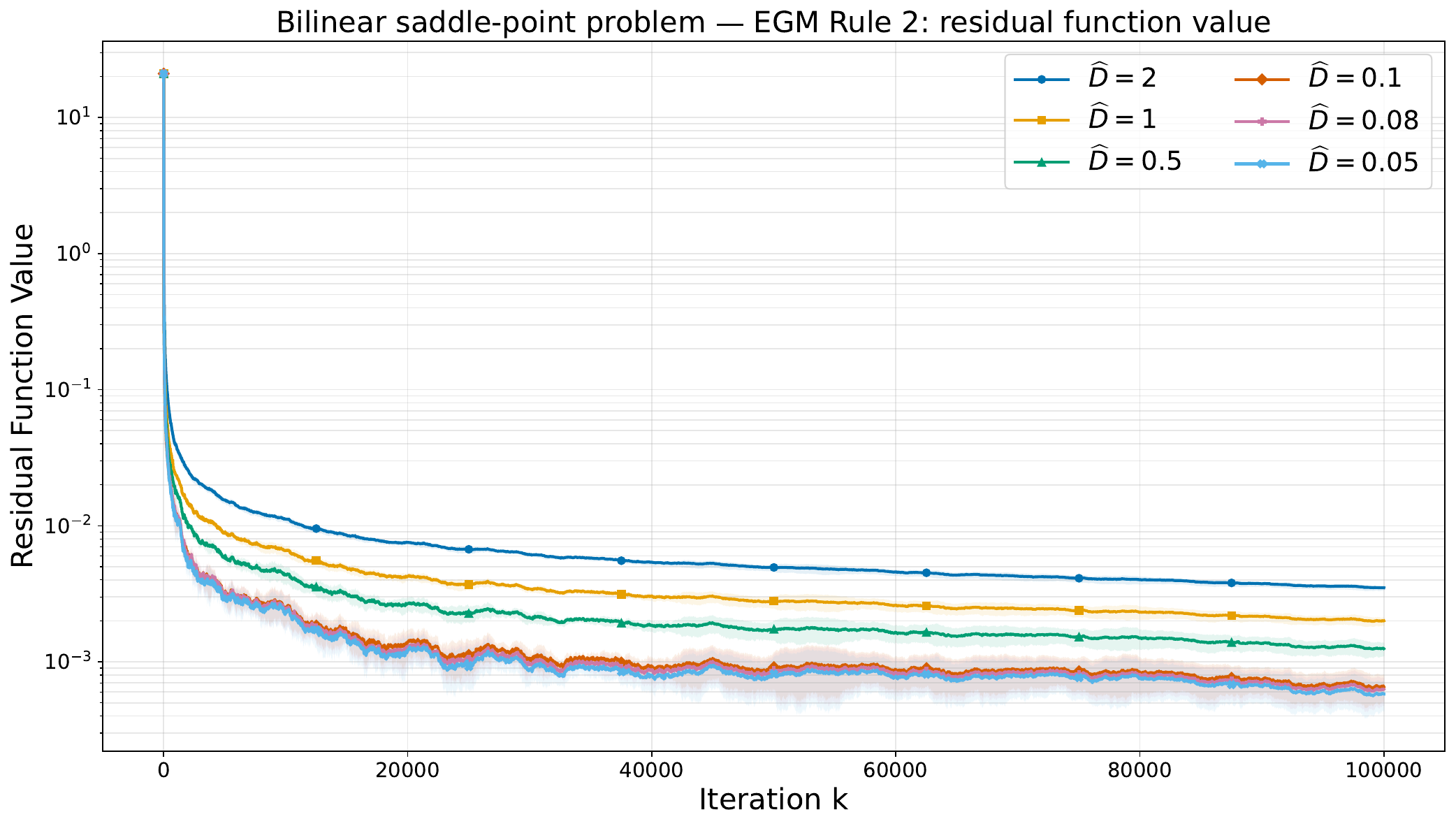}
    \end{subfigure}
    \hfill
    \begin{subfigure}{0.32\textwidth}
        \centering
        \includegraphics[width=\linewidth]
        {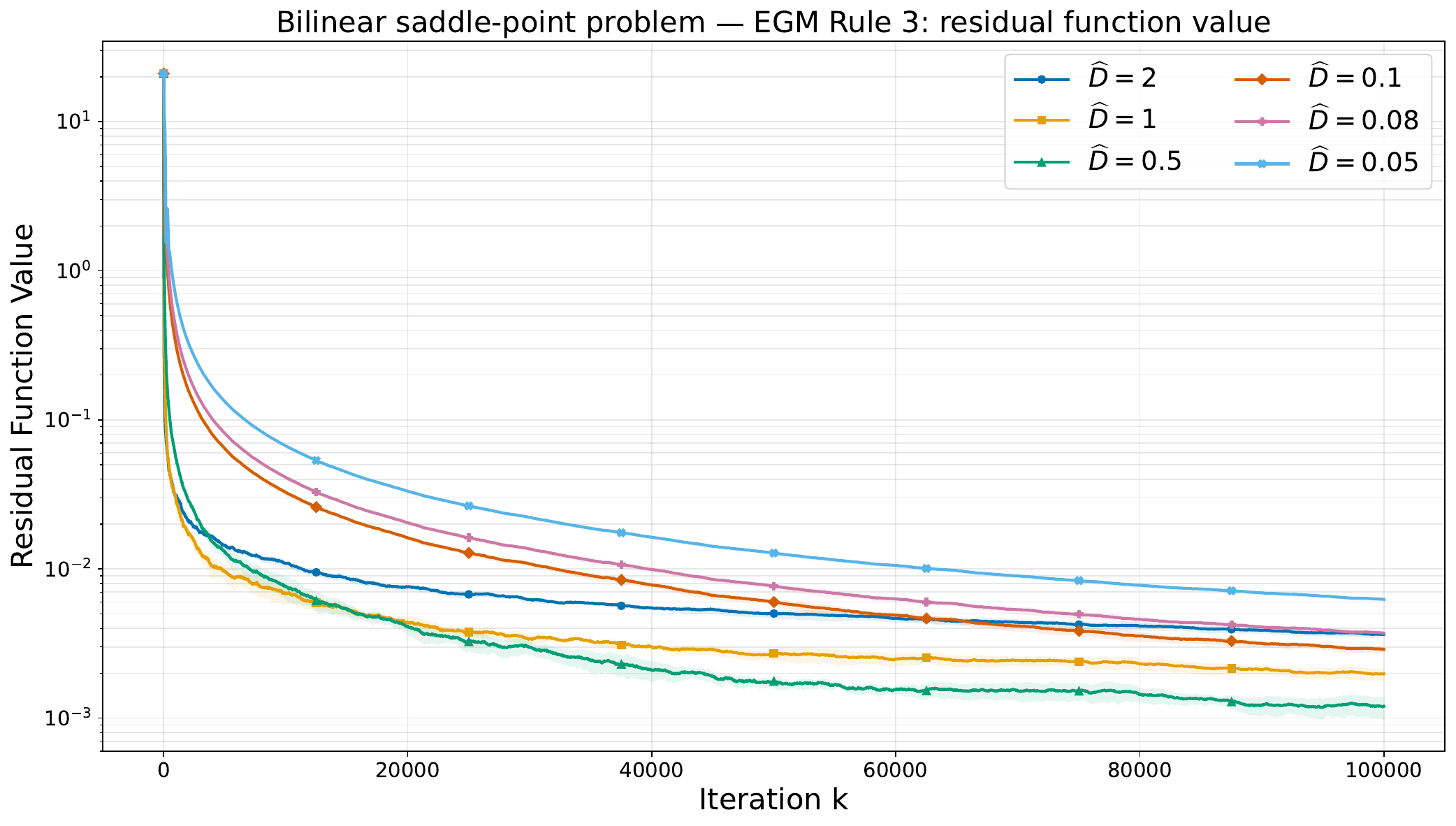}
    \end{subfigure}

    \par\medskip

    \begin{subfigure}{0.32\textwidth}
        \centering
        \includegraphics[width=\linewidth]
        {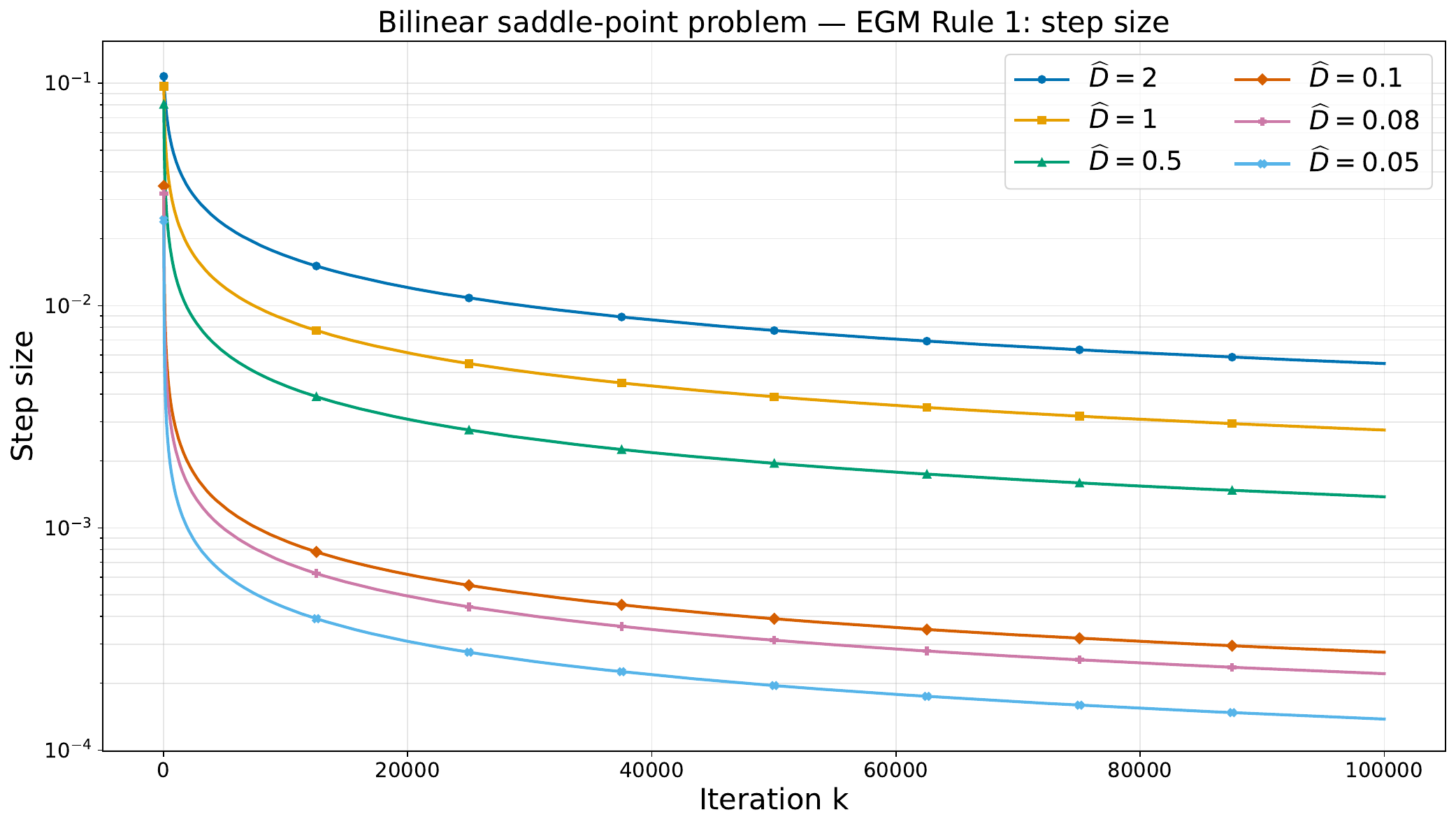}
    \end{subfigure}
    \hfill
    \begin{subfigure}{0.32\textwidth}
        \centering
        \includegraphics[width=\linewidth]
        {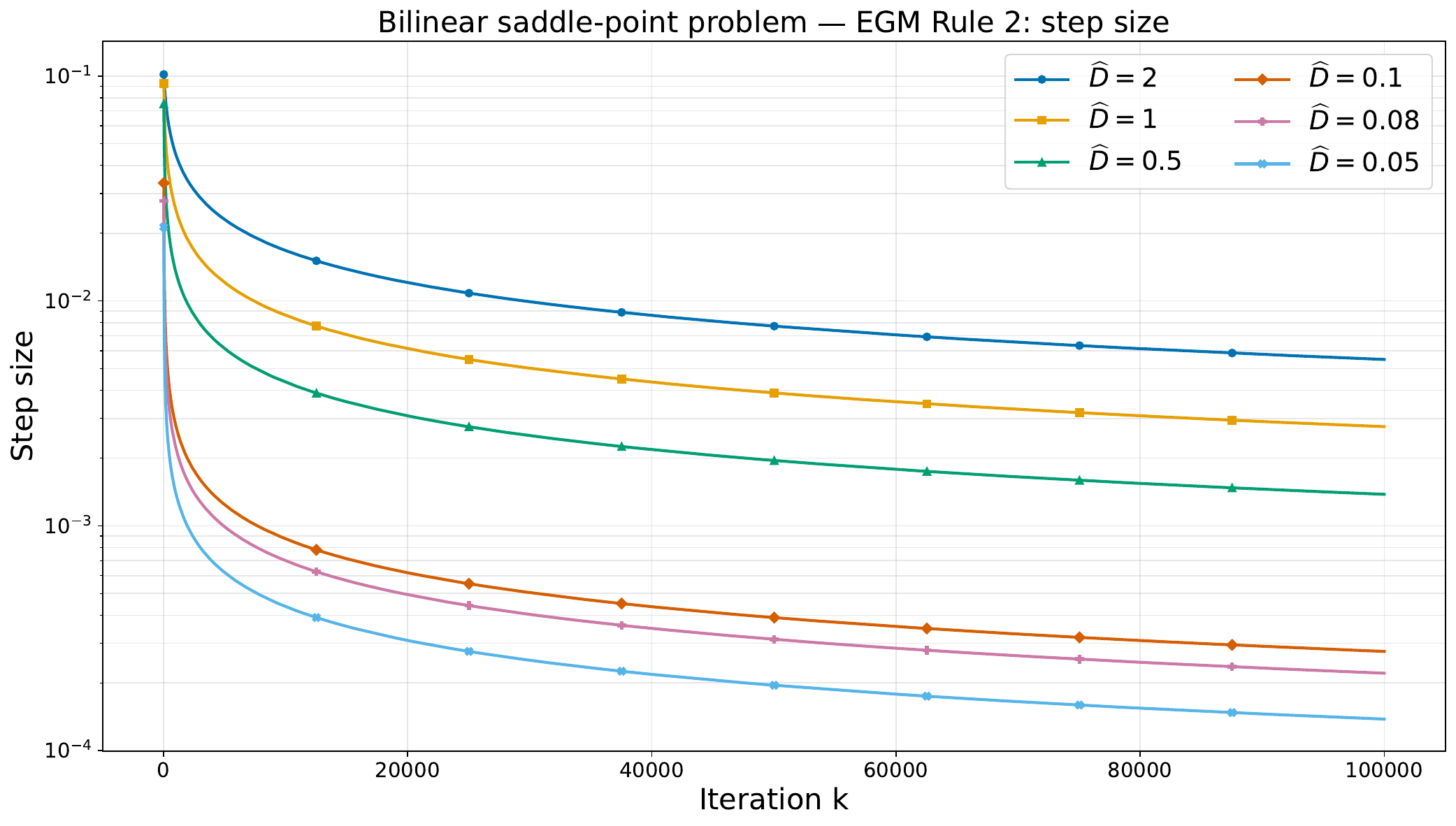}
    \end{subfigure}
    \hfill
    \begin{subfigure}{0.32\textwidth}
        \centering
        \includegraphics[width=\linewidth]
        {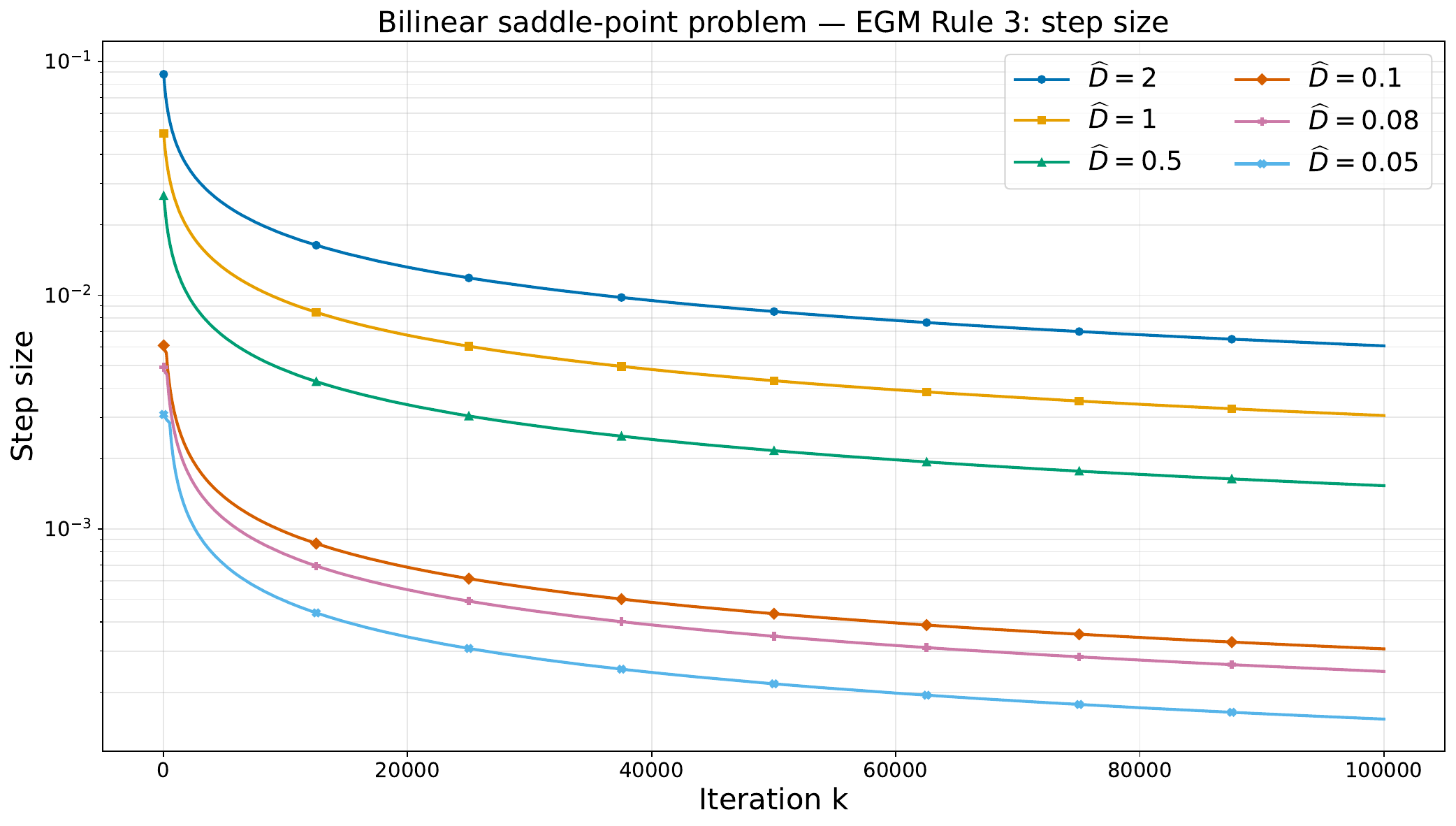}
    \end{subfigure}

    \caption{
        Sensitivity of EGM Rules 1--3 to $\widehat D$.
        The first row shows the residual function values, and the
        second row shows the corresponding step sizes.
    }
    \label{fig:saddle_egm_Dhat_sensitivity}
\end{figure*}



\begin{figure*}[htbp]
    \centering

    \begin{subfigure}{0.48\textwidth}
        \centering
        \includegraphics[width=\linewidth]
        {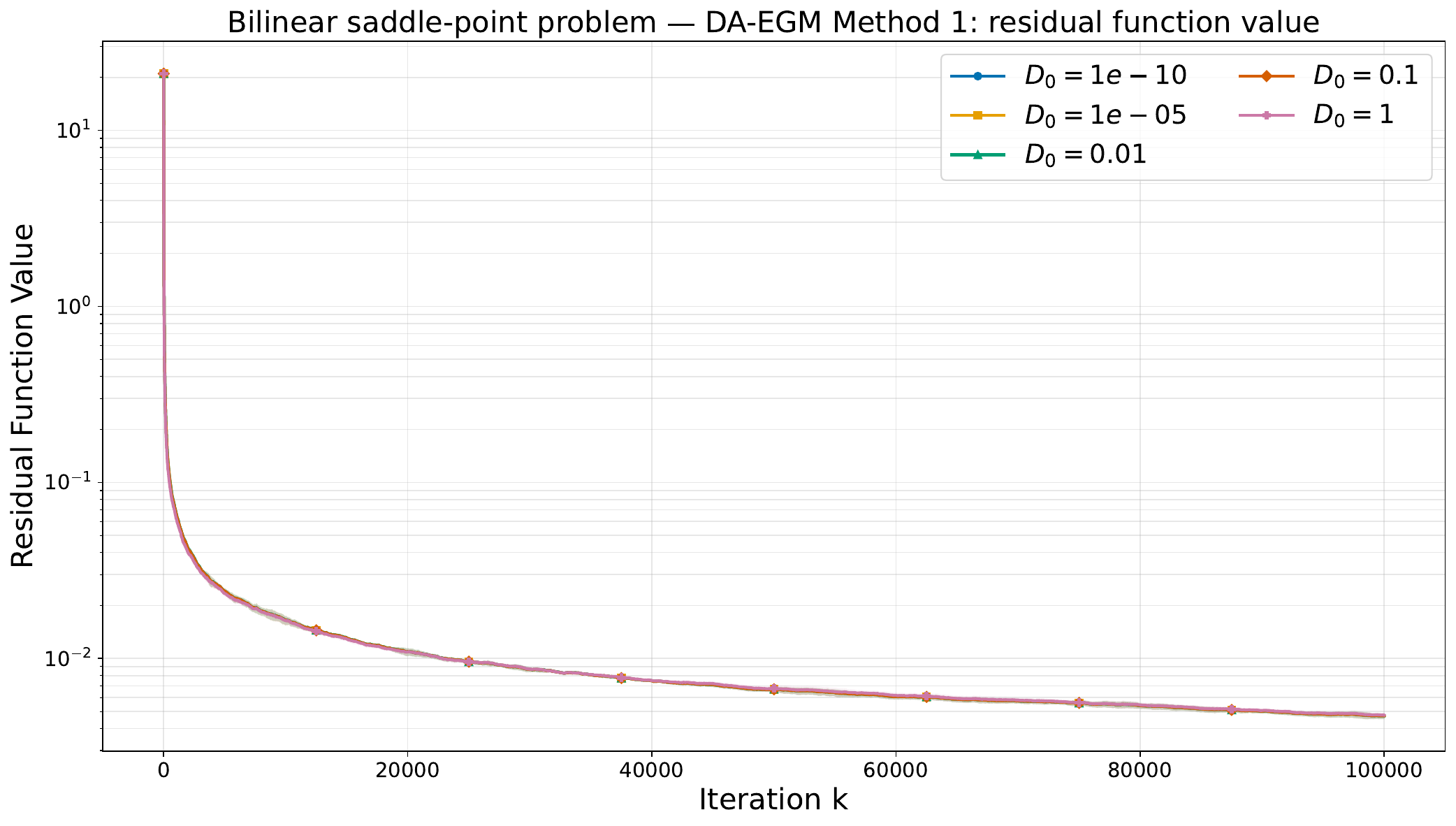}
    \end{subfigure}
    \hfill
    \begin{subfigure}{0.48\textwidth}
        \centering
        \includegraphics[width=\linewidth]
        {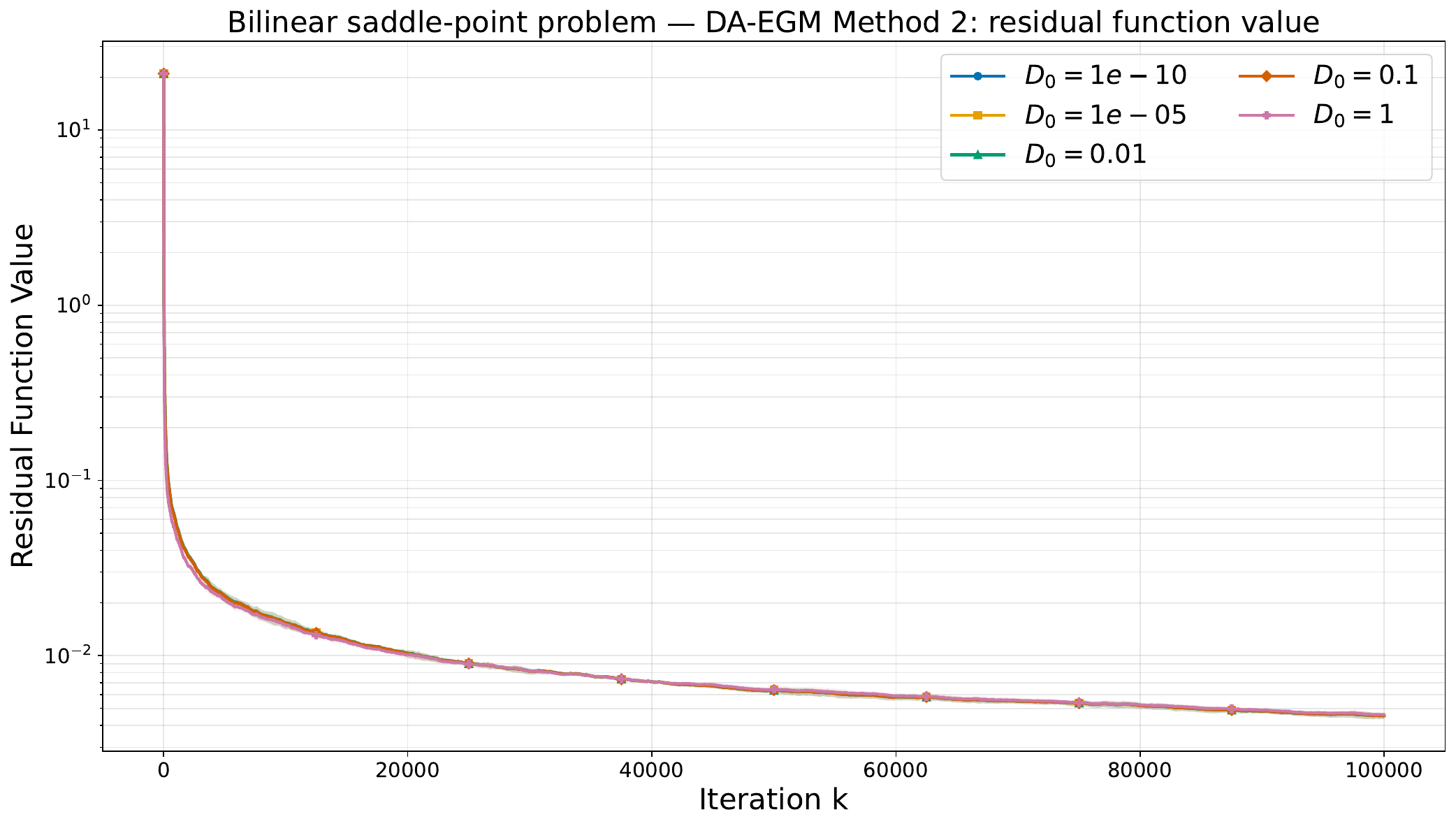}
    \end{subfigure}

    \par\medskip

    \begin{subfigure}{0.48\textwidth}
        \centering
        \includegraphics[width=\linewidth]
        {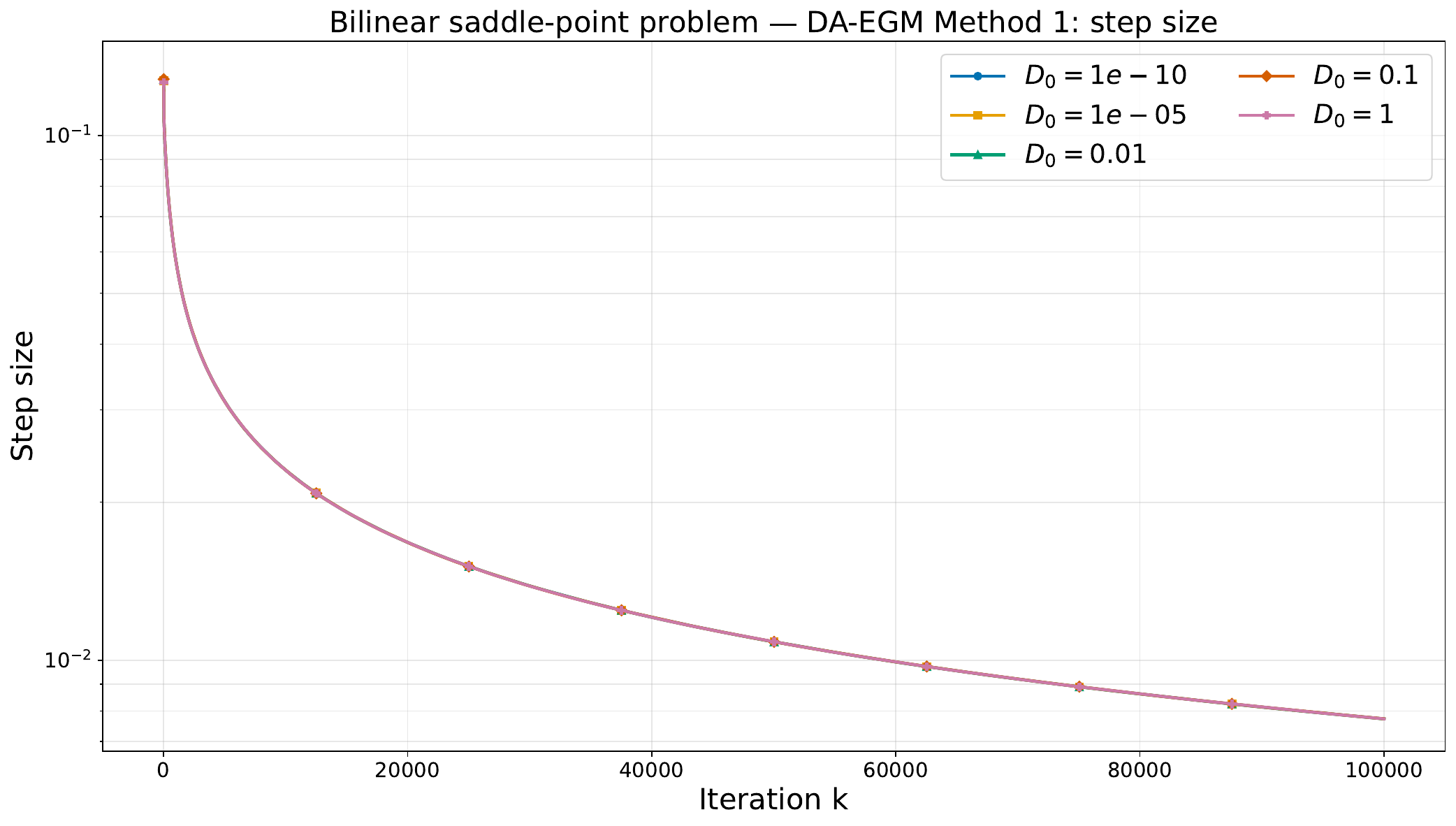}
    \end{subfigure}
    \hfill
    \begin{subfigure}{0.48\textwidth}
        \centering
        \includegraphics[width=\linewidth]
        {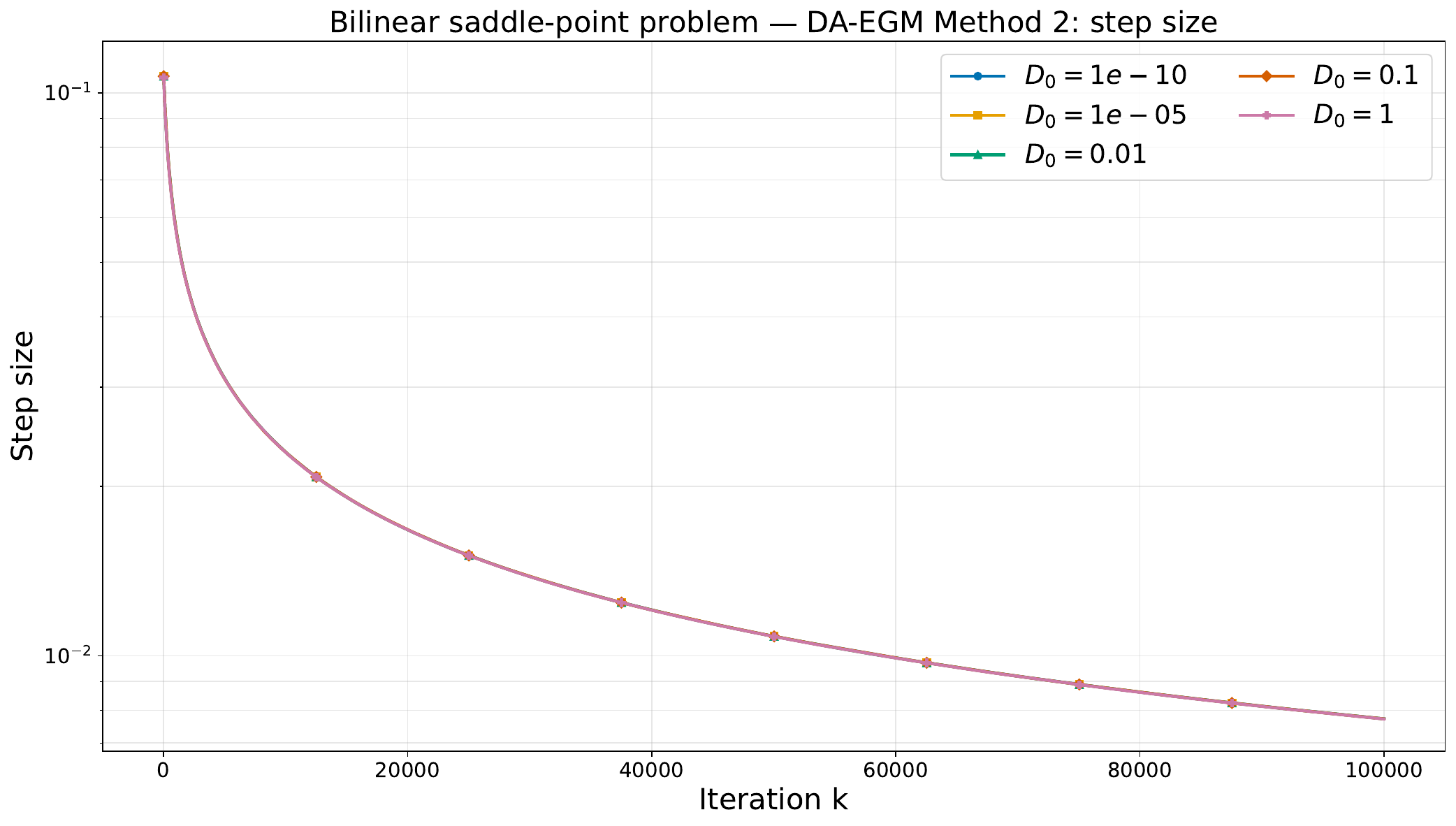}
    \end{subfigure}

    \caption{
        Sensitivity of DA-EGM Methods 1 and 2 to $D_0$.
        The first row shows the residual function values, and the
        second row shows the corresponding step sizes.
    }
    \label{fig:saddle_da_egm_D0_sensitivity}
\end{figure*}


Next, we show how the residual function is computed for the experimental evaluation. Let is denote the iterate at time $t$ as $z_t=[x_t,y_t] \in[-1,1]^2$. By the definition of the residual function (cf. relation~\eqref{residual_val_eq}), we see that
\begin{align*}
    \epsilon_k
    &=
    \max_{z\in[-1,1]^2}
    \frac{1}{\Lambda_k}
    \sum_{t=1}^k
    \lambda_t
    \left\langle F(z_t),z_t-z\right\rangle \\
    &= \frac{1}{\Lambda_k}
    \sum_{t=1}^k
    \lambda_t
    \left\langle F(z_t),z_t\right\rangle
    +
    \max_{z\in[-1,1]^2}
    \left\langle
        -\frac{1}{\Lambda_k}
        \sum_{t=1}^k\lambda_tF(z_t),
        z
    \right\rangle\\
    &=
    \frac{1}{\Lambda_k}
    \sum_{t=1}^k
    \lambda_t
    \left\langle F(z_t),z_t\right\rangle
    +
    \left\|
        \frac{1}{\Lambda_k}
        \sum_{t=1}^k\lambda_tF(z_t)
    \right\|_1,
\end{align*}
where the second equality follows because the support function of
$[-1,1]^2$ is the $\ell_1$ norm. The first term on the right hand side of the preceding relation can be simplified as
\begin{align*}
    \left\langle F(z_t),z_t\right\rangle
    =
    \left(
        My_t+\beta\operatorname{sign}(x_t)
    \right)x_t
    +
    \left(
        -Mx_t+\beta\operatorname{sign}(y_t)
    \right)y_t = \beta\left(|x_t|+|y_t|\right).
\end{align*}
Hence, we obtain the final expression of the residual as
\begin{align*}
    \epsilon_k \!=\! \frac{\beta}{\Lambda_k} \sum_{t=1}^k \lambda_t\left(|x_t| \!+\! |y_t|\right) + \left| \frac{1}{\Lambda_k} \sum_{t=1}^k \lambda_t \left(My_t \!+\! \beta\operatorname{sign}(x_t) \right) \right| + \left| \frac{1}{\Lambda_k}  \sum_{t=1}^k \lambda_t \left(-Mx_t \!+\! \beta\operatorname{sign}(y_t) \right) \right|.
\end{align*}

Next, Figures~\ref{fig:saddle_egm_Dhat_sensitivity} and~\ref{fig:saddle_da_egm_D0_sensitivity} present the sensitivity of EGM and DA EGM to the initialization parameters $\widehat D$ and $D_0$, respectively. We see that $\widehat D = 0.05$ work well for EGM Rules~1 and 2, whereas $\widehat D = 0.5$ works well for Rules~3. DA EGM experiments are pretty insensitive to the choice of the initialization $D_0$ for this problem.


\end{document}